\documentclass[12pt,thmsa]{article}%
\usepackage{amsmath}
\usepackage{sw20lart}
\usepackage{graphicx}%
\usepackage{amsfonts}%
\usepackage{amssymb}

\begin{document}

\author{Richard A. Smith}
\title{Introduction to Vector Spaces, Vector Algebras, and Vector Geometries}
\date{October 14, 2011}
\maketitle

\begin{abstract}
An introductory overview of vector spaces, algebras, and linear geometries
over an arbitrary commutative field is given. Quotient spaces are emphasized
and used in constructing the exterior and the symmetric algebras of a vector
space. The exterior algebra of a vector space and that of its dual are used in
treating linear geometry. Scalar product spaces, orthogonality, and the Hodge
star based on a general basis are treated.

\end{abstract}
\tableofcontents

\newpage

\section{Fundamentals of Structure}

\subsection{What is a Vector Space?}

A \textbf{vector space} is a structured set of elements called
\textbf{vectors}. The structure required for this set of vectors is that of an
\emph{Abelian group} with operators from a specified \emph{field}. (For us,
multiplication in a field is \emph{commutative}, and a number of our results
will depend on this in an essential way.) If the specified field is
$\mathcal{F}$, we say that the vector space is \textbf{over} $\mathcal{F}$. It
is standard for the Abelian group to be written additively, and for its
identity element to be denoted by $0$. The symbol $0$ is also used to denote
the field's zero element, but using the same symbol for both things will not
be confusing. The field elements will operate on the vectors by multiplication
on the left, the result always being a vector. The field elements will be
called \textbf{scalars} and be said to \textbf{scale} the vectors that they
multiply. The scalars are required to scale the vectors in a harmonious
fashion, by which is meant that the following three rules always hold for any
vectors and scalars. (Here $a$ and $b$ denote scalars, and $v$ and $w$ denote vectors.)

\begin{enumerate}
\item Multiplication by the field's unity element is the identity operation:
$1\cdot v=v$.

\item Multiplication is associative: $a\cdot(b\cdot v)=(ab)\cdot v$.

\item Multiplication distributes over addition: $(a+b)\cdot v=a\cdot v+b\cdot
v$ and $a\cdot(v+w)=a\cdot v+a\cdot w$.
\end{enumerate}

Four things which one might reasonably expect to be true are true.

\begin{proposition}
If $a$ is a scalar and $v$ a vector, $a\cdot0=0$, $0\cdot v=0$, $(-a)\cdot
v=-(a\cdot v)$, and given that $a\cdot v=0$, if $a\neq0$ then $v=0$.
\end{proposition}

Proof: $a\cdot0=a\cdot(0+0)=a\cdot0+a\cdot0$; adding $-\left(  a\cdot0\right)
$ to both ends yields the first result. $0\cdot v=(0+0)\cdot v=0\cdot v+0\cdot
v$; adding $-\left(  0\cdot v\right)  $ to both ends yields the second result.
Also, $0=0\cdot v=(a+(-a))\cdot v=a\cdot v+(-a)\cdot v$; adding $-\left(
a\cdot v\right)  $ on the front of both ends yields the third result. If
$a\cdot v=0$ with $a\neq0$, then $a^{-1}\cdot(a\cdot v)=a^{-1}\cdot0$, so that
$(a^{-1}a)\cdot v=0$ and hence $v=0$ as claimed. $\blacksquare$

\begin{exercise}
Given that $a\cdot v=0$, if $v\neq0$ then $a=0$.
\end{exercise}

\begin{exercise}
$-(a\cdot v)=a\cdot(-v)$.
\end{exercise}

Because of the harmonious fashion in which the scalars operate on the vectors,
the elementary arithmetic of vectors contains no surprises.

\subsection{Some Vector Space Examples}

There are many familiar spaces that are vector spaces, although in their usual
appearances they may have additional structure as well. A familiar vector
space over the real numbers $\mathbb{R}$ is $\mathbb{R}^{n}$, the space of
$n$-tuples of real numbers with component-wise addition and scalar
multiplication by elements of $\mathbb{R}$. Closely related is another vector
space over $\mathbb{R}$, the space $\Delta\mathbb{R}^{n}$ of all translations
of $\mathbb{R}^{n}$, the general element of which is the function from
$\mathbb{R}^{n}$ to itself that sends the general element ($\xi_{1}%
,...,\xi_{n})$ of $\mathbb{R}^{n}$ to ($\xi_{1}+\delta_{1},...,\xi_{n}%
+\delta_{n})$ for a fixed real $n$-tuple $(\delta_{1},...,\delta_{n})$. Thus
an element of $\Delta\mathbb{R}^{n}$ uniformly \textit{vectors} the elements
of $\mathbb{R}^{n}$ to new positions, displacing each element of
$\mathbb{R}^{n}$ by the same vector, and in this sense $\Delta\mathbb{R}^{n}$
is truly a space of vectors. In the same fashion as with $\mathbb{R}^{n}$ over
$\mathbb{R}$, over any field $\mathcal{F}$, the space $\mathcal{F}^{n}$ of all
$n$-tuples of elements of $\mathcal{F}$ constitutes a vector space. The space
$\mathbb{C}^{n}$ of all $n$-tuples of complex numbers in addition to being
considered to be a vector space over $\mathbb{C}$ may also be considered to be
a (different, because it is taken over a different field) vector space over
$\mathbb{R}$. Another familiar algebraic structure that is a vector space, but
which also has additional structure, is the space of all polynomials with
coefficients in the field $\mathcal{F}$. To conclude this initial collection
of examples, the set consisting of the single vector $0$ may be considered to
be a vector space over any field; these are the \textbf{trivial} vector spaces.

\subsection{Subspace, Linear Combination and Span}

Given a vector space $\mathcal{V}$ over the field $\mathcal{F}$, $\mathcal{U}
$ is a \textbf{subspace} of $\mathcal{V}$, denoted $\mathcal{U}%
\vartriangleleft\mathcal{V}$, if it is a subgroup of $\mathcal{V}$ that is
itself a vector space over $\mathcal{F}$. To show that a subset $\mathcal{U}$
of a vector space is a subspace, it suffices to show that $\mathcal{U}$ is
closed under sum and under product by any scalar. We call $\mathcal{V}$ and
$\left\{  0\right\}  $ \textbf{improper} subspaces of the vector space
$\mathcal{V}$, and we call all other subspaces \textbf{proper}.

\begin{exercise}
As a vector space over itself, $\mathbb{R}$ has no proper subspace. The set of
all integers is a subgroup, but not a subspace.
\end{exercise}

\begin{exercise}
The intersection of any number of subspaces is a subspace.
\end{exercise}

A \textbf{linear combination} of a finite set $\mathcal{S}$ of vectors is any
sum $\sum_{s\in\mathcal{S}}c_{s}\cdot s$ obtained by adding together exactly
one multiple, by some scalar coefficient, of each vector of $\mathcal{S}$. A
\textbf{linear combination} of an infinite set $\mathcal{S}$ of vectors is a
linear combination of any finite subset of $\mathcal{S}$. The empty sum which
results in forming the linear combination of the empty set is taken to be $0$,
by convention. A subspace must contain the value of each of its linear
combinations. If $\mathcal{S}$ is any subset of a vector space, by the
\textbf{span} of $\mathcal{S}$, denoted $\left\langle \mathcal{S}\right\rangle
$, is meant the set of the values of all linear combinations of $\mathcal{S}$.
$\left\langle \mathcal{S}\right\rangle $ is a subspace. If $\mathcal{U}%
=\left\langle \mathcal{S}\right\rangle $, then we say that $\mathcal{S}$
\textbf{spans} $\mathcal{U}$ and that it is a \textbf{spanning set} for
$\mathcal{U}$.

\begin{exercise}
For any set $\mathcal{S}$ of vectors, $\left\langle \mathcal{S}\right\rangle $
is the intersection of all subspaces that contain $\mathcal{S}$, and therefore
$\left\langle \mathcal{S}\right\rangle $ itself is the only subspace of
$\left\langle \mathcal{S}\right\rangle $ that contains $\mathcal{S}$.
\end{exercise}

\subsection{Independent Set, Basis and Dimension}

Within a given vector space, a \textbf{dependency} is said to exist in a given
set of vectors if the zero vector is the value of a nontrivial (scalar
coefficients not all zero) linear combination of some finite nonempty subset
of the given set. A set in which a dependency exists is (\textbf{linearly})
\textbf{dependent} and otherwise is (\textbf{linearly}) \textbf{independent}.
The ``linearly'' can safely be omitted and we will do this to simplify our writing.

For example, if $v\neq0$, then $\left\{  v\right\}  $ is an independent set.
The empty set is an independent set. Any subset of an independent set is
independent. On the other hand, $\left\{  0\right\}  $ is dependent, as is any
set that contains $0$.

\begin{exercise}
A set $\mathcal{S}$ of vectors is independent if and only if no member is in
the span of the other members: for all $v\in\mathcal{S}$, $v\notin\left\langle
\mathcal{S}\smallsetminus\left\{  v\right\}  \right\rangle $. What would the
similar statement need to be if the space were not a vector space over a
field, but only the similar type of space over the \emph{ring} $\mathbb{Z}$ of
ordinary integers?
\end{exercise}

Similar ideas of dependency and independence may be introduced for vector
sequences. A vector sequence is \textbf{dependent} if the set of its terms is
a dependent set \emph{or} if the sequence has any repeated terms. Otherwise
the sequence is \textbf{independent}.

\begin{exercise}
The terms of the finite vector sequence $v_{1},\ldots,v_{n}$ satisfy a linear
relation $a_{1}\cdot v_{1}+\cdots+a_{n}\cdot v_{n}=0$ with at least one of the
scalar coefficients $a_{i}$ nonzero if and only if the sequence is dependent.
If the sequence is dependent, then it has a term that is in the span of the
set of preceding terms and this set of preceding terms is nonempty if
$v_{1}\neq0$.
\end{exercise}

A \textbf{maximal independent set} in a vector space, i. e., an independent
set that is not contained in any other independent set, is said to be a
\textbf{basis set}, or, for short, a \textbf{basis} (plural: \textbf{bases}),
for the vector space. For example, $\left\{  (1)\right\}  $ is a basis for the
vector space $\mathcal{F}^{1}$, where $\mathcal{F}$ is any field. The empty
set is the unique basis for a trivial vector space.

\begin{exercise}
Let $\Phi$ be a family of independent sets which is linearly ordered by set
inclusion. Then the union of all the sets in $\Phi$ is an independent set.
\end{exercise}

\smallskip\ From the result of the exercise above and \emph{Zorn's Lemma}
applied to independent sets which contain a given independent set, we obtain
the following important result.

\begin{theorem}
\label{Basis}Every vector space has a basis, and, more generally, every
independent set is contained in some basis. $\blacksquare$
\end{theorem}

A basis may be characterized in various other ways besides its definition as a
maximal independent set.

\begin{proposition}
[Independent Spanning Set]A basis for the vector space $\mathcal{V}$ is
precisely an independent set that spans $\mathcal{V}$.
\end{proposition}

Proof: In the vector space $\mathcal{V}$, let $\mathcal{B}$ be an independent
set such that $\left\langle \mathcal{B}\right\rangle =\mathcal{V}$. Suppose
that $\mathcal{B}$ is contained in the larger independent set $\mathcal{C}$.
Then, choosing any $v\in\mathcal{C}\smallsetminus\mathcal{B}$, the set
$\mathcal{B}^{\prime}=\mathcal{B}\cup\left\{  v\right\}  $ is independent
because it is a subset of the independent set $\mathcal{C}$. But then
$v\notin\left\langle \mathcal{B}^{\prime}\smallsetminus\left\{  v\right\}
\right\rangle $, i. e., $v\notin\left\langle \mathcal{B}\right\rangle $,
contradicting $\left\langle \mathcal{B}\right\rangle =\mathcal{V}$. Hence
$\mathcal{B}$ is a basis for $\mathcal{V}$.

On the other hand, let $\mathcal{B}$ be a basis, i. e., a maximal independent
set, but suppose that $\mathcal{B}$ does not span $\mathcal{V}$. There must
then exist $v\in\mathcal{V}\smallsetminus\left\langle \mathcal{B}\right\rangle
$. However, $\mathcal{B}\cup\left\{  v\right\}  $ would then be independent,
since a nontrivial dependency relation in $\mathcal{B}\cup\left\{  v\right\}
$ would have to involve $v$ with a nonzero coefficient and this would put
$v\in\left\langle \mathcal{B}\right\rangle $. But the independence of
$\mathcal{B}\cup\left\{  v\right\}  $ contradicts the hypothesis that
$\mathcal{B}$ is a maximal independent set. Hence, the basis $\mathcal{B}$ is
an independent set that spans $\mathcal{V}$. $\blacksquare$

\begin{corollary}
An independent set is a basis for its own span. $\blacksquare$
\end{corollary}

\begin{proposition}
[Minimal Spanning Set]$\mathcal{B}$ is a basis for the vector space
$\mathcal{V}$ if and only if $\mathcal{B}$ is a \emph{\textbf{minimal spanning
set}} for $\mathcal{V}$, i. e., $\mathcal{B}$ spans $\mathcal{V}$, but no
subset of $\mathcal{B}$ unequal to $\mathcal{B}$ also spans $\mathcal{V}$.
\end{proposition}

Proof: Suppose $\mathcal{B}$ spans $\mathcal{V}$ and no subset of
$\mathcal{B}$ unequal to $\mathcal{B}$ also spans $\mathcal{V}$. If
$\mathcal{B}$ were not independent, there would exist $b\in\mathcal{B}$ for
which $b\in\left\langle \mathcal{B}\smallsetminus\left\{  b\right\}
\right\rangle $. But then the span of $\mathcal{B}\smallsetminus\left\{
b\right\}  $ would be the same as the span of $\mathcal{B}$, contrary to what
we have supposed. Hence $\mathcal{B}$ is independent, and by the previous
proposition, $\mathcal{B}$ must therefore be a basis for $\mathcal{V}$ since
$\mathcal{B}$ also spans $\mathcal{V}$.

On the other hand, suppose $\mathcal{B}$ is a basis for the vector space
$\mathcal{V}$. By the previous proposition, $\mathcal{B}$ spans $\mathcal{V}$.
Suppose that $\mathcal{A}$ is a subset of $\mathcal{B}$ that is unequal to
$\mathcal{B}$. Then no element $v\in$ $\mathcal{B}\smallsetminus\mathcal{A}$
can be in $\left\langle \mathcal{B}\smallsetminus\left\{  v\right\}
\right\rangle $ since $\mathcal{B}$ is independent, and certainly, then, $v$
$\notin$ $\left\langle \mathcal{A}\right\rangle $ since $\mathcal{A}%
\subset\mathcal{B}\smallsetminus\left\{  v\right\}  $. Thus $\mathcal{A}$ does
not span $\mathcal{V}$, and $\mathcal{B}$ is therefore a minimal spanning set
for $\mathcal{V}$. $\blacksquare$

\begin{exercise}
Consider a finite spanning set. Among its subsets that also span, there is at
least one of smallest size. Thus a basis must exist for any vector space that
has a finite spanning set, independently verifying what we already know from
Theorem \ref{Basis}.
\end{exercise}

\begin{exercise}
Over any field $\mathcal{F}$, the set of $n$ $n$-tuples that have a $1$ in one
position and $0$ in all other positions is a basis for $\mathcal{F}^{n}$.
\emph{(This basis is called the \textbf{standard basis} for }$\mathcal{F}%
^{n}\emph{\ over}$\emph{\ $\mathcal{F}$.)}
\end{exercise}

\begin{proposition}
[Unique Linear Representation]$\mathcal{B}$\emph{\ }is a basis for the vector
space $\mathcal{V}$ if and only if $\mathcal{B}$ has the \emph{\textbf{unique
linear representation property}}, i. e., each vector of $\mathcal{V}$ has a
unique linear representation as a linear combination $\sum_{x\in\mathcal{X}%
}a_{x}\cdot x$ where $\mathcal{X}$ is some finite subset of $\mathcal{B}$ and
all of the scalars $a_{x}$ are nonzero.
\end{proposition}

Proof: Let $\mathcal{B}$ be a basis for $\mathcal{V}$, and let $v\in
\mathcal{V}$. Since $\mathcal{V}$ is the span of $\mathcal{B}$, $v$ certainly
is a linear combination of a finite subset of $\mathcal{B}$ with all scalar
coefficients nonzero. Suppose that $v$ has two different expressions of this
kind. Subtracting, we find that $0$ is equal to a nontrivial linear
combination of the union $\mathcal{U}$ of the subsets of $\mathcal{B}$
involved in the two expressions. But $\mathcal{U}$ is an independent set since
it is a subset of the independent set $\mathcal{B}$. Hence $v$ has the unique
representation as claimed.

On the other hand, suppose that $\mathcal{B}$ has the unique linear
representation property. Then $\mathcal{B}$ spans $\mathcal{V}$. We complete
the proof by showing that $\mathcal{B}$ must be an independent set. Suppose
not. Then there would exist $b\in\mathcal{B}$ for which $b\in\left\langle
\mathcal{B}\smallsetminus\left\{  b\right\}  \right\rangle .$ But this $b$
would then have two different representations with nonzero coefficients,
contrary to hypothesis, since we always have $1\cdot b$ as one such
representation and there would also be another one not involving $b$ in virtue
of $b\in\left\langle \mathcal{B}\smallsetminus\left\{  b\right\}
\right\rangle $. $\blacksquare$

\begin{exercise}
A finite set of vectors is a basis for a vector space if and only if each
vector in the vector space has a unique representation as a linear combination
of this set: $\left\{  x_{1},\ldots,x_{n}\right\}  $ \emph{(}with distinct
$x_{i}$, of course\emph{)} is a basis if and only if each $v=a_{1}\cdot
x_{1}+\cdots+a_{n}\cdot x_{n}$ for \emph{unique} scalars $a_{1},\ldots,a_{n}$.
\end{exercise}

\begin{exercise}
If $S$ is a finite independent set, $\sum_{s\in\mathcal{S}}c_{s}\cdot
s=\sum_{s\in\mathcal{S}}d_{s}\cdot s$ implies $c_{s}=d_{s}$ for all $s$.
\end{exercise}

\begin{example}
In this brief digression we now apply the preceding two propositions. Let
$v_{0},v_{1},\ldots,v_{n}$ be vectors in a vector space over the field
$\mathcal{F}$, and suppose that $v_{0}$ is in the span $\mathcal{V}$ of the
other $v_{i}$. Then the equation
\[
\xi_{1}\cdot v_{1}+\cdots+\xi_{n}\cdot v_{n}=v_{0}%
\]
for the $\xi_{i}\in\mathcal{F}$ has at least one solution. Renumbering the
$v_{i}$ if necessary, we may assume that the distinct vectors $v_{1}%
,\ldots,v_{m}$ form a basis set for $\mathcal{V}$. If $m=n$ the equation has a
unique solution. Suppose, however, that $1\leqslant m<n.$ Then the equation
\[
\xi_{1}\cdot v_{1}+\cdots+\xi_{m}\cdot v_{m}=v_{0}-\xi_{m+1}\cdot
v_{m+1}-\cdots-\xi_{n}\cdot v_{n}%
\]
has a unique solution for $\xi_{1},\ldots,\xi_{m}$ for each fixed set of
values we give the other $\xi_{i}$, since the right-hand side is always in
$\mathcal{V}$. Thus $\xi_{1},\ldots,\xi_{m}$ are functions of the variables
$\xi_{m+1},\ldots,\xi_{n}$, where each of these variables is allowed to range
freely over the entire field $\mathcal{F}$. When $\mathcal{F}$ is the field
$\mathbb{R}$ of real numbers, we have deduced a special case of the
\emph{Implicit Function Theorem} of multivariable calculus. When the $v_{i}$
are $d$-tuples \emph{(}i. e., elements of $\mathcal{F}^{d}$\emph{)}, the
original vector equation is the same thing as the general set of $d$
consistent numerical linear equations in $n$ unknowns, and our result
describes how, in the general solution, $n-m$ of the $\xi_{i}$ are arbitrary
parameters upon which the remaining $m$ of the $\xi_{i}$ depend.
\end{example}

We have now reached a point where we are able to give the following key result.

\begin{theorem}
[Replacement Theorem]\label{ReplacementTheorem}Let $\mathcal{V}$ be a vector
space with basis $\mathcal{B}$, and let $\mathcal{C}$ be an independent set in
$\mathcal{V}$. Given any finite subset of $\mathcal{C}$ no larger than
$\mathcal{B}$, its elements can replace an equal number of elements of
$\mathcal{B}$ and the resulting set will still be a basis for $\mathcal{V}$.
\end{theorem}

Proof: The theorem holds in the case when no elements are replaced. Suppose
that the elements $y_{1},...,y_{N}$ of $\mathcal{C}$ have replaced $N$
elements of $\mathcal{B}$ and the resulting set $\mathcal{B}_{N}$ is still a
basis for $\mathcal{V}$. An element $y_{N+1}$ of $\mathcal{C}\smallsetminus
\left\{  y_{1},...,y_{N}\right\}  $ has a unique linear representation as a
linear combination with nonzero coefficients of some finite nonempty subset
$\mathcal{X}=\left\{  x_{1},...,x_{K}\right\}  $ of $\mathcal{B}_{N}$. There
must be an element $x^{\ast}\in\mathcal{X}\smallsetminus\left\{
y_{1},...,y_{N}\right\}  $ because $y_{N+1}$ cannot be a linear combination of
$\left\{  y_{1},...,y_{N}\right\}  $ since $\mathcal{C}$ is independent. In
$\mathcal{B}_{N}$ replace $x^{\ast}$ with $y_{N+1}$. Clearly, the result
$\mathcal{B}_{N+1}$ still spans $\mathcal{V}$.

The assumption that $x\in\mathcal{B}_{N+1}$ is a linear combination of the
other elements of $\mathcal{B}_{N+1}$ will be seen to always contradict the
independence of $\mathcal{B}_{N}$, proving that $\mathcal{B}_{N+1}$ is
independent. For if $x=y_{N+1}$, then we can immediately solve for $x^{\ast} $
as a linear combination of the other elements of $\mathcal{B}_{N}$. And if
$x\neq y_{N+1}$, writing $x=a\cdot y_{N+1}+y$, where $y$ is a linear
combination of elements of $\mathcal{B}_{N}\smallsetminus\left\{  x,x^{\ast
}\right\}  $, we find that we can solve for $x^{\ast}$ as a linear combination
of the other elements of $\mathcal{B}_{N}$ unless $a=0$, but then $x$ is a
linear combination of the other elements of $\mathcal{B}_{N}$.

Thus the theorem is shown to hold for $N+1$ elements replaced given that it
holds for $N$. By induction, the theorem holds in general. $\blacksquare$

\smallskip Suppose that $\mathcal{B}$ is a finite basis of $n$ vectors and
$\mathcal{C}$ is an independent set at least as large as $\mathcal{B}$. By the
theorem, we can replace all the vectors of $\mathcal{B}$ by any $n$ vectors
from $\mathcal{C}$, the set of which must be a basis, and therefore a maximal
independent set, so that $\mathcal{C}$ actually contains just $n$ elements.
Hence, when a finite basis exists, no independent set can exceed it in size,
and no basis can then be larger than any other basis. We therefore have the
following important result.

\begin{corollary}
\label{BasesEquipotent}If a vector space has a finite basis, then each of its
bases is finite, and each has the same number of elements. $\blacksquare
$\smallskip\ 
\end{corollary}

A vector space is called \textbf{finite}-\textbf{dimensional} if it has a
finite basis, and otherwise is called \textbf{infinite}-\textbf{dimensional}.
The \textbf{dimension} of a vector space is the number of elements in any
basis for it if it is finite-dimensional, and $\infty$ otherwise. We write
$\dim\mathcal{V}$ for the dimension of $\mathcal{V}$. If a finite-dimensional
vector space has dimension $n$, then we say that it is $\emph{n}%
$-\textbf{dimensional}. Over any field $\mathcal{F}$, $\mathcal{F}^{n}$ is
$\emph{n}$-dimensional.

\begin{exercise}
In an $n$-dimensional vector space, any vector sequence with $n+1$ terms is dependent.
\end{exercise}

\begin{exercise}
Let the vector space $\mathcal{V}$ have the same finite dimension as its
subspace $\mathcal{U}$. Then $\mathcal{U}=\mathcal{V}$.
\end{exercise}

The corollary above has the following analog for infinite-dimensional spaces.

\begin{theorem}
All bases of an infinite-dimensional vector space have the same cardinality.
\end{theorem}

Proof: Let $\mathcal{B}$ be an infinite basis for a vector space and let
$\mathcal{C}$ be another basis for the same space. For each $y\in\mathcal{C}$
let $\mathcal{X}_{y}$ be the finite nonempty subset of $\mathcal{B}$ such that
$y$ is the linear combination of its elements with nonzero scalar
coefficients. Then each $x\in\mathcal{B}$ must appear in some $\mathcal{X}%
_{y}$, for if some $x\in\mathcal{B}$ appears in no $\mathcal{X}_{y}$, then
that $x$ could not be in the span of $\mathcal{C}$. Thus $\mathcal{B}=$
$\bigcup_{y\in\mathcal{C}}\mathcal{X}_{y}$. Using $\left|  \mathcal{S}\right|
$ to denote cardinality of the set $\mathcal{S}$, we then have $\left|
\mathcal{B}\right|  =\left|  \bigcup_{y\in\mathcal{C}}\mathcal{X}_{y}\right|
$ and $\mathcal{C}$ must be infinite for this to hold. Since $\mathcal{C}$ is
infinite, for each $y$ we have $|\mathcal{X}_{y}|\leqslant|\mathcal{C}|$ and
therefore $\left|  \mathcal{B}\right|  \leqslant\left|  \mathcal{C\times
C}\right|  =|\mathcal{C}|$, where the last equality is a well-known, but not
so easily proved, result found in books on set theory and in online references
such as the Wikipedia entry for \textit{Cardinal number}. Similarly, we find
$\left|  \mathcal{C}\right|  \leqslant\left|  \mathcal{B}\right|  $. The
\emph{Schroeder-Bernstein theorem} then gives $\left|  \mathcal{B}\right|
=\left|  \mathcal{C}\right|  $. $\blacksquare$

\subsection{Sum and Direct Sum of Subspaces}

By the \textbf{sum} of any number of subspaces is meant the span of their
union. We use ordinary additive notation for subspace sums.

\begin{exercise}
For any sets of vectors, the sum of their spans is the span of their union.
\end{exercise}

\begin{exercise}
Each sum of subspaces is equal to the set of all finite sums of vectors
selected from the various subspace summands, no two vector summands being
selected from the same subspace summand.
\end{exercise}

The trivial subspace $\left\{  0\right\}  $ acts as a neutral element in
subspace sums. Because of this, we will henceforth often write $0$ for
$\left\{  0\right\}  $, particularly when subspace sums are involved. By
convention, the empty subspace sum is $0$.

With \emph{partial order} $\subset$, \emph{meet} $\cap$, and \emph{join} $+ $,
the subspaces of a given vector space form a \emph{complete lattice}, which is
\emph{modular}, but generally not \emph{distributive}.

\begin{exercise}
For subspaces, $(\mathcal{U}_{1}\cap\mathcal{U}_{3})+(\mathcal{U}_{2}%
\cap\mathcal{U}_{3})\subset(\mathcal{U}_{1}+\mathcal{U}_{2})\cap
\mathcal{U}_{3}$, and \emph{(Modular Law)} $\mathcal{U}_{2}\vartriangleleft
\mathcal{U}_{3}\Rightarrow(\mathcal{U}_{1}\cap\mathcal{U}_{3})+\mathcal{U}%
_{2}=(\mathcal{U}_{1}+\mathcal{U}_{2})\cap\mathcal{U}_{3}$.
\end{exercise}

\begin{example}
In $\mathbb{R}^{2}$, let $\mathcal{X}$ denote the $x$-axis, let $\mathcal{Y}$
denote the $y$-axis, and let $\mathcal{D}$ denote the line where $y=x$. Then
$\left(  \mathcal{X}+\mathcal{Y}\right)  \cap\mathcal{D}=\mathcal{D}$ while
$\left(  \mathcal{X}\cap\mathcal{D}\right)  +\left(  \mathcal{Y}%
\cap\mathcal{D}\right)  =0$.
\end{example}

A sum of subspaces is called \textbf{direct} if the intersection of each
summand with the sum of the other summands is $0$. We often replace $+$ with
$\oplus,$ or $\sum$ with $\bigoplus$, to indicate that a subspace sum is direct.

The summands in a direct sum may be viewed as being independent, in some
sense. Notice, though, that $0$ may appear as a summand any number of times in
a direct sum. Nonzero direct summands must be unique, however, and more
generally, we have the following result.

\begin{lemma}
If a subspace sum is direct, then, for each summand $\mathcal{U}$,
$\mathcal{U}\cap\Sigma=0$, where $\Sigma$ denotes the sum of any selection of
the summands which does not include $\mathcal{U}$. In particular, two nonzero
summands must be distinct.
\end{lemma}

Proof: $\Sigma\subset\overline{\Sigma}$ where $\overline{\Sigma}$ denotes the
sum of all the summands with the exception of $\mathcal{U}$. By the definition
of direct sum, $\mathcal{U}\cap\overline{\Sigma}=0$, and hence $\mathcal{U}%
\cap\Sigma=0$ as was to be shown. $\blacksquare$\smallskip\ 

This lemma leads immediately to the following useful result.

\begin{theorem}
\label{SelectionofDirectSummands}If a subspace sum is direct, then the sum of
any selection of the summands is direct. $\blacksquare$
\end{theorem}

There are a number of alternative characterizations of directness of sums of
finitely many subspaces. (The interested reader will readily extend the first
of these, and its corollary, to similar results for sums of infinitely many subspaces.)

\begin{theorem}
$u\in\mathcal{U}_{1}\oplus\cdots\oplus\,\mathcal{U}_{n}\Leftrightarrow
u=u_{1}+\cdots+u_{n}$ for \emph{unique} $u_{1}\in\mathcal{U}_{1},\ldots
,u_{n}\in\mathcal{U}_{n}$.
\end{theorem}

Proof: Suppose that $u\in\mathcal{U}_{1}\oplus\cdots\oplus\,\mathcal{U}_{n}$,
and that $u=u_{i}+\overline{u}_{i}=v_{i}+\overline{v}_{i}$ where $u_{i}$ and
$v_{i}$ are in $\mathcal{U}_{i}$ and $\overline{u}_{i}$ and $\overline{v}_{i}$
are in $\sum_{j\neq i}\mathcal{U}_{j}$. Then $w=u_{i}-v_{i}=\overline{v}%
_{i}-\overline{u}_{i}$ must be both in $\mathcal{U}_{i}$ and in $\sum_{j\neq
i}\mathcal{U}_{j}$ so that $w=0$ and therefore $u_{i}=v_{i}$. Hence the
representation of $u$ as $u=u_{1}+\cdots+u_{n}$, where $u_{1}\in
\mathcal{U}_{1},\ldots,u_{n}\in\mathcal{U}_{n}$, is unique.

On the other hand, suppose that each $u\in\mathcal{U}_{1}+\cdots
+\,\mathcal{U}_{n}$ has the unique representation $u=u_{1}+\cdots+u_{n}$ where
$u_{1}\in\mathcal{U}_{1},\ldots,u_{n}\in\mathcal{U}_{n}$. Let $u\in
\mathcal{U}_{i}\cap\sum_{j\neq i}\mathcal{U}_{j}$. Then, since $u\in
\mathcal{U}_{i}$, it has the representation $u=u_{i}$ for some $u_{i}%
\in\mathcal{U}_{i}$, and the uniqueness of the representation of $u$ as
$u=u_{1}+\cdots+u_{n}$ where $u_{1}\in\mathcal{U}_{1},\ldots,u_{n}%
\in\mathcal{U}_{n}$ then implies that $u_{j}=0$ for all $j\neq i$. Also, since
$u\in\sum_{j\neq i}\mathcal{U}_{j}$, it has the representation $u=\sum_{j\neq
i}u_{j}$ for some $u_{j}\in\mathcal{U}_{j}$, and the uniqueness of the
representation of $u$ as $u=u_{1}+\cdots+u_{n}$ where $u_{1}\in\mathcal{U}%
_{1},\ldots,u_{n}\in\mathcal{U}_{n}$ implies that $u_{i}=0$. Hence $u=0$ and
therefore for each $i$, $\mathcal{U}_{i}\cap\sum_{j\neq i}\mathcal{U}_{j}=0$.
The subspace sum therefore is direct. $\blacksquare$

\begin{corollary}
$\mathcal{U}_{1}+\cdots+\,\mathcal{U}_{n}$ is direct if and only if $0$ has a
unique representation as $\sum_{1\leqslant j\leqslant n}u_{j}$, where each
$u_{j}\in\mathcal{U}_{j}$.
\end{corollary}

Proof: If $\mathcal{U}_{1}+\cdots+\,\mathcal{U}_{n}$ is direct, then $0$ has a
unique representation of the form $\sum_{1\leqslant j\leqslant n}u_{j}$, where
each $u_{j}\in\mathcal{U}_{j}$, namely that one where each $u_{j}=0$. On the
other hand, suppose that $0$ has that unique representation. Then given
$\sum_{1\leqslant j\leqslant n}u_{j}^{^{\prime}}=\sum_{1\leqslant j\leqslant
n}u_{j}^{^{\prime\prime}}$, where each $u_{j}^{^{\prime}}\in\mathcal{U}_{j}$
and each $u_{j}^{^{\prime\prime}}\in\mathcal{U}_{j}$, $0=\sum_{1\leqslant
j\leqslant n}u_{j}$, where each $u_{j}=u_{j}^{^{\prime}}-u_{j}^{^{\prime
\prime}}\in\mathcal{U}_{j}$. Hence $u_{j}^{^{\prime}}-u_{j}^{^{\prime\prime}%
}=0$ for each $j$, and the subspace sum is direct. $\blacksquare$

\begin{proposition}
$\mathcal{U}_{1}+\cdots+\,\mathcal{U}_{n}$ is direct if and only if both
$\overline{\mathcal{U}}_{n}=\mathcal{U}_{1}+\cdots+\,\mathcal{U}_{n-1}$ and
$\overline{\mathcal{U}}_{n}+\,\mathcal{U}_{n}$ are direct, or, more
informally, $\mathcal{U}_{1}\oplus\cdots\oplus\,\mathcal{U}_{n}=\left(
\mathcal{U}_{1}\oplus\cdots\oplus\,\mathcal{U}_{n-1}\right)  \oplus
\,\mathcal{U}_{n}$.
\end{proposition}

Proof: Suppose that $\mathcal{U}_{1}+\cdots+\,\mathcal{U}_{n}$ is direct. By
Theorem \ref{SelectionofDirectSummands}, $\overline{\mathcal{U}}_{n}$ is
direct, and the definition of direct sum immediately implies that
$\overline{\mathcal{U}}_{n}+\,\mathcal{U}_{n}$ is direct.

On the other hand, suppose that both $\overline{\mathcal{U}}_{n}%
=\mathcal{U}_{1}+\cdots+\,\mathcal{U}_{n-1}$ and $\overline{\mathcal{U}}%
_{n}+\,\mathcal{U}_{n}$ are direct. Then $0$ has a unique representation of
the form $u_{1}+\cdots+u_{n-1}$, where $u_{j}\in\mathcal{U}_{j}$, and $0$ has
a unique representation of the form $\overline{u}_{n}+u_{n}$, where
$\overline{u}_{n}\in\overline{\mathcal{U}}_{n}$ and $u_{n}\in\,\mathcal{U}%
_{n}$. Let $0$ be represented in the form $u_{1}^{^{\prime}}+\cdots
+u_{n}^{^{\prime}}$, where $u_{j}^{^{\prime}}\in\mathcal{U}_{j}$. Then
$0=\overline{u}_{n}^{^{\prime}}+u_{n}^{^{\prime}}$ where $\overline{u}%
_{n}^{^{\prime}}=u_{1}^{^{\prime}}+\cdots+u_{n-1}^{^{\prime}}\in
\overline{\mathcal{U}}_{n}$. Hence $\overline{u}_{n}^{^{\prime}}=0$ and
$u_{n}^{^{\prime}}=0$ by the corollary above, but then, by the same corollary,
$u_{1}^{^{\prime}}=\cdots=u_{n-1}^{^{\prime}}=0$. Using the corollary once
more, $\mathcal{U}_{1}+\cdots+\,\mathcal{U}_{n}$ is proved to be direct.
$\blacksquare$\smallskip\ \ 

Considering in turn $\left(  \mathcal{U}_{1}\oplus\mathcal{U}_{2}\right)
\oplus\mathcal{U}_{3}$, $\left(  \left(  \mathcal{U}_{1}\oplus\mathcal{U}%
_{2}\right)  \oplus\mathcal{U}_{3}\right)  \oplus\mathcal{U}_{4}$, etc., we
get the following result.

\begin{corollary}
$\mathcal{U}_{1}+\cdots+\mathcal{U}_{n}$ is direct if and only if
$(\mathcal{U}_{1}+\cdots+\,\mathcal{U}_{i})\cap\mathcal{U}_{i+1}=0$ for
$i=1,\ldots,n-1$. $\blacksquare$
\end{corollary}

\begin{exercise}
Let $\Sigma_{k}=\mathcal{U}_{i_{k-1}+1}+\cdots+\mathcal{U}_{i_{k}}$ for
$k=1,\ldots,K$, where $0=i_{0}<i_{1}<\cdots<i_{K}=n$. Then $\mathcal{U}%
_{1}+\cdots+\,\mathcal{U}_{n}$ is direct if and only if each $\Sigma_{k}$ sum
is direct and $\Sigma_{1}+\cdots+\Sigma_{n}$ is direct.
\end{exercise}

\begin{exercise}
$\mathcal{U}_{1}\oplus\cdots\oplus\mathcal{U}_{6}=((\mathcal{U}_{1}%
\oplus(\mathcal{U}_{2}\oplus\mathcal{U}_{3}))\oplus(\mathcal{U}_{4}%
\oplus\mathcal{U}_{5}))\oplus\mathcal{U}_{6}$ in the fullest sense.
\end{exercise}

In the finite-dimensional case, additivity of dimension is characteristic.

\begin{theorem}
The finite-dimensional subspace sum $\mathcal{U}_{1}+\cdots+\mathcal{U}_{n}$
is direct if and only if
\[
\dim\left(  \mathcal{U}_{1}+\cdots+\mathcal{U}_{n}\right)  =\dim
\mathcal{U}_{1}+\cdots+\dim\mathcal{U}_{n}.
\]
\end{theorem}

Proof: Let each $\mathcal{U}_{j}$ have the basis $\mathcal{B}_{j}$. It is
clear that $\mathcal{B}_{1}\cup\cdots\cup\mathcal{B}_{n}$ spans the subspace sum.

Suppose now that $\mathcal{U}_{1}+\cdots+\mathcal{U}_{n}$ is direct.
$\mathcal{B}_{1}\cup\cdots\cup\mathcal{B}_{n}$ is independent, and therefore a
basis, because, applying the definition of direct sum, a dependency in it
immediately implies one in one of the $\mathcal{B}_{j}$. And, of course, no
element of any $\mathcal{B}_{i}$ can be in any of the other $\mathcal{B}_{j}$,
or it would be in the sum of their spans. Hence the $\mathcal{B}_{j}$ are
disjoint and $\dim\left(  \mathcal{U}_{1}+\cdots+\mathcal{U}_{n}\right)
=\dim\mathcal{U}_{1}+\cdots+\dim\mathcal{U}_{n}$.

On the other hand, suppose that $\dim\left(  \mathcal{U}_{1}+\cdots
+\mathcal{U}_{n}\right)  =\dim\mathcal{U}_{1}+\cdots+\dim\mathcal{U}_{n}$.
Because a minimal spanning set is a basis, $\mathcal{B}_{1}\cup\cdots
\cup\mathcal{B}_{n}$ must contain at least $\dim\mathcal{U}_{1}+\cdots
+\,\dim\mathcal{U}_{n}$ distinct elements, but clearly cannot contain more,
and must therefore be a basis for the subspace sum. A nonzero element in
$\mathcal{U}_{i}=\left\langle \mathcal{B}_{i}\right\rangle $ and
simultaneously in the span of the other $\mathcal{B}_{j}$, would entail a
dependency in the basis $\mathcal{B}_{1}\cup\cdots\cup\mathcal{B}_{n}$, which
of course is not possible. Hence no nonzero element of any $\mathcal{U}_{i}$
can be in the sum of the other $\mathcal{U}_{j}$ and therefore the subspace
sum is direct. $\blacksquare$

\begin{exercise}
Let $\mathcal{U}$ be any subspace of the vector space $\mathcal{V}$. As we
already know, $\mathcal{U}$ has a basis $\mathcal{A}$ which is part of a basis
$\mathcal{B}$ for $\mathcal{V}$. Let $\overline{\mathcal{U}}=\left\langle
\mathcal{B}\smallsetminus\mathcal{A}\right\rangle $. Then $\mathcal{V}%
=\mathcal{U}\oplus\overline{\mathcal{U}}$.
\end{exercise}

\begin{exercise}
Let $\mathcal{B}$ be a basis for $\mathcal{V}$. Then $\mathcal{V}%
=\bigoplus_{x\in\mathcal{B}}\left\langle \left\{  x\right\}  \right\rangle $.
\end{exercise}

\subsection{Problems}

\begin{enumerate}
\item Assuming only that $+$ is a group operation, not necessarily that of an
\emph{Abelian group}, deduce from the rules for scalars operating on vectors
that $u+v=v+u$ for any vectors $v$ and $u$.

\item Let $\mathcal{T}$ be a subspace that does not contain the vectors $u$
and $v$. Then
\[
v\in\left\langle \mathcal{T}\cup\left\{  u\right\}  \right\rangle
\Leftrightarrow u\in\left\langle \mathcal{T}\cup\left\{  v\right\}
\right\rangle \text{.}%
\]

\item Let $\mathcal{S}$ and $\mathcal{T}$ be independent subsets of a vector
space and suppose that $\mathcal{S}$ is finite and $\mathcal{T}$ is larger
than $\mathcal{S}$. Then $\mathcal{T}$ contains a vector $x$ such that
$\mathcal{S}\cup\left\{  x\right\}  $ is independent. Deduce Corollary
\ref{BasesEquipotent} from this.

\item Let $\mathcal{B}$ and $\mathcal{B}^{\prime}$ be bases for the same
vector space. Then for every $x\in\mathcal{B}$ there exists $x^{\prime}%
\in\mathcal{B}^{\prime}$ such that both $\left(  \mathcal{B}\smallsetminus
\left\{  x\right\}  \right)  \cup\left\{  x^{\prime}\right\}  $ and $\left(
\mathcal{B}^{\prime}\smallsetminus\left\{  x^{\prime}\right\}  \right)
\cup\left\{  x\right\}  $ are also bases.

\item Let $\mathcal{W}$ and $\mathcal{X}$ be finite-dimensional subspaces of a
vector space. Then
\[
\dim\left(  \mathcal{W}\cap\mathcal{X}\right)  \leqslant\dim\mathcal{X}%
\text{\qquad and\qquad}\dim\left(  \mathcal{W}+\mathcal{X}\right)
\geqslant\dim\mathcal{W}\text{.}%
\]
When does equality hold?

\item Let $\mathcal{W}$ and $\mathcal{X}$ be subspaces of a vector space
$\mathcal{V}$. Then $\mathcal{V}$ has a basis $\mathcal{B}$ such that each of
$\mathcal{W}$ and $\mathcal{X}$ are spanned by subsets of $\mathcal{B}$.

\item Instead of elements of a field $\mathcal{F}$, let the set of scalars be
the integers $\mathbb{Z}$, so the result is not a vector space, but what is
known as a $\mathbb{Z}$-\emph{module}. $\mathbb{Z}_{3}$, the integers modulo
3, is a $\mathbb{Z}$-module that has no linearly independent element. Thus
$\mathbb{Z}_{3}$ has $\varnothing$ as its only maximal linearly independent
set, but $\left\langle \varnothing\right\rangle \neq\mathbb{Z}_{3}$.
$\mathbb{Z}$ itself is also a $\mathbb{Z}$-module, and has $\{2\}$ as a
maximal linearly independent set, but $\left\langle \{2\}\right\rangle
\neq\mathbb{Z}$. On the other hand, $\{1\}$ is also a maximal linearly
independent set in $\mathbb{Z}$, and $\left\langle \{1\}\right\rangle
=\mathbb{Z}$, so we are prone to declare that $\{1\}$ is a basis for
$\mathbb{Z}$. Every $\mathbb{Z}$-module, in fact, has a maximal linearly
independent set, but the span of such a set may very well not be the whole
$\mathbb{Z}$-module. We evidently should reject \emph{maximal linearly
independent set} as a definition of \emph{basis} for $\mathbb{Z}$-modules, but
\emph{linearly independent spanning set}, equivalent for vector spaces, does
give us a reasonable definition for those $\mathbb{Z}$-modules that are indeed
spanned by an independent set. Let us adopt the latter as our definition for
basis of a $\mathbb{Z}$-module. Then, is it possible for a $\mathbb{Z}$-module
to have finite bases with differing numbers of elements? Also, is it then true
that no finite $\mathbb{Z}$-module has a basis?
\end{enumerate}

\newpage

\section{Fundamentals of Maps}

\subsection{Structure Preservation and Isomorphism}

A structure-preserving function from a structured set to another with the same
kind of structure is known in general as a \textbf{homomorphism}. The shorter
term \textbf{map} is often used instead, as we shall do here. When the
specific kind of structure is an issue, a qualifier may also be included, and
so one speaks of group maps, ring maps, etc. We shall be concerned here mainly
with vector space maps, so for us a map, unqualified, is a vector space map.
(A vector space map is also commonly known as a \textbf{linear transformation}.)

For two vector spaces to be said to have the same kind of structure, we
require them to be over the same field. Then for $f$ to be a map between them,
we require $f(x+y)=f(x)+f(y)$ and $f\left(  a\cdot x\right)  =a\cdot f\left(
x\right)  $ for all vectors $x,y$ and any scalar $a$, or, equivalently, that
$f\left(  a\cdot x+b\cdot y\right)  =a\cdot f\left(  x\right)  +b\cdot
f\left(  y\right)  $ for all vectors $x,y$ and all scalars $a$,$b$. It is easy
to show that a map sends $0$ to $0$, and a constant map must then send
everything to $0$. Under a map, the image of a subspace is always a subspace,
and the inverse image of a subspace is always a subspace. The identity
function on a vector space is a map. Maps which compose have a map as their composite.

\begin{exercise}
Let $\mathcal{V}$ and $\mathcal{W}$ be vector spaces over the same field and
let $\mathcal{S}$ be a spanning set for $\mathcal{V}$. Then any map $f:$
$\mathcal{V}\rightarrow\mathcal{W}$ is already completely determined by its
values on $\mathcal{S}$.
\end{exercise}

The values of a map may be arbitrarily assigned on a minimal spanning set, but
not on any larger spanning set.

\begin{theorem}
\label{BasisFreedom}Let $\mathcal{V}$ and $\mathcal{W}$ be vector spaces over
the same field and let $\mathcal{B}$ be a basis for $\mathcal{V}$. Then given
any function $f_{0}:\mathcal{B}\rightarrow\mathcal{W}$, there is a unique map
$f:\mathcal{V}\rightarrow\mathcal{W}$ such that $f$ agrees with $f_{0}$ on
$\mathcal{B}$.
\end{theorem}

Proof: By the unique linear representation property of a basis, given
$v\in\mathcal{V}$, there is a unique subset $\mathcal{X}$ of $\mathcal{B}$ and
unique nonzero scalars $a_{x}$ such that $v=\sum_{x\in\mathcal{X}}a_{x}\cdot
x$. Because a map preserves linear combinations, any map $f$ that agrees with
$f_{0}$ on $\mathcal{B}$ can only have the value $f\left(  v\right)
=\sum_{x\in\mathcal{X}}a_{x}\cdot f(x)=\sum_{x\in\mathcal{X}}a_{x}\cdot
f_{0}(x)$. Setting $f\left(  v\right)  =\sum_{x\in\mathcal{X}}a_{x}\cdot
f_{0}(x)$ does define a function $f:\mathcal{V}\rightarrow\mathcal{W}$ and
this function $f$ clearly agrees with $f_{0}$ on $\mathcal{B}$. Moreover, one
immediately verifies that this $f$ is a map. $\blacksquare$

It is standard terminology to refer to a function that is both one-to-one
(\emph{injective}) and onto (\emph{surjective}) as \emph{bijective}, or
\emph{invertible}. Each invertible $f:\mathcal{X}\rightarrow\mathcal{Y}$ has a
(unique) \emph{inverse }$f^{-1}$ such that $f^{-1}\circ f$ and $f\circ f^{-1}
$are the respective identity functions on $\mathcal{X}$ and $\mathcal{Y}$, and
an $f$ that has such an $f^{-1}$ is invertible. The composite of functions is
bijective if and only if each individual function is bijective.

For some types of structure, the inverse of a bijective map need not always be
a map itself. However, for vector spaces, inverting does always yield another map.

\begin{theorem}
The inverse of a bijective map is a map.
\end{theorem}

Proof: $f^{-1}\left(  a\cdot v+b\cdot w\right)  =a\cdot f^{-1}\left(
v\right)  +b\cdot f^{-1}\left(  w\right)  $ means precisely the same as
$f\left(  a\cdot f^{-1}\left(  v\right)  +b\cdot f^{-1}\left(  w\right)
\right)  =a\cdot v+b\cdot w$, which is clearly true. $\blacksquare$

\begin{corollary}
A one-to-one map preserves independent sets, and conversely, a map that
preserves independent sets is one-to-one.
\end{corollary}

Proof: A one-to-one map sends its domain bijectively onto its image and this
image is a subspace of the codomain. Suppose a set in the one-to-one map's
image is dependent. Then clearly the inverse image of this set is also
dependent. An independent set therefore cannot be sent to a dependent set by a
one-to-one map.

On the other hand, suppose that a map sends the distinct vectors $u$ and $v$
to the same image vector. Then it sends the nonzero vector $v-u$ to $0$, and
hence it sends the independent set $\left\{  v-u\right\}  $ to the dependent
set $\left\{  0\right\}  $. $\blacksquare$

\smallskip Because their inverses are also maps, the bijective maps are the
\textbf{isomorphisms} of vector spaces. If there is a bijective map from the
vector space $\mathcal{V}$ onto the vector space $\mathcal{W}$, we say that
$\mathcal{W}$ is \textbf{isomorphic} to $\mathcal{V}$. The notation
$\mathcal{V}\cong\mathcal{W}$ is commonly used to indicate that $\mathcal{V}$
and $\mathcal{W}$ are isomorphic. Viewed as a relation between vector spaces,
isomorphism is \emph{reflexive}, \emph{symmetric} and \emph{transitive}, hence
is an \emph{equivalence relation}. If two spaces are isomorphic to each other
we say that each is an \textbf{alias} of the other.

\begin{theorem}
Let two vector spaces be over the same field. Then they are isomorphic if and
only if they have bases of the same cardinality.
\end{theorem}

Proof: Applying Theorem \ref{BasisFreedom}, the one-to-one correspondence
between their bases extends to an isomorphism. On the other hand, an
isomorphism restricts to a one-to-one correspondence between bases.
$\blacksquare$

\begin{corollary}
Any $n$-dimensional vector space over the field $\mathcal{F}$ is isomorphic to
$\mathcal{F}^{n}$. $\blacksquare$
\end{corollary}

\begin{corollary}
\emph{(\textbf{Fundamental Theorem of Finite-Dimensional Vector Spaces})} Two
finite-dimensional vector spaces over the same field are isomorphic if and
only if they have the same dimension. $\blacksquare$
\end{corollary}

\smallskip\ The following well-known and often useful result may seem somewhat
surprising on first encounter.

\begin{theorem}
\label{Injmap}A map of a finite-dimensional vector space into an alias of
itself is one-to-one if and only if it is onto.
\end{theorem}

Proof: Suppose the map is one-to-one. Since one-to-one maps preserve
independent sets, the image of a basis is an independent set with the same
number of elements. The image of a basis must then be a maximal independent
set, and hence a basis for the codomain, since the domain and codomain have
equal dimension. Because a basis and its image are in one-to-one
correspondence under the map, expressing elements in terms of these bases, it
is clear that each element has a preimage element.

On the other hand, suppose the map is onto. The image of a basis is a spanning
set for the codomain and must contain at least as many distinct elements as
the dimension of the codomain. Since the domain and codomain have equal
dimension, the image of a basis must in fact then be a minimal spanning set
for the codomain, and therefore a basis for it. Expressing elements in terms
of a basis and its image, we find that, due to unique representation, if the
map sends $u$ and $v$ to the same element, then $u=v$. $\blacksquare$

\smallskip\ For infinite-dimensional spaces, the preceding theorem fails to
hold. It is a pleasant dividend of finite dimension, analogous to the result
that a function between equinumerous finite sets is one-to-one if and only if
it is onto.

\subsection{Kernel, Level Sets and Quotient Space\label{LevelSets}}

The \textbf{kernel} of a map is the inverse image of $\left\{  0\right\}  $.
The kernel of $f$ is a subspace of the domain of $f$. A \textbf{level set} of
a map is a nonempty inverse image of a singleton set. The kernel of a map, of
course, is one of its level sets, the only one that contains $0$ and is a
subspace. The level sets of a particular map make up a \emph{partition} of the
domain of the map. The following proposition internally characterizes the
level sets as the \emph{cosets} of the kernel.

\begin{proposition}
Let $\mathcal{L}$ be a level set and $\mathcal{K}$ be the kernel of the map
$f$. Then for any $v$ in $\mathcal{L}$, $\mathcal{L}=v+\mathcal{K}=\left\{
v+x\mid x\in\mathcal{K}\right\}  $.
\end{proposition}

Proof: Let $v\in\mathcal{L}$ and $x\in\mathcal{K}$. Then
$f(v+x)=f(v)+f(x)=f(v)+0$ so that $v+x\in\mathcal{L}$. Hence $v+\mathcal{K}%
\subset\mathcal{L}$. On the other hand, suppose that $u\in\mathcal{L}$ and
write $u=v+(u-v)$. Then since both $u$ and $v$ are in $\mathcal{L}$,
$f(u)=f(v)$, and hence $f(u-v)=0$, so that $u-v\in\mathcal{K}$. Hence
$\mathcal{L}\subset v+\mathcal{K}$. $\blacksquare$

\begin{corollary}
A map is one-to-one if and only if its kernel is $\left\{  0\right\}  $.
$\blacksquare$
\end{corollary}

\smallskip\ Given any vector $v$ in the domain of $f$, it must be in some
level set, and now we know that this level set is $v+\mathcal{K}$, independent
of the $f$ which has domain $\mathcal{V}$ and kernel $\mathcal{K} $.

\begin{corollary}
Maps with the same domain and the same kernel have identical level set
families. For maps with the same domain, then, the kernel uniquely determines
the level sets. $\blacksquare$
\end{corollary}

\smallskip\ The level sets of a map $f$ are in one-to-one correspondence with
the elements of the image of the map. There is a very simple way in which the
level sets can be made into a vector space isomorphic to the image of the map.
Just use the one-to-one correspondence between the level sets and the image
elements to give the image elements new labels according to the
correspondence. Thus $z$ gets the new label $f^{-1}\left(  \left\{  z\right\}
\right)  $. To see what $f^{-1}\left(  \left\{  z\right\}  \right)
+f^{-1}\left(  \left\{  w\right\}  \right)  $ is, figure out what $z+w$ is and
then take $f^{-1}(\left\{  z+w\right\}  )$, and similarly to see what $a\cdot
f^{-1}(\left\{  z\right\}  )$ is, figure out what $a\cdot z$ is and take
$f^{-1}(\left\{  a\cdot z\right\}  )$. All that has been done here is to give
different names to the elements of the image of $f$. This level-set alias of
the image of $f$ is the \textbf{vector space quotient (}or\textbf{\ quotient
space}) of the domain $\mathcal{V}$ of $f$ modulo the subspace $\mathcal{K}%
=\operatorname*{Kernel}(f)$ (which is the $0$ of the space), and is denoted by
$\mathcal{V}/\mathcal{K}$. The following result ``internalizes'' the
operations in $\mathcal{V}/\mathcal{K}$.

\begin{proposition}
Let $f$ be a map with domain $\mathcal{V}$ and kernel $\mathcal{K}$. Let
$\mathcal{L}=u+\mathcal{K}$ and $\mathcal{M}=v+\mathcal{K}$. Then in
$\mathcal{V}/\mathcal{K}$, $\mathcal{L}+\mathcal{M}=(u+v)+\mathcal{K}$, and
for any scalar $a$, $a\cdot\mathcal{L}=a\cdot u+\mathcal{K}$.
\end{proposition}

Proof: $\mathcal{L}=f^{-1}(\left\{  f\left(  u\right)  \right\}  )$ and
$\mathcal{M}=f^{-1}(\left\{  f\left(  v\right)  \right\}  )$. Hence
$\mathcal{L}+\mathcal{M}=f^{-1}(\left\{  f\left(  u\right)  +f\left(
v\right)  \right\}  )=f^{-1}(\left\{  f\left(  u+v\right)  \right\}
)=(u+v)+\mathcal{K}$. Also, $a\cdot\mathcal{L}=f^{-1}(\left\{  a\cdot f\left(
u\right)  \right\}  )=f^{-1}(\left\{  f(a\cdot u)\right\}  )=a\cdot
u+\mathcal{K}$. $\blacksquare$

\medskip\ The following is immediate.

\begin{corollary}
$\mathcal{V}/\mathcal{K}$ depends only on $\mathcal{V}$ and $\mathcal{K}$, and
not on the map $f$ that has domain $\mathcal{V}$ and kernel $\mathcal{K}$.
$\blacksquare$
\end{corollary}

\smallskip\ The next result tells us, for one thing, that kernel subspaces are
not special, and hence the quotient space exists for any subspace
$\mathcal{K}$. (In the vector space $\mathcal{V}$, the subspace
$\mathcal{\overline{%
K%
}}$ is a \textbf{complementary subspace, }or\textbf{\ complement}, of the
subspace $\mathcal{K}$ if $\mathcal{K}$ and $\mathcal{\overline{%
K%
}}$ have disjoint bases that together form a basis for $\mathcal{V}$, or what
is the same, $\mathcal{V}=\mathcal{K}\bigoplus\mathcal{\overline{%
K%
}}$. Every subspace $\mathcal{K} $ of $\mathcal{V}$ has at least one
complement $\mathcal{\overline{%
K%
}}$, because any basis $\mathcal{A}$ of $\mathcal{K}$ is contained in some
basis $\mathcal{B}$ of $\mathcal{V}$, and we may then take $\mathcal{\overline
{%
K%
}}=\left\langle \mathcal{B}\smallsetminus\mathcal{A}\right\rangle $.)

\begin{proposition}
\label{Complementary}Let $\mathcal{K}$ be a subspace of the vector space
$\mathcal{V}$ and let $\mathcal{\overline{%
K%
}}$ be a complementary subspace of $\mathcal{K}$. Then there is a map $\phi$
from $\mathcal{V}$ into itself which has $\mathcal{K}$ as its kernel and has
$\mathcal{\overline{%
K%
}}$ as its image. Hence given any subspace $\mathcal{K}$ of $\mathcal{V}$, the
quotient space $\mathcal{V}/\mathcal{K}$ exists and is an alias of any
complementary subspace of $\mathcal{K}$.
\end{proposition}

Proof: Let $\mathcal{A}$ be a basis for $\mathcal{K}$, and let $\overline
{\mathcal{A}}$ be a basis for $\mathcal{\overline{%
K%
}}$, so that $\mathcal{A}\cup\overline{\mathcal{A}}$ is a basis for
$\mathcal{V}$. Define $\phi$ as the map of $\mathcal{V}$ into itself which
sends each element of $\mathcal{A}$ to $0$, and each element of $\overline
{\mathcal{A}}$ to itself. Then $\phi$ has kernel $\mathcal{K}$ and image
$\mathcal{\overline{%
K%
}}$. $\blacksquare$

\smallskip Supposing that $\mathcal{K}$ is the kernel of some map
$f:\mathcal{V}\rightarrow\mathcal{W}$, the image of $f$ is an alias of the
quotient space, as is the image of the map $\phi$ of the proposition above.
Hence any complementary subspace of the kernel of a map is an alias of the
map's image. For maps with finite-dimensional domains, we then immediately
deduce the following important and useful result. (The \textbf{rank} of a map
is the dimension of its image, and the \textbf{nullity} of a map is the
dimension of its kernel.)

\begin{theorem}
[Rank Plus Nullity Theorem]If the domain of a map has finite dimension, the
sum of its rank and its nullity equals the dimension of its domain.
$\blacksquare$
\end{theorem}

\begin{exercise}
Let $f:\mathbb{R}^{3}\rightarrow\mathbb{R}^{2}$ be the map that sends
$(x,y,z)$ to $\left(  x-y,0\right)  $. Draw a picture that illustrates for
this particular map the concepts of kernel and level set, and how
complementary subspaces of the kernel are aliases of the image.
\end{exercise}

\smallskip\ If we specify domain $\mathcal{V}$ and kernel $\mathcal{K}$, we
have the level sets turned into the vector space $\mathcal{V}/\mathcal{K}$ as
an alias of the image of any map with that domain and kernel. Thus
$\mathcal{V}/\mathcal{K}$ can be viewed as a generic image for maps with
domain $\mathcal{V}$ and kernel $\mathcal{K}$. To round out the picture, we
now single out a map that sends $\mathcal{V}$ onto $\mathcal{V}/\mathcal{K}$.

\begin{proposition}
\label{NaturalProjection}The function $p:\mathcal{V}\rightarrow\mathcal{V}%
/\mathcal{K}$ which sends each $v\in\mathcal{V}$ to $v+\mathcal{K}$ is a map
with kernel $\mathcal{K}$.
\end{proposition}

Proof: Supposing, as we may, that $\mathcal{K}$ is the kernel of some map
$f:\mathcal{V}\rightarrow\mathcal{W}$, let the map $f^{\flat}:\mathcal{V}%
\rightarrow\operatorname*{Image}(f)$ be defined by $f^{\flat}\left(  v\right)
=f\left(  v\right)  $ for all $v\in\mathcal{V}$. Let $\Theta
:\operatorname*{Image}(f)\rightarrow\mathcal{V}/\mathcal{K}$ be the
isomorphism from the image of $f$ to $\mathcal{V}/\mathcal{K}$ obtained
through the correspondence of image elements with their level sets. Notice
that $\Theta$ sends $w=f\left(  v\right)  $ to $v+\mathcal{K}$. Then
$p=\Theta\circ f^{\flat}$ and is obviously a map with kernel $\mathcal{K}$.
$\blacksquare$

\smallskip\ We call $p$ the \textbf{natural projection}. It serves as a
generic map for obtaining the generic image. For any map $f:\mathcal{V}%
\rightarrow\mathcal{W}$ with kernel $\mathcal{K}$, we always have the composition%

\[
\mathcal{V}\rightarrow\mathcal{V}/\mathcal{K}\longleftrightarrow
\operatorname*{Image}(f)\hookrightarrow\mathcal{W}%
\]
where $\mathcal{V}\rightarrow\mathcal{V}/\mathcal{K}$ is generic, denoting the
natural projection $p$, and the remainder is the nongeneric specialization to
$f$ via an isomorphism and an inclusion map. Our next result details how the
generic map $p$ in fact serves in a more general way as a universal factor of
maps with kernel containing $\mathcal{K}$.

\begin{theorem}
\label{InducedMap}Let $p:\mathcal{V}\rightarrow\mathcal{V}/\mathcal{K}$ be the
natural projection and let $f:\mathcal{V}\rightarrow\mathcal{W}$. If
$\mathcal{K}\subset\operatorname*{Kernel}(f)$ then there is a unique
\textbf{\emph{induced map}} $f_{\mathcal{V}/\mathcal{K}}:$ $\mathcal{V}%
/\mathcal{K}\rightarrow\mathcal{W}$ such that $f_{\mathcal{V}/\mathcal{K}%
}\circ p=f$.
\end{theorem}

Proof: The prescription $f_{\mathcal{V}/\mathcal{K}}\circ p=f$ determines
exactly what $f_{\mathcal{V}/\mathcal{K}}$ must do: for each $v\in\mathcal{V}%
$, $f_{\mathcal{V}/\mathcal{K}}\left(  p(v)\right)  =f_{\mathcal{V}%
/\mathcal{K}}\left(  v+\mathcal{K}\right)  =f\left(  v\right)  $. For this to
unambiguously define the value of $f_{\mathcal{V}/\mathcal{K}}$ at each
element of $\mathcal{V}/\mathcal{K}$, it must be the case that if
$u+\mathcal{K}=v+\mathcal{K}$, then $f\left(  u\right)  =f\left(  v\right)  .$
But this is true because $\mathcal{K}\subset\operatorname*{Kernel}(f)$ implies
that $f\left(  k\right)  =0$ for each $k\in\mathcal{K}$. The unique function
$f_{\mathcal{V}/\mathcal{K}}$ so determined is readily shown to be a map.
$\blacksquare$

\begin{exercise}
Referring to the theorem above, $f_{\mathcal{V}/\mathcal{K}}$ is one-to-one if
and only if $\operatorname*{Kernel}(f)=\mathcal{K}$, and $f_{\mathcal{V}%
/\mathcal{K}}$ is onto if and only if $f$ is onto.
\end{exercise}

The rank of the natural map $p:\mathcal{V}\rightarrow\mathcal{V}/\mathcal{K} $
is the dimension of $\mathcal{V}/\mathcal{K}$, also known as the
\textbf{codimension} of the subspace $\mathcal{K}$ and denoted by
$\operatorname*{codim}\mathcal{K}$. Thus the Rank Plus Nullity Theorem may
also be expressed as
\[
\dim\mathcal{K}+\operatorname*{codim}\mathcal{K}=\dim\mathcal{V}%
\]
when $\dim\mathcal{V}$ is finite.

\begin{exercise}
Let $\mathcal{K}$ be a subspace of finite codimension in the
infinite-dimensional vector space $\mathcal{V}$. Then $\mathcal{K}%
\cong\mathcal{V}$.
\end{exercise}

\subsection{Short Exact Sequences}

There is a special type of map sequence that is a way to view a quotient. The
map sequence $\mathcal{X}\rightarrow\mathcal{Y}\rightarrow\mathcal{Z}$ is said
to be \textbf{exact} at $\mathcal{Y}$ if the image of the map going into
$\mathcal{Y}$ is exactly the kernel of the map coming out of $\mathcal{Y}$. A
\textbf{short exact sequence} is a map sequence of the form $0\rightarrow
\mathcal{K}\rightarrow\mathcal{V}\rightarrow\mathcal{W}\rightarrow0$ which is
exact at $\mathcal{K}$, $\mathcal{V}$, and $\mathcal{W}$. Notice that here the
map $\mathcal{K}\rightarrow\mathcal{V}$ is one-to-one, and the map
$\mathcal{V}\rightarrow\mathcal{W}$ is onto, in virtue of the exactness
requirements and the restricted nature of maps to or from $0=\left\{
0\right\}  $. Hence if $\mathcal{V}\rightarrow\mathcal{W}$ corresponds to the
function $f$, then it is easy to see that $\mathcal{K}$ is an alias of
$\operatorname*{Kernel}(f)$, and $\mathcal{W}$ is an alias of $\mathcal{V}%
/\operatorname*{Kernel}(f)$. Thus the short exact sequence captures the idea
of quotient modulo any subspace. It is true that with the same $\mathcal{K}$,
$\mathcal{V}$, and $\mathcal{W}$, there are maps such that the original
sequence with all its arrows reversed is also a short exact sequence, one
complementary to the original, it might be said. This is due to the fact that
$\mathcal{K}$ and $\mathcal{W}$ are always aliases of complementary subspaces
of $\mathcal{V}$, which is something special that vector space structure makes
possible. Thus it is said that, for vector spaces, a short exact sequence
always \textbf{splits}, as illustrated by the short exact sequence
$0\rightarrow\mathcal{K}\rightarrow\mathcal{K}\oplus\overline{\mathcal{K}%
}\rightarrow\overline{\mathcal{K}}\rightarrow0$.

An example of the use of the exact sequence view of quotients stems from
considering the following diagram where the rows and columns are short exact
sequences. Each second map in a row or column is assumed to be an inclusion,
and each third map is assumed to be a natural projection.
\[%
\begin{tabular}
[c]{ccccccccc}
&  &  &  & $0$ &  & $0$ &  & \\
&  &  &  & $\downarrow$ &  & $\downarrow$ &  & \\
$0$ & $\rightarrow$ & $\mathcal{K}$ & $\rightarrow$ & $\mathcal{H}$ &
$\rightarrow$ & $\mathcal{H}/\mathcal{K}$ & $\rightarrow$ & $0$\\
&  &  &  & $\downarrow$ &  & $\downarrow$ &  & \\
$0$ & $\rightarrow$ & $\mathcal{K}$ & $\rightarrow$ & $\mathcal{V}$ &
$\rightarrow$ & $\mathcal{V}/\mathcal{K}$ & $\rightarrow$ & $0$\\
&  &  &  & $\downarrow$ &  & $\downarrow$ &  & \\
&  &  &  & $\mathcal{V}/\mathcal{H}$ &  & $(\mathcal{V}/\mathcal{K}%
)/(\mathcal{H}/\mathcal{K})$ &  & \\
&  &  &  & $\downarrow$ &  & $\downarrow$ &  & \\
&  &  &  & $0$ &  & $0$ &  &
\end{tabular}
\]
The sequence $0\rightarrow\mathcal{H}\rightarrow\mathcal{V}\rightarrow
\mathcal{V}/\mathcal{K}\rightarrow(\mathcal{V}/\mathcal{K})/(\mathcal{H}%
/\mathcal{K})\rightarrow0$ is a subdiagram. If we compose the pair of onto
maps of $\mathcal{V}\rightarrow\mathcal{V}/\mathcal{K}\rightarrow
(\mathcal{V}/\mathcal{K})/(\mathcal{H}/\mathcal{K})$ to yield the composite
onto map $\mathcal{V}\rightarrow(\mathcal{V}/\mathcal{K})/(\mathcal{H}%
/\mathcal{K})$, we then have the sequence $0\rightarrow\mathcal{H}%
\rightarrow\mathcal{V}\rightarrow(\mathcal{V}/\mathcal{K})/(\mathcal{H}%
/\mathcal{K})\rightarrow0$ which is exact at $\mathcal{H}$ and at
$(\mathcal{V}/\mathcal{K})/(\mathcal{H}/\mathcal{K})$, and would also be exact
at $\mathcal{V}$ if $\mathcal{H}$ were the kernel of the composite. But it is
precisely the $h\in\mathcal{H}$ that map to the $h+\mathcal{K}\in
\mathcal{H}/\mathcal{K\subset V}/\mathcal{K}$ which then are precisely the
elements that map to $\mathcal{H}/\mathcal{K} $ which is the $0$ of
$(\mathcal{V}/\mathcal{K})/(\mathcal{H}/\mathcal{K})$. Thus we have obtained
the following isomorphism theorem.

\begin{theorem}
Let $\mathcal{K}$ be a subspace of $\mathcal{H}$, and let $\mathcal{H}$ be a
subspace of $\mathcal{V}$. Then $\mathcal{V}/\mathcal{H}$ is isomorphic to
$(\mathcal{V}/\mathcal{K})/(\mathcal{H}/\mathcal{K})$. $\blacksquare$
\end{theorem}

\begin{exercise}
Let $\mathcal{V}$ be $\mathbb{R}^{3}$, let $\mathcal{H}$ be the $(x,y)$-plane,
and let $\mathcal{K}$ be the $x$-axis. Interpret the theorem above for this case.
\end{exercise}

Now consider the diagram below where $\mathcal{X}$ and $\mathcal{Y}$ are
subspaces and again the rows and columns are short exact sequences and the
second maps are inclusions while the third maps are natural projections.
\[%
\begin{tabular}
[c]{ccccccccc}
&  & $0$ &  & $0$ &  &  &  & \\
&  & $\downarrow$ &  & $\downarrow$ &  &  &  & \\
$0$ & $\rightarrow$ & $\mathcal{%
X%
}\cap\mathcal{%
Y%
}$ & $\rightarrow$ & $\mathcal{%
Y%
}$ & $\rightarrow$ & $\mathcal{%
Y%
}/(\mathcal{%
X%
}\cap\mathcal{%
Y%
})$ & $\rightarrow$ & $0$\\
&  & $\downarrow$ &  & $\downarrow$ &  &  &  & \\
$0$ & $\rightarrow$ & $\mathcal{%
X%
}$ & $\rightarrow$ & $\mathcal{%
X%
}+\mathcal{%
Y%
}$ & $\rightarrow$ & $(\mathcal{%
X%
}+\mathcal{%
Y%
})/\mathcal{%
X%
}$ & $\rightarrow$ & $0$\\
&  & $\downarrow$ &  & $\downarrow$ &  &  &  & \\
&  & $\mathcal{%
X%
}/(\mathcal{%
X%
}\cap\mathcal{%
Y%
}) $ &  & $(\mathcal{%
X%
}+\mathcal{%
Y%
})/\mathcal{%
Y%
}$ &  &  &  & \\
&  & $\downarrow$ &  & $\downarrow$ &  &  &  & \\
&  & $0$ &  & $0$ &  &  &  &
\end{tabular}
\]
As a subdiagram we have the sequence $0\rightarrow\mathcal{%
X%
}\cap\mathcal{%
Y%
}\rightarrow\mathcal{%
X%
}\rightarrow\mathcal{%
X%
}+\mathcal{%
Y%
}\rightarrow(\mathcal{%
X%
}+\mathcal{%
Y%
})/\mathcal{%
Y%
}\rightarrow0$. Replacing the sequence $\mathcal{%
X%
}\rightarrow\mathcal{%
X%
}+\mathcal{%
Y%
}\rightarrow(\mathcal{%
X%
}+\mathcal{%
Y%
})/\mathcal{%
Y%
}$ with its composite, the result would be a short exact sequence if the
composite $\mathcal{%
X%
}\rightarrow(\mathcal{%
X%
}+\mathcal{%
Y%
})/\mathcal{%
Y%
}$ were an onto map with kernel $\mathcal{%
X%
}\cap\mathcal{%
Y%
}$.

To see if the composite is onto, we must see if each element of $(\mathcal{%
X%
}+\mathcal{%
Y%
})/\mathcal{%
Y%
}$ is the image of some element of $\mathcal{%
X%
}$. Now each element of $(\mathcal{%
X%
}+\mathcal{%
Y%
})/\mathcal{%
Y%
}$ is the image of some element $w\in\mathcal{%
X%
}+\mathcal{%
Y%
}$ where $w=u+v$ for some $u\in\mathcal{%
X%
}$ and some $v\in\mathcal{%
Y%
}$. But clearly the element $u$ by itself has the same image modulo $\mathcal{%
Y%
}$. Hence the composite is onto. Also, the elements of $\mathcal{%
X%
}$ that map to the $0$ of $(\mathcal{%
X%
}+\mathcal{%
Y%
})/\mathcal{%
Y%
}$, namely $\mathcal{%
Y%
}$, are precisely those that then are also in $\mathcal{%
Y%
}$. Hence the kernel of the composite is $\mathcal{%
X%
}\cap\mathcal{%
Y%
}$.

Thus we have obtained another isomorphism theorem and an immediate corollary.

\begin{theorem}
Let $\mathcal{%
X%
}$ and $\mathcal{%
Y%
}$ be any subspaces of some vector space. Then $\mathcal{%
X%
}/\left(  \mathcal{X}\cap\mathcal{Y}\right)  $ is isomorphic to $(\mathcal{%
X%
}+\mathcal{%
Y%
})/\mathcal{%
Y%
}$. $\blacksquare$
\end{theorem}

\begin{corollary}
[Grassmann's Relation]\label{GrassmannRel}Let $\mathcal{%
X%
}$ and $\mathcal{%
Y%
}$ be finite-dimensional subspaces of a vector space. Then
\[
\dim\mathcal{%
X%
}+\dim\mathcal{%
Y%
}=\dim(\mathcal{%
X%
}+\mathcal{%
Y%
})+\dim(\mathcal{%
X%
}\cap\mathcal{%
Y%
})\text{.}:EndProof
\]
\end{corollary}

\begin{exercise}
Two \emph{2}-dimensional subspaces of $\mathbb{R}^{3}$ are either identical or
intersect in a \emph{1}-dimensional subspace.
\end{exercise}

\begin{exercise}
Let $\mathcal{N}$, $\mathcal{N}^{\,^{\prime}}$, $\mathcal{P}$, and
$\mathcal{P}^{\,^{\prime}}$ be subspaces of the same vector space and let
$\mathcal{N}\subset\mathcal{P}$ and $\mathcal{N}^{\,^{\prime}}\subset
\mathcal{P}^{\,^{\prime}}$. Set $\mathcal{%
X%
}=\mathcal{P}\cap\mathcal{P}^{\,^{\prime}}$ and $\mathcal{Y}=(\mathcal{P}%
\cap\mathcal{N}^{\,^{\prime}})+\mathcal{N}$. Verify that $\mathcal{X}%
\cap\mathcal{Y}=(\mathcal{P}^{\,^{\prime}}\cap\mathcal{N})+(\mathcal{P}%
\cap\mathcal{N}^{\,^{\prime}})$ and $\mathcal{X}+\mathcal{Y}=\left(
\mathcal{P}\cap\mathcal{P}^{\,^{\prime}}\right)  +\mathcal{N}$, so that
\[
\frac{\left(  \mathcal{P}^{\,^{\prime}}\cap\mathcal{P}\right)  +\mathcal{N}%
^{\,^{\prime}}}{(\mathcal{P}^{\,^{\prime}}\cap\mathcal{N})+\mathcal{N}%
^{\,^{\prime}}}\cong\frac{\mathcal{P}\cap\mathcal{P}^{\,^{\prime}}%
}{(\mathcal{P}^{\,^{\prime}}\cap\mathcal{N})+(\mathcal{P}\cap\mathcal{N}%
^{\,^{\prime}})}\cong\frac{\left(  \mathcal{P}\cap\mathcal{P}^{\,^{\prime}%
}\right)  +\mathcal{N}}{(\mathcal{P}\cap\mathcal{N}^{\,^{\prime}}%
)+\mathcal{N}}\;\;\text{\textbf{.}}%
\]
\end{exercise}

\subsection{\label{Projections}Projections and Reflections on $\mathcal{V}%
=\mathcal{W}\oplus\mathcal{X}$}

Relative to the decomposition of the vector space $\mathcal{V}$ by a pair
$\mathcal{W},\mathcal{X}$ of direct summands, we may define some noteworthy
self-maps on $\mathcal{V}$. Expressing the general vector $v\in\mathcal{V}%
=\mathcal{W}\oplus\mathcal{X}$ uniquely as $v=w+x$ where $w\in\mathcal{W}$ and
$x\in\mathcal{X}$, two types of such are defined by
\[
P_{\mathcal{%
W%
\mid%
X%
}}\left(  v\right)  =w\text{ \ and\ }R_{\mathcal{%
W%
\mid%
X%
}}\left(  v\right)  =w-x\text{.}%
\]
That these really are maps is easy to establish. We call $P_{\mathcal{%
W%
\mid%
X%
}}$ the \textbf{projection onto }$\mathcal{W}$ \textbf{along }$\mathcal{X}$
and we call $R_{\mathcal{%
W%
\mid%
X%
}}$ the \textbf{reflection in }$\mathcal{W}$ \textbf{along }$\mathcal{X}$.
(The function $\phi$ in the proof of Proposition \ref{Complementary} was the
projection onto $\mathcal{\overline{%
K%
}}$ along $\mathcal{K}$, so we have already employed the projection type to
advantage.) Denoting the identity map on $\mathcal{V}$ by $I$, we have
\[
P_{\mathcal{X}\mid\mathcal{W}}=I-P_{\mathcal{%
W%
\mid%
X%
}}%
\]
and
\[
R_{\mathcal{%
W%
\mid%
X%
}}=P_{\mathcal{%
W%
\mid%
X%
}}-P_{\mathcal{%
X%
\mid%
W%
}}=I-2P_{\mathcal{%
X%
\mid%
W%
}}=2P_{\mathcal{%
W%
\mid%
X%
}}-I\text{.}%
\]
It bears mention that a given $\mathcal{W}$ generally has many different
complements, and if $\mathcal{X}$ and $\mathcal{Y}$ are two such,
$P_{\mathcal{%
W%
\mid%
X%
}}$ and $P_{\mathcal{W}\mid\mathcal{Y}}$ will generally differ.

The image of $P_{\mathcal{%
W%
\mid%
X%
}}$ is $\mathcal{W}$ and its kernel is $\mathcal{X}$. Thus $P_{\mathcal{%
W%
\mid%
X%
}}$ is a self-map with the special property that its image and its kernel are
complements. The image and kernel of $R_{\mathcal{%
W%
\mid%
X%
}}$ are also complements, but trivially, as the kernel of $R_{\mathcal{%
W%
\mid%
X%
}}$ is $0=\left\{  0\right\}  $. The kernel of $R_{\mathcal{%
W%
\mid%
X%
}}$ is $0$ because $w$ and $x$ are equal only if they are both $0$ since
$\mathcal{W}\cap\mathcal{X}=0$ from the definition of direct sum.

A double application of $P_{\mathcal{%
W%
\mid%
X%
}}$ has the same effect as a single application: $P_{\mathcal{%
W%
\mid%
X%
}}^{2}=P_{\mathcal{%
W%
\mid%
X%
}}\circ P_{\mathcal{%
W%
\mid%
X%
}}=P_{\mathcal{%
W%
\mid%
X%
}}$ ($P_{\mathcal{%
W%
\mid%
X%
}}$ is \emph{idempotent}). We then have $P_{\mathcal{%
W%
\mid%
X%
}}\circ P_{\mathcal{%
X%
\mid%
W%
}}=0$ (the constant map on $\mathcal{V}$ that sends everything to $0$). We
also find that $P_{\mathcal{%
W%
\mid%
X%
}}^{2}=P_{\mathcal{%
W%
\mid%
X%
}}$ implies that $R_{\mathcal{%
W%
\mid%
X%
}}^{2}=I$ ($R_{\mathcal{%
W%
\mid%
X%
}}$ is \emph{involutary}). An idempotent map always turns out to be some
$P_{\mathcal{%
W%
\mid%
X%
}}$, and an involutary map is \emph{almost} always some $R_{\mathcal{%
W%
\mid%
X%
}}$, as we now explain.

\begin{proposition}
Let $P:\mathcal{V}\rightarrow\mathcal{V}$ be a map such that $P^{2}=P.$ Then
$P$ is the projection onto $\operatorname*{Image}\left(  P\right)  $ along
$\operatorname*{Kernel}\left(  P\right)  $.
\end{proposition}

Proof: Let $w$ be in $\operatorname*{Image}\left(  P\right)  $, say
$w=P\left(  v\right)  $. Then
\[
w=P^{2}\left(  v\right)  =P\left(  P\left(  v\right)  \right)  =P\left(
w\right)
\]
so that if $w$ is also in $\operatorname*{Kernel}\left(  P\right)  $, it must
be that $w=0$. Hence
\[
\operatorname*{Image}\left(  P\right)  \cap\operatorname*{Kernel}\left(
P\right)  =0.
\]

Let $Q=I-P$ so that $P\left(  v\right)  +Q\left(  v\right)  =I\left(
v\right)  =v$ and thus $\operatorname*{Image}\left(  P\right)
+\operatorname*{Image}\left(  Q\right)  =\mathcal{V}.$ But
$\operatorname*{Image}\left(  Q\right)  =\operatorname*{Kernel}\left(
P\right)  $ since if $x=Q\left(  v\right)  =I\left(  v\right)  -P\left(
v\right)  $, then $P\left(  x\right)  =P\left(  v\right)  -P^{2}\left(
v\right)  =0$. Therefore
\[
\mathcal{V}=\operatorname*{Image}\left(  P\right)  \oplus
\operatorname*{Kernel}\left(  P\right)  .
\]
Now if $v=w+x$, $w\in\operatorname*{Image}\left(  P\right)  $, $x\in
\operatorname*{Kernel}\left(  P\right)  $, then $P\left(  v\right)  =P\left(
w+x\right)  =P\left(  w\right)  +P\left(  x\right)  =w+0=w$, and $P$ is indeed
the projection onto $\operatorname*{Image}\left(  P\right)  $ along
$\operatorname*{Kernel}\left(  P\right)  $. $\blacksquare$\smallskip\ 

When $\mathcal{V}$ is over a field in which $1+1=0$, the only possibility for
$R_{\mathcal{%
W%
\mid%
X%
}}$ is the identity $I$. However, when $\mathcal{V}$ is over such a field
there can be an involutary map on $\mathcal{V}$ which is not the identity, for
example the map of the following exercise.

\begin{exercise}
Over the $2$-element field $\mathcal{F}=\left\{  0,1\right\}  $ \emph{(}%
wherein $1+1=0$\emph{)}, let $\mathcal{V}$ be the vector space $\mathcal{F}%
^{2}$ of $2$-tuples of elements of $\mathcal{F}$. Let $F$ be the map on
$\mathcal{V}$ that interchanges $\left(  0,1\right)  $ and $\left(
1,0\right)  $. Then $F\circ F$ is the identity map.
\end{exercise}

The $F$ of the exercise can reasonably be called a reflection, but it is not
some $R_{\mathcal{%
W%
\mid%
X%
}}$. If we say that any involutary map is a reflection, then reflections of
the form $R_{\mathcal{%
W%
\mid%
X%
}}$ do not quite cover all the possibilities. However, when $1+1\neq0$, they
do. (For any function $f:\mathcal{V}\rightarrow\mathcal{V}$, we denote by
$\operatorname*{Fix}\left(  f\right)  $ the set of all elements in
$\mathcal{V}$ that are fixed by $f$, i. e., $\operatorname*{Fix}\left(
f\right)  =\left\{  v\in\mathcal{V\ }|\mathcal{\ }f\left(  v\right)
=v\right\}  $.)

\begin{exercise}
Let $\mathcal{V}$ be over a field such that $1+1\neq0$. Let $R:\mathcal{V}%
\rightarrow\mathcal{V}$ be a map such that $R^{2}=I.$ Then $R=2P-I$, where
$P=\frac{1}{2}\left(  R+I\right)  $ is the projection onto
$\operatorname*{Fix}\left(  R\right)  $ along $\operatorname*{Fix}\left(
-R\right)  $.
\end{exercise}

\subsection{Problems}

\ 

1. Let $f:S\rightarrow T$ be a function. Then for all subsets $A,B$ of $S$,
\[
f\left(  A\cap B\right)  \subset f\left(  A\right)  \cap f\left(  B\right)
\text{,}%
\]
but
\[
f\left(  A\cap B\right)  =f\left(  A\right)  \cap f\left(  B\right)
\]
for all subsets $A,B$ of $S$ if and only if $f$ is one-to-one. What about for
$\cup$ instead of $\cap$? And what about for $\cap$ and for $\cup$ when $f$ is
replaced by $f^{-1}$?

\medskip\ 

2. \textbf{(}Enhancing Theorem \ref{BasisFreedom} for nontrivial
codomains\textbf{)} Let $\mathcal{B}$ be a subset of the vector space
$\mathcal{V}$ and let $\mathcal{W}\neq\left\{  0\right\}  $ be a vector space
over the same field. Then $\mathcal{B}$ is independent if and only if every
function $f_{0}:\mathcal{B}\rightarrow\mathcal{W}$ extends to at least one map
$f:\mathcal{V}\rightarrow\mathcal{W}$, and $\mathcal{B}$ spans $\mathcal{V}$
if and only if every function $f_{0}:\mathcal{B}\rightarrow\mathcal{W}$
extends to at most one map $f:\mathcal{V}\rightarrow\mathcal{W}$. Hence
$\mathcal{B}$ is a basis for $\mathcal{V}$ if and only if every function
$f_{0}:\mathcal{B}\rightarrow\mathcal{W}$ extends uniquely to a map
$f:\mathcal{V}\rightarrow\mathcal{W}$.

\medskip\ 

3. Let $\mathcal{B}$ be a basis for the vector space $\mathcal{V}$ and let
$f:\mathcal{V}\rightarrow\mathcal{W}$ be a map. Then $f$ is one-to-one if and
only if $f\left(  \mathcal{B}\right)  $ is an independent set, and $f$ is onto
if and only if $f\left(  \mathcal{B}\right)  $ spans $\mathcal{W}$.

\medskip\ 

4. Deduce Theorem \ref{Injmap} as a corollary of the Rank Plus Nullity Theorem.

\medskip\ 

5. Let $\mathcal{V}$ be a finite-dimensional vector space of dimension $n$
over the finite field of $q$ elements. Then the number of vectors in
$\mathcal{V}$ is $q^{n}$, the number of nonzero vectors that are not in a
given $m$-dimensional subspace is
\[
\left(  q^{n}-1\right)  -\left(  q^{m}-1\right)  =\left(  q^{n}-q^{m}\right)
\text{,}%
\]
and the number of different basis sets that can be extracted from
$\mathcal{V}$ is
\[
\frac{\left(  q^{n}-1\right)  \left(  q^{n}-q\right)  \cdots\left(
q^{n}-q^{n-1}\right)  }{n!}%
\]
since there are
\[
\left(  q^{n}-1\right)  \left(  q^{n}-q\right)  \cdots\left(  q^{n}%
-q^{n-1}\right)
\]
different sequences of $n$ distinct basis vectors that can be drawn from
$\mathcal{V}$.

\medskip\ 

6. (Correspondence Theorem) Let $\mathcal{U}\lhd\mathcal{V}$. The natural
projection $p:\mathcal{V}\rightarrow\mathcal{V}/\mathcal{U}$ puts the
subspaces of $\mathcal{V}$ that contain $\mathcal{U}$ in one-to-one
correspondence with the subspaces of $\mathcal{V}/\mathcal{U}$, and every
subspace of $\mathcal{V}/\mathcal{U}$ is of the form $\mathcal{W}/\mathcal{U}
$ for some subspace $\mathcal{W}$ that contains $\mathcal{U}$.

\medskip\ 

7. $P_{\mathcal{W}_{1}\mathcal{\mid%
X%
}_{1}}+P_{\mathcal{%
W%
}_{2}\mathcal{\mid%
X%
}_{2}}=P_{\mathcal{%
W%
\mid%
X%
}}$ if and only if $P_{\mathcal{W}_{1}\mathcal{\mid%
X%
}_{1}}\circ P_{\mathcal{%
W%
}_{2}\mathcal{\mid%
X%
}_{2}}=P_{\mathcal{%
W%
}_{2}\mathcal{\mid%
X%
}_{2}}\circ P_{\mathcal{W}_{1}\mathcal{\mid%
X%
}_{1}}=0$, in which case $\mathcal{W}=\mathcal{W}_{1}\oplus\mathcal{W}_{2}$
and $\mathcal{X}=\mathcal{X}_{1}\cap\mathcal{X}_{2}$. What happens when
$P_{\mathcal{W}_{1}\mathcal{\mid%
X%
}_{1}}-P_{\mathcal{%
W%
}_{2}\mathcal{\mid%
X%
}_{2}}=P_{\mathcal{%
W%
\mid%
X%
}}$ is considered instead?

\newpage

\section{More on Maps and Structure}

\subsection{Function Spaces and Map Spaces}

The\textbf{\ }set $\mathcal{W}^{\mathcal{D}}$ of all functions from a set
$\mathcal{D}$ into the vector space $\mathcal{W}$ may be made into a vector
space by defining the operations in a pointwise fashion. That is, suppose that
$f$ and $g$ are in $\mathcal{W}^{\mathcal{D}}$, and define $f+g$ by
$(f+g)\left(  x\right)  =f\left(  x\right)  +g\left(  x\right)  $, and for any
scalar $a$ define $a\cdot f$ by $(a\cdot f)\left(  x\right)  =a\cdot\left(
f\left(  x\right)  \right)  $. This makes $\mathcal{W}^{\mathcal{D}}$ into a
vector space, and it is this vector space that we mean whenever we refer to
$\mathcal{W}^{\mathcal{D}}$ as a \textbf{function space}. Considering
$n$-tuples as functions with domain $\left\{  1,\ldots,n\right\}  $ and the
field $\mathcal{F}$ as a vector space over itself, we see that the familiar
$\mathcal{F}^{n}$ is an example of a function space.

Subspaces of a function space give us more examples of vector spaces. For
instance, consider the set $\left\{  \mathcal{V}\rightarrow\mathcal{W}%
\right\}  $ of all maps from the vector space $\mathcal{V}$ into the vector
space $\mathcal{W}$. It is easy to verify that it is a subspace of the
function space $\mathcal{W}^{\mathcal{V}}$. It is this subspace that we mean
when we refer to $\left\{  \mathcal{V}\rightarrow\mathcal{W}\right\}  $ as a
\textbf{map space}.

\subsection{Dual of a Vector Space, Dual of a Basis}

By a \textbf{linear functional} on a vector space $\mathcal{V}$ over a field
$\mathcal{F}$ is meant an element of the map space $\left\{  \mathcal{V}%
\rightarrow\mathcal{F}^{1}\right\}  $. (Since $\mathcal{F}^{1}$ is just
$\mathcal{F}$ viewed as a vector space over itself, a linear functional is
also commonly described as a map into $\mathcal{F}$.) For a vector space
$\mathcal{V}$ over a field $\mathcal{F}$, the map space $\left\{
\mathcal{V}\rightarrow\mathcal{F}^{1}\right\}  $ will be called the
\textbf{dual} of $\mathcal{V}$ and will be denoted by $\mathcal{V}^{\top}$.

Let the vector space $\mathcal{V}$ have the basis $\mathcal{B}$. For each
$x\in\mathcal{B}$, let $x^{\top}$ be the linear functional on $\mathcal{V}$
that sends $x$ to $1$ and sends all other elements of $\mathcal{B}$ to $0$. We
call $x^{\top}$ the \textbf{coordinate function} corresponding to the basis
vector $x$. The set $\mathcal{B}^{\,\#}$ of all such coordinate functions is
an independent set. For if for some finite $\mathcal{X}\subset\mathcal{B}$ we
have $\sum_{x\in\mathcal{X}}a_{x}\cdot x^{\top}\left(  v\right)  =0$ for all
$v\in\mathcal{V}$, then for any $v\in\mathcal{X}$ we find that $a_{v}=0$ so
there can be no dependency in $\mathcal{B}^{\,\#}$.

When the elements of $\mathcal{B}^{\,\#}$ span $\mathcal{V}^{\top}$, we say
that $\mathcal{B}$ has $\mathcal{B}^{\#}$ as its \textbf{dual} and we may then
replace $\mathcal{B}^{\#}$ with the more descriptive notation $\mathcal{B}%
^{\top}$. When $\mathcal{V}$ is finite-dimensional, $\mathcal{B}^{\,\#}$ does
span $\mathcal{V}^{\top}$, we claim. For the general linear functional $f$ is
defined by prescribing its value on each element of $\mathcal{B}$, and given
any $v\in\mathcal{B}$ we have $f\left(  v\right)  =\sum_{x\in\mathcal{B}%
}f\left(  x\right)  \cdot x^{\top}\left(  v\right)  $. Therefore any linear
functional on $\mathcal{V}$ may be expressed as the linear combination
$f=\sum_{x\in\mathcal{B}}f\left(  x\right)  \cdot x^{\top}$, so that
$\mathcal{B}^{\,\#}$ spans $\mathcal{V}^{\top}$ as we claimed.

We thus have the following result.

\begin{theorem}
Let $\mathcal{V}$ be a finite-dimensional vector space. Then any basis of
$\mathcal{V}$ has a dual, and therefore $\mathcal{V}^{\top}$ has the same
dimension as $\mathcal{V}$. $\blacksquare$
\end{theorem}

When $\mathcal{V}$ is infinite-dimensional, $\mathcal{B}^{\,\#}$ need not span
$\mathcal{V}^{\top}$.

\begin{example}
\label{DualofUltZeroSeqs}Let $\mathbb{N}$ denote the natural numbers $\left\{
0,1,\ldots\right\}  $. Let $\mathcal{V}$ be the subspace of the function space
$\mathbb{R}^{\mathbb{N}}$ consisting of only those real sequences $a_{0}%
,a_{1},\ldots$ that are ultimately zero, i. e., for which there exists a
smallest $N\in\mathbb{N}$ \emph{(}called the \emph{length} of the ultimately
zero sequence\emph{)}, such that $a_{n}=0$ for all $n\geqslant N$. Then
$\phi:\mathcal{V}\rightarrow\mathbb{R}$ defined by $\phi(a_{0},a_{1}%
,\ldots)=\sum_{n\in\mathbb{N}}\alpha_{n}a_{n}$ is a linear functional given
\emph{any} \emph{(}not necessarily ultimately zero\emph{)} sequence of real
coefficients $\alpha_{0},\alpha_{1},\ldots$. \emph{(}Given $\alpha_{0}%
,\alpha_{1},\ldots$ that is not ultimately zero, it cannot be replaced by an
ultimately zero sequence that produces the same functional, since there would
always be sequences $a_{0},a_{1},\ldots$ of greater length for which the
results would differ.\emph{)} The $\mathcal{B}^{\,\#}$ derived from the basis
$\mathcal{B}=\left\{  \left(  1,0,0,\ldots\right)  ,\left(  0,1,0,0,\ldots
\right)  ,\ldots\right\}  $ for $\mathcal{V}$ does not span $\mathcal{V}%
^{\top}$ because $\left\langle \mathcal{B}^{\,\#}\right\rangle $ only contains
elements corresponding to sequences $\alpha_{0},\alpha_{1},\ldots$ that are
ultimately zero. Also note that by Theorem \ref{BasisFreedom}, each element of
$\mathcal{V}^{\top}$ is determined by an unrestricted choice of value at each
and every element of $\mathcal{B}$.
\end{example}

Fixing an element $v\in\mathcal{V}$, the scalar-valued function $F_{v}$
defined on $\mathcal{V}^{\top}$ by $F_{v}(f)=f\left(  v\right)  $ is a linear
functional on $\mathcal{V}^{\top}$, i. e., $F_{v}$ is an element of
$\mathcal{V}^{\top\top}=\left(  \mathcal{V}^{\top}\right)  ^{\top}$. The
association of each such $v\in\mathcal{V}$ with the corresponding $F_{v}%
\in\mathcal{V}^{\top\top}$ gives a map $\Theta$ from $\mathcal{V}$ into
$\mathcal{V}^{\top\top}$. $\Theta$ is one-to-one, we claim. For suppose that
$v$ and $w$ are in $\mathcal{V}$ and that $F_{v}=F_{w}$ so that for all
$f\in\mathcal{V}^{\top}$, $f\left(  v\right)  =f\left(  w\right)  $. Then
$f(v-w)=0$ for every $f\in\mathcal{V}^{\top}$. Let $\mathcal{B}$ be a basis
for $\mathcal{V}$, and for each $x\in\mathcal{B}$ let $x^{\top}$ be the
corresponding coordinate function. Let $\mathcal{X}$ be a finite subset of
$\mathcal{B}$ such that $v-w=\sum_{x\in\mathcal{X}}a_{x}\cdot x$. Applying
$x^{\top}$ to $v-w$ for each $x\in\mathcal{X}$, we conclude that there is no
nonzero $a_{x}$. Hence $v=w$, and $\Theta$ is one-to-one as claimed.

Thus it has been shown that each vector $v\in\mathcal{V}$ defines a unique
element $F_{v}\in\mathcal{V}^{\top\top}$. But the question still remains as to
whether $\Theta$ is onto, that is whether every $F\in\mathcal{V}^{\top\top}$
is an $F_{v}$ for some $v\in\mathcal{V}$. We claim that $\Theta$ is onto if
$\mathcal{V}$ is finite-dimensional. This is because if we have bases
$\mathcal{B}$ and $\mathcal{B}^{\top}$ as above and we form the basis
$\mathcal{B}^{\top\top}=\left\{  x_{1}^{\top\top},\ldots,x_{n}^{\top\top
}\right\}  $ dual to $\mathcal{B}^{\top}$, we readily verify that $x_{i}%
^{\top\top}=F_{x_{i}}$ for each $i$. Thus for $F=a_{1}\cdot x_{1}^{\top\top
}+\cdots+a_{n}\cdot x_{n}^{\top\top}$, and $v=a_{1}\cdot x_{1}+\cdots
+a_{n}\cdot x_{n}$, we have $F=F_{v}$.

We have proved the following.

\begin{theorem}
The map $\Theta:\mathcal{V\rightarrow V}^{\top\top}$ which sends
$v\in\mathcal{V}$ to the linear functional $F_{v}$ on $\mathcal{V}^{\top}$
defined by $F_{v}\left(  f\right)  =f\left(  v\right)  $ is always one-to-one,
and moreover is onto if $\mathcal{V}$ is finite-dimensional. $\blacksquare$
\end{theorem}

$\Theta$ does not depend on any choice of basis, so $\mathcal{V}$ and
$\Theta\left(  \mathcal{V}\right)  \subset\mathcal{V}^{\top\top}$ are
isomorphic in a ``natural'' basis-invariant manner. We say that $\Theta$
\textbf{naturally embeds} $\mathcal{V}$ in $\mathcal{V}^{\top\top}$ and we
call $\Theta$ the \textbf{natural injection} of $\mathcal{V}$ into
$\mathcal{V}^{\top\top}$, or, when it is onto, the \textbf{natural
correspondence }between $\mathcal{V}$ and $\mathcal{V}^{\top\top}$.
Independent of basis choices, a vector in $\mathcal{V}$ can be considered to
act directly on $\mathcal{V}^{\top}$ as a linear functional. If one writes
$f\left(  v\right)  $, where $f\in\mathcal{V}^{\top}$ and $v\in\mathcal{V}$,
either of $f$ or $v$ can be the variable. When $f$ is the variable, then $v$
is really playing the r\^{o}le of a linear functional on $\mathcal{V}^{\top}
$, i. e., an element of $\mathcal{V}^{\top\top}$. One often sees notation such
as $\left\langle f,v\right\rangle $ used for $f\left(  v\right)  $ in order to
more clearly indicate that $f$ and $v$ are to viewed as being on similar
footing, with either being a functional of the other. Thus, as convenience
dictates, we may \emph{identify} $\mathcal{V}$ with $\Theta\left(
\mathcal{V}\right)  \subset\mathcal{V}^{\top\top}$, and in particular, we
agree always to do so whenever $\Theta$ is onto.

\subsection{Annihilators\label{Annihilator}}

Let $\mathcal{S}$ be a subset of the vector space $\mathcal{V}$. By the
\textbf{annihilator} $\mathcal{S}^{0}$ of $\mathcal{S}$ is meant the subset of
$\mathcal{V}^{\top}$ which contains all the linear functionals $f$ such that
$f\left(  s\right)  =0$ for all $s\in\mathcal{S}$. It is easy to see that
$\mathcal{S}^{0}$ is actually a subspace of $\mathcal{V}^{\top}$ and that
$\mathcal{S}^{0}=\left\langle \mathcal{S}\right\rangle ^{0}$. Obviously,
$\left\{  0\right\}  ^{0}=\mathcal{V}^{\top}$. Also, since a linear functional
that is not the zero functional, i. e., which is not identically zero, must
have at least one nonzero value, $\mathcal{V}^{0}=\left\{  0\right\}  $.

\begin{exercise}
Because of the natural injection of $\mathcal{V}$ into $\mathcal{V}^{\top\top
}$, the only $x\in\mathcal{V}$ such that $f\left(  x\right)  =0$ for all
$f\in\mathcal{V}^{\top}$ is $x=0$.
\end{exercise}

When $\mathcal{V}$ is finite-dimensional, the dimension of a subspace
determines that of its annihilator.

\begin{theorem}
\label{AnnihilatorDimension}Let $\mathcal{V}$ be an $n$-dimensional vector
space and let $\mathcal{U}\vartriangleleft\mathcal{V}$ have dimension $m$.
Then $\mathcal{U}^{0}$ has dimension $n-m$.
\end{theorem}

Proof: Let $\mathcal{B}=\left\{  x_{1},\ldots,x_{n}\right\}  $ be a basis for
$\mathcal{V}$ such that $\left\{  x_{1},\ldots,x_{m}\right\}  $ is a basis for
$\mathcal{U}$. We claim that $\left\{  x_{m+1}^{\top},\ldots,x_{n}^{\top
}\right\}  \subset\mathcal{B}^{\top}$ spans $\mathcal{U}^{0}$. Clearly,
$\left\langle \left\{  x_{m+1}^{\top},\ldots,x_{n}^{\top}\right\}
\right\rangle \subset\mathcal{U}^{0}$. On the other hand, if $f=b_{1}\cdot
x_{1}^{\top}+\cdots+b_{n}\cdot x_{n}^{\top}\in\mathcal{U}^{0}$ then $f\left(
x\right)  =0$ for each $x\in\left\{  x_{1},\ldots,x_{m}\right\}  $ so that
$b_{1}=\cdots=b_{m}=0$ and therefore $\mathcal{U}^{0}\subset\left\langle
\left\{  x_{m+1}^{\top},\ldots,x_{n}^{\top}\right\}  \right\rangle $.
$\blacksquare$

\begin{exercise}
Let $\mathcal{V}$ be a finite-dimensional vector space. For any $\mathcal{U}%
\vartriangleleft\mathcal{V}$, $\left(  \mathcal{U}^{0}\right)  ^{0}%
=\mathcal{U}$, remembering that we identify $\mathcal{V}$ with $\mathcal{V}%
^{\top\top}$.
\end{exercise}

\subsection{Dual of a Map}

Let $f:\mathcal{V}\rightarrow\mathcal{W}$ be a map. Then for $\varphi
\in\mathcal{W}^{\top}$, the prescription $f^{\top}\left(  \varphi\right)
=\varphi\circ f$, or $(f^{\top}\left(  \varphi\right)  )\left(  v\right)
=\varphi(f\left(  v\right)  )$ for all $v\in\mathcal{V}$, defines a function
$f^{\top}:\mathcal{W}^{\top}\rightarrow\mathcal{V}^{\top}$. It is easy to show
that $f^{\top}$ is a map. $f^{\top}$ is called the \textbf{dual} of $f$. The
dual of an identity map is an identity map. The dual of the composite map
$f\circ g$ is $g^{\top}\circ f^{\top}$. Thus, if $f$ is a bijective map, the
dual of its inverse is the inverse of its dual.

\begin{exercise}
\label{CoordinatesofDualMap}Let $\mathcal{V}$ and $\mathcal{W}$ be
finite-dimensional vector spaces with respective bases $\mathcal{B}$ and
$\mathcal{C}$. Let the map $f:\mathcal{V}\rightarrow\mathcal{W}$ satisfy
$f\left(  x\right)  =\sum_{y\in\mathcal{C}}a_{y,x}\cdot y$ for each
$x\in\mathcal{B}$. Then $f^{\top}$ satisfies $f^{\top}\left(  y^{\top}\right)
=\sum_{x\in\mathcal{B}}a_{y,x}\cdot x^{\top}$ for each $y\in\mathcal{C}$
\emph{(}where $x^{\top}\in\mathcal{B}^{\top}$ and $y^{\top}\in\mathcal{C}%
^{\top}$ are coordinate functions corresponding respectively to $x\in
\mathcal{B}$ and $y\in\mathcal{C}$, of course.\emph{)}
\end{exercise}

The kernel of $f^{\top}$ is found to have a special property.

\begin{theorem}
\label{KernelofDualvsImage}The kernel of $f^{\top}$ is the annihilator of the
image of the map $f:\mathcal{V}\rightarrow\mathcal{W}$.
\end{theorem}

Proof: $\varphi$ is in the kernel of $f^{\top}\Leftrightarrow(f^{\top}\left(
\varphi\right)  )\left(  v\right)  =0$ for all $v\in\mathcal{V}\Leftrightarrow
\varphi\left(  f\left(  v\right)  \right)  =0$ for all $v\in\mathcal{V}%
\Leftrightarrow\varphi\in\left(  f\left(  \mathcal{V}\right)  \right)  ^{0}$.
$\blacksquare$

\medskip\ 

When the codomain of $f$ has finite dimension $m$, say, then the image of $f$
has finite dimension $r$, say. Hence, the kernel of $f^{\top}$ (which is the
annihilator of the image of $f$) must have dimension $m-r$. By the Rank Plus
Nullity Theorem, the rank of $f^{\top}$ must then be $r$. Thus we have the
following result.

\begin{theorem}
[Rank Theorem]\emph{\ }When the map $f$ has a finite-dimensional codomain, the
rank of $f^{\top}$ equals the rank of $f$. $\blacksquare$
\end{theorem}

\begin{exercise}
Let $f$ be a map from one finite-dimensional vector space to another. Then,
making the standard identifications, $\left(  f^{\top}\right)  ^{\top}=f$ and
the kernel of $f$ is the annihilator of the image of $f^{\top}$.
\end{exercise}

\begin{exercise}
For $\Phi\subset\mathcal{V}^{\top}$, let $\Phi^{\Box}$ denote the subset of
$\mathcal{V}$ that contains all the vectors $v$ such that $\varphi\left(
v\right)  =0$ for all $\varphi\in\Phi$. Prove analogs of Theorems
\ref{AnnihilatorDimension} and \ref{KernelofDualvsImage}. Finally prove that
when the map $f$ has a finite-dimensional \emph{domain}, the rank of $f^{\top
}$ equals the rank of $f$ \emph{(}another\emph{\ Rank Theorem).}
\end{exercise}

\subsection{The Contragredient\label{Contragredient}}

Suppose that $\mathcal{V}$ is a finite-dimensional vector space with basis
$\mathcal{B}$ and dual basis $\mathcal{B}^{\top}$, and suppose that
$\mathcal{W}$ is an alias of $\mathcal{V}$ via the isomorphism $f:\mathcal{V}%
\rightarrow\mathcal{W}$ which then sends the basis $\mathcal{B}$ to the basis
$f\left(  \mathcal{B}\right)  $. Then the \textbf{contragredient }$f^{-\top
}=\left(  f^{-1}\right)  ^{\top}=\left(  f^{\top}\right)  ^{-1}$ maps
$\mathcal{V}^{\top}$ isomorphically onto $\mathcal{W}^{\top}$ and maps
$\mathcal{B}^{\top}$ onto the dual of the basis $f\left(  \mathcal{B}\right)
$.

\begin{theorem}
Let $f:\mathcal{V}\rightarrow\mathcal{W}$ be an isomorphism between
finite-dimensional vector spaces. Let $\mathcal{B}$ be a basis for
$\mathcal{V}$. Then the dual of the basis $f\left(  \mathcal{B}\right)  $ is
\textbf{\ }$f^{-\top}\left(  \mathcal{B}^{\top}\right)  $.
\end{theorem}

Proof: Let $x,y\in\mathcal{B}$. We then have $\left(  f^{-\top}\left(
x^{\top}\right)  \right)  \left(  f\left(  y\right)  \right)  =\left(
x^{\top}\circ f^{-1}\right)  \left(  f\left(  y\right)  \right)  =x^{\top
}\left(  y\right)  $. $\blacksquare$

This result is most often applied in the case when $\mathcal{V}=\mathcal{W}$
and the isomorphism $f$ (now an \textbf{automorphism}) is specifically being
used to effect a change of basis and the contragredient is being used to
effect the corresponding change of the dual basis.

\subsection{Product Spaces}

The function space concept introduced at the beginning of this chapter can be
profitably generalized. Instead of considering functions each value of which
lies in the same vector space $\mathcal{W}$, we are now going to allow
different vector spaces (all over the same field) for different elements of
the domain $\mathcal{D}$. Thus for each $x\in\mathcal{D}$ there is to be a
vector space $\mathcal{W}_{x}$, with each of these $\mathcal{W}_{x}$ being
over the same field, and we form all functions $f$ such that $f\left(
x\right)  \in\mathcal{W}_{x}$. The set of all such $f$ is the familiar
\emph{Cartesian product} of the $\mathcal{W}_{x}$. Any such Cartesian product
of vector spaces may be made into a vector space using pointwise operations in
the same way that we earlier used them to make $\mathcal{W}^{\mathcal{D}}$
into a function space, and it is this vector space that we mean whenever we
refer to the Cartesian product of vector spaces as a \textbf{product}
\textbf{space}. Sometimes we also refer to a product space or a Cartesian
product simply as a \textbf{product}, when there is no chance of confusion. We
use the general notation $\prod_{x\in\mathcal{D}}\mathcal{W}_{x}$ for a
product, but for the product of $n$ vector spaces $\mathcal{W}_{1}%
,\ldots,\mathcal{W}_{n}$ we often write $\mathcal{W}_{1}\times\cdots
\times\mathcal{W}_{n}$. For a field $\mathcal{F}$, we recognize $\mathcal{F}%
\times\cdots\times\mathcal{F}$ (with $n$ factors) as the familiar
$\mathcal{F}^{n}$.

When there are only a finite number of factors, a product space is the direct
sum of the ``internalized'' factors. For example, $\mathcal{W}_{1}%
\times\mathcal{W}_{2}=\left(  \mathcal{W}_{1}\times\left\{  0\right\}
\right)  \oplus\left(  \left\{  0\right\}  \times\mathcal{W}_{2}\right)  $.

\begin{exercise}
If $\mathcal{U}_{1}$ and $\mathcal{U}_{2}$ are subspaces of the same vector
space and have only $0$ in common, then $\mathcal{U}_{1}\times\mathcal{U}%
_{2}\cong\mathcal{U}_{1}\oplus\mathcal{U}_{2}$ in a manner not requiring any
choice of basis.
\end{exercise}

However, when there are an infinite number of factors, a product space is not
equal to the direct sum of ``internalized'' factors.

\begin{exercise}
Let $\mathbb{N}=\left\{  0,1,2,...\right\}  $ be the set of natural numbers
and let $\mathcal{F}$ be a field. Then the function space $\mathcal{F}%
^{\mathbb{N}}$, which is the same thing as the product space with each factor
equal to $\mathcal{F}$ and with one such factor for each natural number, is
not equal to its subspace $\bigoplus_{i\in\mathbb{N}}\left\langle
x_{i}\right\rangle $ where $x_{i}\left(  j\right)  $ is equal to $1$ if $i=j$
and is equal to $0$ otherwise.
\end{exercise}

Selecting from a product space only those functions that have finitely many
nonzero values gives a subspace called the \textbf{weak product space}. As an
example, the set of all formal power series with coefficients in a field
$\mathcal{F}$ is a product space (isomorphic to the $\mathcal{F}^{\mathbb{N}}$
of the exercise above) for which the corresponding weak product space is the
set of all polynomials with coefficients in $\mathcal{F}$. The weak product
space always equals the direct sum of the ``internalized'' factors, and for
that reason is also called the \textbf{external direct sum}. We will write
$\biguplus_{x\in\mathcal{D}}\mathcal{W}_{x}$ for the weak product space
derived from the product $\prod_{x\in\mathcal{D}}\mathcal{W}_{x}$.

\subsection{Maps Involving Products\label{MapsOnProducts}}

We single out some basic families of maps that involve products and their
factors. The function $\pi_{x}$ that sends $f\in\prod_{x\in\mathcal{D}%
}\mathcal{W}_{x}$ (or $\biguplus_{x\in\mathcal{D}}\mathcal{W}_{x}$) to
$f(x)\in\mathcal{W}_{x}$, is called the \textbf{canonical projection} onto
$\mathcal{W}_{x}$. The function $\eta_{x}$ that sends $v\in\mathcal{W}_{x}$ to
the $f\in\prod_{x\in\mathcal{D}}\mathcal{W}_{x}$ (or $\biguplus_{x\in
\mathcal{D}}\mathcal{W}_{x}$) such that $f\left(  x\right)  =v$ and $f\left(
y\right)  =0$ for $y\neq x$, is called the \textbf{canonical injection} from
$\mathcal{W}_{x} $. It is easy to verify that $\pi_{x}$ and $\eta_{x}$ are
maps. Also, we introduce $\xi_{x}=\eta_{x}\circ\pi_{x},$ called the
\textbf{component projection} for component $x$. Note that $\pi_{x}\circ
\eta_{x}=id_{x}$ (identity on $\mathcal{W}_{x}$), $\pi_{x}\circ\eta_{y}=0$
whenever $x\neq y$, and $\xi_{x}\circ\xi_{x}=\xi_{x}$. Note also that
$\sum_{x\in\mathcal{D}}\xi_{x}\left(  f\right)  =f$ so that $\sum
_{x\in\mathcal{D}}\xi_{x}$ is the identity function, where we have made the
convention that the pointwise sum of any number of $0$s is $0$.

A map $F:\mathcal{V}\rightarrow\prod_{x\in\mathcal{D}}\mathcal{W}_{x}$ from a
vector space $\mathcal{V}$ into the product of the vector spaces
$\mathcal{W}_{x}$ is uniquely determined by the map projections $\pi_{x}\circ
F$.

\begin{theorem}
\label{MapsIntoProducts}There exists exactly one map $F:\mathcal{V}%
\rightarrow\prod_{x\in\mathcal{D}}\mathcal{W}_{x}$ such that $\pi_{x}\circ
F=\varphi_{x}$ where for each $x\in\mathcal{D}$, $\varphi_{x}:\mathcal{V}%
\rightarrow\mathcal{W}_{x}$ is a prescribed map.
\end{theorem}

Proof: The function $F$ that sends $v\in\mathcal{V}$ to the $f\in\prod
_{x\in\mathcal{D}}\mathcal{W}_{x}$ such that $f\left(  x\right)  =\varphi
_{x}\left(  v\right)  $ is readily seen to be a map such that $\pi_{x}\circ
F=\varphi_{x}$. Suppose that the map $G:\mathcal{V}\rightarrow\prod
_{x\in\mathcal{D}}\mathcal{W}_{x}$ satisfies $\pi_{x}\circ G=\varphi_{x}$ for
each $x\in\mathcal{D}$. For $v\in\mathcal{V}$ let $G\left(  v\right)  =g$.
Then $g\left(  x\right)  =\pi_{x}(g)=\pi_{x}(G\left(  v\right)  )=(\pi
_{x}\circ G)\left(  v\right)  =\varphi_{x}\left(  v\right)  =f\left(
x\right)  $ and hence $g=f$. But then $G\left(  v\right)  =f=F\left(
v\right)  $ so that $G=F$, which shows the uniqueness of $F$. $\blacksquare$

\medskip\ 

For weak products, it is the maps from them, rather than the maps into them,
that have the special property.

\begin{theorem}
There exists exactly one map $F:\biguplus_{x\in\mathcal{D}}\mathcal{W}%
_{x}\rightarrow\mathcal{V}$ such that $F\circ\eta_{x}=\varphi_{x}$ where for
each $x\in\mathcal{D}$, $\varphi_{x}:\mathcal{W}_{x}\rightarrow\mathcal{V}$ is
a prescribed map.
\end{theorem}

Proof: The function $F$ that sends $f\in\biguplus_{x\in\mathcal{D}}%
\mathcal{W}_{x}$ to $\sum_{\left\{  x\in\mathcal{D}|f(x)\neq0\right\}
}\varphi_{x}\left(  f\left(  x\right)  \right)  $ is readily seen to be a map
such that $F\circ\eta_{x}=\varphi_{x}$. Suppose that the map $G:\biguplus
_{x\in\mathcal{D}}\mathcal{W}_{x}\rightarrow\mathcal{V}$ satisfies $G\circ
\eta_{x}=\varphi_{x}$ for each $x\in\mathcal{D}$. Let $f\in\biguplus
_{x\in\mathcal{D}}\mathcal{W}_{x}$. Then $f=\sum_{\left\{  x\in\mathcal{D}%
|f(x)\neq0\right\}  }\eta_{x}\left(  f\left(  x\right)  \right)  $ so that
$G\left(  f\right)  =\sum_{\left\{  x\in\mathcal{D}|f(x)\neq0\right\}
}G\left(  \eta_{x}\left(  f\left(  x\right)  \right)  \right)  =\sum_{\left\{
x\in\mathcal{D}|f(x)\neq0\right\}  }\varphi_{x}\left(  f\left(  x\right)
\right)  =F\left(  f\right)  $, which shows the uniqueness of $F$.
$\blacksquare$

\begin{exercise}
$\biguplus_{x\in\mathcal{D}}\mathcal{W}_{x}=\bigoplus_{x\in\mathcal{D}}%
\eta_{x}\left(  \mathcal{W}_{x}\right)  $.
\end{exercise}

\begin{exercise}
$\bigoplus_{x\in\mathcal{D}}\mathcal{U}_{x}\cong\biguplus_{x\in\mathcal{D}%
}\mathcal{U}_{x}$ in a manner not requiring any choice of basis.
\end{exercise}

\begin{exercise}
Let $\mathcal{W}_{1},\ldots,\mathcal{W}_{n}$ be finite dimensional vector
spaces over the same field. Then
\[
\dim\left(  \mathcal{W}_{1}\times\cdots\times\mathcal{W}_{n}\right)
=\dim\mathcal{W}_{1}+\cdots+\dim\mathcal{W}_{n}\text{.}%
\]
\end{exercise}

The last two theorems lay a foundation for two important isomorphisms. As
before, with all vector spaces over the same field, let the vector space
$\mathcal{V}$ be given, and let a vector space $\mathcal{W}_{x}$ be given for
each $x$ in some set $\mathcal{D}$. Then, utilizing map spaces, we may form
the pairs of vector spaces
\[
\mathcal{M}=\prod_{x\in\mathcal{D}}\left\{  \mathcal{V}\rightarrow
\mathcal{W}_{x}\right\}  \text{, }\mathcal{N}=\left\{  \mathcal{V}%
\rightarrow\prod_{x\in\mathcal{D}}\mathcal{W}_{x}\right\}
\]
and
\[
\mathcal{M}^{^{\prime}}=\prod_{x\in\mathcal{D}}\left\{  \mathcal{W}%
_{x}\rightarrow\mathcal{V}\right\}  \text{, }\mathcal{N}^{^{\,\prime}%
}=\left\{  \biguplus_{x\in\mathcal{D}}\mathcal{W}_{x}\rightarrow
\mathcal{V}\right\}  \text{.}%
\]

\begin{theorem}
$\mathcal{M}\cong\mathcal{N}$ and $\mathcal{M}^{^{\prime}}\cong\mathcal{N}%
^{^{\,\prime}}$.
\end{theorem}

Proof: 

$f\in\mathcal{M}$ means that for each $x\in\mathcal{D}$, $f\left(  x\right)
=\varphi_{x}$ for some map $\varphi_{x}:\mathcal{V}\rightarrow\mathcal{W}_{x}%
$, and hence there exists exactly one map $F\in\mathcal{N}$ such that $\pi
_{x}\circ F=f(x)$ for each $x\in\mathcal{D}$. This association of $F$ with $f$
therefore constitutes a well-defined function $\Lambda:\mathcal{M}%
\rightarrow\mathcal{N}$ which is easily seen to be a map. The map $\Lambda$ is
one-to-one. For if $\Lambda$ sends both $f$ and $g$ to $F$, then for each
$x\in\mathcal{D}$, $\pi_{x}\circ F=f(x)=g\left(  x\right)  $, and $f$ and $g$
must therefore be equal. Also, each $F\in\mathcal{N}$ is the image under
$\Lambda$ of the $f\in\mathcal{M}$ such that $f\left(  x\right)  =\pi_{x}\circ
F$ so that $\Lambda$ is also onto. Hence $\Lambda$ is an isomorphism from
$\mathcal{M}$ onto $\mathcal{N}$.

A similar argument shows that $\mathcal{M}^{^{\prime}}\cong\mathcal{N}%
^{^{\,\prime}}$. $\blacksquare$

\begin{exercise}
$\left(  \biguplus_{x\in\mathcal{D}}\mathcal{W}_{x}\right)  ^{\top}\cong
\prod_{x\in\mathcal{D}}\mathcal{W}_{x}^{\top}$. \emph{(}Compare with Example
\ref{DualofUltZeroSeqs} at the beginning of this chapter.\emph{)}
\end{exercise}

\begin{exercise}
$(\mathcal{W}_{1}\times\cdots\times\mathcal{W}_{n})^{\top}\cong\mathcal{W}%
_{1}^{\top}\times\cdots\times\mathcal{W}_{n}^{\top}$.
\end{exercise}

We note that the isomorphisms of the two exercises and theorem above do not
depend on any choice of basis. This will generally be the case for the
isomorphisms we will be establishing. From now on, we will usually skip
pointing out when an isomorphism does not depend on specific choices, but will
endeavor to point out any opposite cases.\medskip\ 

Up to isomorphism, the order of the factors in a product does not matter. This
can be proved quite readily directly from the product definition, but the
theorem above on maps into a product gives an interesting approach from
another perspective. We also incorporate the ``obvious'' result that using
aliases of the factors does not affect the product, up to isomorphism.

\begin{theorem}
$\prod_{x\in\mathcal{D}}\mathcal{W}_{x}\cong\prod_{x\in\mathcal{D}}%
\mathcal{W}_{\sigma(x)}^{^{\prime}}$ where $\sigma$ is a bijection of
$\mathcal{D}$ onto itself, and for each $x\in\mathcal{D}$ there is an
isomorphism $\theta_{x}:\mathcal{W}_{x}^{^{\prime}}\rightarrow\mathcal{W}_{x}$.
\end{theorem}

Proof: In Theorem \ref{MapsIntoProducts}, take the product to be $\prod
_{x\in\mathcal{D}}\mathcal{W}_{x}$, take $\mathcal{V}=\prod_{x\in\mathcal{D}%
}\mathcal{W}_{\sigma(x)}^{^{\prime}}$, and take $\varphi_{x}=\theta_{x}%
\circ\pi_{\sigma^{-1}(x)}$ so that there exists a map $\Psi$ from $\prod
_{x\in\mathcal{D}}\mathcal{W}_{\sigma(x)}^{^{\prime}}$ to $\prod
_{x\in\mathcal{D}}\mathcal{W}_{x}$ such that $\pi_{x}\circ\Psi=\theta_{x}%
\circ\pi_{\sigma^{-1}(x)}$. Interchanging the product spaces and applying the
theorem again, there exists a map $\Phi$ from $\prod_{x\in\mathcal{D}%
}\mathcal{W}_{x}$ to $\prod_{x\in\mathcal{D}}\mathcal{W}_{\sigma(x)}%
^{^{\prime}}$ such that $\pi_{\sigma^{-1}(x)}\circ\Phi=\theta_{x}^{-1}\circ
\pi_{x}$. Then $\pi_{x}\circ\Psi\circ\Phi=\theta_{x}\circ\pi_{\sigma^{-1}%
(x)}\circ\Phi=\pi_{x}$, and $\pi_{\sigma^{-1}(x)}\circ\Phi\circ\Psi=\theta
_{x}^{-1}\circ\pi_{x}\circ\Psi=\pi_{\sigma^{-1}(x)}$. Now applying the theorem
to the two cases where the $\mathcal{V}$ and the product space are the same,
the unique map determined in each case must be the identity. Hence $\Psi
\circ\Phi$ and $\Phi\circ\Psi$ must each be identity maps, thus making $\Phi$
and $\Psi$ a pair of mutually inverse maps that link $\prod_{x\in\mathcal{D}%
}\mathcal{W}_{x}$ and $\prod_{x\in\mathcal{D}}\mathcal{W}_{\sigma
(x)}^{^{\prime}}$ isomorphically.  $\blacksquare$

\subsection{Problems}

\ 

1. Let $\mathcal{V}$ be a vector space over the field $\mathcal{F}$. Let
$f\in\mathcal{V}^{\top}$. Then $f^{\top}$ is the map that sends the element
$\varphi$ of $\left(  \mathcal{F}^{1}\right)  ^{\top}$to $(\varphi\left(
1\right)  )\cdot f\in\mathcal{V}^{\top}$.

\medskip\ \ 

2.\ Let $\mathcal{V}\vartriangleleft\mathcal{W}$. Let $j:\mathcal{V}%
\rightarrow\mathcal{W}$ be the inclusion map of $\mathcal{V}$ into
$\mathcal{W}$ which sends elements in $\mathcal{V}$ to themselves. Then
$j^{\top}$ is the corresponding ``restriction of domain'' operation on the
functionals of $\mathcal{W}^{\top}$, restricting each of their domains to
$\mathcal{V}$ from the original $\mathcal{W}$. (Thus, using a frequently-seen
notation, $j^{\top}\left(  \varphi\right)  =\varphi_{\mid\mathcal{V}}$ for any
$\varphi\in\mathcal{W}^{\top}$.)

\medskip\ 

3. Let $f:\mathcal{V}\rightarrow\mathcal{W}$ be a map. Then for all
$w\in\mathcal{W}$ the equation
\[
f\left(  x\right)  =w
\]
has a unique solution for $x$ if and only if the equation
\[
f\left(  x\right)  =0
\]
has only the trivial solution $x=0$. (This result is sometimes called the
\emph{Alternative Theorem}. However, see the next problem.)

\medskip\ 

4. (Alternative Theorem) Let $\mathcal{V}$ and $\mathcal{W}$ be
finite-dimensional and let $f:\mathcal{V}\rightarrow\mathcal{W}$ be a map.
Then for a given $w$ there exists $x$ such that
\[
f\left(  x\right)  =w
\]
if and only if
\[
w\in\left(  \operatorname*{Kernel}\left(  f^{\top}\right)  \right)
^{0}\text{,}%
\]
i. e., if and only if $\varphi\left(  w\right)  =0$ for all $\varphi
\in\operatorname*{Kernel}\left(  f^{\top}\right)  $. (The reason this is
called the \emph{Alternative} Theorem is because it is usually phrased in the
equivalent form: \emph{Either} the equation $f\left(  x\right)  =w$ can be
solved for $x$, \emph{or} there is a $\varphi\in\operatorname*{Kernel}\left(
f^{\top}\right)  $ such that $\varphi\left(  w\right)  \neq0$.)

\medskip\ 

5. Let $\mathcal{V}$ and $\mathcal{W}$ be finite-dimensional and let
$f:\mathcal{V}\rightarrow\mathcal{W}$ be a map. Then $f$ is onto if and only
if $f^{\top}$ is one-to-one and $f^{\top}$ is onto if and only if $f$ is one-to-one.

\medskip\ 

6. Let $\mathcal{V}$ and $\mathcal{W}$ have the respective finite dimensions
$m$ and $n$ over the finite field $\mathcal{F}$ of $q$ elements. How many maps
$\mathcal{V}\rightarrow\mathcal{W}$ are there? If $m\geqslant n$ how many of
these are onto, if $m\leqslant n$ how many are one-to-one, and if $m=n$ how
many are isomorphisms?

\newpage

\section{Multilinearity and Tensoring}

\subsection{Multilinear Transformations}

Maps, or linear transformations, are functions of a single vector variable.
Also important are the multilinear transformations. These are the functions of
several vector variables (i. e., functions on the product space $\mathcal{V}%
_{1}\times\cdots\times\mathcal{V}_{n}$) which are linear in each variable
separately. In detail, we call the function $f:\mathcal{V}_{1}\times
\cdots\times\mathcal{V}_{n}\rightarrow\mathcal{W}$, where the $\mathcal{V}%
_{i}$ and $\mathcal{W}$ are vector spaces over the same field $\mathcal{F}$, a
\textbf{multilinear transformation}, if for each $k$%
\[
f\left(  v_{1},\ldots,v_{k-1},a\cdot u+b\cdot v,\ldots,v_{n}\right)  =
\]
\[
=a\cdot f\left(  v_{1},\ldots,v_{k-1},u,\ldots,v_{n}\right)  +b\cdot f\left(
v_{1},\ldots,v_{k-1},v,\ldots,v_{n}\right)
\]
for all vectors $u,v,v_{1},\ldots,v_{n}$ and all scalars $a,b$. We also say
that such an $f$ is $n$-\textbf{linear} (\textbf{bilinear}, \textbf{trilinear}
when $n=$ $2$, $3$). When $\mathcal{W}=\mathcal{F}$, we call a multilinear
transformation a \textbf{multilinear functional}. In general, a multilinear
transformation is not a map on $\mathcal{V}_{1}\times\cdots\times
\mathcal{V}_{n}$. However, we will soon construct a vector space on which the
linear transformations ``are'' the multilinear transformations on
$\mathcal{V}_{1}\times\cdots\times\mathcal{V}_{n}$.\medskip\ 

A unique multilinear transformation is determined by the arbitrary assignment
of values on all tuples of basis vectors.

\begin{theorem}
Let $\mathcal{V}_{1},\ldots,\mathcal{V}_{n},\mathcal{W}$ be vector spaces over
the same field and for each $\mathcal{V}_{i}$ let the basis $\mathcal{B}_{i}$
be given. Then given the function $f_{0}:\mathcal{B}_{1}\times\cdots
\times\mathcal{B}_{n}\rightarrow\mathcal{W}$, there is a unique multilinear
transformation $f:\mathcal{V}_{1}\times\cdots\times\mathcal{V}_{n}%
\rightarrow\mathcal{W}$ such that $f$ agrees with $f_{0}$ on $\mathcal{B}%
_{1}\times\cdots\times\mathcal{B}_{n}$.
\end{theorem}

Proof: Given $v_{i}\in\mathcal{V}_{i}$ there is a unique finite subset
$\mathcal{X}_{i}$ of $\mathcal{B}_{i}$ and unique nonzero scalars $a_{x_{i}}$
such that $v_{i}=\sum_{x_{i}\in\mathcal{X}_{i}}a_{x_{i}}\cdot x_{i}$. The
$n$-linearity of $f:\mathcal{V}_{1}\times\cdots\times\mathcal{V}%
_{n}\rightarrow\mathcal{W}$ implies
\[
f\left(  v_{1},\ldots,v_{n}\right)  =\sum_{x_{1}\in\mathcal{X}_{1}}\cdots
\sum_{x_{n}\in\mathcal{X}_{n}}a_{x_{1}}\cdots a_{x_{n}}\cdot f\left(
x_{1},\ldots,x_{n}\right)  \text{.}%
\]
Hence for $f$ to agree with $f_{0}$ on $\mathcal{B}_{1}\times\cdots
\times\mathcal{B}_{n}$, it can only have the value
\[
f\left(  v_{1},\ldots,v_{n}\right)  =\sum_{x_{1}\in\mathcal{X}_{1}}\cdots
\sum_{x_{n}\in\mathcal{X}_{n}}a_{x_{1}}\cdots a_{x_{n}}\cdot f_{0}\left(
x_{1},\ldots,x_{n}\right)  \text{.}%
\]
On the other hand, setting
\[
f\left(  v_{1},\ldots,v_{n}\right)  =\sum_{x_{1}\in\mathcal{X}_{1}}\cdots
\sum_{x_{n}\in\mathcal{X}_{n}}a_{x_{1}}\cdots a_{x_{n}}\cdot f_{0}\left(
x_{1},\ldots,x_{n}\right)
\]
does define a function $f:\mathcal{V}_{1}\times\cdots\times\mathcal{V}%
_{n}\rightarrow\mathcal{W}$ which clearly agrees with $f_{0}$ on
$\mathcal{B}_{1}\times\cdots\times\mathcal{B}_{n}$, and this $f$ is readily
verified to be $n$-linear. $\blacksquare$

\subsection{Functional Product Spaces}

Let $\mathcal{S}_{1},\ldots,\mathcal{S}_{n}$ be sets of linear functionals on
the respective vector spaces $\mathcal{V}_{1},\ldots,\mathcal{V}_{n}$ all over
the same field $\mathcal{F}$. By the \textbf{functional product }%
$\mathcal{S}_{1}\cdots\mathcal{S}_{n}$ is meant the set of all products
$f_{1}\cdots f_{n}$ with each $f_{i}\in\mathcal{S}_{i}$, where by such a
product of linear functionals is meant the function $f_{1}\cdots
f_{n}:\mathcal{V}_{1}\times\cdots\times\mathcal{V}_{n}\rightarrow\mathcal{F}$
such that $f_{1}\cdots f_{n}\left(  v_{1},\ldots,v_{n}\right)  =f_{1}\left(
v_{1}\right)  \cdots f_{n}\left(  v_{n}\right)  $. The function-space subspace
$\left\langle \mathcal{S}_{1}\cdots\mathcal{S}_{n}\right\rangle $ obtained by
taking all linear combinations of finite subsets of the functional product
$\mathcal{S}_{1}\cdots\mathcal{S}_{n}$ is called the \textbf{functional
product space} generated by $\mathcal{S}_{1}\cdots\mathcal{S}_{n}$.

\begin{exercise}
Each element of $\left\langle \mathcal{S}_{1}\cdots\mathcal{S}_{n}%
\right\rangle $ is a multilinear functional on $\mathcal{V}_{1}\times
\cdots\times\mathcal{V}_{n}$.
\end{exercise}

\begin{exercise}
$\left\langle \left\langle \mathcal{S}_{1}\right\rangle \cdots\left\langle
\mathcal{S}_{n}\right\rangle \right\rangle =\left\langle \mathcal{S}_{1}%
\cdots\mathcal{S}_{n}\right\rangle $.
\end{exercise}

\begin{lemma}
Let $\mathcal{W}_{1},\ldots,\mathcal{W}_{n}$ be vector spaces of linear
functionals, all over the same field, with respective bases $\mathcal{B}%
_{1},\ldots,\mathcal{B}_{n}$. Then $\mathcal{B}_{1}\cdots\mathcal{B}_{n}$ is
an independent set.
\end{lemma}

Proof: Induction will be employed. The result is clearly true for $n=1$, and
suppose that it is true for $n=N-1$. Let $\mathcal{X}$ be a finite subset of
$\mathcal{B}_{1}\cdots\mathcal{B}_{N}$ such that $\sum_{x\in\mathcal{X}}%
a_{x}\cdot x=0$. Each $x\in\mathcal{X}$ has the form $yz$ where $y\in
\mathcal{B}_{1}\cdots\mathcal{B}_{N-1}$ and $z$ $\in\mathcal{B}_{N}$.
Collecting terms on the distinct $z$, $\sum_{x\in\mathcal{X}}a_{x}\cdot
x=\sum_{z}(\sum_{x=yz}a_{x}\cdot y)z$. The independence of $\mathcal{B}_{N}$
implies that for each $z$, $\sum_{x=yz}a_{x}\cdot y=0$, and the assumed
independence of $\mathcal{B}_{1}\cdots\mathcal{B}_{N-1}$ then implies that
each $a_{x}$ is zero. Hence $\mathcal{B}_{1}\cdots\mathcal{B}_{N}$ is an
independent set, and by induction, $\mathcal{B}_{1}\cdots\mathcal{B}_{n}$ is
an independent set for each $n$. $\blacksquare$

\begin{proposition}
Let $\mathcal{W}_{1},\ldots,\mathcal{W}_{n}$ be vector spaces of linear
functionals, all over the same field, with respective bases $\mathcal{B}%
_{1},\ldots,\mathcal{B}_{n}$. Then $\mathcal{B}_{1}\cdots\mathcal{B}_{n}$ is a
basis for $\left\langle \mathcal{W}_{1}\cdots\mathcal{W}_{n}\right\rangle $.
\end{proposition}

Proof: By the exercise above, $\left\langle \mathcal{W}_{1}\cdots
\mathcal{W}_{n}\right\rangle =\left\langle \mathcal{B}_{1}\cdots
\mathcal{B}_{n}\right\rangle $ and hence $\mathcal{B}_{1}\cdots\mathcal{B}%
_{n}$ spans $\left\langle \mathcal{W}_{1}\cdots\mathcal{W}_{n}\right\rangle $.
By the lemma above, $\mathcal{B}_{1}\cdots\mathcal{B}_{n}$ is an independent
set. Thus $\mathcal{B}_{1}\cdots\mathcal{B}_{n}$ is a basis for $\left\langle
\mathcal{W}_{1}\cdots\mathcal{W}_{n}\right\rangle $. $\blacksquare$

\subsection{Functional Tensor Products}

In the previous chapter, it was shown that the vector space $\mathcal{V}$ may
be embedded in its double dual $\mathcal{V}^{\top\top}$ independently of any
choice of basis, effectively making each vector in $\mathcal{V}$ into a linear
functional on $\mathcal{V}^{\top}$. Thus given the vector spaces
$\mathcal{V}_{1},\ldots,\mathcal{V}_{n}$, all over the same field
$\mathcal{F}$, and the natural injection $\Theta_{i}$ that embeds each
$\mathcal{V}_{i}$ in its double dual, we may form the functional product space
$\left\langle \Theta_{1}(\mathcal{V}_{1})\cdots\Theta_{n}(\mathcal{V}%
_{n})\right\rangle $ which will be called the \textbf{functional tensor
product} of the $\mathcal{V}_{i}$ and which will be denoted by $\mathcal{V}%
_{1}\bigotimes\cdots\bigotimes\mathcal{V}_{n}$. Similarly, an element
$\Theta_{1}(v_{1})\cdots\Theta_{n}(v_{n})$ of $\Theta_{1}(\mathcal{V}%
_{1})\cdots\Theta_{n}(\mathcal{V}_{n})$ will be called the \textbf{functional
tensor product} of the $v_{i}$ and will be denoted by $v_{1}\otimes
\cdots\otimes v_{n}$.

If $\mathcal{W}$ is a vector space over the same field as the $\mathcal{V}%
_{i}$, the ``universal'' multilinear transformation $\Xi:\mathcal{V}_{1}%
\times\cdots\times\mathcal{V}_{n}\rightarrow\mathcal{V}_{1}\bigotimes
\cdots\bigotimes\mathcal{V}_{n}$ which sends $(v_{1},\ldots,v_{n})$ to
$v_{1}\otimes\cdots\otimes v_{n}$ exchanges the multilinear transformations
$f:\mathcal{V}_{1}\times\cdots\times\mathcal{V}_{n}\rightarrow\mathcal{W}$
with the linear transformations $\varphi:\mathcal{V}_{1}\bigotimes
\cdots\bigotimes\mathcal{V}_{n}\rightarrow\mathcal{W}$ via the composition of
functions $f=\varphi\circ\Xi$.

\begin{theorem}
\label{TensorProperty}Let $\mathcal{V}_{1},\ldots,\mathcal{V}_{n},\mathcal{W}$
be vector spaces over the same field. Then for each multilinear transformation
$f:\mathcal{V}_{1}\times\cdots\times\mathcal{V}_{n}\rightarrow\mathcal{W}$
there exists a unique linear transformation $\varphi:\mathcal{V}_{1}%
\bigotimes\cdots\bigotimes\mathcal{V}_{n}\rightarrow\mathcal{W}$ such that
$f=\varphi\circ\Xi$, where $\Xi:\mathcal{V}_{1}\times\cdots\times
\mathcal{V}_{n}\rightarrow\mathcal{V}_{1}\bigotimes\cdots\bigotimes
\mathcal{V}_{n}$ is the tensor product function that sends $(v_{1}%
,\ldots,v_{n})$ to $v_{1}\otimes\cdots\otimes v_{n}$.
\end{theorem}

Proof: For each $\mathcal{V}_{i}$, choose some basis $\mathcal{B}_{i}$. Then
by the proposition above, $\mathcal{B}=\left\{  x_{1}\otimes\cdots\otimes
x_{n}\,|\,x_{1}\in\mathcal{B}_{1},\ldots,x_{n}\in\mathcal{B}_{n}\right\}  $ is
a basis for $\mathcal{V}_{1}\bigotimes\cdots\bigotimes\mathcal{V}_{n}$. Given
the multilinear transformation $f:\mathcal{V}_{1}\times\cdots\times
\mathcal{V}_{n}\rightarrow\mathcal{W}$, consider the linear transformation
$\varphi:\mathcal{V}_{1}\bigotimes\cdots\bigotimes\mathcal{V}_{n}%
\rightarrow\mathcal{W}$ such that for $x_{1}\in\mathcal{B}_{1},\ldots,x_{n}%
\in\mathcal{B}_{n}$
\[
\varphi\left(  x_{1}\otimes\cdots\otimes x_{n}\right)  =\varphi\left(
\Xi(x_{1},\ldots,x_{n})\right)  =f\left(  x_{1},\ldots,x_{n}\right)  \text{.}%
\]
Since $\varphi\circ\Xi$ is clearly $n$-linear, it must then equal $f$ by the
theorem above. Every $\varphi$ determines an $n$-linear $f$ via $f=\varphi
\circ\Xi$, but since the values of $\varphi$ on $\mathcal{B}$ determine it
completely, and the values of $f$ on $\mathcal{B}_{1}\times\cdots
\times\mathcal{B}_{n}$ determine it completely, a different $\varphi$
determines a different $f$. $\blacksquare$\medskip\ 

Thus, while multilinear transformations themselves are generally not maps on
$\mathcal{V}_{1}\times\cdots\times\mathcal{V}_{n}$, they do correspond
one-to-one to the maps on the related vector space $\mathcal{V}_{1}%
\bigotimes\cdots\bigotimes\mathcal{V}_{n}=\left\langle \left\{  v_{1}%
\otimes\cdots\otimes v_{n}\right\}  \right\rangle $. In fact, we can interpret
this correspondence as an isomorphism if we introduce the function space
subspace $\left\{  \mathcal{V}_{1}\times\cdots\times\mathcal{V}_{n}%
\overset{(n-linear)}{\rightarrow}\mathcal{W}\right\}  $ of $n$-linear
functions, which is then clearly isomorphic to the map space $\left\{
\mathcal{V}_{1}\bigotimes\cdots\bigotimes\mathcal{V}_{n}\rightarrow
\mathcal{W}\right\}  $. The vector space $\left\{  \mathcal{V}_{1}\times
\cdots\times\mathcal{V}_{n}\overset{(n-linear)}{\rightarrow}\mathcal{F}%
\right\}  $ of multilinear functionals is therefore isomorphic to
$(\mathcal{V}_{1}\bigotimes\cdots\bigotimes\mathcal{V}_{n})^{\top}$. If all
the $\mathcal{V}_{i}$ are finite-dimensional, we then have $\mathcal{V}%
_{1}\bigotimes\cdots\bigotimes\mathcal{V}_{n}\cong\left\{  \mathcal{V}%
_{1}\times\cdots\times\mathcal{V}_{n}\overset{(n-linear)}{\rightarrow
}\mathcal{F}\right\}  ^{\top}$.

When the $\mathcal{V}_{i}$ are finite-dimensional, we may discard any even
number of consecutive dualization operators, since we identify a
finite-dimensional space with its double dual. Therefore, $\mathcal{V}%
_{1}^{\top}\bigotimes\cdots\bigotimes\mathcal{V}_{n}^{\top}$ is just
$\left\langle \mathcal{V}_{1}^{\top}\cdots\mathcal{V}_{n}^{\top}\right\rangle
$. Thus $\mathcal{V}_{1}^{\top}\bigotimes\cdots\bigotimes\mathcal{V}_{n}%
^{\top}$ is a subspace of $\left\{  \mathcal{V}_{1}\times\cdots\times
\mathcal{V}_{n}\overset{(n-linear)}{\rightarrow}\mathcal{F}\right\}  $, and it
is easy to see that both have the same dimension, so they are actually equal.
This and our observations immediately above then give us the following result.

\begin{theorem}
\label{DualofTensorProduct}Let $\mathcal{V}_{1},\ldots,\mathcal{V}_{n}$ be
vector spaces over the field $\mathcal{F}$. Then
\[
(\mathcal{V}_{1}%
{\textstyle\bigotimes}
\cdots%
{\textstyle\bigotimes}
\mathcal{V}_{n})^{\top}\cong\left\{  \mathcal{V}_{1}\times\cdots
\times\mathcal{V}_{n}\overset{(n-linear)}{\rightarrow}\mathcal{F}\right\}  .
\]
If the $\mathcal{V}_{i}$ are all finite-dimensional, we also have
\[
\mathcal{V}_{1}%
{\textstyle\bigotimes}
\cdots%
{\textstyle\bigotimes}
\mathcal{V}_{n}\cong\left\{  \mathcal{V}_{1}\times\cdots\times\mathcal{V}%
_{n}\overset{(n-linear)}{\rightarrow}\mathcal{F}\right\}  ^{\top}%
\]
and
\[
(\mathcal{V}_{1}%
{\textstyle\bigotimes}
\cdots%
{\textstyle\bigotimes}
\mathcal{V}_{n})^{\top}\cong\mathcal{V}_{1}^{\top}%
{\textstyle\bigotimes}
\cdots%
{\textstyle\bigotimes}
\mathcal{V}_{n}^{\top}=\left\{  \mathcal{V}_{1}\times\cdots\times
\mathcal{V}_{n}\overset{(n-linear)}{\rightarrow}\mathcal{F}\right\}
.:EndProof
\]
\end{theorem}

\subsection{Tensor Products in General}

Theorem \ref{TensorProperty} contains the essence of the tensor product
concept. By making a definition out of the \emph{universal property} that this
theorem ascribes to $\mathcal{V}_{1}\bigotimes\cdots\bigotimes\mathcal{V}_{n}$
and $\Xi$, a general tensor product concept results. Let the vector spaces
$\mathcal{V}_{1},\ldots,\mathcal{V}_{n}$ each over the same field
$\mathcal{F}$ be given. A \textbf{tensor product} of $\mathcal{V}_{1}%
,\ldots,\mathcal{V}_{n}$ (in that order) is a vector space $\Pi$ (a
\textbf{tensor product space}) over $\mathcal{F}$ along with an $n$-linear
function $\Upsilon:\mathcal{V}_{1}\times\cdots\times\mathcal{V}_{n}%
\rightarrow\Pi$ (a \textbf{tensor product function}) such that given any
vector space $\mathcal{W}$, every $n$-linear function $f:\mathcal{V}_{1}%
\times\cdots\times\mathcal{V}_{n}\rightarrow\mathcal{W}$ is equal to
$\varphi\circ\Upsilon$ for a unique map $\varphi:\Pi\rightarrow\mathcal{W}$.
The functional tensor product $\mathcal{V}_{1}\bigotimes\cdots\bigotimes
\mathcal{V}_{n}$ along with the functional tensor product function $\Xi$ are
thus an example of a tensor product of $\mathcal{V}_{1},\ldots,\mathcal{V}%
_{n}$.

Another tensor product of the same factor spaces may be defined using any
alias of a tensor product space along with the appropriate tensor product function.

\begin{theorem}
\label{IsoToTensorProd}Let $\mathcal{V}_{1},\ldots,\mathcal{V}_{n}$ be vector
spaces over the same field. Then if $\Pi$ along with the $n$-linear function
$\Upsilon:\mathcal{V}_{1}\times\cdots\times\mathcal{V}_{n}\rightarrow\Pi$ is a
tensor product of $\mathcal{V}_{1},\ldots,\mathcal{V}_{n}$, and $\Theta
:\Pi\rightarrow\Pi^{\prime}$ is an isomorphism, then $\Pi^{\prime}$ along with
$\Theta\circ\Upsilon$ is also a tensor product of $\mathcal{V}_{1}%
,\ldots,\mathcal{V}_{n}$.
\end{theorem}

Proof: Given the vector space $\mathcal{W}$, the $n$-linear function
$f:\mathcal{V}_{1}\times\cdots\times\mathcal{V}_{n}\rightarrow\mathcal{W}$ is
equal to $\varphi\circ\Upsilon$ for a unique map $\varphi:\Pi\rightarrow
\mathcal{W}$. The map $\varphi\circ\Theta^{-1}$ therefore satisfies $\left(
\varphi\circ\Theta^{-1}\right)  \circ\left(  \Theta\circ\Upsilon\right)  =f$.
On the other hand, if $\varphi^{\prime}\circ\left(  \Theta\circ\Upsilon
\right)  =(\varphi^{\prime}\circ\Theta)\circ\Upsilon=f$, then by the
uniqueness of $\varphi$, $\varphi^{\prime}\circ\Theta=\varphi$ and hence
$\varphi^{\prime}=\varphi\circ\Theta^{-1}$. $\blacksquare$

\begin{exercise}
Make sense of the tensor product concept in the case where $n=1$, so that we
have $1$-linearity and the tensor product of a single factor.\ \ 
\end{exercise}

If in the manner of the previous theorem, we have isomorphic tensor product
spaces with the tensor product functions related by the isomorphism as in the
theorem, we say that the tensor products are \textbf{linked} by the
isomorphism of their tensor product spaces. We now show that tensor products
of the same factors always have isomorphic tensor product spaces, and in fact
there is a unique isomorphism between them that links them.

\begin{theorem}
Let $\mathcal{V}_{1},\ldots,\mathcal{V}_{n}$ be vector spaces over the same
field. Then if $\Pi$ along with the $n$-linear function $\Upsilon
:\mathcal{V}_{1}\times\cdots\times\mathcal{V}_{n}\rightarrow\Pi$, and
$\Pi^{\prime}$ along with the $n$-linear function $\Upsilon^{\prime
}:\mathcal{V}_{1}\times\cdots\times\mathcal{V}_{n}\rightarrow\Pi^{\prime}$,
are each tensor products of $\mathcal{V}_{1},\ldots,\mathcal{V}_{n}$, then
there is a unique isomorphism\textbf{\ }$\Theta:\Pi\rightarrow\Pi^{\prime}$
such that $\Theta\circ\Upsilon=\Upsilon^{\prime}$.
\end{theorem}

Proof: In the tensor product definition above, take $\mathcal{W}=\Pi^{\prime}$
and $f=\Upsilon^{\prime}$. Then there is a unique map $\Theta:\Pi
\rightarrow\Pi^{\prime}$ such that $\Upsilon^{\prime}=\Theta\circ\Upsilon$.
Similarly, there is a unique map $\Theta^{\prime}:\Pi^{\prime}\rightarrow\Pi$
such that $\Upsilon=\Theta^{\prime}\circ\Upsilon^{\prime}$. Hence
$\Upsilon^{\prime}=\Theta\circ\Theta^{\prime}\circ\Upsilon^{\prime}$ and
$\Upsilon=\Theta^{\prime}\circ\Theta\circ\Upsilon$. Applying the tensor
product definition again, the unique map $\varphi:\Pi\rightarrow\Pi$ such that
$\Upsilon=\varphi\circ\Upsilon$ must be the identity. We conclude that
$\Theta^{\prime}\circ\Theta$ is the identity and the same for $\Theta
\circ\Theta^{\prime}$. Hence the unique map $\Theta$ has the inverse
$\Theta^{\prime}$ and is therefore an isomorphism. $\blacksquare$\smallskip\ 

We will write any tensor product space as $\mathcal{V}_{1}\bigotimes
\cdots\bigotimes\mathcal{V}_{n}$ just as we did for the functional tensor
product. Similarly we will write $v_{1}\otimes\cdots\otimes v_{n}$ for
$\Upsilon\left(  v_{1},\ldots,v_{n}\right)  $, and call it the \textbf{tensor
product} of the vectors $v_{1},\ldots,v_{n}$, just as we did for the
functional tensor product $\Xi\left(  v_{1},\ldots,v_{n}\right)  $. With
isomorphic tensor product spaces, we will always assume that the tensor
products are linked by the isomorphism. By assuming this, tensor products of
the same vectors will then correspond among all tensor products of the same
factor spaces:$\;v_{1}\otimes\cdots\otimes v_{n}$ in one of these spaces is
always an isomorph of $v_{1}\otimes\cdots\otimes v_{n}$ in any of the others.

\begin{exercise}
A tensor product space is spanned by the image of its associated tensor
product function: $\left\langle \left\{  v_{1}\otimes\cdots\otimes
v_{n}\,|\,v_{1}\in\mathcal{V}_{1},\ldots,v_{n}\in\mathcal{V}_{n}\right\}
\right\rangle =\mathcal{V}_{1}\bigotimes\cdots\bigotimes\mathcal{V}_{n}$.
\end{exercise}

Up to isomorphism, the order of the factors in a tensor product does not
matter. We also incorporate the ``obvious'' result that using aliases of the
factors does not affect the result, up to isomorphism.

\begin{theorem}
Let the vector spaces $\mathcal{V}_{1}^{^{\prime}},\ldots,\mathcal{V}%
_{n}^{^{\prime}}$, all over the same field, be \emph{(}possibly
reordered\emph{)} aliases of $\mathcal{V}_{1},\ldots,\mathcal{V}_{n}$, and let
$\mathcal{V}_{\otimes}=\mathcal{V}_{1}\bigotimes\cdots\bigotimes
\mathcal{V}_{n}$ along with the tensor product function $\Upsilon
:\mathcal{V}_{\times}=\mathcal{V}_{1}\times\cdots\times\mathcal{V}%
_{n}\rightarrow\mathcal{V}_{\otimes}$ be a tensor product, and let
$\mathcal{V}_{\otimes}^{^{\prime}}=\mathcal{V}_{1}^{^{\prime}}\bigotimes
\cdots\bigotimes\mathcal{V}_{n}^{^{\prime}}$ along with the tensor product
function $\Upsilon^{^{\prime}}:\mathcal{V}_{\times}^{^{\prime}}=\mathcal{V}%
_{1}^{^{\prime}}\times\cdots\times\mathcal{V}_{n}^{^{\prime}}\rightarrow
\mathcal{V}_{\otimes}^{^{\prime}}$ be a tensor product. Then $\mathcal{V}%
_{\otimes}\cong\mathcal{V}_{\otimes}^{^{\prime}}$.
\end{theorem}

Proof: Let $\Phi:\mathcal{V}_{\times}\rightarrow\mathcal{V}_{\times}%
^{^{\prime}}$ denote the isomorphism that exists from the product space
$\mathcal{V}_{\times}$ to the product space $\mathcal{V}_{\times}^{^{\prime}}%
$. It is easy to see that $\Upsilon^{^{\prime}}\circ\Phi$ and $\Upsilon
\circ\Phi^{-1}$are $n$-linear. Hence there is a (unique) map $\Theta
:\mathcal{V}_{\otimes}\rightarrow\mathcal{V}_{\otimes}^{^{\prime}}$ such that
$\Theta\circ\Upsilon=\Upsilon^{^{\prime}}\circ\Phi$ and there is a (unique)
map $\Theta^{^{\prime}}:\mathcal{V}_{\otimes}^{^{\prime}}\rightarrow
\mathcal{V}_{\otimes}$ such that $\Theta^{^{\prime}}\circ\Upsilon^{^{\prime}%
}=\Upsilon\circ\Phi^{-1}$. From this we readily deduce that $\Upsilon
^{^{\prime}}=\Theta\circ\Theta^{^{\prime}}\circ\Upsilon^{^{\prime}}$ and
$\Upsilon=\Theta^{^{\prime}}\circ\Theta\circ\Upsilon$. Hence the map $\Theta$
has the inverse $\Theta^{^{\prime}}$ and is therefore an isomorphism.
$\blacksquare$\smallskip\ 

Up to isomorphism, tensor multiplication of vector spaces may be performed iteratively.

\begin{theorem}
\label{Associative}Let $\mathcal{V}_{1},\ldots,\mathcal{V}_{n}$ be vector
spaces over the same field. Then for any integer $k$ such that $1\leqslant
k<n$ there is an isomorphism
\[
\Theta:(\mathcal{V}_{1}%
{\textstyle\bigotimes}
\cdots%
{\textstyle\bigotimes}
\mathcal{V}_{k})%
{\textstyle\bigotimes}
(\mathcal{V}_{k+1}%
{\textstyle\bigotimes}
\cdots%
{\textstyle\bigotimes}
\mathcal{V}_{n})\rightarrow\mathcal{V}_{1}%
{\textstyle\bigotimes}
\cdots%
{\textstyle\bigotimes}
\mathcal{V}_{n}%
\]
such that
\[
\Theta\left(  \left(  v_{1}\otimes\cdots\otimes v_{k}\right)  \otimes
(v_{k+1}\otimes\cdots\otimes v_{n})\right)  =v_{1}\otimes\cdots\otimes v_{n}.
\]
\end{theorem}

Proof: Set $\mathcal{V}_{\times}=\mathcal{V}_{1}\times\cdots\times
\mathcal{V}_{n}$, $\mathcal{V}_{\otimes}=\mathcal{V}_{1}\bigotimes
\cdots\bigotimes\mathcal{V}_{n}$, $\overline{\mathcal{V}}_{\times}%
=\mathcal{V}_{1}\times\cdots\times\mathcal{V}_{k}$, $\overline{\mathcal{V}%
}_{\otimes}=\mathcal{V}_{1}\bigotimes\cdots\bigotimes\mathcal{V}_{k}$,
$\overline{\overline{\mathcal{V}}}_{\times}=\mathcal{V}_{k+1}\times
\cdots\times\mathcal{V}_{n}$, and $\overline{\overline{\mathcal{V}}}_{\otimes
}=\mathcal{V}_{k+1}\bigotimes\cdots\bigotimes\mathcal{V}_{n}$.

For fixed $v_{k+1}\otimes\cdots\otimes v_{n}\in\overline{\overline
{\mathcal{V}}}_{\otimes}$ we define the $k$-linear function $\overline
{f}_{v_{k+1}\otimes\cdots\otimes v_{n}}:\overline{\mathcal{V}}_{\times
}\rightarrow\mathcal{V}_{\otimes}$ by $\overline{f}_{v_{k+1}\otimes
\cdots\otimes v_{n}}\left(  v_{1},\ldots,v_{k}\right)  =v_{1}\otimes
\cdots\otimes v_{n}$. Corresponding to $\overline{f}_{v_{k+1}\otimes
\cdots\otimes v_{n}}$ is the (unique) map $\overline{\Theta}_{v_{k+1}%
\otimes\cdots\otimes v_{n}}:\overline{\mathcal{V}}_{\otimes}\rightarrow
\mathcal{V}_{\otimes}$ such that
\[
\overline{\Theta}_{v_{k+1}\otimes\cdots\otimes v_{n}}\left(  v_{1}%
\otimes\cdots\otimes v_{k}\right)  =v_{1}\otimes\cdots\otimes v_{n}.
\]
Similarly, for fixed $v_{1}\otimes\cdots\otimes v_{k}\in\overline{\mathcal{V}%
}_{\otimes}$ we define the $\left(  n-k\right)  $-linear function
$\overline{\overline{f}}_{v_{1}\otimes\cdots\otimes v_{k}}:\overline
{\overline{\mathcal{V}}}_{\times}\rightarrow\mathcal{V}_{\otimes}$ by
$\overline{\overline{f}}_{v_{1}\otimes\cdots\otimes v_{k}}\left(
v_{k+1},\ldots,v_{n}\right)  =v_{1}\otimes\cdots\otimes v_{n}$. Corresponding
to $\overline{\overline{f}}_{v_{1}\otimes\cdots\otimes v_{k}}$ is the (unique)
map $\overline{\overline{\Theta}}_{v_{1}\otimes\cdots\otimes v_{k}}%
:\overline{\overline{\mathcal{V}}}_{\otimes}\rightarrow\mathcal{V}_{\otimes}$
such that
\[
\overline{\overline{\Theta}}_{v_{1}\otimes\cdots\otimes v_{k}}\left(
v_{k+1}\otimes\cdots\otimes v_{n}\right)  =v_{1}\otimes\cdots\otimes v_{n}.
\]

We then define the function $f:\overline{\mathcal{V}}_{\otimes}\times
\overline{\overline{\mathcal{V}}}_{\otimes}\rightarrow\mathcal{V}_{\otimes}$
by the formula
\[
f\left(  x,y\right)  =\sum a_{v_{1}\otimes\cdots\otimes v_{k}}\overline
{\overline{\Theta}}_{v_{1}\otimes\cdots\otimes v_{k}}\left(  y\right)
\]
when $x=\sum a_{v_{1}\otimes\cdots\otimes v_{k}}\cdot v_{1}\otimes
\cdots\otimes v_{k}$. We claim this formula gives the same result no matter
how $x$ is expressed as a linear combination of elements of the form
$v_{1}\otimes\cdots\otimes v_{k}$. To aid in verifying this, let $y$ be
expressed as $y=\sum b_{v_{k+1}\otimes\cdots\otimes v_{n}}\cdot v_{k+1}%
\otimes\cdots\otimes v_{n}$. The formula may then be written
\[
f\left(  x,y\right)  =\sum a_{v_{1}\otimes\cdots\otimes v_{k}}\cdot\sum
b_{v_{k+1}\otimes\cdots\otimes v_{n}}\cdot\overline{\overline{\Theta}}%
_{v_{1}\otimes\cdots\otimes v_{k}}\left(  v_{k+1}\otimes\cdots\otimes
v_{n}\right)  \text{.}%
\]
But $\overline{\overline{\Theta}}_{v_{1}\otimes\cdots\otimes v_{k}}\left(
v_{k+1}\otimes\cdots\otimes v_{n}\right)  =\overline{\Theta}_{v_{k+1}%
\otimes\cdots\otimes v_{n}}\left(  v_{1}\otimes\cdots\otimes v_{k}\right)  $
so that the formula is equivalent to
\[
f\left(  x,y\right)  =\sum a_{v_{1}\otimes\cdots\otimes v_{k}}\cdot\sum
b_{v_{k+1}\otimes\cdots\otimes v_{n}}\cdot\overline{\Theta}_{v_{k+1}%
\otimes\cdots\otimes v_{n}}\left(  v_{1}\otimes\cdots\otimes v_{k}\right)
\text{.}%
\]
But this is the same as
\[
f\left(  x,y\right)  =\sum b_{v_{k+1}\otimes\cdots\otimes v_{n}}\cdot
\overline{\Theta}_{v_{k+1}\otimes\cdots\otimes v_{n}}\left(  x\right)
\]
which of course has no dependence on how $x$ is expressed as a linear
combination of elements of the form $v_{1}\otimes\cdots\otimes v_{k}$.

$f$ is clearly linear in either argument when the other argument is held
fixed. Corresponding to $f$ is the (unique) map $\Theta:\overline{\mathcal{V}%
}_{\otimes}\bigotimes\overline{\overline{\mathcal{V}}}_{\otimes}%
\rightarrow\mathcal{V}_{\otimes}$ such that $\Theta\left(  x\otimes y\right)
=f\left(  x,y\right)  $ so that, in particular,
\[
\Theta\left(  \left(  v_{1}\otimes\cdots\otimes v_{k}\right)  \otimes\left(
v_{k+1}\otimes\cdots\otimes v_{n}\right)  \right)  =v_{1}\otimes\cdots\otimes
v_{n}.
\]

We define an $n$-linear function $f^{^{\prime}}:\mathcal{V}_{\times
}\rightarrow\overline{\mathcal{V}}_{\otimes}\bigotimes\overline{\overline
{\mathcal{V}}}_{\otimes}$ by setting
\[
f^{^{\prime}}\left(  v_{1},\ldots,v_{n}\right)  =\left(  v_{1}\otimes
\cdots\otimes v_{k}\right)  \otimes\left(  v_{k+1}\otimes\cdots\otimes
v_{n}\right)  .
\]
Corresponding to $f^{^{\prime}}$ is the (unique) map $\Theta^{^{\prime}%
}:\mathcal{V}_{\otimes}\rightarrow\overline{\mathcal{V}}_{\otimes}%
\bigotimes\overline{\overline{\mathcal{V}}}_{\otimes}$ such that
\[
\Theta^{^{\prime}}(v_{1}\otimes\cdots\otimes v_{n})=\left(  v_{1}\otimes
\cdots\otimes v_{k}\right)  \otimes\left(  v_{k+1}\otimes\cdots\otimes
v_{n}\right)  .
\]

Thus each of $\Theta^{^{\prime}}\circ\Theta$ and $\Theta\circ\Theta^{^{\prime
}}$ coincides with the identity on a spanning set for its domain, so each is
in fact the identity map. Hence the map $\Theta$ has the inverse
$\Theta^{^{\prime}}$ and therefore is an isomorphism. $\blacksquare$

\begin{exercise}
$\left(  \cdots\left(  \left(  \mathcal{V}_{1}\bigotimes\mathcal{V}%
_{2}\right)  \bigotimes\mathcal{V}_{3}\right)  \bigotimes\cdots\bigotimes
\mathcal{V}_{n-1}\right)  \bigotimes\mathcal{V}_{n}\cong\mathcal{V}%
_{1}\bigotimes\cdots\bigotimes\mathcal{V}_{n}$.
\end{exercise}

\begin{corollary}
$\left(  \mathcal{V}_{1}\bigotimes\mathcal{V}_{2}\right)  \bigotimes
\mathcal{V}_{3}\cong\mathcal{V}_{1}\bigotimes\left(  \mathcal{V}_{2}%
\bigotimes\mathcal{V}_{3}\right)  \cong\mathcal{V}_{1}\bigotimes
\mathcal{V}_{2}\bigotimes\mathcal{V}_{3}$. $\blacksquare$\medskip\ 
\end{corollary}

The two preceding theorems form the foundation for the following general
associativity result.

\begin{theorem}
Tensor products involving the same vector spaces are isomorphic no matter how
the factors are grouped.
\end{theorem}

Proof: Complete induction will be employed. The result is trivially true for
one or two spaces involved. Suppose the result is true for $r$ spaces for
every $r<n$. Consider a grouped (perhaps nested) tensor product $\Pi$ of the
$n>2$ vector spaces $\mathcal{V}_{1},...,\mathcal{V}_{n}$ in that order. At
the outermost essential level of grouping, we have an expression of the form
$\Pi=\mathcal{W}_{1}\bigotimes\cdots\bigotimes\mathcal{W}_{k}$, where
$1<k\leqslant n$ and the $n$ factor spaces $\mathcal{V}_{1},...,\mathcal{V}%
_{n}$ are distributed in (perhaps nested) groupings among the $\mathcal{W}%
_{i}$. Then $\Pi\cong(\mathcal{W}_{1}\bigotimes\cdots\bigotimes\mathcal{W}%
_{k-1})\bigotimes\mathcal{W}_{k}$ and by the induction hypothesis we may
assume that, for some $i$ $<n$, $\mathcal{W}_{1}\bigotimes\cdots
\bigotimes\mathcal{W}_{k-1}\cong\mathcal{V}_{1}\bigotimes\cdots\bigotimes
\mathcal{V}_{i}$ and $\mathcal{W}_{k}\cong$ $\mathcal{V}_{i+1}\bigotimes
\cdots\bigotimes\mathcal{V}_{n}$. Hence $\Pi\cong\left(  \mathcal{V}%
_{1}\bigotimes\cdots\bigotimes\mathcal{V}_{i}\right)  \bigotimes\left(
\mathcal{V}_{i+1}\bigotimes\cdots\bigotimes\mathcal{V}_{n}\right)  $ and
therefore $\Pi\cong\mathcal{V}_{1}\bigotimes\cdots\bigotimes\mathcal{V}_{n}$
no matter how the $n$ factors $\mathcal{V}_{1},...,\mathcal{V}_{n}$ are
grouped. The theorem thus holds by induction. $\blacksquare$\medskip\ 

Not surprisingly, when we use the field $\mathcal{F}$ (considered as a vector
space over itself) as a tensor product factor, it does not really do anything.

\begin{theorem}
\label{Scalar}Let $\mathcal{V}$ be a vector space over the field $\mathcal{F}
$, and let $\mathcal{F}$ be considered also to be a vector space over itself.
Then there is an isomorphism from $\mathcal{F}\bigotimes\mathcal{V}$
\emph{(}and another from $\mathcal{V}\bigotimes\mathcal{F}$\emph{)} to
$\mathcal{V}$ sending $a\otimes v$ \emph{(}and the other sending $v\otimes
a$\emph{)} to $a\cdot v$.
\end{theorem}

Proof: The function from $\mathcal{F}\times\mathcal{V}$ to $\mathcal{V}$ that
sends $\left(  a,v\right)  $ to $a\cdot v$ is clearly bilinear. Hence there
exists a (unique) map $\Theta:\mathcal{F}\bigotimes\mathcal{V}\rightarrow
\mathcal{V}$ that sends $a\otimes v$ to $a\cdot v$. Let $\Theta^{^{\prime}%
}:\mathcal{V}\rightarrow\mathcal{F}\bigotimes\mathcal{V}$ be the map that
sends $v$ to $1\otimes v$. Then $\Theta\circ\Theta^{^{\prime}}$ is clearly the
identity on $\mathcal{V}$. On the other hand, $\Theta^{^{\prime}}\circ\Theta$
sends each element in $\mathcal{F}\bigotimes\mathcal{V}$ of the form $a\otimes
v$ to itself, and so is the identity on $\mathcal{F}\bigotimes\mathcal{V}$.

The other part is proved similarly. $\blacksquare$\newpage

\subsection{Problems}

\begin{enumerate}
\item Suppose that in $\mathcal{V}_{1}\bigotimes\mathcal{V}_{2}$ we have
$u\otimes v=u\otimes w$ for some nonzero $u\in\mathcal{V}_{1}$. Is it then
possible for the vectors $v$ and $w$ in $\mathcal{V}_{2}$ to be different?

\item When all $\mathcal{V}_{i}$ are finite-dimensional, $\mathcal{V}%
_{1}^{\top}%
{\textstyle\bigotimes}
\cdots%
{\textstyle\bigotimes}
\mathcal{V}_{n}^{\top}\cong(\mathcal{V}_{1}%
{\textstyle\bigotimes}
\cdots%
{\textstyle\bigotimes}
\mathcal{V}_{n})^{\top}$ via the unique map $\Phi:\mathcal{V}_{1}^{\top}%
{\textstyle\bigotimes}
\cdots%
{\textstyle\bigotimes}
\mathcal{V}_{n}^{\top}\rightarrow(\mathcal{V}_{1}%
{\textstyle\bigotimes}
\cdots%
{\textstyle\bigotimes}
\mathcal{V}_{n})^{\top}$ for which%
\[
\Phi\left(  \varphi_{1}\otimes\cdots\otimes\varphi_{n}\right)  \left(
v_{1}\otimes\cdots\otimes v_{n}\right)  =\varphi_{1}\left(  v_{1}\right)
\cdots\varphi_{n}\left(  v_{n}\right)  ,
\]
or, in other words,%
\[
\Phi\left(  \varphi_{1}\otimes\cdots\otimes\varphi_{n}\right)  =\varphi
_{1}\cdots\varphi_{n}.
\]
Thus, for example, if the $\mathcal{V}_{i}$ are all equal to $\mathcal{V}$,
and $\mathcal{V}$ has the basis $\mathcal{B}$, then $\Phi\left(  x_{1}^{\top
}\otimes\cdots\otimes x_{n}^{\top}\right)  =x_{1}^{\top}\cdots x_{n}^{\top
}=\left(  x_{1}\otimes\cdots\otimes x_{n}\right)  ^{\top}$, where each
$x_{i}\in\mathcal{B}$.
\end{enumerate}

\newpage

\section{Vector Algebras}

\subsection{The New Element: Vector Multiplication}

A vector space \emph{per se }has but one type of multiplication, namely
multiplication (scaling) of vectors by scalars. However, additional structure
may be added to a vector space by defining other multiplications. Defining a
vector multiplication such that any two vectors may be multiplied to yield
another vector, and such that this multiplication acts in a harmonious fashion
with respect to vector addition and scaling, gives a new kind of structure
when this vector multiplication is included as a structural element. The maps
of such an enhanced structure, that is, the functions that preserve the entire
structure, including the vector multiplication, are still vector space maps,
but some vector space maps may no longer qualify due to the requirement that
the vector multiplication must also be preserved. Adding a new structural
element, giving whatever benefits it may, also then may give new burdens in
the logical development of the theory involved if we are to exploit structural
aspects through use of function images. We will see, though, that at the cost
of a little added complication, some excellent constructs will be obtained as
a result of including various vector multiplications.

Thus we will say that the vector space $\mathcal{V}$ over the field
$\mathcal{F}$ becomes a \textbf{vector algebra} (\textbf{over} $\mathcal{F}$)
when there is defined on it a \textbf{vector multiplication }$\mu
:\mathcal{V}\times\mathcal{V}\rightarrow\mathcal{V}$\textbf{,} written
$\mu\left(  u,v\right)  =u\ast v,$ which is required to be bilinear so that it
satisfies both the distributive laws
\[
(t+u)\ast v=t\ast v+u\ast v\text{, and }u\ast(v+w)=u\ast v+u\ast w
\]
and also the law
\[
a\cdot\left(  u\ast v\right)  =(a\cdot u)\ast v=u\ast(a\cdot v)
\]
for all $a\in\mathcal{F}$ and all $t,u,v,w\in\mathcal{V}$. A vector algebra is
\textbf{associative} or \textbf{commutative} according as its vector
multiplication is associative or commutative. If there is a multiplicative
neutral element (a \textbf{unit element) }the vector algebra is called
\textbf{unital}. Since vector algebras are the only algebras treated here, we
will often omit the ``vector'' qualifier and simply refer to a vector algebra
as an \textbf{algebra}.

\begin{exercise}
$v\ast0=0\ast v=0$.
\end{exercise}

\begin{exercise}
In a unital algebra with unit element $1$, $(-1)\ast v=-v$.
\end{exercise}

A (\textbf{vector}) \textbf{subalgebra} of an algebra is a vector space
subspace that is closed under the vector multiplication of the algebra, so
that it is an algebra in its own right with a vector multiplication that has
the same effect as the vector multiplication of the original algebra, but is
actually its restriction to the subspace.

An \textbf{algebra map} (or sometimes simply \textbf{map} when the kind of map
is clear) is a function between algebras over the same field, which function
is a vector space map that also preserves vector multiplication. The image of
a map is a subalgebra of the codomain and is associative (resp. commutative)
when the domain algebra is. Under a map, the image of a unit element acts as a
unit element in the image algebra.

\begin{exercise}
The inverse of a bijective algebra map is an algebra map. \emph{(The
isomorphisms of algebras are thus the bijective algebra maps.)}
\end{exercise}

\subsection{Quotient Algebras}

\label{QuotientAlgebras}

The \textbf{kernel} of an algebra map is defined to be the kernel of the same
function viewed as a vector space map, namely those vectors of the domain that
map to the zero vector of the codomain. The kernel of an algebra map is a
subalgebra and more. It is an \textbf{ideal} in the algebra, i. e., a
subalgebra closed under multiplication on the left, and on the right, by all
the elements of the whole original algebra. Thus for an ideal $\mathcal{I}$ in
the algebra $\mathcal{V}$ we have $x*v\in\mathcal{I}$ and $v*x\in\mathcal{I}$
for all $x\in\mathcal{I}$ and for all $v\in\mathcal{V}$. Every ideal is a
subalgebra, but not every subalgebra is an ideal. Hence there can be
subalgebras that are not kernels of any algebra map. However every ideal does
turn out to be the kernel of some algebra map, as we shall soon see.

The level sets of an algebra map are the level sets of a vector space map, and
they may be made into a vector space, and indeed into a vector algebra, in the
same manner as was previously used in Section \ref{LevelSets}. This vector
algebra made by mirroring the image of an algebra map in the level sets is the
\textbf{quotient algebra} $\mathcal{V}/\mathcal{%
K%
} $ of the vector algebra $\mathcal{V}$ by the algebra map's kernel $\mathcal{%
K%
}$. The theme remains the same: a quotient ``is'' the image of a map, up to
isomorphism, and all maps with the same domain and the same kernel have
isomorphic images. In a quotient algebra, the product of cosets is
\textbf{modular,} as the following exercise describes.

\begin{exercise}
Let $\mathcal{V}$ be a vector algebra and let $\mathcal{%
K%
}$ be the kernel of an algebra map from $\mathcal{V}$. Then in $\mathcal{V}%
/\mathcal{%
K%
}$, $\left(  u+\mathcal{K}\right)  *\left(  v+\mathcal{K}\right)
=u*v+\mathcal{%
K%
}$ for all $u,v\in\mathcal{V}$.
\end{exercise}

The cosets of any ideal make up a vector space quotient. It is always possible
to impose the modular product on the cosets and thereby make the vector space
quotient into an algebra.

\begin{exercise}
Let $\mathcal{V}$ be a vector algebra over $\mathcal{F}$ and let $\mathcal{I}
$ be an ideal in $\mathcal{V}$. On the vector space $\mathcal{V}/\mathcal{I}$
, prescribing the modular product
\[
\left(  u+\mathcal{I}\right)  *\left(  v+\mathcal{I}\right)  =\left(
u*v\right)  +\mathcal{I}%
\]
for all $u,v\in\mathcal{V}$ gives a well-defined multiplication that makes
$\mathcal{V}/\mathcal{I}$ into a vector algebra over $\mathcal{F}$.
\end{exercise}

For any vector algebra $\mathcal{V}$ and ideal $\mathcal{I}\lhd\mathcal{V}$,
the natural projection function $p:\mathcal{V}\rightarrow\mathcal{V}%
/\mathcal{I}$ that sends $v$ to the coset $v+\mathcal{I}$ is a vector space
map with kernel $\mathcal{I}$, according to Proposition
\ref{NaturalProjection}. $p$ is also an algebra map if we employ the modular
product, for then $p\left(  u\right)  *p\left(  v\right)  =\left(
u+\mathcal{I}\right)  *\left(  v+\mathcal{I}\right)  =u*v+\mathcal{I}=p(u*v)$.
Thus we have the following result.

\begin{theorem}
\label{QuotientAlgebra}Let $\mathcal{V}$ be a vector algebra and let
$\mathcal{I}$ be an ideal in $\mathcal{V}$. Then the natural map
$p:\mathcal{V}\rightarrow\mathcal{V}/\mathcal{I}$ that sends $v$ to the coset
$v+\mathcal{I}$ is an algebra map with kernel $\mathcal{I}$, if the modular
product is used on $\mathcal{V}/\mathcal{I}$. Hence each ideal $\mathcal{I}$
in $\mathcal{V}$ is the kernel of some algebra map, and thus a quotient
algebra $\mathcal{V}/\mathcal{I}$ exists for every ideal $\mathcal{I}$ in
$\mathcal{V}$. $\blacksquare$
\end{theorem}

\subsection{Algebra Tensor Product}

\label{AlgebraTensorProduct}

For vector algebras $\mathcal{V}$ and $\mathcal{W}$ over the field
$\mathcal{F}$, the vector space tensor product $\mathcal{V}\bigotimes
\mathcal{W}$ may be made into an algebra by requiring
\[
\left(  v\otimes w\right)  \ast\left(  v\,^{\prime}\otimes w\,^{\prime
}\right)  =\left(  v\ast v\,^{\prime}\right)  \otimes(w\ast w\,^{\prime})
\]
and extending as a bilinear function on $\left(  \mathcal{V}\bigotimes
\mathcal{W}\right)  \times\left(  \mathcal{V}\bigotimes\mathcal{W}\right)  $
so that
\[
\left(  \sum_{i}v_{i}\otimes w_{i}\right)  \ast\left(  \sum_{j}v_{j}{}%
^{\prime}\otimes w_{j}{}^{\prime}\right)  =\sum_{i,j}\left(  v_{i}\ast v_{j}%
{}^{\prime}\right)  \otimes\left(  w_{i}\ast w_{j}{}^{\prime}\right)  \text{.}%
\]
This does give a well-defined vector multiplication for $\mathcal{V}%
\bigotimes\mathcal{W}$.

\begin{proposition}
Given vectors algebras $\mathcal{V}$ and $\mathcal{W}$ over the field
$\mathcal{F}$, there exits a bilinear function $\mu:(\mathcal{V}%
\bigotimes\mathcal{W})\times(\mathcal{V}\bigotimes\mathcal{W})\rightarrow
\mathcal{V}\bigotimes\mathcal{W}$ such that for all $v,v\,^{\prime}%
\in\mathcal{V}$ and all $w,w\,^{\prime}\in\mathcal{W}$, $\mu(v\otimes
w,v\,^{\prime}\otimes w\,^{\prime})=\left(  v*v\,^{\prime}\right)
\otimes\left(  w*w\,^{\prime}\right)  $.
\end{proposition}

Proof: Consider the function from $\mathcal{V}\times\mathcal{W}\times
\mathcal{V}\times\mathcal{W}$ to $\mathcal{V}\bigotimes\mathcal{W}$ that sends
$(v,w,v\,^{\prime},w\,^{\prime})$ to $\left(  v\ast v\,^{\prime}\right)
\otimes\left(  w\ast w\,^{\prime}\right)  $. This function is linear in each
variable separately and therefore there is a vector space map from
$\mathcal{V}\bigotimes\mathcal{W}\bigotimes\mathcal{V}\bigotimes\mathcal{W}$
to $\mathcal{V}\bigotimes\mathcal{W}$ that sends $v\otimes w\otimes
v\,^{\prime}\otimes w\,^{\prime}$ to $\left(  v\ast v\,^{\prime}\right)
\otimes\left(  w\ast w\,^{\prime}\right)  $. Since there is an isomorphism
from $(\mathcal{V}\bigotimes\mathcal{W})\bigotimes(\mathcal{V}\bigotimes
\mathcal{W})$ to $\mathcal{V}\bigotimes\mathcal{W}\bigotimes\mathcal{V}%
\bigotimes\mathcal{W}$ that sends $(v\otimes w)\otimes(v\,^{\prime}\otimes
w\,^{\prime})$ to $v\otimes w\otimes v\,^{\prime}\otimes w\,^{\prime}$, there
is therefore a vector space map from $(\mathcal{V}\bigotimes\mathcal{W}%
)\bigotimes(\mathcal{V}\bigotimes\mathcal{W})$ to $\mathcal{V}\bigotimes
\mathcal{W}$ that sends $(v\otimes w)\otimes(v\,^{\prime}\otimes w\,^{\prime
})$ to $\left(  v\ast v\,^{\prime}\right)  \otimes\left(  w\ast w\,^{\prime
}\right)  $, and corresponding to this vector space map is a bilinear function
$\mu:(\mathcal{V}\bigotimes\mathcal{W})\times(\mathcal{V}\bigotimes
\mathcal{W})\rightarrow\mathcal{V}\bigotimes\mathcal{W}$ that sends $(v\otimes
w,v\,^{\prime}\otimes w\,^{\prime})$ to $\left(  v\ast v\,^{\prime}\right)
\otimes\left(  w\ast w\,^{\prime}\right)  $. $\blacksquare$

\smallskip\ 

$\mathcal{V}\bigotimes\mathcal{W}$ with this vector multiplication function is
the \textbf{algebra tensor product} of the vector algebras $\mathcal{V}$ and
$\mathcal{W}$ over the field $\mathcal{F}$.

\begin{exercise}
If the vector algebras $\mathcal{V}$ and $\mathcal{W}$ over the field
$\mathcal{F}$ are both commutative, or both associative, or both unital, then
so is their algebra tensor product.
\end{exercise}

\subsection{The Tensor Algebras of a Vector Space}

The \textbf{contravariant} \textbf{tensor algebra} $\bigotimes\mathcal{V}$ of
a vector space $\mathcal{V}$ over the field $\mathcal{F}$ is formed from the
\textbf{tensor powers} $\bigotimes^{0}\mathcal{V}=\mathcal{F}$, and
$\bigotimes^{k}\mathcal{V}=\mathcal{V}\bigotimes\cdots\bigotimes\mathcal{V} $
(with $k$ factors $\mathcal{V}$, $k=1,2,\ldots$), by taking the weak product
$\biguplus_{k\geqslant0}\bigotimes^{k}\mathcal{V}$ and defining on it a vector
multiplication. To simplify notation, we will not distinguish $\bigotimes
^{k}\mathcal{V}$ from its alias $\xi_{k}\left(  \bigotimes\mathcal{V}\right)
=\eta_{k}\left(  \bigotimes^{k}\mathcal{V}\right)  $ (using the functions
introduced in section \ref{MapsOnProducts}), so that, for example,
$\bigotimes\mathcal{V}$ may also be written as the direct sum $\bigotimes
\mathcal{V}=\bigoplus_{k\geqslant0}\bigotimes^{k}\mathcal{V}$. For
$u,v\in\bigotimes\mathcal{V}$ such that $u\in\bigotimes^{i}\mathcal{V} $ and
$v\in\bigotimes^{j}\mathcal{V}$, $u\ast v$ will be defined as the element of
$\bigotimes^{i+j}\mathcal{V}$ that corresponds to $u\otimes v$ under the
isomorphism of Theorem \ref{Associative} or Theorem \ref{Scalar}. In general,
$u,v\in\bigotimes\mathcal{V}=$ $\bigoplus_{k\geqslant0}\bigotimes
^{k}\mathcal{V}$ can be written uniquely as finite sums $u=\sum u_{i}$ and
$v=\sum v_{j}$ where $u_{i}\in\bigotimes^{i}\mathcal{V}$ and $v_{j}%
\in\bigotimes^{j}\mathcal{V}$ and their product will then be defined as $u\ast
v=\sum w_{k}$ where $w_{k}=\sum_{i+j=k}u_{i}\ast v_{j}$. With this bilinear
vector multiplication, $\bigotimes\mathcal{V}$ is an associative, unital
algebra. For most $\mathcal{V}$, this algebra is not commutative, however. As
is customary for this algebra, we will henceforth use $\otimes$ instead of
$\ast$ as the vector multiplication symbol.

$\bigotimes\mathcal{V}^{\top}$is the \textbf{covariant tensor algebra} of a
vector space $\mathcal{V}$. Taking the algebra tensor product of the
contravariant algebra with the covariant algebra gives $\mathsf{T}%
\mathcal{V}=(\bigotimes\mathcal{V})\bigotimes(\bigotimes\mathcal{V}^{\top})$,
the \textbf{(full) tensor algebra} of $\mathcal{V}$. The elements of
$\mathsf{T}\mathcal{V}$ are called \textbf{tensor combinations}. In
$\mathsf{T}\mathcal{V}$, the algebra product of $r\in$ $\bigotimes\mathcal{V}$
and $s\in\bigotimes\mathcal{V}^{\top}$ is denoted $r\otimes s$ and when the
$r$ is in $\bigotimes^{p}\mathcal{V}$ and the $s$ is in $\bigotimes
^{q}\mathcal{V}^{\top}$, the product $r\otimes s$ is said to be a
\textbf{homogeneous tensor combination}, or simply a \textbf{tensor}, of
\textbf{contravariant degree} $p$, of \textbf{covariant degree} $q$, and of
\textbf{total degree} $p+q$. It is not hard to see that $\mathsf{T}%
\mathcal{V}=\bigoplus_{p,q\geqslant0}\mathsf{T}_{q}^{p}\mathcal{V}$ where
$\mathsf{T}_{q}^{p}\mathcal{V}=(\bigotimes^{p}\mathcal{V})\bigotimes
(\bigotimes^{q}\mathcal{V}^{\top})$.

The elements of the form $v_{1}\otimes\cdots\otimes v_{p}$ in $\bigotimes
^{p}\mathcal{V}$, $p>0$, where each factor is in $\bigotimes^{1}\mathcal{V}$,
are the \textbf{e-products} (\textbf{elementary products)} of \textbf{degree}
$p$ in $\bigotimes\mathcal{V}$. The nonzero elements of $\bigotimes
^{0}\mathcal{V} $ are the \textbf{e-products} of \textbf{degree} $0$.
$\bigotimes\mathcal{V} $ thus is the vector space of linear combinations of
e-products of various degrees, and $\mathsf{T}\mathcal{V}$ is the vector space
of linear combinations of \textbf{mixed e-products} $r\otimes s$ where the
various $r$ and $s$ are e-products of various degrees in $\ \bigotimes
\mathcal{V}$ and $\bigotimes\mathcal{V}^{\top}$, \emph{respectively}. This way
of defining the tensor algebra of $\mathcal{V}$ always keeps the elements of
$\bigotimes^{1}\mathcal{V}$ all together on the left, and the elements of
$\bigotimes^{1}\mathcal{V}^{\top}$all together on the right, in each mixed e-product.

$\bigotimes\mathcal{V}$ is a basic construct from which two important algebras
may be derived as quotients, as we do in the next two sections.

\subsection{The Exterior Algebra of a Vector Space}

Given an ideal $\mathcal{I}$ in an algebra $\mathcal{%
W%
}$, the quotient algebra $\mathcal{%
W%
}/\mathcal{I}$ is a kind of scaled-down image of $\mathcal{%
W%
}$ that suppresses $\mathcal{I}$ and glues together into one image element all
the elements of each coset $v+\mathcal{I}$ in $\mathcal{%
W%
}$, effectively treating the elements of $\mathcal{%
W%
}$ \textbf{modulo} $\mathcal{I}$. Elements of a particular type in $\mathcal{%
W%
}$ that can be expressed as the members of an ideal $\mathcal{I}$ in
$\mathcal{%
W%
}$ can be ``factored out'' by passing to the quotient $\mathcal{%
W%
}/\mathcal{I}$. We now proceed along these lines by defining in $\bigotimes
\mathcal{V}$ an ideal that when factored out will give us a very useful new algebra.

We say that the e-product $v_{1}\otimes\cdots\otimes v_{p}$ is
\textbf{dependent} if the sequence $v_{1},\ldots,v_{p}$ is dependent, that is
if the set $\left\{  v_{1},\ldots,v_{p}\right\}  $ is dependent \emph{or} if
there are any equal vectors among $v_{1},\ldots,v_{p}$. We similarly call the
tuple $\left(  v_{1},\ldots,v_{p}\right)  $, \textbf{dependent} under the same
circumstances. Such that are not dependent are \textbf{independent}. The set
of values of all linear combinations of the set of dependent e-products in
$\bigotimes\mathcal{V}$ will be denoted by $\mathcal{D}$. It is clear that
$\mathcal{D}$ is an ideal in $\bigotimes\mathcal{V}$. The quotient algebra
$\bigwedge\mathcal{V}=\left(  \bigotimes\mathcal{V}\right)  /\mathcal{D}$ is
the \textbf{exterior algebra} of $\mathcal{V}$. The \textbf{exterior e-product
} $v_{1}\wedge\cdots\wedge v_{p}$ is the image of the e-product $v_{1}%
\otimes\cdots\otimes v_{p}$ under the natural projection $\pi_{\wedge}:$
$\bigotimes\mathcal{V}\rightarrow\bigwedge\mathcal{V}$. $\bigwedge
^{p}\mathcal{V}=$ $\pi_{\wedge}\left(  \bigotimes^{p}\mathcal{V}\right)  $,
the $p$th \textbf{exterior power} of $\mathcal{V}$, is the subspace of
elements of \textbf{degree} $p$, and is spanned by all the exterior e-products
of degree $p$. Noting that $\bigwedge^{0}\mathcal{V}$ is an alias of
$\mathcal{F}$, and $\bigwedge^{1}\mathcal{V}$ is an alias of $\mathcal{V}$, we
will identify $\bigwedge^{0}\mathcal{V}$ with $\mathcal{F}$ and $\bigwedge
^{1}\mathcal{V}$ with $\mathcal{V}$. $\bigwedge\mathcal{V}$ is an associative
algebra because it is the algebra map image of an associative algebra. It is
customary to use $\wedge$ as the multiplication symbol in $\bigwedge
\mathcal{V}$.

Suppose that the set\emph{\ }$\left\{  v_{1},\ldots,v_{p}\right\}  $ of
nonzero vectors in $\mathcal{V}$ is dependent. Then for some $i$, $v_{i}%
=\sum_{j\neq i}a_{j}\cdot v_{j}$ and then
\[
v_{1}\otimes\cdots\otimes v_{p}=\sum_{j\neq i}a_{j}\cdot v_{1}\otimes
\cdots\otimes v_{i-1}\otimes v_{j}\otimes v_{i+1}\otimes\cdots\otimes v_{p}%
\]
so that such a dependent e-product is always a linear combination of
e-products each of which has equal vectors among its factors. Hence the ideal
$\mathcal{D}$ in $\bigotimes\mathcal{V}$ may be more primitively described as
the set of all linear combinations of e-products each of which has at least
one pair of equal vectors among its factors.

Let $u,v\in\bigotimes^{1}\mathcal{V}$. Then $0=\left(  u+v\right)
\wedge\left(  u+v\right)  =u\wedge u+u\wedge v+v\wedge u+v\wedge v=u\wedge
v+v\wedge u$ and therefore $u\wedge v=-\left(  v\wedge u\right)  $.

\begin{exercise}
Let $v_{1}\wedge\cdots\wedge v_{p}$ be an exterior e-product. Then
\[
v_{\sigma\left(  1\right)  }\wedge\cdots\wedge v_{\sigma(p)}=(-1)^{\sigma
}\cdot\left(  v_{1}\wedge\cdots\wedge v_{p}\right)
\]
where $(-1)^{\sigma}=+1$ or $-1$ according as the permutation $\sigma$ of
$\left\{  1,\ldots,p\right\}  $ is even or odd.
\end{exercise}

\begin{exercise}
Let $r\in\ \bigwedge^{p}\mathcal{V}$ and let $s\in\bigwedge^{q}\mathcal{V}$.
Then $s\wedge r=\left(  -1\right)  ^{pq}\cdot r\wedge s$.\emph{\ (The sign
changes only when both }$p$\emph{\ and }$q$\emph{\ are odd.)}
\end{exercise}

\begin{exercise}
Let $r\in\ \bigwedge^{p}\mathcal{V}$. Then it follows at once from the above
exercise that $r\wedge r=0$ if $p$ is odd and if $1+1\neq0$ in $\mathcal{F}$.
What about when $p$ is odd and $1+1=0$ in $\mathcal{F}$?
\end{exercise}

Let $\mathcal{B}$ be a basis for $\mathcal{V}$. Based on $\mathcal{B}$, the
set of \textbf{basis monomials} of degree $p$ (all the e-products of degree
$p$ with factors chosen only from $\mathcal{B}$), is a basis for
$\bigotimes^{p}\mathcal{V}$. Let $r$ denote the e-product $v_{1}\otimes
\cdots\otimes v_{p}$. Then $r$ has an expansion of the form
\[
r=\sum_{i_{1}}\cdots\sum_{i_{p}}a_{1,i_{1}}\cdots a_{p,i_{p}}\cdot x_{i_{1}%
}\otimes\cdots\otimes x_{i_{p}}%
\]
in terms of some $x_{1},\ldots,x_{N}\in\mathcal{B}$. Let us write $r=$
$r_{=}+r_{\neq}$, where $r_{=}$ is the sum of those terms that have some equal
subscripts and $r_{\neq}$ is the sum of the remaining terms that have no equal
subscripts. We have
\[
r_{\neq}=\sum_{i_{1}<\cdots<i_{p}}\sum_{\sigma\in\mathcal{S}_{p}%
}a_{1,i_{\sigma(1)}}\cdots a_{p,i_{\sigma(p)}}\cdot x_{i_{\sigma(1)}}%
\otimes\cdots\otimes x_{i_{\sigma(p)}}%
\]
where $\mathcal{S}_{p}$ denotes the set of all permutations of $\left\{
1,\ldots,p\right\}  $. In the event that two of the factors, say $v_{l}$ and
$v_{m}$, $l<m$, are equal in $r$, the terms in the summation above occur in
pairs with the same coefficient since
\[
a_{1,i_{\sigma(1)}}\cdots a_{p,i_{\sigma(p)}}=a_{1,i_{\sigma(1)}}\cdots
a_{l,i_{\sigma(m)}}\cdots a_{m,i_{\sigma(l)}}\cdots a_{p,i_{\sigma(p)}%
}\text{.}%
\]
When $v_{l}$ $=v_{m}$, $l<m$, we may thus write
\[
r_{\neq}=\sum_{i_{1}<\cdots<i_{p}}\sum_{\sigma\in\mathcal{A}_{p}%
}a_{1,i_{\sigma(1)}}\cdots a_{p,i_{\sigma(p)}}\cdot
\]
\[
\cdot(x_{i_{\sigma(1)}}\otimes\cdots\otimes x_{i_{\sigma(p)}}+x_{i_{\sigma
(1)}}\otimes\cdots\otimes x_{i_{\sigma(m)}}\otimes\cdots\otimes x_{i_{\sigma
(l)}}\otimes\cdots\otimes x_{i_{\sigma(p)}})
\]
where $\mathcal{A}_{p}$ is the set of all \emph{even} permutations of
$\left\{  1,\ldots,p\right\}  $. We say that the two e-products $v_{1}%
\otimes\cdots\otimes v_{p}$ and $v_{\sigma\left(  1\right)  }\otimes
\cdots\otimes v_{\sigma(p)}$ are of the \textbf{same parity}, or of
\textbf{opposite} \textbf{parity}, according as $\sigma$ is an even, or an
odd, permutation of $\left\{  1,\ldots,p\right\}  $. We see therefore that the
dependent basis monomials, along with the sums of pairs of independent basis
monomials of opposite parity, span $\mathcal{D}$.

Let $\mathcal{X}=\left\{  x_{1},\ldots,x_{p}\right\}  $ be a subset of the
basis $\mathcal{B}$ for $\mathcal{V}$. From the elements of $\mathcal{X}$,
$K=p!$ independent basis monomials may be formed by multiplying together the
elements of $\mathcal{X}$ in the various possible orders. Let $\mathcal{T}%
_{1}=\left\{  t_{1},t_{3,}\ldots,t_{K-1}\right\}  $ be the independent basis
monomials of degree $p$, with factors from $\mathcal{X}$, and of the same
parity as $t_{1}=x_{1}\otimes\cdots\otimes x_{p}$, and let $\mathcal{T}%
_{2}=\left\{  t_{2},t_{4},\ldots,t_{K}\right\}  $ be those of the opposite
parity. Then
\[
\mathcal{%
T%
}=\left\{  t_{1},t_{1}+t_{2},t_{2}+t_{3},\ldots,t_{K-1}+t_{K}\right\}
\]
is a set of $K$ independent elements with the same span as $\mathcal{T}%
_{1}\cup\mathcal{T}_{2}$. Moreover, we claim that for any $s\in\mathcal{T}%
_{1}$ and any $t\in\mathcal{T}_{2}$, $s+t$ is in the span of $\mathcal{%
T%
}\smallsetminus\left\{  t_{1}\right\}  $. It suffices to show this for $s+t$
of the form $t_{i}+t_{k}$ where $i<k$ and exactly one of $i,k$ is odd, or,
what is the same, for $s+t=$ $t_{i}+t_{i+2j+1}$. But $t_{i}+t_{i+2j+1}=\left(
t_{i}+t_{i+1}\right)  -\left(  t_{i+1}+t_{i+2}\right)  +\left(  t_{i+2}%
+t_{i+3}\right)  -\cdots+\left(  t_{i+2j}+t_{i+2j+1}\right)  $, verifying the
claim. The following result is now clear.

\begin{proposition}
Let $\mathcal{B}$ be a basis for $\mathcal{V}$. From each nonempty subset
$\mathcal{X}=\left\{  x_{1},\ldots,x_{p}\right\}  $ of $p$ elements of
$\mathcal{B}$, let the sets $\mathcal{T}_{1}=\left\{  t_{1},t_{3,}%
\ldots,t_{K-1}\right\}  $ and $\mathcal{T}_{2}=\left\{  t_{2},t_{4}%
,\ldots,t_{K}\right\}  $ comprising in total the $K=p!$ degree $p$ independent
basis monomials be formed, $\mathcal{T}_{1}$ being those that are of the same
parity as $t_{1}=x_{1}\otimes\cdots\otimes x_{p}$, and $\mathcal{T}_{2}$ being
those of parity opposite to $t_{1}$, and let
\[
\mathcal{%
T%
}=\left\{  t_{1},t_{1}+t_{2},t_{2}+t_{3},\ldots,t_{K-1}+t_{K}\right\}
\text{.}%
\]
Then $\left\langle \mathcal{T}\right\rangle =\left\langle \mathcal{T}_{1}%
\cup\mathcal{T}_{2}\right\rangle $ and if $s\in\mathcal{T}_{1}$ and
$t\in\mathcal{T}_{2}$, then $s+t\in\left\langle \mathcal{T}\smallsetminus
\left\{  t_{1}\right\}  \right\rangle $. Let $\mathcal{A}_{0}$ denote the set
of all dependent basis monomials based on $\mathcal{B}$, let $\mathcal{A}_{1}$
denote the union of all the sets $\mathcal{%
T%
}\smallsetminus\left\{  t_{1}\right\}  $ for all $p$, and let $\mathcal{E}$
denote the union of all the singleton sets $\left\{  t_{1}\right\}  $ for all
$p$, and $\left\{  1\right\}  $. Then $\mathcal{A}_{0}\cup\mathcal{A}_{1}%
\cup\mathcal{E}$ is a basis for $\bigotimes\mathcal{V}$, $\mathcal{A}_{0}%
\cup\mathcal{A}_{1}$ is a basis for the ideal $\mathcal{D}$ of all linear
combinations of dependent e-products in $\bigotimes\mathcal{V}$, and
$\mathcal{E}$ is a basis for a complementary subspace of $\mathcal{D}$.
$\blacksquare$
\end{proposition}

Inasmuch as any independent set of vectors is part of some basis, in the light
of Proposition \ref{Complementary} we then immediately infer the following
useful result which also assures us that in suppressing the dependent
e-products we have not suppressed any independent ones.

\begin{corollary}
\label{CriterionForIndependence}\emph{(\textbf{Criterion for Independence})
}The sequence $v_{1},\ldots,v_{p}$ of vectors of $\mathcal{V}$ is independent
if and only if $v_{1}\wedge\cdots\wedge v_{p}\neq0$. $\blacksquare$
\end{corollary}

\begin{exercise}
\label{FiniteDimensionalExterior}Let $\mathcal{V}$ be an $n$-dimensional
vector space with basis $\mathcal{B}=\left\{  x_{1},\ldots,x_{n}\right\}  $.
Then for $1\leqslant p\leqslant n$, the set $\left\{  x_{i_{1}}\wedge
\cdots\wedge x_{i_{p}}\,|\,i_{1}<\cdots<i_{p}\right\}  $ of $\binom{n}{p}$
elements is a basis for $\bigwedge^{p}\mathcal{V}$, while for $p>n$,
$\bigwedge^{p}\mathcal{V}=0$. In particular, then, the singleton set $\left\{
x_{1}\wedge\cdots\wedge x_{n}\right\}  $ is a basis for $\bigwedge
^{n}\mathcal{V}$.
\end{exercise}

We may restrict the domain of the natural projection $\pi_{\wedge}:$
$\bigotimes\mathcal{V}\rightarrow\bigwedge\mathcal{V}$ to $\bigotimes
^{p}\mathcal{V}$ and the codomain to $\bigwedge^{p}\mathcal{V}$ and thereby
obtain the vector space map $\pi_{\wedge^{p}}:\bigotimes^{p}\mathcal{V}%
\rightarrow\bigwedge^{p}\mathcal{V}$ which has kernel $\mathcal{D}%
_{p}=\mathcal{D}\cap\bigotimes^{p}\mathcal{V}$. It is easy to see that
$\bigwedge^{p}\mathcal{V}=\bigotimes^{p}\mathcal{V}/\mathcal{D}_{p}$ so that
$\pi_{\wedge^{p}}$ is the natural projection of $\bigwedge^{p}\mathcal{V}$
onto $\bigotimes^{p}\mathcal{V}/\mathcal{D}_{p}$. Let $\Upsilon_{p}%
:\mathcal{V}^{p}\rightarrow\bigotimes^{p}\mathcal{V}$ be a tensor product
function, and let $f$ be a $p$-linear function from $\mathcal{V}^{p}$ to a
vector space $\mathcal{W}$ over the same field as $\mathcal{V}$. Then by the
universal property of the tensor product, there is a unique vector space map
$f_{\otimes}:\bigotimes^{p}\mathcal{V}\rightarrow\mathcal{W}$ such that
$f=f_{\otimes}\circ\Upsilon_{p}$. Assume now that $f$ is zero on dependent
$p$-tuples, so that the corresponding $f_{\otimes}$ is zero on dependent
e-products of degree $p$. Applying Theorem \ref{InducedMap}, we infer the
existence of a unique map $f_{\wedge}:\bigwedge^{p}\mathcal{V}\rightarrow
\mathcal{W}$ such that $f_{\otimes}=f_{\wedge}\circ\pi_{\wedge^{p}}$. Thus
there is a universal property for the $p$th exterior power, which we now
officially record as a theorem. (An \textbf{alternating} $p$-linear function
is one that vanishes on dependent $p$-tuples.)

\begin{theorem}
\label{UnivExt}\emph{(\textbf{Universal Property of Exterior Powers})} Let
$\mathcal{V}$ and $\mathcal{W}$ be vector spaces over the same field, and let
$f:$ $\mathcal{V}^{p}\rightarrow\mathcal{W}$ be an alternating $p$-linear
function. Then there is a unique vector space map $f_{\wedge}:\bigwedge
^{p}\mathcal{V}\rightarrow\mathcal{W}$ such that%
\[
f_{\wedge}(v_{1}\wedge\cdots\wedge v_{p})=f(v_{1},\ldots,v_{p}).
\]
$\blacksquare$
\end{theorem}

\begin{exercise}
\label{MapExtPower}Let $f:\mathcal{V}\rightarrow\mathcal{W}$ be a vector space
map. Then there is a unique map $\bigwedge^{p}f:\bigwedge^{p}\mathcal{V}%
\rightarrow\bigwedge^{p}\mathcal{W}$ that satisfies%
\[%
{\textstyle\bigwedge\nolimits^{p}}
f\left(  v_{1}\wedge\cdots\wedge v_{p}\right)  =f\left(  v_{1}\right)
\wedge\cdots\wedge f\left(  v_{p}\right)  .
\]
\end{exercise}

The map $\bigwedge^{p}f$ of the above exercise is commonly called the $p$th
\textbf{exterior power} of $f$.

Let $f:\mathcal{V}\rightarrow$ $\mathcal{V}$ be a vector space map from the
$n$-dimensional vector space $\mathcal{V}$ to itself. Then by the above
exercise there is a unique map $\bigwedge^{n}f:\bigwedge^{n}\mathcal{V}%
\rightarrow\bigwedge^{n}\mathcal{V}$ such that $\bigwedge^{n}f\left(
v_{1}\wedge\cdots\wedge v_{n}\right)  =f\left(  v_{1}\right)  \wedge
\cdots\wedge f\left(  v_{n}\right)  $. Since $\bigwedge^{n}\mathcal{V}$ is
$1$-dimensional, $\bigwedge^{n}f\left(  t\right)  =a\cdot t$ for some scalar
$a=\det f$, the \textbf{determinant} of $f$. Note that the determinant of $f$
is independent of any basis choice. However, the determinant is only defined
for self-maps on finite-dimensional spaces.

\begin{theorem}
[Product Theorem]\label{ProductTheorem}Let $f,g$ be vector space maps on an
$n$-dimensional vector space. Then $\det g\circ f=\left(  \det g\right)
\left(  \det f\right)  $.
\end{theorem}

Proof: $\bigwedge^{n}\left(  g\circ f\right)  \left(  v_{1}\wedge\cdots\wedge
v_{n}\right)  =g\left(  f\left(  v_{1}\right)  \right)  \wedge\cdots\wedge
g\left(  f\left(  v_{n}\right)  \right)  =\left(  \det g\right)  \cdot
f\left(  v_{1}\right)  \wedge\cdots\wedge f\left(  v_{n}\right)  =\left(  \det
g\right)  \left(  \det f\right)  \cdot v_{1}\wedge\cdots\wedge v_{n}$.
$\blacksquare$

\begin{proposition}
\bigskip Let $f:\mathcal{V}\rightarrow\mathcal{W}$ be a vector space map from
the $n$-dimensional vector space $\mathcal{V}$ to its alias $\mathcal{W}.$ Let
$\left\{  x_{1},\ldots,x_{n}\right\}  $ be a basis for $\mathcal{V}$. Then $f$
is invertible if and only if $f\left(  x_{1})\wedge\cdots\wedge f(x_{n}%
\right)  \neq0$.
\end{proposition}

Proof: Suppose that $f\left(  x_{1}\right)  \wedge\cdots\wedge f\left(
x_{n}\right)  =0$. Then, by Corollary \ref{CriterionForIndependence}, the
sequence $f\left(  x_{1}\right)  ,\ldots,f\left(  x_{n}\right)  $ is
dependent, so there exist scalars $a_{1},\ldots,a_{n}$, not all zero, such
that
\[
0=a_{1}\cdot f\left(  x_{1}\right)  +\ldots+a_{n}\cdot f\left(  x_{n}\right)
=f\left(  a_{1}\cdot x_{1}+\ldots+a_{n}\cdot x_{n}\right)  \text{.}%
\]
Hence $f$ sends the nonzero vector $a_{1}\cdot x_{1}+\ldots+a_{n}\cdot x_{n}$
to $0$, and thus the kernel of $f$ is not $\left\{  0\right\}  $. Therefore
$f$ fails to be one-to-one and hence is not invertible.

On the other hand, if $f\left(  x_{1}\right)  \wedge\cdots\wedge f\left(
x_{n}\right)  \neq0$, the sequence $f\left(  x_{1}\right)  ,\ldots,f\left(
x_{n}\right)  $ is independent and makes up a basis of $\mathcal{W}$. Thus $f$
sends each nonzero vector in $\mathcal{V}$ to a nonzero vector in
$\mathcal{W}$, and thus the kernel of $f$ is $\left\{  0\right\}  $. Therefore
$f$ is one-to-one and hence invertible. $\blacksquare$

\smallskip

\begin{corollary}
Let $f:\mathcal{V}\rightarrow\mathcal{V}$ be a vector space map from the
$n$-dimensional vector space $\mathcal{V}$ to itself. Then $f$ is invertible
if and only if $\det f\neq0$. $\blacksquare$
\end{corollary}

Suppose the $n$-dimensional vector space aliases $\mathcal{V}$ and
$\mathcal{W}$ have the respective bases $\left\{  x_{1},\ldots,x_{n}\right\}
$, and $\left\{  y_{1},\ldots,y_{n}\right\}  $. The map $f:\mathcal{V}%
\rightarrow$ $\mathcal{W}$ then sends $x_{j}$ to $\sum_{i}a_{i,j}\cdot y_{i}$
for some scalars $a_{i,j}$. We have
\[
f\left(  x_{1}\right)  \wedge\cdots\wedge f\left(  x_{n}\right)  =\left(
\sum_{i}a_{i,1}\cdot y_{i}\right)  \wedge\cdots\wedge\left(  \sum_{i}%
a_{i,n}\cdot y_{i}\right)  =
\]%
\[
=\sum_{\sigma\in\mathcal{S}_{n}}a_{\sigma\left(  1\right)  ,1}\cdots
a_{\sigma\left(  n\right)  ,n}\cdot y_{\sigma\left(  1\right)  }\wedge
\cdots\wedge y_{\sigma\left(  n\right)  }=
\]%
\[
=\left(  \sum_{\sigma\in\mathcal{S}_{n}}\left(  -1\right)  ^{\sigma}%
a_{\sigma\left(  1\right)  ,1}\cdots a_{\sigma\left(  n\right)  ,n}\right)
\cdot y_{1}\wedge\cdots\wedge y_{n}%
\]
where $\mathcal{S}_{n}$ is the set of all permutations of $\{1,\ldots,n\}$ and
$(-1)^{\sigma}=+1$ or $-1$ according as the permutation $\sigma$ is even or
odd. Therefore $f$ is invertible if and only if%
\[
0\neq\sum_{\sigma\in\mathcal{S}_{n}}\left(  -1\right)  ^{\sigma}%
a_{\sigma\left(  1\right)  ,1}\cdots a_{\sigma\left(  n\right)  ,n}%
=\det\left[  a_{i,j}\right]  \text{,}%
\]
i. e., if and only if the familiar \emph{determinant} of the $n\times n$
\emph{matrix} $\left[  a_{i,j}\right]  $ is nonzero. We observe that when
$\mathcal{V}=\mathcal{W}$ and $x_{i}=y_{i}$ for all $i$, $\det\left[
a_{i,j}\right]  =\det f$. We record this result as follows.

\begin{proposition}
\label{Matrixdet}Let the $n$-dimensional vector space aliases $\mathcal{V}$
and $\mathcal{W}$ have the respective bases $\left\{  x_{1},\ldots
,x_{n}\right\}  $, and $\left\{  y_{1},\ldots,y_{n}\right\}  $. Let the map
$f:\mathcal{V}\rightarrow$ $\mathcal{W}$ send each $x_{j}$ to $\sum_{i}%
a_{i,j}\cdot y_{i}$. Then $f$ is invertible if and only if the determinant of
the matrix $\left[  a_{i,j}\right]  $, denoted $\det\left[  a_{i,j}\right]  $,
is nonzero. When $\mathcal{V}=\mathcal{W}$ and $x_{i}=y_{i} $ for all $i$,
$\det\left[  a_{i,j}\right]  =\det f$.
\end{proposition}

\begin{exercise}
The determinant of the $n\times n$ matrix $\left[  a_{i,j}\right]  $ is an
alternating $n$-linear functional of its columns \emph{(}and also of its rows,
due to the familiar result that a matrix and its \emph{transpose} have the
same determinant\emph{)}.
\end{exercise}

The universal property of exterior powers leads at once to the conclusion that
if $\mathcal{V}$ is over the field $\mathcal{F}$, the map space subspace
$\left\{  \mathcal{V}^{p}\overset{(alt\ p-lin)}{\rightarrow}\mathcal{F}%
\right\}  $ of alternating $p$-linear functionals is isomorphic to $\left(
\bigwedge^{p}\mathcal{V}\right)  ^{\top}$. There is also an analog of Theorem
\ref{DualofTensorProduct}, which we now start to develop by exploring the
coordinate functions on $\bigwedge^{p}\mathcal{V}$ relative to a basis
$\mathcal{B}$ for $\mathcal{V}$. Let $t=v_{1}\wedge\cdots\wedge v_{p}$ be an
exterior e-product in $\bigwedge^{p}\mathcal{V}$. Let $x_{1}\wedge\cdots\wedge
x_{p}$ be a typical basis vector of $\bigwedge^{p}\mathcal{V}$, formed from
elements $x_{i}\in\mathcal{B}$, and let $t$ be expanded in terms of such basis
vectors. Then the coordinate function $\left(  x_{1}\wedge\cdots\wedge
x_{p}\right)  ^{\top}$ corresponding to $x_{1}\wedge\cdots\wedge x_{p}$ is the
element of $\left(  \bigwedge^{p}\mathcal{V}\right)  ^{\top}$ that gives the
coefficient of $x_{1}\wedge\cdots\wedge x_{p}$ in this expansion of $t$. This
coefficient is readily seen to be $\det\left[  x_{i}^{\top}\left(
v_{j}\right)  \right]  $, where $\ x_{i}^{\top}$ is the coordinate function on
$\mathcal{V}$ that corresponds to $x_{i}\in\mathcal{B}$. Thus we are led to
consider what are easily seen to be alternating $p$-linear functionals
$f_{\phi_{1},\ldots,\phi_{p}}$ of the form $f_{\phi_{1},\ldots,\phi_{p}%
}\left(  v_{1},\ldots,v_{p}\right)  =\det\left[  \phi_{i}\left(  v_{j}\right)
\right]  $ where the $\phi_{i}$ are linear functionals on $\mathcal{V}$. By
the universal property, $f_{\phi_{1},\ldots,\phi_{p}}\left(  v_{1}%
,\ldots,v_{p}\right)  =\widehat{f}_{\phi_{1},\ldots,\phi_{p}}\left(
v_{1}\wedge\cdots\wedge v_{p}\right)  $ for a unique linear functional
$\widehat{f}_{\phi_{1},\ldots,\phi_{p}}\in\left(  \bigwedge^{p}\mathcal{V}%
\right)  ^{\top}$. Hence there is a function $\Phi:\left(  \mathcal{V}^{\top
}\right)  ^{p}\rightarrow\left(  \bigwedge^{p}\mathcal{V}\right)  ^{\top}$
such that $\Phi\left(  \phi_{1},\ldots,\phi_{p}\right)  =\widehat{f}_{\phi
_{1},\ldots,\phi_{p}}$. But $\Phi$ is clearly an alternating $p$-linear
function on $\left(  \mathcal{V}^{\top}\right)  ^{p}$ and by the universal
property there is a unique vector space map $\widehat{\Phi}:\bigwedge
^{p}\mathcal{V}^{\top}\rightarrow\left(  \bigwedge^{p}\mathcal{V}\right)
^{\top}$ such that $\widehat{\Phi}\left(  \phi_{1}\wedge\cdots\wedge\phi
_{p}\right)  =\Phi\left(  \phi_{1},\ldots,\phi_{p}\right)  $. Based on our
study of coordinate functions above, we see that $\widehat{\Phi}(x_{1}^{\top
}\wedge\cdots\wedge x_{p}^{\top})=\left(  x_{1}\wedge\cdots\wedge
x_{p}\right)  ^{\top}$. Therefore, for finite-dimensional $\mathcal{V}$,
$\widehat{\Phi}$ is an isomorphism because it sends a basis to a basis which
makes it onto, and it is therefore also one-to-one by Theorem \ref{Injmap}.
The following statement summarizes our results.

\begin{theorem}
\label{DualExtIso}Let $\mathcal{V}$ be a vector space over the field
$\mathcal{F}$. Then $\left(  \bigwedge^{p}\mathcal{V}\right)  ^{\top}%
\cong\left\{  \mathcal{V}^{p}\overset{(alt\ p-lin)}{\rightarrow}%
\mathcal{F}\right\}  $. If $\mathcal{V}$ is finite-dimensional, we also have
$\bigwedge^{p}\mathcal{V}^{\top}\cong\left(  \bigwedge^{p}\mathcal{V}\right)
^{\top}$ via the map $\widehat{\Phi}:\bigwedge^{p}\mathcal{V}^{\top
}\rightarrow\left(  \bigwedge^{p}\mathcal{V}\right)  ^{\top}$ for which
$\widehat{\Phi}\left(  \phi_{1}\wedge\cdots\wedge\phi_{p}\right)  \left(
v_{1}\wedge\cdots\wedge v_{p}\right)  =\det\left[  \phi_{i}\left(
v_{j}\right)  \right]  $. $\blacksquare$
\end{theorem}

\begin{exercise}
What is $\widehat{\Phi}$ in the case $p=1$?
\end{exercise}

Exterior algebra is an important subject that will receive further attention
in later chapters. Now we consider another example of a vector algebra that is
an algebra map image of the contravariant tensor algebra $\bigotimes
\mathcal{V}$.

\subsection{The Symmetric Algebra of a Vector Space}

Noncommutativity in $\bigotimes\mathcal{V}$ can be suppressed by passing to
the quotient $\mathsf{S}\mathcal{V}=\bigotimes\mathcal{V}/\mathcal{N}$ where
$\mathcal{N}$ is an ideal that expresses noncommutativity. The result is a
commutative algebra that is customarily called the \textbf{symmetric} algebra
of $\mathcal{V}$. The noncommutativity ideal $\mathcal{N}$ may be taken to be
the set of all linear combinations of differences of pairs of e-products of
the same degree $p$, $p\geqslant2$, which contain the same factors with the
same multiplicities, but in different orders. The effect of passing to the
quotient then is to identify all the e-products of the same degree which have
the same factors with the same multiplicities without regard to the order in
which these factors are being multiplied. The \textbf{symmetric e-product
}$v_{1}\cdots v_{p}$ is the image of the e-product $v_{1}\otimes\cdots\otimes
v_{p}$ under the natural projection $\pi_{\mathsf{S}}:$ $\bigotimes
\mathcal{V}\rightarrow\mathsf{S}\mathcal{V}$. $\mathsf{S}^{p}\mathcal{V}=$
$\pi_{\mathsf{S}}\left(  \bigotimes^{p}\mathcal{V}\right)  $, the $p$th
\textbf{symmetric power}$\mathsf{\ }$of $\mathcal{V}$, is the subspace of
elements of \textbf{degree} $p$. Noting that $\mathsf{S}^{0}\mathcal{V}$ is an
alias of $\mathcal{F}$, and $\mathsf{S}^{1}\mathcal{V} $ is an alias of
$\mathcal{V}$, we will identify $\mathsf{S}^{0}\mathcal{V}$ with $\mathcal{F}$
and $\mathsf{S}^{1}\mathcal{V}$ with $\mathcal{V}$. We will write products in
$\mathsf{S}\mathcal{V}$ with no sign to indicate the product operation, and we
will indicate repeated adjacent factors by the use of an exponent. Thus
$r_{3}r_{2}^{2}r_{1}$, $r_{2}r_{3}r_{1}r_{2}$, and $r_{1}r_{2}^{2}r_{3}$ are
typical (and equal) product expressions in $\mathsf{S}\mathcal{V}$.

Following a plan similar to that used with the exterior algebra, we will
produce a basis for the ideal $\mathcal{N}$ and a basis for a complementary
subspace of $\mathcal{N}$. Let $\mathcal{B}$ be a basis for $\mathcal{V}$.
Consider the difference of the pair of e-products $v_{1}\otimes\cdots\otimes
v_{p}$ and $v_{\sigma(1)}\otimes\cdots\otimes v_{\sigma(p)}$, where $\sigma$
is some permutation of $\left\{  1,\ldots,p\right\}  $. In terms of some
$x_{1},\ldots,x_{N}\in\mathcal{B}$ we have for this difference
\[
\sum_{i_{1}}\cdots\sum_{i_{p}}(a_{1,i_{1}}\cdots a_{p,i_{p}}\cdot x_{i_{1}%
}\otimes\cdots\otimes x_{i_{p}}-a_{\sigma(1),i_{\sigma(1)}}\cdots
a_{\sigma(p),i_{\sigma(p)}}\cdot x_{i_{\sigma(1)}}\otimes\cdots\otimes
x_{i_{\sigma(p)}})
\]
\[
=\sum_{i_{1}}\cdots\sum_{i_{p}}a_{1,i_{1}}\cdots a_{p,i_{p}}\cdot(x_{i_{1}%
}\otimes\cdots\otimes x_{i_{p}}-x_{i_{\sigma(1)}}\otimes\cdots\otimes
x_{i_{\sigma(p)}})\text{,}%
\]
where the multiplicative commutativity of the scalars has been exploited. Thus
$\mathcal{N}$ is spanned by the differences of pairs of basis monomials that
contain the same factors with the same multiplicities, but in different orders.

To each basis monomial $t=x_{1}\otimes\cdots\otimes x_{p}$, where the $x_{i}$
are (not necessarily distinct) elements of $\mathcal{B}$, there corresponds
the \textbf{multiset}
\[
\mathcal{M}=\left\{  (\xi_{1},\nu_{1}),\ldots,(\xi_{m},\nu_{m})\right\}
\]
where the $\xi_{j}$ are the distinct elements of $\left\{  x_{1},\ldots
,x_{p}\right\}  $ and $\nu_{j}$ is the multiplicity with which $\xi_{j}$
appears as a factor in $t$. Putting $\nu=\left(  \nu_{1},\ldots,\nu
_{m}\right)  $ and $\left|  \nu\right|  =\nu_{1}+\cdots+\nu_{m}$, we then have
$\left|  \nu\right|  =p$. Given a particular multiset $\mathcal{M}$, let
$\mathcal{T}_{0}=\left\{  t_{1},\ldots,t_{K}\right\}  $ be the set of all
basis monomials to which $\mathcal{M}$ corresponds. Here $K=\binom{\left|
\nu\right|  }\nu$ is the \emph{multinomial coefficient} given by\emph{\ }
$\binom{\left|  \nu\right|  }\nu=\left|  \nu\right|  !/\nu!$ where
$\nu!=\left(  \nu_{1}!\right)  \cdots\left(  \nu_{m}!\right)  $. From each
such $\mathcal{T}_{0}$ we may form the related set $\mathcal{%
T%
}=\left\{  t_{1},t_{1}-t_{2},\ldots,t_{K-1}-t_{K}\right\}  $ with the same
span. Now $t_{i}-t_{i+j}=(t_{i}-t_{i+1})+\cdots+(t_{i+j-1}-t_{i+j})$ so that
the difference of any two elements of $\mathcal{T}_{0}$ is in $\left\langle
\mathcal{%
T%
}\smallsetminus\left\{  t_{1}\right\}  \right\rangle $. We therefore have the
following result.

\begin{proposition}
Let $\mathcal{B}$ be a basis for $\mathcal{V}$. From each nonempty multiset
$\mathcal{M}=\left\{  (\xi_{1},\nu_{1}),\ldots,(\xi_{k},\nu_{m})\right\}  $
where the $\xi_{j}$ are distinct elements of $\mathcal{B}$, let the set
$\mathcal{T}_{0}=\left\{  t_{1},\ldots,t_{K}\right\}  $ of basis monomials
that correspond be formed, where $K=\binom{\left|  \nu\right|  }{\nu}$. Let
\[
\mathcal{%
T%
}=\left\{  t_{1},t_{1}-t_{2},t_{2}-t_{3},\ldots,t_{K-1}-t_{K}\right\}
\text{.}%
\]
Then $\left\langle \mathcal{T}\right\rangle =\left\langle \mathcal{T}%
_{0}\right\rangle $ and if $s,t$ $\in\mathcal{T}_{0}$, then $s-t\in
\left\langle \mathcal{T}\smallsetminus\left\{  t_{1}\right\}  \right\rangle $.
Let $\mathcal{A}$ denote the union of all the sets $\mathcal{%
T%
}\smallsetminus\left\{  t_{1}\right\}  $ for all $\nu$, and let $\mathcal{S}$
denote the union of all the singleton sets $\left\{  t_{1}\right\}  $ for all
$\nu$, and $\left\{  1\right\}  $. Then $\mathcal{A}\cup\mathcal{S}$ is a
basis for $\bigotimes\mathcal{V}$, $\mathcal{A}$ is a basis for the
noncommutativity ideal $\mathcal{N}$, and $\mathcal{S}$ is a basis for a
complementary subspace of $\mathcal{N}$. $\blacksquare$
\end{proposition}

\begin{exercise}
Let $\mathcal{V}$ be an $n$-dimensional vector space with basis $\mathcal{B}%
=\left\{  x_{1},\ldots,x_{n}\right\}  $. Then for $p>0$, the set $\left\{
x_{i_{1}}\cdots x_{i_{p}}\,|\,i_{1}\leqslant\cdots\leqslant i_{p}\right\}  $
of $\binom{n+p-1}{p}$ elements is a basis for $\mathsf{S}^{p}\mathcal{V}$.
\end{exercise}

\begin{exercise}
\emph{(\textbf{Universal Property of Symmetric Powers})} Let $\mathcal{V}$ and
$\mathcal{W}$ be vector spaces over the same field, and let $f:$
$\mathcal{V}^{p}\rightarrow\mathcal{W}$ be a \emph{\textbf{symmetric}%
}\textbf{\ }$p$-linear function, i. e., a $p$-linear function such that
$f(v_{1},\ldots,v_{p})=f\left(  v_{\sigma(1)},\ldots,v_{\sigma(p)}\right)  $
for every permutation $\sigma$ of $\left\{  1,\ldots,p\right\}  $. Then there
is a unique vector space map $f_{\mathsf{S}}:\mathsf{S}^{p}\mathcal{V}%
\rightarrow\mathcal{W}$ such that $f(v_{1},\ldots,v_{p})=f_{\mathsf{S}}\left(
v_{1}\cdots v_{p}\right)  $.
\end{exercise}

\begin{exercise}
$\left(  \mathsf{S}^{p}\mathcal{V}\right)  ^{\top}$ is isomorphic to the
vector space $\left\{  \mathcal{V}^{p}\overset{(sym\ p-lin)}{\rightarrow
}\mathcal{F}\right\}  $ of symmetric $p$-linear functionals.
\end{exercise}

\subsection{Null Pairs and Cancellation}

By a \textbf{null pair} in an algebra we mean a pair $u,v$ of nonzero
elements, not necessarily distinct, such that $u\ast v=0$. (We avoid the
rather misleading term \emph{zero divisor}, or the similar \emph{divisor of
zero},\emph{\ }usually used in describing this notion.) It is the same for an
algebra to have no null pair as it is for it to support cancellation of
nonzero factors, as the following proposition details.

\begin{proposition}
An algebra has no null pair if and only if it supports \textbf{cancellation}:
the equations $u\ast x=u\ast y$ and $x\ast v=y\ast v$, where $u,v$ are any
nonzero elements, each imply $x=y$.
\end{proposition}

Proof: Suppose that the algebra has no null pair. The equation $u\ast x=u\ast
y$ with $u\neq0$ implies that $u\ast(x-y)=0$ which then implies that $x-y=0$
since the algebra has no null pair. Similarly $x\ast v=y\ast v$ with $v\neq0$
implies $x-y=0$. Thus in each case, $x=y$ as was to be shown.

On the other hand, assume that cancellation is supported. Suppose that $u\ast
v=0$ with $u\neq0$. Then $u\ast v=u\ast0$, and canceling $u$ gives $v=0$.
Hence $u\ast v=0$ with $u\neq0$ implies that $v=0$ and it is therefore not
possible that $u\ast v=0$ with both $u$ and $v$ nonzero. $\blacksquare$

\smallskip An element $u$ is called \textbf{left regular} if $\ $there is no
nonzero element $v$ such that $u\ast v=0$, and similarly $v$ is called
\textbf{right regular} if there is no nonzero $u$ such that $u\ast v=0$. An
element is \textbf{regular} if it is both left regular and right regular.

\begin{exercise}
$u$ may be canceled from the left of $u\ast x=u\ast y$ if and only if it is
left regular, and $v$ may be canceled from the right of $x\ast v=y\ast v$ if
and only if it is right regular.
\end{exercise}

\subsection{Problems}

\begin{enumerate}
\item Give an example of an algebra map $f:A\rightarrow B$ where $A$ and $B$
are both unital with unit elements $1_{A}$ and $1_{B}$ respectively, but
$f\left(  1_{A}\right)  \neq1_{B}$.

\item Use exterior algebra to derive the \emph{Laplace expansion} of the
determinant of an $n\times n$ matrix.

\item Use exterior algebra to derive \emph{Cramer's rule} for solving a set of
$n$ linear equations in $n$ unknowns when the coefficient matrix has nonzero
determinant. Then derive the well-known formula for the inverse of the
coefficient matrix in terms of \emph{cofactors}.

\item If $\mathcal{V}$ is finite-dimensional over a field of
\emph{characteristic} 0, then $\mathsf{S}^{p}\mathcal{V}^{\top}\cong\left(
\mathsf{S}^{p}\mathcal{V}\right)  ^{\top}$ via the vector space map
$\widehat{\Psi}:\mathsf{S}^{p}\mathcal{V}^{\top}\rightarrow\left(
\mathsf{S}^{p}\mathcal{V}\right)  ^{\top}$ for which $\widehat{\Psi}\left(
\psi_{1}\cdots\psi_{p}\right)  \left(  r_{1}\cdots r_{p}\right)
=\operatorname*{per}\left[  \psi_{i}\left(  r_{j}\right)  \right]  $ where
$\operatorname*{per}\left[  a_{i,j}\right]  =\sum_{\sigma\in\mathcal{S}_{n}%
}a_{\sigma\left(  1\right)  ,1}\cdots a_{\sigma\left(  n\right)  ,n}$ is the
\emph{permanent} of the matrix $\left[  a_{i,j}\right]  .$

\item If $u\ast v=0$ in an algebra, must then $v\ast u=0$ also?

\item An algebra supports \textbf{left-cancellation} if $v\ast x=v\ast y$
implies $x=y$ whenever $v\neq0$, and \textbf{right-cancellation} is defined
similarly. An algebra supports left-cancellation if and only if it supports right-cancellation.
\end{enumerate}

\newpage

\section{Vector Affine Geometry}

\label{AffineGeo}

\subsection{Basing a Geometry on a Vector Space}

One may base a geometry on a vector space by defining the fundamental
geometric objects to be the members of some specified family of subsets of
vectors, and by declaring that two such objects are \emph{incident} if one is
included in the other by set inclusion. It is natural to take as fundamental
geometric objects a collection of ``flat'' sets, such as subspaces or cosets
of subspaces, and this is what we shall do here. By basing a geometry in the
same way on each vector space of a family of vector spaces over the same
field, we will, upon also specifying suitable structure-preserving maps, be
able to view the entire collection of fundamental geometric objects of all the
vector spaces of the family as a \emph{category} of geometric spaces having
the same kind of structure.

In this chapter we derive from an underlying vector space over a general field
a type of geometry called \emph{affine}, which embodies the more primitive
notions of what is no doubt the most familiar type of geometry, namely
Euclidean geometry. In affine geometry, just like in Euclidean geometry, there
is uniformity in the sense that the view from each point is the same. The
points all have equal standing, and any one of them may be considered to be
the origin. In addition to incidence, we also have parallelism as a
fundamental concept in affine geometry. Working, as we will, over a general,
possibly unordered, field, the availability of geometrically useful quantities
that can be expressed as field values will be quite limited, and even over the
real numbers, far fewer such quantities are available to us than in ordinary
Euclidean geometry. Nevertheless, many important concepts and results are
embraced within affine geometry, and it provides a suitable starting point for
our exploration of vector geometry.

\subsection{Affine Flats}

Let $\mathcal{V}$ be a vector space $\mathcal{V}$ over the field $\mathcal{F}
$. We denote by $\mathsf{V}\left(  \mathcal{V}\right)  $ the set of all
subspaces of $\mathcal{V}$. We utilize $\mathsf{V}\left(  \mathcal{V}\right)
$ to obtain an affine geometry by taking as fundamental geometric objects the
set $\mathcal{V}+\mathsf{V}\left(  \mathcal{V}\right)  $ of all \textbf{affine
flats} of $\mathcal{V}$, which are just the cosets of the vector space
subspaces of $\mathcal{V}$. Thus we introduce $\mathsf{A}\left(
\mathcal{V}\right)  =\mathcal{V}+\mathsf{V}\left(  \mathcal{V}\right)  $ as
the \textbf{affine structure} on $\mathcal{V}$. In $\mathsf{A}\left(
\mathcal{V}\right)  $, the \textbf{points} are the cosets of $\left\{
0\right\}  $ (i. e., the singleton subsets of $\mathcal{V}$, or essentially
the vectors of $\mathcal{V}$), the \textbf{lines} are the cosets of the
one-dimensional subspaces, the \textbf{planes} are the cosets of the
two-dimensional subspaces, and the $n$\textbf{-planes} are the cosets of the
$n$-dimensional subspaces. The \textbf{hyperplanes} are the cosets of the
subspaces of codimension 1. The \textbf{dimension} of an affine flat is the
dimension of the subspace of which it is a coset. If $\Phi\in\mathsf{A}\left(
\mathcal{V}\right)  $ is given by $\Phi=v+\mathcal{W}$ where $\mathcal{W}%
\in\mathsf{V}\left(  \mathcal{V}\right)  $, we call $\mathcal{W}$ the
\textbf{directional subspace,} or \textbf{underlying subspace}, of $\Phi$. The
notation $\Phi^{\mathsf{v}}$ will be used to signify the unique directional
subspace of the affine flat $\Phi$.

\subsection{Parallelism in Affine Structures}

We say that two affine flats, distinct or not, in the same affine structure
are \textbf{parallel} if their directional subspaces are incident (one is
contained entirely within the other). They are \textbf{strictly parallel} if
their directional subspaces are equal. They are \textbf{nontrivially parallel}
if they are parallel and neither is a point. They are \textbf{properly
parallel} if they are parallel and unequal. Any of these terms applies to a
set of affine flats if it applies pairwise within the set. Thus, for example,
to say that $\mathcal{S}$ is a set of strictly parallel affine flats means
that for every $\Phi,\Psi\in\mathcal{S}$, $\Phi$ and $\Psi$ are strictly parallel.

\begin{exercise}
Two strictly parallel affine flats that intersect must be equal. Hence, two
parallel affine flats that intersect must be incident.
\end{exercise}

\begin{exercise}
The relation of being strictly parallel is an equivalence relation. However,
in general, the relation of being parallel is not transitive.
\end{exercise}

\newpage\ 

\subsection{Translates and Other Expressions}

For vectors $v,w,x$ and scalar $a,$ a \textbf{translate} of $v$ (in the
direction of $x-w$) is an expression of the form $v+a\cdot(x-w)$. This is a
fundamental affine concept. If $x,w\in\mathcal{S}$, we call $v+a\cdot(x-w)$ an
$\mathcal{S}$\textbf{-translate} of $v$.

\begin{proposition}
The nonempty subset $\mathcal{S}$ of the vector space $\mathcal{V}$ is an
affine flat if and only if it contains the value of each $\mathcal{S}%
$-translate of each of its elements.
\end{proposition}

Proof: Let $\mathcal{S}\subset\mathcal{V}$ be such that it contains each
$\mathcal{S}$-translate of each of its elements and let $v\in\mathcal{S}$.
Then we wish to show that $\mathcal{S}=v+\mathcal{W}$ where $\mathcal{W}%
\lhd\mathcal{V}$. This is the same as showing that 1) $v+a\cdot\left(
w-v\right)  $ $\in\mathcal{S}$ for all $w\in\mathcal{S}$ and all scalars $a$,
and that 2) $v+(\left(  w-v\right)  +\left(  x-v\right)  )\in\mathcal{S}$ for
all $w,x\in\mathcal{S}$. The expression in 1) is clearly an $\mathcal{S}%
$-translate of $v$ and hence is in $\mathcal{S}$ as required. The expression
in 2) may be written as $v+\left(  \left(  w+\left(  x-v\right)  \right)
-v\right)  $ and is seen to be in $\mathcal{S}$ because $z=w+\left(
x-v\right)  $ is an $\mathcal{S}$-translate of $w$ and is therefore in
$\mathcal{S}$ so that our expression $v+\left(  z-v\right)  $ is also an
$\mathcal{S}$-translate of $v$ and must be in $\mathcal{S}$.

On the other hand, suppose that $\mathcal{S}$ is an affine flat in
$\mathcal{V}$ and let $v\in\mathcal{S}$. Let $x=v+x^{\mathsf{v}}$ and
$w=v+w^{\mathsf{v}}$ where $x^{\mathsf{v}},w^{\mathsf{v}}\in\mathcal{S}%
^{\mathsf{v}}$. Then $a\cdot\left(  x-w\right)  =a\cdot\left(  x^{\mathsf{v}%
}-w^{\mathsf{v}}\right)  \in\mathcal{S}^{\mathsf{v}}$ and hence $v+a\cdot
(x-w)\in\mathcal{S}$. $\blacksquare$

\medskip\ 

\begin{exercise}
The intersection of affine flats is an affine flat or empty.
\end{exercise}

Iterating the process of forming translates gives the
\textbf{multi-translate}
\[
v+b_{1}\cdot\left(  x_{1}-w_{1}\right)  +\cdots+b_{n}\cdot\left(  x_{n}%
-w_{n}\right)
\]
of the vector $v$ in the directions of the $x_{i}-w_{i}$ based on the vectors
$w_{i},x_{i}$ and the scalars $b_{i}$.

An expression of the form $a_{1}\cdot v_{1}+\cdots+a_{n}\cdot v_{n}$ is a
\textbf{linear expression} in the vectors $v_{1},\ldots,v_{n}$, and
$a_{1}+\cdots+a_{n}$ is its \textbf{weight}. Two important cases are specially
designated: A linear expression of weight $1$ is an \textbf{affine
expression}, and a linear expression of weight $0$ is a \textbf{directional
expression}. If we expand out each of the terms $b_{i}\cdot\left(  x_{i}%
-w_{i}\right)  $ as $b_{i}\cdot x_{i}+\left(  -b_{i}\right)  \cdot w_{i}$, the
multi-translate above becomes an affine expression in the $w_{i}$, the $x_{i}%
$, and $v$, and it becomes $v$ plus a directional expression in the $w_{i}$
and the $x_{i}$. On the other hand, since
\[
a_{1}\cdot v_{1}+\cdots+a_{n}\cdot v_{n}=\left(  a_{1}+\cdots+a_{n}\right)
\cdot v_{1}+a_{2}\cdot\left(  v_{2}-v_{1}\right)  +\cdots+a_{n}\cdot\left(
v_{n}-v_{1}\right)  \text{,}%
\]
our typical affine expression may be written as the multi-translate
\[
v_{1}+a_{2}\cdot\left(  v_{2}-v_{1}\right)  +\cdots+a_{n}\cdot\left(
v_{n}-v_{1}\right)  \text{,}%
\]
and our typical directional expression may be written as the linear
expression
\[
a_{2}\cdot\left(  v_{2}-v_{1}\right)  +\cdots+a_{n}\cdot\left(  v_{n}%
-v_{1}\right)
\]
in the differences $\left(  v_{2}-v_{1}\right)  ,\ldots,\left(  v_{n}%
-v_{1}\right)  $.

\begin{exercise}
\label{AffineDependence}A directional expression in the vectors $v_{1}%
,\ldots,v_{n}$ is nontrivially equal to $0$ \emph{(}some scalar coefficient
not $0$\emph{)} if and only if some $v_{j}$ is equal to an affine expression
in the other $v_{i}$.
\end{exercise}

\subsection{Affine Span and Affine Sum}

An \textbf{affine combination} of a finite nonempty set of vectors is any
linear combination of the set such that the coefficients sum to $1$, and thus
is an affine expression equal to a multi-translate of any one of the vectors
based on all of them. If $\mathcal{X}$ is a nonempty set of vectors, its
\textbf{affine span} $\left\langle \mathcal{X}\right\rangle _{\mathsf{A}} $ is
the set of all affine combinations of all its finite nonempty subsets. We say
that $\mathcal{X}$ \textbf{affine-spans} $\left\langle \mathcal{X}%
\right\rangle _{\mathsf{A}}$, and is an \textbf{affine spanning set} for
$\left\langle \mathcal{X}\right\rangle _{\mathsf{A}}$.

\begin{exercise}
Let $\mathcal{X}$ be a nonempty subset of the vector space $\mathcal{V}$. Then
$\left\langle \mathcal{X}\right\rangle _{\mathsf{A}}$ is the smallest affine
flat containing $\mathcal{X}$, i. e., it is the intersection of all affine
flats containing $\mathcal{X}$.
\end{exercise}

Let $\Phi,\Psi\in\mathsf{A}\left(  \mathcal{V}\right)  $. Then $\Phi
+_{\mathsf{A}}\Psi=\left\langle \Phi\cup\Psi\right\rangle _{\mathsf{A}}$ is
their \textbf{affine} \textbf{sum,} and similarly for more summands. With
partial order $\subset$, meet $\cap$, and join $+_{\mathsf{A}}$,
$\mathsf{A}\left(  \mathcal{V}\right)  \cup\left\{  \O\right\}  $ is a
\emph{complete lattice} which is \emph{upper semi-modular} but in general is
neither \emph{modular} nor \emph{distributive} (see, e. g., Hall's \emph{The
Theory of Groups} for definitions and characterizations).

\begin{exercise}
Let $\Phi,\Psi\in\mathsf{A}\left(  \mathcal{V}\right)  $. Suppose first that
$\Phi\cap\Psi\neq\O$. Then $\left(  \Phi+_{\mathsf{A}}\Psi\right)
^{\mathsf{v}}=\Phi^{\mathsf{v}}+\Psi^{\mathsf{v}}$ and $\left(  \Phi\cap
\Psi\right)  ^{\mathsf{v}}=\Phi^{\mathsf{v}}\cap\Psi^{\mathsf{v}}$, so that
\[
\dim\left(  \Phi+_{\mathsf{A}}\Psi\right)  +\dim\left(  \Phi\cap\Psi\right)
=\dim\Phi+\dim\Psi\text{.}%
\]
On the other hand, supposing that $\Phi\cap\Psi=\O$, we have
\[
\dim\left(  \Phi+_{\mathsf{A}}\Psi\right)  =\dim\Phi+\dim\Psi-\dim\left(
\Phi^{\mathsf{v}}\cap\Psi^{\mathsf{v}}\right)  +1\text{.}%
\]
\emph{(}Assigning $\O$ the standard dimension of $-1$, the formula above may
be written
\[
\dim\left(  \Phi+_{\mathsf{A}}\Psi\right)  +\dim\left(  \Phi\cap\Psi\right)
=\dim\Phi+\dim\Psi-\dim\left(  \Phi^{\mathsf{v}}\cap\Psi^{\mathsf{v}}\right)
\text{.}%
\]
Hence, no matter whether $\Phi$ and $\Psi$ intersect or not, we have
\[
\dim\left(  \Phi+_{\mathsf{A}}\Psi\right)  +\dim\left(  \Phi\cap\Psi\right)
\leqslant\dim\Phi+\dim\Psi
\]
which is a characterizing inequality for upper semi-modularity.\emph{)} Note
also that we have $\dim\left(  \Phi+_{\mathsf{A}}\Psi\right)  -\dim\left(
\Phi^{\mathsf{v}}+\Psi^{\mathsf{v}}\right)  =0$ or $1$ according as $\Phi$ and
$\Psi$ intersect or not.
\end{exercise}

\subsection{Affine Independence and Affine Frames}

Using an affine structure on a vector space allows any vector to serve as
origin. The idea of a dependency in a set of vectors needs to be extended in
this light. Here the empty set plays no significant r\^ole, and it is tacitly
assumed henceforth that only nonempty sets and nonempty subsets are being
addressed relative to these considerations. A \textbf{directional combination}
of a finite (nonempty) set of vectors is any linear combination of the set
such that the coefficients sum to $0$. An \textbf{affine dependency} is said
to exist in a set of vectors if the zero vector is a nontrivial (scalar
coefficients not all zero) directional combination of one of its finite
subsets (see also Exercise \ref{AffineDependence}). A set in which an affine
dependency exists is \textbf{affine dependent} and otherwise is \textbf{affine
independent}. A point is affine independent.

\begin{exercise}
If $\mathcal{X}$ is a nonempty set of vectors and $w\in\mathcal{X}$, then
$\mathcal{X}$ is \textbf{\emph{independent relative to}} $w$ if $\left\{
v-w\,|\,v\in\mathcal{X},v\neq w\right\}  $ is an independent set of vectors
\emph{(}in the original sense\emph{)}. A nonempty set of vectors is affine
independent if and only if it is independent relative to an arbitrarily
selected one of its elements.
\end{exercise}

\begin{exercise}
In an $n$-dimensional vector space, the number $m$ of elements in an affine
independent subset must satisfy $1\leqslant m\leqslant n+1$ and for each such
value of $m$ there exists an affine independent subset having that value of
$m$ as its number of elements.
\end{exercise}

An \textbf{affine frame} for the affine flat $\Phi$ is an affine independent
subset that affine-spans the flat.

\begin{exercise}
An affine frame for the affine flat $\Phi$ is a subset of the form $\left\{
v\right\}  \cup\left(  v+\mathcal{B}\right)  $ where $\mathcal{B}$ is a basis
for the directional subspace $\Phi^{\mathsf{v}}$.
\end{exercise}

\begin{exercise}
Let $\mathcal{A}$ be an affine frame for the affine flat $\Phi.$ Then each
vector of $\Phi$ has a unique expression as an affine combination $\sum
_{x\in\mathcal{X}}a_{x}\cdot x$ where $\mathcal{X}$ is some finite subset of
$\mathcal{A}$ and all of the scalars $a_{x}$ are nonzero.
\end{exercise}

\begin{exercise}
Let the affine flat $\Phi$ have the finite affine frame $\mathcal{X}$. Then
each vector of $\Phi$ has a unique expression as an affine combination
$\sum_{x\in\mathcal{X}}a_{x}\cdot x$. The scalars $a_{x}$ are the
\emph{\textbf{barycentric coordinates}} of the vector relative to the frame
$\mathcal{X}$.
\end{exercise}

\subsection{Affine Maps}

Let $\mathcal{V}$ and $\mathcal{W}$ be vector spaces over the same field,
which we will always understand to be the case whenever there are to be maps
between their affine flats. For $\Phi\in\mathsf{A}\left(  \mathcal{V}\right)
$ and $\Psi\in\mathsf{A}\left(  \mathcal{W}\right)  $, we call the function
$\mathcal{\alpha}:\Phi\rightarrow\Psi$ an \textbf{affine map} if it preserves
translates:
\[
\mathcal{\alpha}(v+a\cdot(x-w))=\mathcal{\alpha}(v)+a\cdot(\mathcal{\alpha
}(x)-\mathcal{\alpha}(w))
\]
for all $v,w,x\in\Phi$ and all scalars $a$.

In a series of exercises we now point out a number of important results
concerning affine maps. We lead off with two basic results.

\begin{exercise}
The composite of affine maps is an affine map.
\end{exercise}

\begin{exercise}
A vector space map is an affine map.
\end{exercise}

If $\Phi\,^{\prime},\Phi$ are elements of the same affine structure and
$\Phi\,^{\prime}\subset\Phi$, we say that $\Phi\,^{\prime}$ is a
\textbf{subflat} of $\Phi$. Thus two elements of the same affine structure are
incident whenever one is a subflat of the other. The image of an affine map
should be an affine flat, if we are even to begin to say that affine structure
is preserved. This is true, as the following exercise points out, and as a
consequence subflats are also sent to subflats with the result that incidence
is completely preserved.

\begin{exercise}
The image of an affine flat under an affine map is an affine flat. Since the
restriction of an affine map to a subflat is clearly an affine map, an affine
map therefore maps subflats to subflats and thereby preserves incidence.
\end{exercise}

Affine expressions are preserved as well, and the preservation of affine spans
and affine sums then follows.

\begin{exercise}
An affine map preserves all affine expressions of vectors of its domain. Hence
under an affine map, the image of the affine span of a subset of the domain is
the affine span of its image, and affine maps preserve affine sums.
\end{exercise}

A one-to-one function clearly preserves all intersections. However, in general
it can only be said that the intersection of subset images contains the image
of the intersection of the subsets. Thus it should not be too suprising that
affine maps do not necessarily even preserve the intersection of flats as the
following exercise demonstrates.

\begin{exercise}
As usual, let $\mathcal{F}$ be a field. Let $\alpha:\mathcal{F}^{2}%
\rightarrow\mathcal{F}^{2}$ be the vector space map that sends $\left(
a,b\right)  $ to $\left(  0,a+b\right)  .$ Let $\Phi=\left\{  \left(
a,0\right)  \mid a\in\mathcal{F}\right\}  $ and $\Psi=\left\{  \left(
0,a\right)  \mid a\in\mathcal{F}\right\}  $. Then $\alpha\left(  \Phi\cap
\Psi\right)  \neq\alpha\left(  \Phi\right)  \cap\alpha\left(  \Psi\right)  $.
\end{exercise}

A basic result about vector space maps generalizes to affine maps.

\begin{exercise}
The inverse of a bijective affine map is an affine map.
\end{exercise}

The bijective affine maps are therefore isomorphisms, so that $\Phi
\in\mathsf{A}\left(  \mathcal{V}\right)  $ and $\Psi\in\mathsf{A}\left(
\mathcal{W}\right)  $ are isomorphic (as affine flats) when there is a
bijective affine map from one to the other. Isomorphic objects, such as
isomorphic affine flats, are always \textbf{aliases} of each other, with
context serving to clarify the type of object.

Affine maps and vector space maps are closely related, as the following
theorem details.

\begin{theorem}
Let $\alpha$ be an affine map from $\Phi\in\mathsf{A}\left(  \mathcal{V}%
\right)  $ to $\Psi\in\mathsf{A}\left(  \mathcal{W}\right)  $and fix an
element $v\in\Phi$. Let $f:\Phi^{\mathsf{v}}\rightarrow\Psi^{\mathsf{v}} $ be
the function defined by
\[
f\left(  h\right)  =\alpha\left(  v+h\right)  -\alpha\left(  v\right)
\]
for all $h\in\Phi^{\mathsf{v}}$. Then $f$ is a vector space map that is
independent of the choice of $v$.
\end{theorem}

Proof:  Let $x,w\in\Phi$ be arbitrary, so that $h=x-w$ is an arbitrary element
of $\Phi^{\mathsf{v}}$. Then
\begin{align*}
f\left(  a\cdot h\right)   &  =\alpha\left(  v+a\cdot\left(  x-w\right)
\right)  -\alpha\left(  v\right)  =a\cdot\alpha\left(  x\right)  -a\cdot
\alpha\left(  w\right)  =\\
&  =a\cdot\left(  \alpha\left(  v\right)  +\alpha\left(  x\right)
-\alpha\left(  w\right)  \right)  -a\cdot\alpha\left(  v\right)  =\\
&  =a\cdot\left(  \alpha\left(  v+h\right)  -\alpha\left(  v\right)  \right)
=a\cdot f\left(  h\right)  \text{.}%
\end{align*}
Similarly, we find that $f\left(  h+k\right)  =f\left(  h\right)  +f\left(
k\right)  $ for arbitrary $h,k\in\Phi^{\mathsf{v}}$. Hence $f$ is a vector
space map. Replacing $h$ with $x-w$, we find that $f\left(  h\right)
=\alpha\left(  v+x-w\right)  -\alpha\left(  v\right)  =\alpha\left(  x\right)
-\alpha\left(  w\right)  $, showing that $f\left(  h\right)  $ is actually
independent of $v$. $\blacksquare$

\medskip\ 

The vector space map $f$ that corresponds to the affine map $\alpha$ as in the
theorem above is the \textbf{underlying vector space map} of $\alpha$, which
we will denote by $\alpha^{\mathsf{v}}$. Thus for any affine map $\alpha$ we
have
\[
\alpha\left(  v\right)  =\alpha^{\mathsf{v}}\left(  v-u\right)  +\alpha\left(
u\right)
\]
for all $u,v$ in the domain of $\alpha$.

\begin{exercise}
Let $\Phi\in\mathsf{A}\left(  \mathcal{V}\right)  $ and $\Psi\in
\mathsf{A}\left(  \mathcal{W}\right)  $ and let $f:\Phi^{\mathsf{v}%
}\rightarrow\Psi^{\mathsf{v}}$ be a given vector space map. Fix an element
$v\in\Phi$ and let $\alpha:\Phi\rightarrow\Psi$ be the function such that for
each $h\in\Phi^{\mathsf{v}}$, $\alpha\left(  v+h\right)  =\alpha\left(
v\right)  +f\left(  h\right)  $. Then $\alpha$ is an affine map such that
$\alpha^{\mathsf{v}}=f$.
\end{exercise}

Other nice properties of affine maps are now quite apparent, as the following
exercises detail.

\begin{exercise}
Under an affine map, the images of strictly parallel affine flats are strictly
parallel affine flats.
\end{exercise}

\begin{exercise}
An affine map is completely determined by its values on any affine spanning
set, and if that affine spanning set is an affine frame, the values may be
arbitrarily assigned.
\end{exercise}

\begin{exercise}
Let $\alpha:\Phi\rightarrow\Psi$ be an affine map from $\Phi\in\mathsf{A}%
\left(  \mathcal{V}\right)  $ to $\Psi\in\mathsf{A}\left(  \mathcal{W}\right)
$. Then there exists the affine map $\alpha^{\sharp}:\mathcal{V}%
\rightarrow\Psi$ that agrees with $\alpha$ on $\Phi$.
\end{exercise}

\subsection{Some Affine Self-Maps}

An affine self-map $\delta_{t,u;b}:\Phi\rightarrow\Phi$ of the flat $\Phi$
that sends $v\in\Phi$ to
\[
\delta_{t,u;b}\left(  v\right)  =t+b\cdot\left(  v-u\right)  \text{,}%
\]
where $t,u$ are fixed in $\Phi$ and $b$ is a fixed scalar, is known as a
\textbf{dilation} (or \textbf{dilatation}, a term we will not use, but which
means the same). Every \textbf{proper }dilation ($b\neq0$) is invertible:
$\delta_{t,u;b}^{-1}=\delta_{u,t;b^{-1}}$. The following result tells us that
the proper dilations are the direction-preserving automorphisms of affine flats.

\begin{proposition}
A function $\delta:\Phi\rightarrow\Phi$ on the affine flat $\Phi$ is a
dilation $\delta_{t,u;b}$ if and only if $\delta\left(  v\right)
-\delta\left(  w\right)  =b\cdot\left(  v-w\right)  $ for all $v,w\in\Phi$.
\end{proposition}

Proof: A dilation clearly satisfies the condition. On the other hand, suppose
that the condition holds. Fix $u,w\in\Phi$ and define $t\in\Phi$ by
$t=\delta\left(  w\right)  -b\cdot\left(  w-u\right)  $. Then, employing the
condition, we find that for any $v\in\Phi$, $\delta\left(  v\right)
=\delta\left(  w\right)  +b\cdot\left(  v-w\right)  =t+b\cdot\left(
w-u\right)  +b\cdot\left(  v-w\right)  =t+b\cdot\left(  v-u\right)  $.
$\blacksquare$

\medskip\ 

Letting $h=t-u$ and
\[
\tau_{h}\left(  v\right)  =v+h=\delta_{t,u;1}\left(  v\right)
\]
for all $v\in\Phi$ gives the special dilation $\tau_{h}:\Phi\rightarrow\Phi$
which we call a \textbf{translation }by $h$. Every translation is invertible,
$\tau_{h}^{-1}=\tau_{-h}$, and $\tau_{0}$ is the identity map. The
translations on $\Phi$ correspond one-to-one with the vectors of the
directional subspace $\Phi^{\mathsf{v}}$ of which $\Phi$ is a coset. We have
$\tau_{h}\circ\tau_{k}=\tau_{k}\circ\tau_{h}=\tau_{k+h}$ and thus under the
operation of composition, the translations on $\Phi$ are isomorphic to the
additive group of $\Phi^{\mathsf{v}}$. Since every translation is a dilation
$\delta_{t,u;b}$ with $b=1$, from the proposition above we infer the following result.

\begin{corollary}
A function $\tau:\Phi\rightarrow\Phi$ on the affine flat $\Phi$ is a
translation if and only if $\tau\left(  v\right)  -\tau\left(  w\right)  =v-w$
for all $v,w\in\Phi$. $\blacksquare$
\end{corollary}

\begin{exercise}
The translations are the dilations that lack a unique fixed point \emph{(}or
``center''\emph{)}.
\end{exercise}

\subsection{Congruence Under the Affine Group}

Viewing the vector space $\mathcal{V}$ as an affine flat in $\mathsf{A}\left(
\mathcal{V}\right)  $ and using the composition of maps as group operation,
the bijective affine self-maps $\alpha:\mathcal{V}\rightarrow\mathcal{V}$ form
a group $\mathsf{GA}\left(  \mathcal{V}\right)  $ called the \textbf{affine
group} of $\mathcal{V}$. The elements of $\mathsf{GA}\left(  \mathcal{V}%
\right)  $ may be described as the invertible vector space self-maps on
$\mathcal{V}$ composed with the translations on $\mathcal{V}$. Similarly,
there is an affine group $\mathsf{GA}\left(  \Phi\right)  $ for any affine
flat $\Phi\in\mathsf{A}\left(  \mathcal{V}\right)  $. However, $\mathsf{GA}%
\left(  \mathcal{V}\right)  $ already covers the groups $\mathsf{GA}\left(
\Phi\right)  $ in the sense that any such $\mathsf{GA}\left(  \Phi\right)  $
may be viewed as the restriction to $\Phi$ of the subgroup of $\mathsf{GA}%
\left(  \mathcal{V}\right)  $ that fixes $\Phi$ (those maps for which the
image of $\Phi$ remains contained in $\Phi$).

A \textbf{figure} in $\mathcal{V}$ is just geometric language for a subset of
$\mathcal{V}$. Two figures in $\mathcal{V}$ are called \textbf{affine
congruent} if one is the image of the other under an element of $\mathsf{GA}%
\left(  \mathcal{V}\right)  $. Affine congruence is an equivalence relation.

\begin{exercise}
Affine frames in a vector space are affine congruent if and only if they have
the same cardinality.
\end{exercise}

\begin{exercise}
Affine flats in a vector space are affine congruent if and only if their
directional subspaces have bases of the same cardinality.
\end{exercise}

\subsection{Problems}

\ 

\noindent\textbf{1}. When $\mathcal{F}$ is the smallest possible field, namely
the field $\left\{  0,1\right\}  $ of just two elements, \emph{any} nonempty
subset $\mathcal{S}$ of a vector space $\mathcal{V}$ over $\mathcal{F}$
contains the translate $v+a\cdot\left(  x-v\right)  $ for every $v,x\in
\mathcal{S}$ and every $a\in\mathcal{F}$. Not all such $\mathcal{S}$ are
affine flats, however.

On the other hand, if $1+1\neq0$ in $\mathcal{F}$, a nonempty subset
$\mathcal{S}$ of a vector space over $\mathcal{F}$ is an affine flat if it
contains $v+a\cdot\left(  x-v\right)  $ for every $v,x\in\mathcal{S}$ and
every $a\in\mathcal{F}$.

What about when $\mathcal{F}$ is a field of $4$ elements?

\newpage

\section{Basic Affine Results and Methods}

\subsection{Notations}

Throughout this chapter, $\mathcal{V}$ is our usual vector space over the
field $\mathcal{F}$. We generally omit any extra assumptions regarding
$\mathcal{V}$ or $\mathcal{F}$ (such as $\dim\mathcal{V}\neq0$) that each
result might need to make its hypotheses realizable. Without further special
mention, we will frequently use the convenient abuse of notation $P=\left\{
P\right\}  $ so that the point $P\in\mathsf{A}\left(  \mathcal{V}\right)  $ is
the singleton set that contains the vector $P\in\mathcal{V}$.

\subsection{Two Axiomatic Propositions}

The two propositions presented in this section are of an extremely basic
geometric character, and are assumed as axioms in most developments of affine
geometry from the \emph{synthetic} viewpoint. While certainly not universal,
the first of these does apply to many ordinary kinds of geometries. It is
commonly, but somewhat loosely, expressed as ``two points determine a line.''

\begin{proposition}
There is one and only one line that contains the distinct points $P$ and $Q$.
\end{proposition}

Proof: The affine flat $k=\left\{  P+a\cdot\left(  Q-P\right)  \mid
a\in\mathcal{F}\right\}  $ clearly contains both $P$ and $Q$, and is a line
since $k^{\mathsf{v}}=\left\{  a\cdot\left(  Q-P\right)  \mid a\in
\mathcal{F}\right\}  $ is one-dimensional. On the other hand, suppose that $l$
is a line that contains both $P$ and $Q$. Since $l$ is a line that contains
$P$, it has the form $l=P+l^{\mathsf{v}}$ where $l^{\mathsf{v}}$ is
one-dimensional. Since $l$ is an affine flat that contains both $P$ and $Q$,
it contains the translate $P+\left(  Q-P\right)  $ so that $Q-P\in
l^{\mathsf{v}}$, and therefore $l^{\mathsf{v}}=k^{\mathsf{v}}$ and $l=k$.
$\blacksquare$

The line that contains the distinct points $P$ and $Q$ is $P+_{\mathsf{A}}Q$,
of course. We will write $PQ$ to mean the line $P+_{\mathsf{A}}Q$ for
(necessarily distinct) points $P$ and $Q$.

\begin{corollary}
Two distinct intersecting lines intersect in a single point. $\blacksquare$
\end{corollary}

\begin{exercise}
Dimensional considerations then imply that the affine sum of distinct
intersecting lines is a plane.
\end{exercise}

Our second proposition involves parallelism, and can be used in a synthetic
context to prove the transitivity of parallelism for coplanar lines.

\begin{proposition}
Given a point $P$ and a line $l$, there is one and only one line $m$ that
contains $P$ and is parallel to $l$.
\end{proposition}

Proof: In order for $m$ and $l$ to be parallel, we must have $m^{\mathsf{v}%
}=l^{\mathsf{v}}$, since $m$ and $l$ have equal finite dimension. Hence the
unique line sought is $m=P+l^{\mathsf{v}}$. (Recalling that the cosets of a
subspace form a partition of the full space, $P$ must lie in exactly one coset
of $l^{\mathsf{v}}$.) $\blacksquare$

\subsection{A Basic Configuration Theorem}

The following is one of the affine versions of a configuration theorem of
projective geometry attributed to the French geometer Girard Desargues
(1591-1661). The particular affine version treated here is the completely
nonparallel one, the one that looks just like the projective theorem as we
shall later see.

\begin{theorem}
[Desargues]Let $A,A\,^{\prime},B,B\,^{\prime},C,C\,^{\prime}$ be distinct
points and let the distinct lines $AA\,^{\prime}$, $BB\,^{\prime}$, and
$CC\,^{\prime}$ be concurrent in the point $P$. Let the nonparallel lines $AB$
and $A\,^{\prime}B\,^{\prime}$ intersect in the point $C\,^{\prime\prime}$,
the nonparallel lines $BC$ and $B\,^{\prime}C\,^{\prime}$ intersect in the
point $A\,^{\prime\prime}$, and the nonparallel lines $CA$ and $C\,^{\prime
}A\,^{\prime}$ intersect in the point $B\,^{\prime\prime}$. Then
$A\,^{\prime\prime}$, $B\,^{\prime\prime}$, and $C\,^{\prime\prime}$ are
collinear.%
\begin{center}
\includegraphics[
trim=0.000000in 0.000000in -0.355713in -0.003817in,
height=2.2589in,
width=2.7605in
]%
{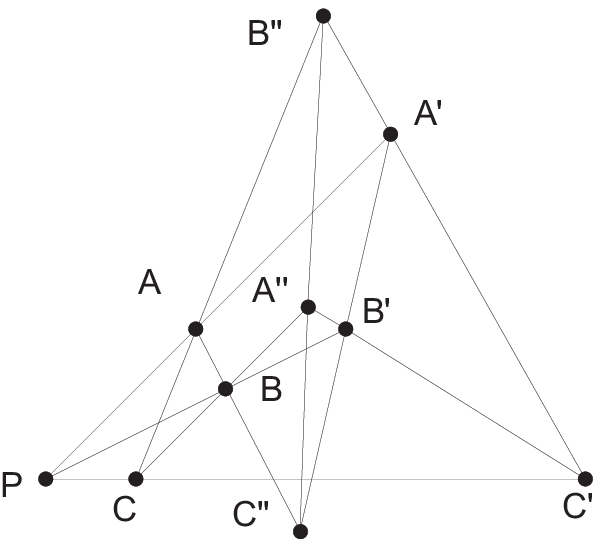}%
\\
Illustrating the Theorem of Desargues
\end{center}
\end{theorem}

Proof: We have $P=a\cdot A+a\,^{\prime}\cdot A\,^{\prime}=b\cdot
B+b\,^{\prime}\cdot B\,^{\prime}=c\cdot C+c\,^{\prime}\cdot C\,^{\prime}$ for
scalars $a,a\,^{\prime},b,b\,^{\prime}c,c\,^{\prime}$ which satisfy
$a+a\,^{\prime}=b+b\,^{\prime}=c+c\,^{\prime}=1$. We then derive the three
equalities
\begin{align*}
a\cdot A-b\cdot B &  =b\,^{\prime}\cdot B\,^{\prime}-a\,^{\prime}\cdot
A\,^{\prime}\text{,}\\
b\cdot B-c\cdot C &  =c\,^{\prime}\cdot C\,^{\prime}-b\,^{\prime}\cdot
B\,^{\prime}\text{,}\\
c\cdot C-a\cdot A &  =a\,^{\prime}\cdot A\,^{\prime}-c\,^{\prime}\cdot
C\,^{\prime}\text{.}%
\end{align*}
Consider the first of these. Suppose that $a=b$. Then $a\,^{\prime
}=b\,^{\prime}$ too, since $a+a\,^{\prime}=b+b\,^{\prime}=1$. But then
$a\neq0$, for if $a=0$, then $a\,^{\prime}=1$, which would give $A\,^{\prime
}=B\,^{\prime}$ contrary to hypothesis. Similarly $a\,^{\prime}\neq0$. Thus
supposing that $a=b$ leads to the nonparallel lines $AB$ and $A\,^{\prime
}B\,^{\prime}$ being parallel. Hence $a-b=$ $b\,^{\prime}-a\,^{\prime}\neq0$,
and similarly, $b-c=c\,^{\prime}-b\,^{\prime}\neq0$ and $c-a=a\,^{\prime
}-c\,^{\prime}\neq0$. Denote these pairs of equal nonzero differences by
$c\,^{\prime\prime}$, $a\,^{\prime\prime}$, and $b\,^{\prime\prime}$,
respectively. Then the equalities above relate to the double-primed
intersection points by
\begin{align*}
a\cdot A-b\cdot B &  =b\,^{\prime}\cdot B\,^{\prime}-a\,^{\prime}\cdot
A\,^{\prime}=c\,^{\prime\prime}\cdot C\,^{\prime\prime}\text{,}\\
b\cdot B-c\cdot C &  =c\,^{\prime}\cdot C\,^{\prime}-b\,^{\prime}\cdot
B\,^{\prime}=a\,^{\prime\prime}\cdot A\,^{\prime\prime}\text{,}\\
c\cdot C-a\cdot A &  =a\,^{\prime}\cdot A\,^{\prime}-c\,^{\prime}\cdot
C\,^{\prime}=b\,^{\prime\prime}\cdot B\,^{\prime\prime}\text{,}%
\end{align*}
since if we divide by the double-primed scalar, each double-primed
intersection point is the equal value of two affine expressions that determine
points on the intersecting lines. Adding up the equalities gives
$0=a\,^{\prime\prime}\cdot A\,^{\prime\prime}+b\,^{\prime\prime}\cdot
B\,^{\prime\prime}+c\,^{\prime\prime}\cdot C\,^{\prime\prime}$ and noting that
$a\,^{\prime\prime}+b\,^{\prime\prime}+c\,^{\prime\prime}=0$ we conclude in
light of Exercise \ref{AffineDependence} that $A\,^{\prime\prime}$,
$B\,^{\prime\prime}$, and $C\,^{\prime\prime}$ do lie on one line as was to be
shown. $\blacksquare$

\subsection{Barycenters}

We now introduce a concept analogous to the physical concept of ``center of
gravity,'' which will help us visualize the result of computing a linear
expression and give us some convenient notation. Given the points
$A_{1},\ldots,A_{n}$, consider the linear expression $a_{1}\cdot A_{1}%
+\cdots+a_{n}\cdot A_{n}$ of weight $a_{1}+\cdots+a_{n}\neq0$ based on
$A_{1},\ldots,A_{n}$. The \textbf{barycenter} based on (the factors of the
terms of) this linear expression is the unique point $X$ such that
\[
\left(  a_{1}+\cdots+a_{n}\right)  \cdot X=a_{1}\cdot A_{1}+\cdots+a_{n}\cdot
A_{n}\text{.}%
\]
(If the weight were zero, there would be no such unique point $X$.) The
barycenter of $a_{1}\cdot A_{1}+\cdots+a_{n}\cdot A_{n}$, which of course lies
in the affine span of $\left\{  A_{1},\ldots,A_{n}\right\}  $, will be denoted
by $\pounds\left[  a_{1}\cdot A_{1}+\cdots+a_{n}\cdot A_{n}\right]  $.
Clearly, for any scalar $m\neq0$, $\pounds\left[  \left(  ma_{1}\right)  \cdot
A_{1}+\cdots+\left(  ma_{n}\right)  \cdot A_{n}\right]  =\pounds\left[
a_{1}\cdot A_{1}+\cdots+a_{n}\cdot A_{n}\right]  $, and this homogeneity
property allows us to use $\pounds\left[  a_{1}\cdot A_{1}+\cdots+a_{n}\cdot
A_{n}\right]  $ as a convenient way to refer to an affine expression by
referring to any one of its nonzero multiples instead. Also, supposing that
$1\leqslant k<n$, $a=a_{1}+\cdots+a_{k}\neq0$, $a\,^{\prime}=a_{k+1}%
+\cdots+a_{n}\neq0$,
\[
X_{k}=\pounds\left[  a_{1}\cdot A_{1}+\cdots+a_{k}\cdot A_{k}\right]  \text{,}%
\]
and
\[
X_{k}{}^{\prime}=\pounds\left[  a_{k+1}\cdot A_{k+1}+\cdots+a_{n}\cdot
A_{n}\right]  \text{,}%
\]
it is easy to see that if $a+a\,^{\prime}\neq0$ then
\[
\pounds\left[  a_{1}\cdot A_{1}+\cdots+a_{n}\cdot A_{n}\right]  =\pounds
\left[  a\cdot X_{k}+a\,^{\prime}\cdot X_{k}{}^{\prime}\right]  \text{.}%
\]
Hence we may compute barycenters in a piecemeal fashion, provided that no zero
weight is encountered along the way.

\subsection{A Second Basic Configuration Theorem}

The following is the affine version of a configuration theorem of projective
geometry which dates back much farther than the theorem of Desargues given
above. Pappus of Alexandria, a Greek who lived in the 4th century, long before
projective geometry was created, used concepts of Euclidean geometry to prove
this same affine version of the projective theorem that now bears his name.

\begin{theorem}
[Theorem of Pappus]\label{Pappus}Let $l,l\,^{\prime}$ be distinct coplanar
lines with distinct points $A,B,C\subset l\smallsetminus l\,^{\prime}$ and
distinct points $A\,^{\prime},B\,^{\prime},C\,^{\prime}\subset l\,^{\prime
}\smallsetminus l$. If $BC\,^{\prime}$ meets $B\,^{\prime}C$ in $A\,^{\prime
\prime}$, $AC\,^{\prime}$ meets $A\,^{\prime}C$ in $B\,^{\prime\prime}$, and
$AB\,^{\prime}$ meets $A\,^{\prime}B$ in $C\,^{\prime\prime}$, then
$A\,^{\prime\prime},B\,^{\prime\prime},C\,^{\prime\prime}$ are collinear.
\end{theorem}%

\begin{center}
\includegraphics[
trim=0.000000in 0.000000in -0.047249in -0.027293in,
height=2.9992in,
width=3.4999in
]%
{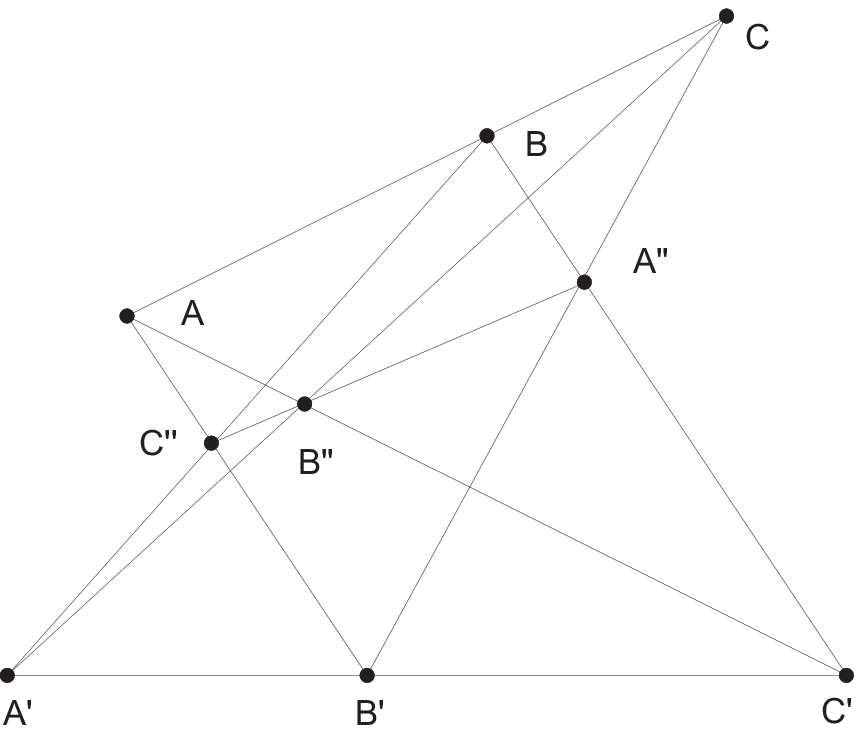}%
\\
Illustrating the Theorem of Pappus
\end{center}

Proof: We will express $A\,^{\prime\prime},B\,^{\prime\prime},C\,^{\prime
\prime}$ in terms of the elements of the affine frame $\left\{  A,B\,^{\prime
},C\right\}  $. First, we may write for some $a,b\,^{\prime},c$ with
$a+b\,^{\prime}+c=1$,
\[
B\,^{\prime\prime}=\pounds \left[  a\cdot A+b\,^{\prime}\cdot B\,^{\prime
}+c\cdot C\right]  \text{.}%
\]
Because $B$ is on $AC$, we may write
\[
B=\pounds \left[  pa\cdot A+c\cdot C\right]  \text{.}%
\]
Because $A\,^{\prime}$ is on $B\,^{\prime\prime}C$, we may write
\begin{align*}
A\,^{\prime} &  =\pounds \left[  B\,^{\prime\prime}+(q-1)c\cdot C\right]  \\
&  =\pounds \left[  a\cdot A+b\,^{\prime}\cdot B\,^{\prime}+qc\cdot C\right]
\text{.}%
\end{align*}
Because $C\,^{\prime}$ is on both $AB\,^{\prime\prime}$ and $A\,^{\prime
}B\,^{\prime}$, we may write $C\,^{\prime}=\pounds \left[  r\cdot A+q\cdot
B\,^{\prime\prime}\right]  =\pounds \left[  A\,^{\prime}+s\cdot B\,^{\prime
}\right]  $ or
\begin{align*}
C\,^{\prime} &  =\pounds \left[  (r+qa)\cdot A+qb\,^{\prime}\cdot B\,^{\prime
}+qc\cdot C\right]  \\
&  =\pounds \left[  a\cdot A+(b\,^{\prime}+s)\cdot B\,^{\prime}+qc\cdot
C\right]
\end{align*}
and selecting between equal coefficients permits us to write
\[
C\,^{\prime}=\pounds \left[  a\cdot A+qb\,^{\prime}\cdot B\,^{\prime}+qc\cdot
C\right]  \text{.}%
\]
Because $A\,^{\prime\prime}$ is on both $BC\,^{\prime}$ and $B\,^{\prime}C$,
we may write $A\,^{\prime\prime}=\pounds \left[  p\cdot C\,^{\prime}-k\cdot
B\right]  =\pounds \left[  w\cdot B\,^{\prime}+x\cdot C\right]  $ or
\begin{align*}
A\,^{\prime\prime} &  =\pounds \left[  (1-k)pa\cdot A+pqb\,^{\prime}\cdot
B\,^{\prime}+(pq-1)c\cdot C\right]  \\
&  =\pounds \left[  w\cdot B\,^{\prime}+x\cdot C\right]
\end{align*}
so that $k=1$ and
\[
A\,^{\prime\prime}=\pounds \left[  pqb\,^{\prime}\cdot B\,^{\prime
}+(pq-1)c\cdot C\right]  \text{.}%
\]
Because $C\,^{\prime\prime}$ is on both $A\,^{\prime}B$ and $AB\,^{\prime}$,
we may write $C\,^{\prime\prime}=\pounds \left[  l\cdot A\,^{\prime}-q\cdot
B\right]  =\pounds \left[  y\cdot A+z\cdot B\,^{\prime}\right]  $ or
\begin{align*}
C\,^{\prime\prime} &  =\pounds \left[  (l-qp)a\cdot A+lb\,^{\prime}\cdot
B\,^{\prime}+(l-1)qc\cdot C\right]  \\
&  =\pounds \left[  y\cdot A+z\cdot B\,^{\prime}\right]
\end{align*}
so that $l=1$ and
\[
C\,^{\prime\prime}=\pounds \left[  (1-qp)a\cdot A+b\,^{\prime}\cdot
B\,^{\prime}\right]  \text{.}%
\]
Letting $a\,^{\prime\prime},b\,^{\prime\prime},c\,^{\prime\prime}$ be the
respective weights of the expressions above of which $A\,^{\prime\prime
},B\,^{\prime\prime},C\,^{\prime\prime}$ are barycenters in terms of
$A,B\,^{\prime},C$, we have
\begin{align*}
a\,^{\prime\prime} &  =pqb\,^{\prime}+(pq-1)c,\\
b\,^{\prime\prime} &  =a+b\,^{\prime}+c,\\
c\,^{\prime\prime} &  =(1-qp)a+b\,^{\prime}.
\end{align*}
We then find that
\[
a\,^{\prime\prime}\cdot A\,^{\prime\prime}+(1-pq)b\,^{\prime\prime}\cdot
B\,^{\prime\prime}+(-c\,^{\prime\prime})\cdot C\,^{\prime\prime}=0
\]
and the left-hand side is a directional expression. However, the coefficients
of this directional expression are not all zero. For if we suppose they are
all zero, then $1-pq=0$, and $b\,^{\prime}=0$, which puts $B\,^{\prime\prime}$
on $AC$ so that $C\,^{\prime}$ would be on $AC$ contrary to hypothesis. We
conclude in light of Exercise \ref{AffineDependence} that $A\,^{\prime\prime}%
$, $B\,^{\prime\prime}$, and $C\,^{\prime\prime}$ are collinear as was to be
shown. $\blacksquare$

\subsection{Vector Ratios}

If two nonzero vectors $u,v$ satisfy $a\cdot u=b\cdot v$ for two nonzero
scalars $a,b$, it makes sense to speak of the proportion $u:v=b:a$. Under
these conditions we say that $u,v$ have the \textbf{same direction} and only
then do we impute a scalar value to the \textbf{vector ratio} $\dfrac uv$ and
write $\dfrac uv=\dfrac ba$.

Let us write $\overrightarrow{PQ}=Q-P$ whenever $P$ and $Q$ are points.
Nonzero vectors $\overrightarrow{PQ}$ and $\overrightarrow{RS}$ have the same
direction precisely when $PQ$ and $RS$ are parallel lines. Hence the ratio
$\overrightarrow{PQ}/\overrightarrow{RS}$ has meaning whenever $PQ$ and $RS$
are parallel lines, and in particular, when they are the same line.

$X$ is a point on the line $AB$, and distinct from $A$ and $B$, if and only if
there are nonzero scalars $a,b$, such that $a+b\neq0$ and such that
$X=\pounds[a\cdot A+b\cdot B]$. But this is the same as saying that there are
nonzero scalars $a$ and $b$ with nonzero sum, such that $\left(  a+b\right)
\cdot\overrightarrow{AX}=b\cdot\overrightarrow{AB}$ and $\left(  a+b\right)
\cdot\overrightarrow{XB}=a\cdot\overrightarrow{AB}$, or equivalently such
that
\[
\frac{\underset{\,}{\overrightarrow{AX}}}{\overset{\,}{\overrightarrow{XB}}%
}=\frac{b}{a}\cdot
\]
We say that $\overrightarrow{AB}=\overrightarrow{AX}+\overrightarrow{XB}$ is
\textbf{divided} by $X$ in the ratio $b:a$.

The result of the following Subcenter Exercise will find useful employment
below as a lemma.

\begin{exercise}
[Subcenter Exercise]Given noncollinear points $A,B,C$ and nonzero scalars
$a,b,c$ such that $X=\pounds\left[  a\cdot A+b\cdot B+c\cdot C\right]  $,
barycentric considerations imply that when $a+b\neq0$, $CX$ and $AB$ meet in
$\pounds\left[  a\cdot A+b\cdot B\right]  $. On the other hand, when $a+b=0$,
$CX$ and $AB$ are parallel. Hence $CX\parallel AB$ if and only if $a+b=0$, and
$CX\cap AB=\pounds\left[  a\cdot A+b\cdot B\right]  $ if and only if $a+b\neq0
$.%
\begin{center}
\includegraphics[
trim=0.000000in 0.000000in -0.276685in 0.000000in,
height=1.2505in,
width=2.0081in
]%
{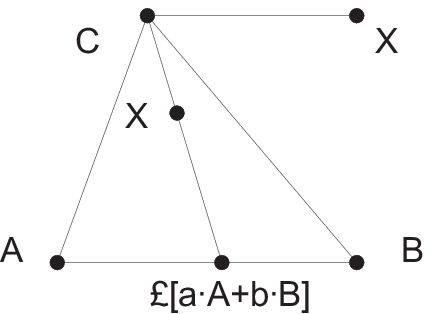}%
\\
Figure for Above Exercise
\label{Figure for Above Exercise}%
\end{center}
\end{exercise}

The next result is the affine version of Proposition VI.2 from Euclid's
\textit{Elements}. By $PQR$ we will always mean the plane that is the affine
sum of three given (necessarily noncollinear) points $P,Q,R$.

\begin{theorem}
[Similarity Theorem]Let $X$ be a point in $ABC$ and not on $AB$, $BC$, or
$CA$, and let $W=BX\cap AC$. Then
\[
CX\parallel AB\Rightarrow\frac{\underset{\,}{\overrightarrow{WC}}}%
{\overset{\,}{\overrightarrow{WA}}}=\frac{\underset{\,}{\overrightarrow{WX}}%
}{\overset{\,}{\overrightarrow{WB}}}=\frac{\underset{\,}{\overrightarrow{CX}}%
}{\overset{\,}{\overrightarrow{AB}}}%
\]
and
\[
\frac{\underset{\,}{\overrightarrow{WC}}}{\overset{\,}{\overrightarrow{WA}}%
}=\frac{\underset{\,}{\overrightarrow{WX}}}{\overset{\,}{\overrightarrow{WB}}%
}\Rightarrow CX\parallel AB\text{.}%
\]
\end{theorem}%

\begin{center}
\includegraphics[
trim=0.000000in 0.000000in -0.004150in -0.356755in,
height=1.4996in,
width=3.4584in
]%
{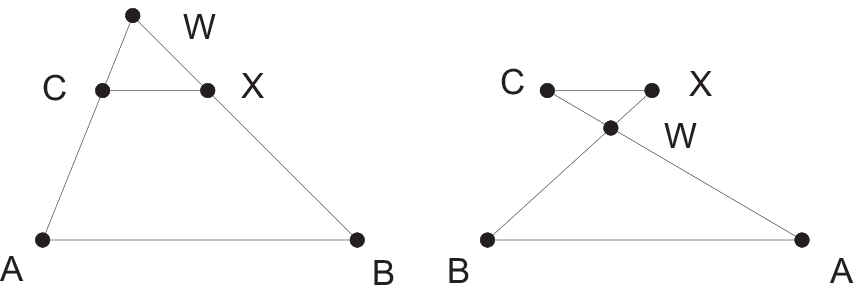}%
\\
Two Possible Figures for the Similarity Theorem
\end{center}

Proof: ``Let $X$ be a point in $ABC$ and not on $AB$, $BC$, or $CA$%
\thinspace'' is just the geometric way to say that there exist nonzero scalars
$a,b,c$ such that $X=\pounds \left[  a\cdot A+b\cdot B+c\cdot C\right]  $.
Suppose first that $CX\parallel AB$. Then by the Subcenter Exercise above,
$a+b=0$ and we may replace $a$ with $-b$ to give $X=\pounds \left[  \left(
-b\right)  \cdot A+b\cdot B+c\cdot C\right]  $ or $c\cdot\left(  X-C\right)
=b\cdot\left(  B-A\right)  $. Hence $\overrightarrow{CX}/\overrightarrow
{AB}=b/c$. Employing the Subcenter Exercise again, we find that $W=BX\cap
AC=\pounds \left[  \left(  -b\right)  \cdot A+c\cdot C\right]  $ so that $W$
divides $\overrightarrow{AC}$ in the ratio $c:-b$ and therefore
$\overrightarrow{WC}/\overrightarrow{WA}=b/c$. $X$ may also be obtained by the
piecemeal calculation $X=\pounds \left[  \left(  c-b\right)  \cdot W+b\cdot
B\right]  $ which is the same as $c\cdot\left(  X-W\right)  =b\cdot\left(
B-W\right)  $ or $\overrightarrow{WX}/\overrightarrow{WB}=b/c$. This proves
the first part.

On the other hand, assume that $\overrightarrow{WC}/\overrightarrow
{WA}=\overrightarrow{WX}/\overrightarrow{WB}$. We may suppose that these two
equal vector ratios both equal $b/c$ for two nonzero scalars $b,c$. Then
$c\cdot\left(  C-W\right)  =b\cdot(A-W)$ and $c\cdot\left(  X-W\right)
=b\cdot(B-W)$. Subtracting the first of these two equations from the second,
we obtain $c\cdot\left(  X-C\right)  =b\cdot\left(  B-A\right)  $, which shows
that $CX\parallel AB$ as desired. $\blacksquare$

\subsection{A Vector Ratio Product and Collinearity}

The following theorem provides a vector ratio product criterion for three
points to be collinear. The Euclidean version is attributed to Menelaus of
Alexandria, a Greek geometer of the late 1st century.

\begin{theorem}
[Theorem of Menelaus]Let $A,B,C$ be noncollinear points, and let $A\,^{\prime
},B\,^{\prime},C\,^{\prime}$ be points none of which is $A$, $B$, or $C$, such
that $C\,^{\prime}$ is on $AB$, $A\,^{\prime}$ is on $BC$, and $B\,^{\prime}$
is on $CA$. Then $A\,^{\prime},B\,^{\prime},C\,^{\prime}$ are collinear if and
only if
\[
\frac{\underset{\,}{\overrightarrow{BA\,^{\prime}}}}{\overset{\,}%
{\overrightarrow{A\,^{\prime}C}}}\cdot\frac{\underset{\,}{\overrightarrow
{CB\,^{\prime}}}}{\overset{\,}{\overrightarrow{B\,^{\prime}A}}}\cdot
\frac{\underset{\,}{\overrightarrow{AC\,^{\prime}}}}{\overset{\,}%
{\overrightarrow{C\,^{\prime}B}}}=-1.
\]
\end{theorem}%

\begin{center}
\includegraphics[
trim=0.119164in 0.000000in -0.600137in -0.001768in,
height=1.6077in,
width=2.5408in
]%
{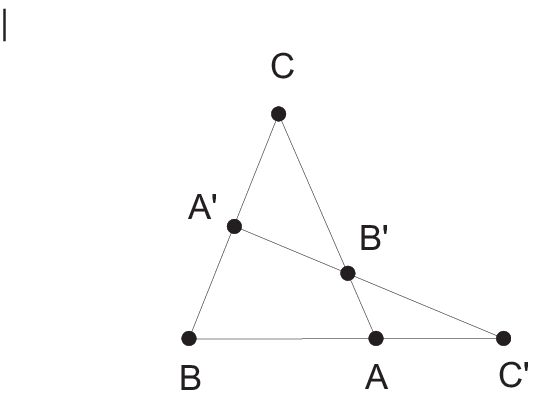}%
\\
Figure for the Theorem of Menelaus
\end{center}

Proof: Let $A\,^{\prime},B\,^{\prime},C\,^{\prime}$ be collinear.
$B,C,C\,^{\prime}$ are noncollinear and we may write $B\,^{\prime
}=\pounds \left[  b\cdot B+c\cdot C+c\,^{\prime}\cdot C\,^{\prime}\right]  $.
Then $A\,^{\prime}=\pounds \left[  b\cdot B+c\cdot C\right]  $ and
$A=\pounds \left[  b\cdot B+c\,^{\prime}\cdot C\,^{\prime}\right]  $ by the
Subcenter Exercise above. $B\,^{\prime}$ may also be obtained by the piecemeal
calculation $B\,^{\prime}=\pounds \left[  c\cdot C+(b+c\,^{\prime})\cdot
A\right]  $. Thus $A\,^{\prime}$ divides $\overrightarrow{BC}$ in the ratio
$c:b$, $B\,^{\prime}$ divides $\overrightarrow{CA}$ in the ratio $\left(
b+c\,^{\prime}\right)  :c$, and $A$ divides $\overrightarrow{BC\,^{\prime}}$
in the ratio $c\,^{\prime}:b$. Hence
\[
\frac{\underset{\,}{\overrightarrow{BA\,^{\prime}}}}{\overset{\,}%
{\overrightarrow{A\,^{\prime}C}}}\cdot\frac{\underset{\,}{\overrightarrow
{CB\,^{\prime}}}}{\overset{\,}{\overrightarrow{B\,^{\prime}A}}}\cdot
\frac{\underset{\,}{\overrightarrow{AC\,^{\prime}}}}{\overset{\,}%
{\overrightarrow{C\,^{\prime}B}}}=\frac{c}{b}\cdot\frac{b+c\,^{\prime}}%
{c}\cdot\frac{b}{-b-c\,^{\prime}}=-1
\]
as required.

On the other hand, suppose that the product of the ratios is $-1$. We may
suppose that
\[
A\,^{\prime}=\pounds \left[  1\cdot B+x\cdot C\right]  ,B\,^{\prime
}=\pounds \left[  y\cdot A+x\cdot C\right]  ,C\,^{\prime}=\pounds \left[
y\cdot A+z\cdot B\right]
\]
where none of $x,y,z,1+x,y+x,y+z$ is zero. Thus
\[
\frac{x}{1}\cdot\frac{y}{x}\cdot\frac{z}{y}=-1
\]
and hence $z=-1$. Then we have
\[
\left(  1+x\right)  \cdot A\,^{\prime}+\left(  -y-x\right)  \cdot B\,^{\prime
}+\left(  y-1\right)  \cdot C\,^{\prime}=0
\]
and $\left(  1+x\right)  +\left(  -y-x\right)  +\left(  y-1\right)  =0.$ We
conclude in light of Exercise \ref{AffineDependence} that $A\,^{\prime
},B\,^{\prime},C\,^{\prime}$ are collinear as required. $\blacksquare$

\subsection{A Vector Ratio Product and Concurrency}

The Italian geometer Giovanni Ceva (1647-1734) is credited with the Euclidean
version of the following theorem which is quite similar to the above theorem
of Menelaus but instead provides a criterion for three lines to be either
concurrent or parallel.

\begin{theorem}
[Theorem of Ceva]Let $A,B,C$ be noncollinear points, and let $A\,^{\prime
},B\,^{\prime},C\,^{\prime}$ be points none of which is $A$, $B$, or $C$, such
that $C\,^{\prime}$ is on $AB$, $A\,^{\prime}$ is on $BC$, and $B\,^{\prime}$
is on $CA$. Then $AA\,^{\prime},BB\,^{\prime},CC\,^{\prime}$ are either
concurrent in a point not lying on $AB$, $BC$, or $CA,$ or are parallel, if
and only if
\[
\frac{\underset{\,}{\overrightarrow{BA\,^{\prime}}}}{\overset{\,}%
{\overrightarrow{A\,^{\prime}C}}}\cdot\frac{\underset{\,}{\overrightarrow
{CB\,^{\prime}}}}{\overset{\,}{\overrightarrow{B\,^{\prime}A}}}\cdot
\frac{\underset{\,}{\overrightarrow{AC\,^{\prime}}}}{\overset{\,}%
{\overrightarrow{C\,^{\prime}B}}}=1.
\]
\end{theorem}%

\begin{center}
\includegraphics[
trim=0.097248in 0.000000in -0.000710in 0.000151in,
height=1.5091in,
width=3.5492in
]%
{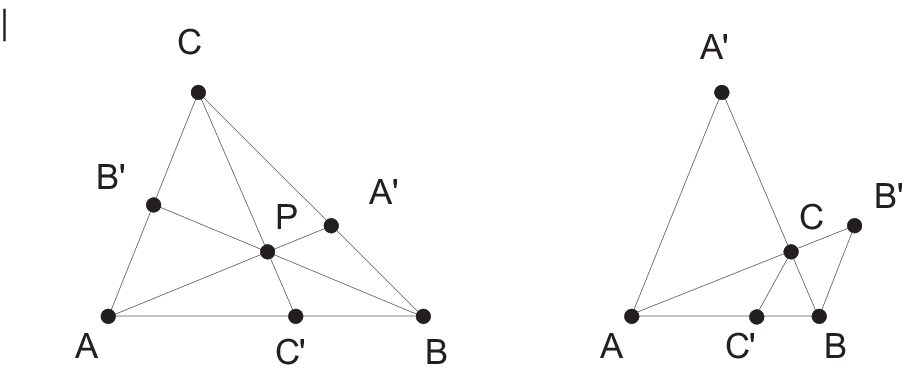}%
\\
Figures for the Theorem of Ceva
\end{center}

Proof: Suppose $AA\,^{\prime},BB\,^{\prime},CC\,^{\prime}$ meet in
$P=\pounds \left[  a\cdot A+b\cdot B+c\cdot C\right]  $ with nonzero scalars
$a,b,c$. By the Subcenter Exercise above, $A\,^{\prime}=\pounds \left[  b\cdot
B+c\cdot C\right]  $ so that $\overrightarrow{BA\,^{\prime}}/\overrightarrow
{A\,^{\prime}C}=c/b$. Similarly, we find that $\overrightarrow{CB\,^{\prime}%
}/\overrightarrow{B\,^{\prime}A}=a/c$ and $\overrightarrow{AC\,^{\prime}%
}/\overrightarrow{C\,^{\prime}B}=b/a$. The desired product result follows at once.

Suppose next that $AA\,^{\prime},BB\,^{\prime},CC\,^{\prime}$ are parallel.
Applying the Similarity Theorem above, $AA\,^{\prime}\parallel CC\,^{\prime}$
implies that $C$ divides $A\,^{\prime}B$ in the same ratio as $C\,^{\prime}$
divides $AB$, and $BB\,^{\prime}\parallel CC\,^{\prime}$ implies that $C$ also
divides $AB\,^{\prime}$ in the same ratio as $C\,^{\prime}$ divides $AB$. The
desired product result follows readily.

Now suppose that the product is as stated. According to the conditions under
which we allow vector ratios to be written, there are nonzero scalars $a,b,c$
such that $\overrightarrow{BA\,^{\prime}}/\overrightarrow{A\,^{\prime}C}=c/b$,
$\overrightarrow{CB\,^{\prime}}/\overrightarrow{B\,^{\prime}A}=a/c$, and
$\overrightarrow{AC\,^{\prime}}/\overrightarrow{C\,^{\prime}B}=b/a$. If
$a+b+c\neq0$, define $P=\pounds \left[  a\cdot A+b\cdot B+c\cdot C\right]  $.
Then $\overrightarrow{BC}$ is divided by $A\,^{\prime}$ in the ratio $c:b$,
$A\,^{\prime}=\pounds \left[  b\cdot B+c\cdot C\right]  $, and barycentric
considerations place $P$ on $AA\,^{\prime}$. Similarly, $P$ is on
$BB\,^{\prime}$ and $CC\,^{\prime}$. If, on the other hand, $a+b+c=0$, then
$c=-a-b$. We see then that $C$ divides $A\,^{\prime}B$ in the same ratio as
$C\,^{\prime}$ divides $AB$, and that $C$ also divides $AB\,^{\prime}$ in the
same ratio as $C\,^{\prime}$ divides $AB$. The Similarity Theorem above then
gives $AA\,^{\prime}\parallel CC\,^{\prime}$ and $BB\,^{\prime}\parallel
CC\,^{\prime}$. $\blacksquare$

\newpage

\section{An Enhanced Affine Environment}

\subsection{Inflating a Flat}

Let $\mathcal{V}$ be our usual vector space over the field $\mathcal{F}$. Any
flat $\Phi\in\mathsf{A}\left(  \mathcal{V}\right)  $ may be embedded in a
vector space $\Phi_{+}$ as a hyperplane not containing $0$, a construction
that we will call \textbf{inflating} $\Phi$. Employing the directional
subspace $\Phi^{\mathsf{v}}$ of $\Phi$, we form $\Phi_{+}=\mathcal{F}%
\times\Phi^{\mathsf{v}}$ and embed $\Phi$ in it as $\left\{  1\right\}
\times\Phi^{\mathsf{v}}$ by the affine isomorphism that first sends $x\in\Phi$
to $x^{\mathsf{v}}=(x-o)\in\Phi^{\mathsf{v}}$ and then sends $x^{\mathsf{v}}$
to $1\times x^{\mathsf{v}}=\left(  1,x^{\mathsf{v}}\right)  $, where $o$ is a
fixed but arbitrary vector in $\Phi$. Forming $\Phi_{+}$ begins a transition
from affine to projective geometry. Benefits of this enhancement of the affine
environment include more freedom in representing points (similar to that
gained by using the $\pounds$ operator) and the ability to completely separate
the vectors representing points from the vectors representing directions.
Having then gained the separation of point-representing vectors and
direction-representing vectors, the directions can be added to the point
community as a new special class of generalized points which from the affine
viewpoint are thought of as lying at infinity. These points at infinity can be
used to deal in a uniform way with cases where an intersection point of lines
disappears to infinity as the lines become parallel. This leads us into the
viewpoint of projective geometry where what were points at infinity become
just ordinary points, and parallel lines cease to exist.

We now assume that any flat $\Phi$ which is of geometric interest to us has
been inflated, so that the r\^{o}le of $\Phi$ will be played by its alias
$\left\{  1\right\}  \times\Phi^{\mathsf{v}}$ in $\mathsf{A}\left(  \Phi
_{+}\right)  $, and the r\^{o}le of $\Phi^{\mathsf{v}}$ then will be played by
its alias $\left\{  0\right\}  \times\Phi^{\mathsf{v}}$. To simplify the
notation, we will refer to $\left\{  a\right\}  \times\Phi^{\mathsf{v}}$ as
$\Phi_{a}$, so that $\Phi_{1}$ is playing the r\^{o}le of $\Phi$ and $\Phi
_{0}$ is playing the r\^{o}le of $\Phi^{\mathsf{v}}$. The r\^{o}le of any
subflat $\Psi\subset\Phi$ is played by the subflat $\Psi_{1}=\left\{
1\right\}  \times\Psi^{\mathsf{v}}\subset\Phi_{1}$, of course. No matter what
holds for $\Phi$ and its directional subspace, $\Phi_{1}$ is disjoint from its
directional subspace $\Phi_{0}$. Lines through the point $O=\left\{
0\right\}  $ in $\Phi_{+}$ are of two fundamentally different types with
regard to $\Phi_{1}$: either such a line meets $\Phi_{1}$ or it does not. The
nonzero vectors of the lines through $O$ that meet $\Phi_{1}$ will be called
\textbf{point vectors}, and the nonzero vectors of the remaining lines through
$O$ (those that lie in $\Phi_{0})$ will be called \textbf{direction vectors}.
There is a unique line through $O$ passing through any given point $P$ of
$\Phi_{1}$, and since any nonzero vector in that line may be used as an
identifier of $P$, that line and each of its point vectors will be said to
\textbf{represent }$P$. Point vectors are homogeneous representatives of their
points in the sense that multiplication by any nonzero scalar gives a result
that still represents exactly the same point. To complete our view, each line
through $O$ lying in $\Phi_{0}$ will be called a \textbf{direction}, and it
and each of its nonzero vectors will be said to \textbf{represent} that
direction. Thus the zero vector stands alone in not representing either a
point or a direction.%

\begin{center}
\includegraphics[
trim=0.000000in 0.174247in -0.004351in -0.580822in,
height=2.6766in,
width=2.4171in
]%
{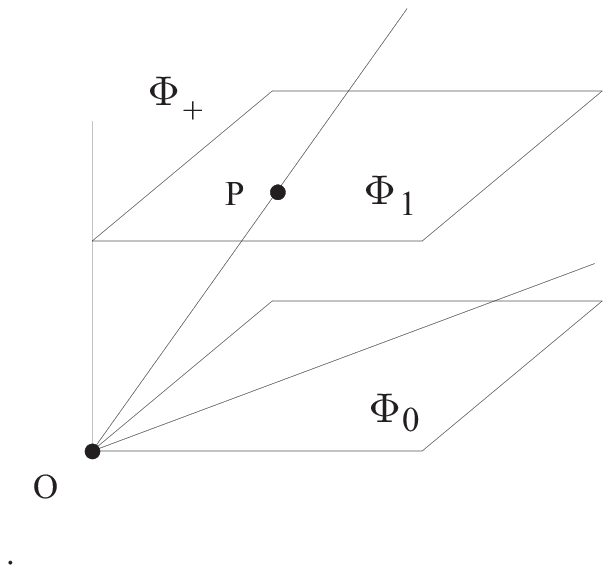}%
\\
Illustrating Inflating $\Phi$%
\end{center}

\newpage\ 

We are now close to the viewpoint of projective geometry where the lines
through $0$ actually \emph{are} the points, and the $n$-dimensional subspaces
are the ($n-1$)-dimensional flats of a projective structure. However, the
inflated-flat viewpoint is still basically an affine viewpoint, but an
enhanced one that is close to the pure projective one. From now on, we assume
that any affine flat has been inflated and freely use the resulting enhanced
affine viewpoint.

Suppose now that $P,Q,R\in\Phi_{1}$ are distinct points represented by the
point vectors $p,q,r$. Our geometric intuition indicates that $P,Q,R$ are
collinear in $\Phi_{1}$ if and only if $p,q,r$ are coplanar in $\Phi_{+}$. The
following proposition confirms this.

\begin{proposition}
Let $P,Q,R\in\Phi_{1}$ be distinct points represented by the point vectors
$p,q,r$. Then $P,Q,R$ are collinear if and only if $\left\{  p,q,r\right\}  $
is dependent.
\end{proposition}

Proof: Let $i$ denote the vector $1\times0$ in $\Phi_{+}=\mathcal{F}\times
\Phi^{\mathsf{v}}$. Then points in $\Phi_{1}$ all have the form $\left\{
i+x\right\}  $ where $x\in\Phi_{0}$, and point vectors all have the form
$d\cdot\left(  i+x\right)  $ for some nonzero scalar $d$. We thus set
$P=\left\{  i+u\right\}  ,Q=\left\{  i+v\right\}  ,R=\left\{  i+w\right\}  $,
and $p=a\cdot\left(  i+u\right)  ,q=b\cdot\left(  i+v\right)  ,r=c\cdot\left(
i+w\right)  $. Supposing that $P,Q,R$ are collinear, there are nonzero scalars
$k,l,m$ such that $k\cdot\left(  i+u\right)  +l\cdot\left(  i+v\right)
+m\cdot\left(  i+w\right)  =0$. Therefore $p,q,r$ are related by
\[
\frac{k}{a}\cdot p+\frac{l}{b}\cdot q+\frac{m}{c}\cdot r=0\text{.}%
\]

On the other hand, suppose that there are scalars $f,g,h$, not all zero, such
that $f\cdot p+g\cdot q+h\cdot r=0$. Then
\[
fa\cdot\left(  i+u\right)  +gb\cdot\left(  i+v\right)  +hc\cdot\left(
i+w\right)  =0\text{.}%
\]
Since $i\notin\Phi_{0}$, $fa+gb+hc=0$ and $P,Q,R$ are collinear by Exercise
\ref{AffineDependence}. $\blacksquare$

\smallskip\ 

Note that since we have complete freedom in scaling point vectors by nonzero
scalars, if we are given distinct collinear $P,Q,R$, there always is a
$\left\{  p,q,r\right\}  $ such that the scalar coefficients of a dependency
linear combination of it are any three nonzero scalars we choose. Moreover,
since a dependency linear combination is unaffected by a nonzero scaling,
besides the three nonzero scalar coefficients, one of $p,q,r$ may also be
picked in advance. This is the freedom in representing points which we get
from having insured that the flat where all the points lie is completely
separated from its directional subspace.

If we are given two distinct points $P,Q$ represented by point vectors $p,q$,
the above proposition implies that any point vector in the span of $\left\{
p,q\right\}  $ represents a point on the line $PQ$, and conversely, that every
point of $PQ$ in $\Phi_{1}$ is represented by a point vector of $\left\langle
\left\{  p,q\right\}  \right\rangle $. But besides point vectors, there are
other nonzero vectors in $\left\langle \left\{  p,q\right\}  \right\rangle $,
namely direction vectors. It is easy to see that the single direction
contained in $\left\langle \left\{  p,q\right\}  \right\rangle $ is
$V=\left\langle \left\{  p,q\right\}  \right\rangle \cap\Phi_{0}$ which we
refer to as the \textbf{direction} of $PQ$. We will consider $V$ to be a
generalized point of $PQ$, and we think of $V$ as lying ``at infinity'' on
$PQ$. Doing this will allow us to treat all nonzero vectors of $\left\langle
\left\{  p,q\right\}  \right\rangle $ in a uniform manner and to say that
$\left\langle \left\{  p,q\right\}  \right\rangle $ represents the line PQ$.$
This will turn out to have the benefit of eliminating the need for separate
consideration of related ``parallel cases'' for many of our results.

\begin{exercise}
Let $P,Q$ be distinct points represented by point vectors $p,q$, and let the
direction $V$ be represented by the direction vector $v$. Then the lines
$\left\langle \left\{  p,v\right\}  \right\rangle \cap\Phi_{1}$ and
$\left\langle \left\{  q,v\right\}  \right\rangle \cap\Phi_{1}$ are parallel
and $\left\langle \left\{  p,v\right\}  \right\rangle \cap\left\langle
\left\{  q,v\right\}  \right\rangle =\left\langle \left\{  v\right\}
\right\rangle $. \emph{(}Thus we speak of the parallel lines $PV$ and $QV$
that meet at $V$.\emph{)}
\end{exercise}

Besides lines like $PV$ and $QV$ of the above exercise, there can be lines of
the form $VW$ where $V$ and $W$ are distinct directions. We speak of such
lines as lines at infinity. If $V$ is represented by $v$ and $W$ by $w$, then
the nonzero vectors of $\left\langle \left\{  v,w\right\}  \right\rangle $
(all of which are direction vectors) are exactly the vectors that represent
the points of $VW$. Thus we view the ``compound direction'' $\left\langle
\left\{  v,w\right\}  \right\rangle $ as a generalized line. Each plane has
its line at infinity, and planes have the same line at infinity if and only if
they are parallel.

\subsection{Desargues Revisited}

We now state and prove another version of the theorem of Desargues treated
previously. The theorem's points and lines are now assumed to exist in our new
generalized framework.

\begin{theorem}
[Desargues]\label{GenEnvDesargues}Let $A,A\,^{\prime},B,B\,^{\prime
},C,C\,^{\prime},P $ be distinct points and let the lines $AA\,^{\prime}$,
$BB\,^{\prime}$, and $CC\,^{\prime}$ be distinct and concurrent in $P$. Let
$AB$ and $A\,^{\prime}B\,^{\prime}$ meet in the point $C\,^{\prime\prime}$,
$BC$ and $B\,^{\prime}C\,^{\prime}$ meet in the point $A\,^{\prime\prime}$,and
$CA$ and $C\,^{\prime}A\,^{\prime}$ meet in the point $B\,^{\prime\prime}$.
Then $A\,^{\prime\prime}$, $B\,^{\prime\prime}$, and $C\,^{\prime\prime}$ all
lie on the same line.
\end{theorem}

Proof: The corresponding small letter will always stand for a representing
vector. We can find representing vectors such that
\[
p=a+a\,^{\prime}=b+b\,^{\prime}=c+c\,^{\prime}%
\]
so that
\begin{align*}
a-b &  =b\,^{\prime}-a\,^{\prime}=c\,^{\prime\prime}\\
b-c &  =c\,^{\prime}-b\,^{\prime}=a\,^{\prime\prime}\\
c-a &  =a\,^{\prime}-c\,^{\prime}=b\,^{\prime\prime}\text{.}%
\end{align*}
Hence $a\,^{\prime\prime}+b\,^{\prime\prime}+c\,^{\prime\prime}=0$.
$\blacksquare$

While we were able to use simpler equations, this proof is not essentially
different from the previous one. What is more important is that new
interpretations now arise from letting various of the generalized points be at
infinity. We depict two of these possible new interpretations in the pair of
figures below.%

\begin{center}
\includegraphics[
trim=0.000000in 0.442636in -0.181993in -0.203793in,
height=2.5953in,
width=5.0194in
]%
{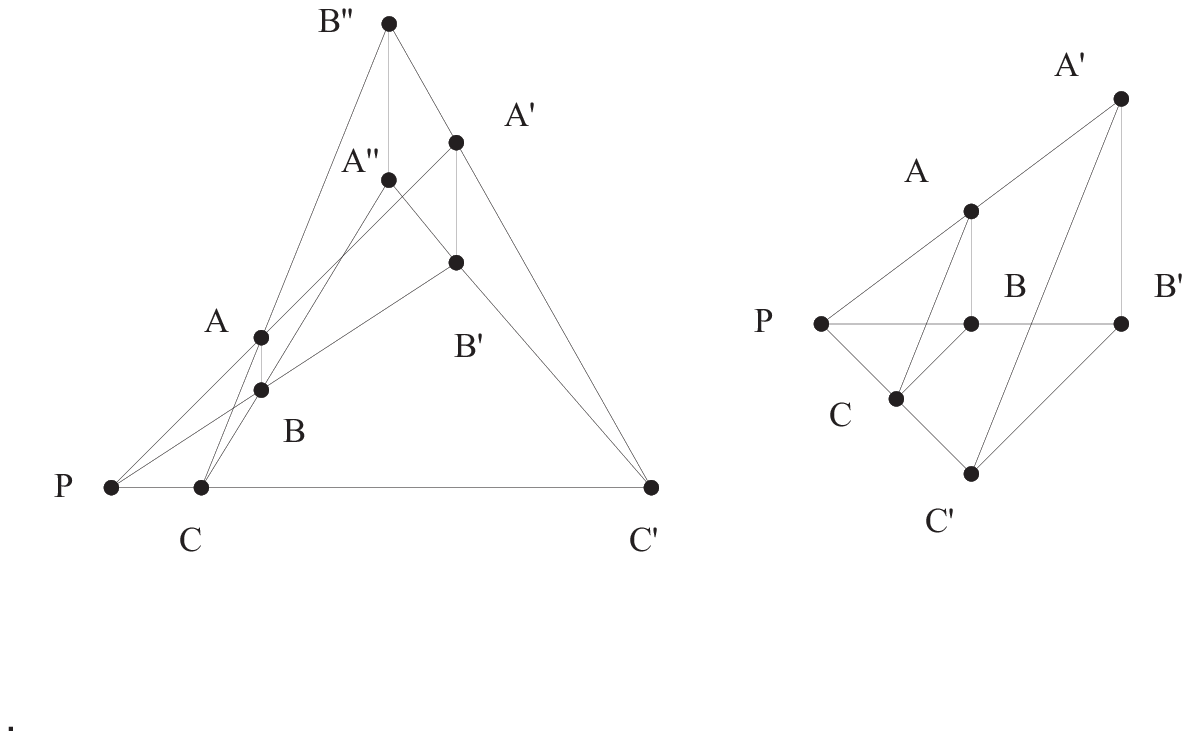}%
\\
New Interpretations Arising from Generalized-Framework Version
\end{center}
\smallskip

In the left-hand figure, notice that $A\,^{\prime\prime}B\,^{\prime\prime}$ is
parallel to the parallel lines $AB$ and $A\,^{\prime}B\,^{\prime}$, as it
contains their intersection point $C\,^{\prime\prime}$ at infinity. In the
right-hand figure, $A\,^{\prime\prime}$, $B\,^{\prime\prime}$, and
$C\,^{\prime\prime}$ all lie on the same line at infinity, and we can say that
the fact that any two of $A\,^{\prime\prime},B\,^{\prime\prime},C\,^{\prime
\prime}$ lie at infinity means that the third must also. Thus the right-hand
figure depicts the result that if any two of the three corresponding pairs
$AB,A\,^{\prime}B\,^{\prime}$, etc., are pairs of parallels, then the third is also.

\subsection{Representing (Generalized) Subflats}

\emph{Notice:} The qualifier ``generalized'' will henceforth usually be
omitted. From now on, points of $\Phi$ or other subflats to which we refer may
conceivably be ``at infinity'' unless specifically stated otherwise. This
includes intersections, so that, for example, parallel lines will be deemed to
intersect in a point just like nonparallel coplanar lines. Similarly, the more
or less obvious generalizations of all encountered affine concepts (affine
span, affine sum, etc.) are to be assumed to be in use whether specifically
defined or not.

\smallskip\ 

Within the context of an inflated flat we can readily show that there is a
one-to-one correspondence between the $(n-1)$-dimensional subflats of $\Phi$
and the $n$-dimensional subspaces of $\Phi_{+}$. Moreover, $n$ points of
$\Phi$ affine-span one of its $(n-1)$-dimensional subflats if and only if
every set of vectors which represents those points is an independent set in
$\Phi_{+}$. (The concept of affine span of points is generalized through use
of the span of the representing vectors to include all points, not just
``finite'' ones.)

The criterion for independence provided by exterior algebra (Corollary
\ref{CriterionForIndependence}) leads to the conclusion that a nonzero
exterior e-product (a \textbf{blade}) represents a subspace of $\Phi_{+}$ and
therefore also represents a subflat of $\Phi$ in the same homogeneous fashion
that a single nonzero vector represents a point. That is, if $\left\{
v_{1},\ldots,v_{n}\right\}  $ is an independent $n$-element set of vectors in
$\Phi_{+}$, then the $n$-\textbf{blade} $v_{1}\wedge\cdots\wedge v_{n}$
homogeneously represents an $\left(  n-1\right)  $-dimensional subflat of
$\Phi$. This is based on the following result.

\begin{proposition}
If $\left\{  v_{1},\ldots,v_{n}\right\}  $ and $\left\{  w_{1},\ldots
,w_{n}\right\}  $ are independent $n$-element sets of vectors that span the
same subspace, then the blades $w_{1}\wedge\cdots\wedge w_{n}$ and
$v_{1}\wedge\cdots\wedge v_{n}$ are proportional, and conversely, if
$w_{1}\wedge\cdots\wedge w_{n}$ and $v_{1}\wedge\cdots\wedge v_{n}$ are
proportional blades made up of the vectors of the independent sets $\left\{
w_{1},\ldots,w_{n}\right\}  $ and $\left\{  v_{1},\ldots,v_{n}\right\}  $,
then those sets of vectors span the same subspace.
\end{proposition}

Proof: Suppose that the independent sets $\left\{  v_{1},\ldots,v_{n}\right\}
$ and $\left\{  w_{1},\ldots,w_{n}\right\}  $ of $n$ vectors span the same
subspace. Then each $w_{i}$ can be expressed as a linear combination of
$\left\{  v_{1},\ldots,v_{n}\right\}  $. Putting these linear combinations in
the blade $w_{1}\wedge\cdots\wedge w_{n}$, after expanding and collecting the
nonzero terms we are left with a nonzero (since $w_{1}\wedge\cdots\wedge
w_{n}$ is nonzero) multiple of the blade $v_{1}\wedge\cdots\wedge v_{n}$.

On the other hand, suppose that $w_{1}\wedge\cdots\wedge w_{n}$ and
$v_{1}\wedge\cdots\wedge v_{n}$ are proportional blades expressed as exterior
products of vectors. Then
\[
0=w_{i}\wedge w_{1}\wedge\cdots\wedge w_{n}=w_{i}\wedge v_{1}\wedge
\cdots\wedge v_{n}%
\]
so that each $w_{i}$ is a linear combination of $\left\{  v_{1},\ldots
,v_{n}\right\}  $. $\blacksquare$

\smallskip\ 

We will find it convenient to class the empty exterior product, which we by
convention set equal to $1$, and all its nonzero multiples, as blades that
represents the empty flat and the subspace $\{0\}$. Thus, there is a
one-to-one correspondence between the sets of proportional blades in
$\bigwedge\Phi_{+}$ and the finite-dimensional subflats of $\Phi$. Given the
two blades $v_{1}\wedge\cdots\wedge v_{k}$ and $v_{k+1}\wedge\cdots\wedge
v_{n}$, then either $v_{1}\wedge\cdots\wedge v_{n}=0$ and the two
corresponding subflats intersect, or $v_{1}\wedge\cdots\wedge v_{n}$
represents the (generalized) affine sum of the two non-intersecting subflats.
By the same token, given the blade $\beta$, the subspace of $\Phi_{+}$ which
it represents is $\left\{  v\in\Phi_{+}\mid v\wedge\beta=0\right\}  $. A blade
that represents a hyperplane will be referred to as a \textbf{hyperplane
blade}, and similarly, the terms \textbf{plane blade}, \textbf{line blade},
and \textbf{point blade} will be used.

\subsection{Homogeneous and Pl\"ucker Coordinates}

\label{coords}

We suppose throughout this section and for the remainder of this chapter that
$\Phi_{+}$ has finite dimension $d$. Fix a basis $\mathcal{B}$ for $\Phi_{+}$.
The coordinates of a nonzero vector $v$ of $\Phi_{+}$ are known as
\textbf{homogenous coordinates} for the point represented by $v$. Bases
$\mathcal{B}^{\wedge}$ for $\bigwedge\Phi_{+}$ may be obtained by choosing
blades that are products of elements of $\mathcal{B}$. The coefficients of the
expansion for the $n$-blade $\beta$ in terms of the $n$-blades of any such
$\mathcal{B}^{\wedge}$ are known, in honor of German geometer Julius Pl\"ucker
(1801-1868), as \textbf{Pl\"ucker coordinates} for the $\left(  n-1\right)
$-dimensional subflat of $\Phi$ represented by $\beta$, and sometimes, by
abuse of language, as Pl\"ucker coordinates for any blade proportional to
$\beta$. Homogenous coordinates are then Pl\"ucker coordinates for a point.
Pl\"ucker coordinates are sometimes called \textbf{Grassmann coordinates} for
German schoolteacher Hermann Grassmann (1809-1877), the first author to
systematically treat vector spaces and vector algebras.

The coefficients of any linear combination of $\mathcal{B}$ are homogenous
coordinates for some point of $\Phi$, but an arbitrary linear combination of
the $n$-blades of a $\mathcal{B}^{\wedge}$ cannot be guaranteed in general to
yield Pl\"ucker coordinates for some $\left(  n-1\right)  $-dimensional
subflat of $\Phi$. That is, it is not true in general that every nonzero
element of $\bigwedge^{n}\Phi_{+}$ is expressible as an $n$-blade. We shall
find, however, that every nonzero element of $\bigwedge^{d-1}\Phi_{+}$ is a
$\left(  d-1\right)  $-blade.

Pl\"{u}cker coordinates of a blade may be expressed in terms of the
coordinates of the vectors that make up the blade by expanding out the
exterior product. If we suppose that we are given the $n$ (possibly dependent)
vectors $v_{1},\ldots,v_{n}$ with each $v_{j}$ given in terms of basis vectors
$x_{i}$ by
\[
v_{j}=\sum_{i=1}^{d}a_{i,j}\cdot x_{i\text{ ,}}%
\]
then
\[
v_{1}\wedge\cdots\wedge v_{n}=\sum_{i_{1}=1}^{d}\cdots\cdot\sum_{i_{n}=1}%
^{d}a_{i_{1},1}\cdots a_{i_{n},n}\cdot x_{i_{1}}\wedge\cdots\wedge x_{i_{n}}%
\]
and collecting terms on a particular set of basis elements (those that have
subscripts that increase from left to right) we then get
\[
\sum_{1\leqslant i_{1}<\cdots<i_{n}\leqslant d}\left(  \sum_{\sigma}\left(
-1\right)  ^{\sigma}a_{\sigma\left(  i_{1}\right)  ,1}\cdots a_{\sigma\left(
i_{n}\right)  ,n}\right)  \cdot x_{i_{1}}\wedge\cdots\wedge x_{i_{n}}%
\]
where the inner sum is over all permutations
\[
\sigma=\left(
\begin{array}
[c]{ccc}%
i_{1} & \cdots & i_{n}\\
\sigma\left(  i_{1}\right)  & \cdots & \sigma\left(  i_{n}\right)
\end{array}
\right)
\]
of each subscript combination that the outer summation has selected. We
recognize the inner sum as a determinant and thus the expansion may be written
as
\[
\sum_{1\leqslant i_{1}<\cdots<i_{n}\leqslant d}\left|
\begin{array}
[c]{ccc}%
a_{i_{1},1} & \cdots & a_{i_{1},n}\\
\vdots & \cdots & \vdots\\
a_{i_{n},1} & \cdots & a_{i_{n},n}%
\end{array}
\right|  \cdot x_{i_{1}}\wedge\cdots\wedge x_{i_{n}}\text{ .}%
\]
We may phrase this result in terms of points and infer the results of the
following exercise which also contains a well-known result about the
independence of the columns of a $d\times n$ matrix and the determinants of
its $n\times n$ minors.

\begin{exercise}
Let $P_{1},\ldots,P_{n}$ be $n\leqslant d$ distinct points of $\Phi$ and for
each $j$ let $A_{j}=\left(  a_{1,j},\ldots,a_{d,j}\right)  $ be a $d$-tuple of
homogenous coordinates for $P_{j}$ in terms of some basis $\mathcal{B}$ for
$\Phi_{+}$. Let $A=\left[  a_{i,j}\right]  $ be the $d\times n$ matrix formed
by using the $A_{j}$ as its columns. Then if $\left\{  P_{1},\ldots
,P_{n}\right\}  $ is affine independent, there is a $\mathcal{B}^{\wedge}$
such that in terms of its $n$-blades the $\binom{d}{n}$ determinants of
$n\times n$ submatrices of $A$ formed in each possible way by deleting all but
$n$ rows of $A$ are Pl\"{u}cker coordinates for the affine span of $\left\{
P_{1},\ldots,P_{n}\right\}  $. On the other hand, if $\left\{  P_{1}%
,\ldots,P_{n}\right\}  $ is affine dependent, these determinants all vanish.
\end{exercise}

\subsection{Dual Description Via Annihilators}

Remembering that we are assuming that $\Phi_{+}$ has finite dimension $d$,
$\Phi_{+}^{\top}$ also has finite dimension $d$. Each subspace $\mathcal{W}$
of dimension $n$ in $\Phi_{+}$ has the annihilator $\mathcal{W}^{0}$ (the set
of elements of $\Phi_{+}^{\top}$ which annihilate every vector of
$\mathcal{W}$) that is a $\left(  d-n\right)  $-dimensional subspace of
$\Phi_{+}^{\top}$. As in Chapter 3, we identify $\left(  \Phi_{+}^{\top
}\right)  ^{\mathbf{\top}}$ with $\Phi_{+}$, and thereby are entitled to write
$\mathcal{W}^{00}=\left(  \mathcal{W}^{0}\right)  ^{0}=\mathcal{W}$. The
assignment of $\mathcal{W}^{0}$ to $\mathcal{W}$ clearly creates a one-to-one
correspondence between the subspaces of $\Phi_{+}$ and the subspaces of
$\Phi_{+}^{\top}$, so that giving $\mathcal{W}^{0}$ is equivalent to giving
$\mathcal{W}$.

$\mathcal{W}=\mathcal{W}^{00}$ describes $\mathcal{W}$ as the annihilator of
$\mathcal{W}^{0}$, which amounts to saying that $\mathcal{W}$ is the set of
those vectors in $\Phi_{+}$ which are annihilated by every element of
$\mathcal{W}^{0}$. This dual description gives $\mathcal{W}$ as the
intersection of all the hyperplanes that contain it. Verification of this is
provided by the following exercise, once we observe that the hyperplanes
through $0$ in $\Phi_{+}$ are precisely the kernels of the linear functionals
on $\Phi_{+}$.

\begin{exercise}
Let $\mathcal{W}$ be a subspace of $\Phi_{+}$. Then $\mathcal{W}^{0}$ is the
set of all elements of $\Phi_{+}^{\top}$ whose kernel contains $\mathcal{W} $.
Let $\mathcal{X}$ be the intersection of all the kernels of the elements of
$\mathcal{W}^{0}$, i. e., the intersection of all hyperplanes containing
$\mathcal{W}$. Show that $\mathcal{W}=\mathcal{X}$ by verifying that
$\mathcal{W}\subset\mathcal{X}$ and $\mathcal{X}\subset\mathcal{W}^{00}.$
\end{exercise}

We now point out and justify what is from a geometric standpoint perhaps the
most important result concerning the annihilator, namely that for any
subspaces $\mathcal{W}$ and $\mathcal{X}$ of $\Phi_{+}$ we have
\[
\left(  \mathcal{W}+\mathcal{X}\right)  ^{0}=\mathcal{W}^{0}\cap
\mathcal{X}^{\,0}\text{.}%
\]
This is the same as saying that $f\in\Phi_{+}^{\top}$ satisfies $f\left(
w+x\right)  =0$ for all $w\in\mathcal{W}$ and for all $x\in\mathcal{X}$ if and
only if $f\left(  w\right)  =0$ for all $w\in\mathcal{W}$ and $f\left(
x\right)  =0$ for all $x\in\mathcal{X}$. The truth of the ``if'' part is
obvious, and that of the ``only if'' part is nearly so. For, since $0$ is in
any subspace, it follows at once that $f\left(  w+0\right)  =0$ for all
$w\in\mathcal{W}$ and $f\left(  0+x\right)  =0$ for all $x\in\mathcal{X}$.
Interchanging $\Phi_{+}$ with $\Phi_{+}^{\top}$ then gives us also%

\[
\left(  \mathcal{W}^{0}+\mathcal{X}^{\,0}\right)  ^{0}=\mathcal{W}%
\cap\mathcal{X}\text{ .}%
\]

\begin{exercise}
Let $f_{1},\ldots,f_{n}$ be elements of $\Phi_{+}^{\top}$. Then
\[
\left\langle \left\{  f_{1},\ldots,f_{n}\right\}  \right\rangle ^{0}=\left(
\left\langle \left\{  f_{1}\right\}  \right\rangle +\cdots+\left\langle
\left\{  f_{n}\right\}  \right\rangle \right)  ^{0}=\left\{  f_{1}\right\}
^{0}\cap\cdots\cap\left\{  f_{n}\right\}  ^{0}\text{,}%
\]
and thus any subspace of $\Phi_{+}$ is the intersection of finitely many hyperplanes.
\end{exercise}

\subsection{Direct and Dual Pl\"ucker Coordinates}

Let $\left\{  t_{1},\ldots,t_{n}\right\}  $ be a basis for the subspace
$\mathcal{W}$ of $\Phi_{+}$, and let $\mathcal{A}=\left\{  t_{1},\ldots
,t_{d}\right\}  $ be any extension of this basis for $\mathcal{W}$ to a basis
for $\Phi_{+}$. Then the basis $\mathcal{A}$ has the dual $\mathcal{A}%
^{\mathbf{\top}}=\left\{  t_{1}^{\top},\ldots,t_{d}^{\top}\right\}  $, and
$\left\{  t_{n+1}^{\top},\ldots,t_{d}^{\top}\right\}  $ is a basis for
$\mathcal{W}^{0}$. Thus we have the blades proportional to $\omega=t_{1}%
\wedge\cdots\wedge t_{n}$ representing $\mathcal{W}$, and the corresponding
\textbf{annihilator blades} proportional to $\omega^{0}=t_{n+1}^{\top}%
\wedge\cdots\wedge t_{d}^{\top}$ representing $\mathcal{W}^{0}$ (and therefore
also representing $\mathcal{W}$ in a dual sense).

\begin{exercise}
Given $\omega$ and $\omega^{0}$ as above, the subspace of $\Phi_{+}$ which
$\omega$ represents is the set of all vectors $v\in\Phi_{+}$ that
simultaneously satisfy the $d-n$ equations $t_{n+1}^{\top}\left(  v\right)
=0,\ldots,t_{d}^{\top}\left(  v\right)  =0$.
\end{exercise}

For coordinatization purposes, fix a basis $\mathcal{B}=\left\{  x_{1}%
,\ldots,x_{d}\right\}  $ for $\Phi_{+}$. Construct a basis $\mathcal{B}%
^{\wedge}$ for $\bigwedge\Phi_{+}$ as we did previously, and similarly
construct a basis $\mathcal{B}^{\top\wedge}$ for $\bigwedge\Phi_{+}^{\top}$
based on the elements of the dual basis $\mathcal{B}^{\top}$. In terms of
$\mathcal{B}^{\wedge}$, Pl\"ucker coordinates stemming from blades
proportional to $\omega$ will now be called \textbf{direct (Pl\"ucker)
coordinates}, while coordinates in terms of $\mathcal{B}^{\top\wedge}$ and
stemming from blades proportional to $\omega^{0}$ will be called \textbf{dual
(Pl\"ucker) coordinates}.

Each $n$-blade $\chi$ of $\mathcal{B}^{\wedge}$ has a matching annihilator
blade $\chi^{0}\in\mathcal{B}^{\top\wedge}$ which is the product of the $d-n$
coordinate functions of the basis vectors not in $\chi$. This gives a
one-to-one match between the basis elements of $\mathcal{B}^{\wedge}$ and
$\mathcal{B}^{\top\wedge}$. We will find that, apart from an easily-determined
sign factor, direct and dual coordinates corresponding to matching basis
elements may be taken to be equal. Using the bases $\mathcal{B}^{\wedge}$ and
$\mathcal{B}^{\top\wedge}$ therefore makes it easy to convert back and forth
between direct and dual coordinates since at worst we need only multiply by
$-1$. To see how this works, let us suppose, as we did just above, that we
have a blade $\omega=t_{1}\wedge\cdots\wedge t_{n}$ and its annihilator blade
$\omega^{0}=t_{n+1}^{\top}\wedge\cdots\wedge t_{d}^{\top}$ where
$\mathcal{A}=\left\{  t_{1},\ldots,t_{d}\right\}  $ is a basis for $\Phi_{+}$.
The bases $\mathcal{A}$ and $\mathcal{A}^{\mathbf{\top}}$ are respectively
related to the coordinatization bases $\mathcal{B}$ and $\mathcal{B}^{\top}$
by an automorphism $f$ and (Section \ref{Contragredient}) its contragredient
$f^{-\top}=\left(  f^{-1}\right)  ^{\top}$ according to
\[
t_{j}=f\left(  x_{j}\right)  \text{\quad and\quad}t_{j}^{\top}=f^{-\top
}\left(  x_{j}^{\top}\right)  \text{.}%
\]

Now suppose that $f$ and $f^{-1}$ are given in terms of the elements of the
basis $\mathcal{B}$ by
\[
f\left(  x_{j}\right)  =\sum_{i=1}^{d}a_{i,j}\cdot x_{i}\text{\quad and\quad
}f^{-1}\left(  x_{j}\right)  =\sum_{i=1}^{d}\alpha_{i,j}\cdot x_{i}%
\]
and thus (Exercise \ref{CoordinatesofDualMap})
\[
f^{-\top}\left(  x_{j}^{\top}\right)  =\sum_{i=1}^{d}\alpha_{j,i}\cdot
x_{i}^{\top}\text{.}%
\]
Therefore
\[
\omega=t_{1}\wedge\cdots\wedge t_{n}=\left(  \sum_{i=1}^{d}a_{i,1}\cdot
x_{i}\right)  \wedge\cdots\wedge\left(  \sum_{i=1}^{d}a_{i,n}\cdot
x_{i}\right)
\]
and
\[
\omega^{0}=t_{n+1}^{\top}\wedge\cdots\wedge t_{d}^{\top}=\left(  \sum
_{i=1}^{d}\alpha_{n+1,i}\cdot x_{i}^{\top}\right)  \wedge\cdots\wedge\left(
\sum_{i=1}^{d}\alpha_{d,i}\cdot x_{i}^{\top}\right)  \text{.}%
\]
For distinct $i_{1},\ldots,i_{n}$ we wish to compare a coordinate
$p_{i_{1},\ldots,i_{n}}$ related to the basis element $x_{i_{1}}\wedge
\cdots\wedge x_{i_{n}}$ in $\omega$ to that of a dual coordinate $\pi
_{i_{n+1},\ldots,i_{d}}$ related to the basis element $x_{i_{n+1}}^{\top
}\wedge\cdots\wedge x_{i_{d}}^{\top}$ in $\omega^{0}$, where the two sets of
subscripts are complements in $\left\{  1,..,d\right\}  $. As our coordinates
we take the coefficients of the relevant blades in the expressions for
$\omega$ and $\omega^{0}$ above, namely the determinants
\[
p_{i_{1},\ldots,i_{n}}=\sum_{\sigma}\left(  -1\right)  ^{\sigma}%
a_{\sigma\left(  i_{1}\right)  ,1}\cdots a_{\sigma\left(  i_{n}\right)  ,n}%
\]
and
\[
\pi_{i_{n+1},\ldots,i_{d}}=\sum_{\sigma}\left(  -1\right)  ^{\sigma}%
\alpha_{n+1,\sigma\left(  i_{n+1}\right)  }\cdots\alpha_{d,\sigma\left(
i_{d}\right)  }%
\]
where each sum is over all permutations $\sigma$ of the indicated subscript
set, and as usual, $(-1)^{\sigma}=+1$ or $-1$ according as the permutation
$\sigma$ is even or odd.

The following result attributed to the noted German mathematician Carl G. J.
Jacobi (1804-1851) provides the final key.

\begin{lemma}
[Jacobi's Determinant Identity]\label{Jacobi}Fix a basis $\mathcal{B}=\left\{
x_{1},\ldots,x_{d}\right\}  $ for $\Phi_{+}$ and let $g$ be the automorphism
of $\Phi_{+}$ such that for each $j$%
\[
g\left(  x_{j}\right)  =\sum_{i=1}^{d}b_{i,j}\cdot x_{i}\text{\quad and\quad
}g^{-1}\left(  x_{j}\right)  =\sum_{i=1}^{d}\beta_{i,j}\cdot x_{i}\text{.}%
\]
Then
\[
\det\left[  c_{i,j}\right]  =\left(  \det g\right)  \det\left[  \gamma
_{i,j}\right]
\]
where $\left[  c_{i,j}\right]  $ is the $n\times n$ matrix with elements
$c_{i,j}=b_{i,j}$, $1\leqslant i,j\leqslant n$, and $\left[  \gamma
_{i,j}\right]  $ is the $\left(  d-n\right)  \times\left(  d-n\right)  $
matrix with elements $\gamma_{i,j}=\beta_{i+n,j+n}$, $1\leqslant i,j\leqslant
d-n$.
\end{lemma}

Proof: Define the self map $h$ of $\Phi_{+}$ by specifying its values on
$\mathcal{B}$ as follows:
\[
\ \ \ \ \ \ \ \ h\left(  x_{j}\right)  =\left\{
\begin{array}
[c]{lll}%
x_{j}, &  & j=1,\ldots,n,\\
g^{-1}\left(  x_{j}\right)  , &  & j=n+1,\ldots,d.
\end{array}
\right.
\]
Then
\[
g\left(  h\left(  x_{j}\right)  \right)  =\left\{
\begin{array}
[c]{lll}%
g\left(  x_{j}\right)  , &  & j=1,\ldots,n,\\
x_{j}, &  & j=n+1,\ldots,d.
\end{array}
\right.
\]
We have
\[
x_{1}\wedge\cdots\wedge x_{n}\wedge g^{-1}(x_{n+1})\wedge\cdots\wedge
g^{-1}(x_{d})=\left(  \det h\right)  \cdot x_{1}\wedge\cdots\wedge x_{d}%
\]
and it is readily established that $\det h=\det\left[  \gamma_{i,j}\right]  $.
Similarly, $\det g\circ h=\det\left[  c_{i,j}\right]  $. The lemma now follows
from the Product Theorem (Theorem \ref{ProductTheorem}). $\blacksquare$

\medskip\ 

We now construct a $g$ that we will use in the lemma to produce the result we
seek. Let $\rho$ be the permutation such that $\rho(k)=i_{k}$. Define the
automorphism $q$ of $\Phi_{+}$ by
\[
q\left(  x_{i}\right)  =x_{\rho^{-1}\left(  i\right)  }\text{.}%
\]
We apply the lemma to $g=q\circ f$ so that
\begin{align*}
g\left(  x_{j}\right)   &  =q(f(x_{j}))=\sum_{i=1}^{d}a_{i,j}\cdot q(x_{i})\\
\  &  =\sum_{i=1}^{d}a_{i,j}\cdot x_{\rho^{-1}\left(  i\right)  }=\sum
_{k=1}^{d}a_{\rho(k),j}\cdot x_{k}\\
\  &  =\sum_{k=1}^{d}a_{i_{k},j}\cdot x_{k}=\sum_{k=1}^{d}b_{k,j}\cdot x_{k}%
\end{align*}
and
\begin{align*}
g^{-1}\left(  x_{j}\right)   &  =f^{-1}\left(  q^{-1}(x_{j})\right)
=f^{-1}\left(  x_{\rho(j)}\right)  =\sum_{i=1}^{d}\alpha_{i,\rho(j)}\cdot
x_{i}\\
\  &  =\sum_{i=1}^{d}\alpha_{i,i_{j}}\cdot x_{i}=\sum_{i=1}^{d}\beta
_{i,j}\cdot x_{i}\text{.}%
\end{align*}
We observe that for this particular $g$%
\[
\det\left[  c_{i,j}\right]  =p_{i_{1},\ldots,i_{n}}\text{, }\det\left[
\gamma_{i,j}\right]  =\pi_{i_{n+1},\ldots,i_{d}}\text{.}%
\]
We also see that
\[
\det g=\left(  \det q\right)  \left(  \det f\right)  =\left(  -1\right)
^{\rho}\det f\text{.}%
\]
Hence, by the lemma
\[
p_{i_{1},\ldots,i_{n}}=\left(  -1\right)  ^{\rho}(\det f)\pi_{i_{n+1}%
,\ldots,i_{d}}%
\]
and this is the result we have been seeking. Up to the factor $\left(
-1\right)  ^{\rho}$, $\omega$ and $\omega^{0}$ may then be taken to have the
same coordinate with respect to matching basis elements, since the factor
$\det f$ may be ignored due to homogeneity.

Thus we may easily obtain a corresponding annihilator blade for any blade that
is expressed as a linear combination of basis blades, and we then have the
means to define a vector space isomorphism that sends each blade in
$\bigwedge\Phi_{+}$ to a corresponding annihilator blade in $\bigwedge\Phi
_{+}^{\top}$, as we now record.

\begin{theorem}
\label{Hodge}Let $\Phi_{+}$ have the basis $\mathcal{B}=\left\{  x_{1}%
,\ldots,x_{d}\right\}  $. Choose a basis $\mathcal{B}^{\wedge}$ for
$\bigwedge\Phi_{+}$ made up of exterior products \emph{(}including the empty
product which we take to equal $1$\emph{)} of the elements of $\mathcal{B}$.
Similarly choose a basis $\mathcal{B}^{\top\wedge}$ for $\bigwedge\Phi
_{+}^{\top}$ made up of exterior products of the elements of $\mathcal{B}%
^{\top}$. Let $H:\bigwedge\Phi_{+}\rightarrow\bigwedge\Phi_{+}^{\top}$ be the
vector space map such that for each $x_{i_{1}}\wedge\cdots\wedge x_{i_{n}}%
\in\mathcal{B}^{\wedge}$%
\[
H\left(  x_{i_{1}}\wedge\cdots\wedge x_{i_{n}}\right)  =\left(  -1\right)
^{\rho}\cdot x_{i_{n+1}}^{\top}\wedge\cdots\wedge x_{i_{d}}^{\top}%
\]
where $x_{i_{n+1}}^{\top}\wedge\cdots\wedge x_{i_{d}}^{\top}\in\mathcal{B}%
^{\top\wedge}$, $\left\{  i_{1},\ldots,i_{n}\right\}  \cup\left\{
i_{n+1},\ldots,i_{d}\right\}  =\left\{  1,\ldots,d\right\}  $, and
\[
\rho=\left(
\begin{array}
[c]{ccc}%
1 & \cdots & d\\
i_{1} & \cdots & i_{d}%
\end{array}
\right)  .
\]
Then $H$ is a vector space isomorphism that sends each blade to a
corresponding annihilator blade. $\blacksquare$
\end{theorem}

\begin{exercise}
The $H$ of the theorem above does not depend on the particular order of the
factors used to make up each element of the chosen bases $\mathcal{B}^{\wedge
}$ and $\mathcal{B}^{\top\wedge}$. That is, for any $n$ and any permutation
\[
\rho=\left(
\begin{array}
[c]{ccc}%
1 & \cdots & d\\
i_{1} & \cdots & i_{d}%
\end{array}
\right)
\]
we have
\[
H\left(  x_{i_{1}}\wedge\cdots\wedge x_{i_{n}}\right)  =\left(  -1\right)
^{\rho}\cdot x_{i_{n+1}}^{\top}\wedge\cdots\wedge x_{i_{d}}^{\top}%
\]
and therefore
\[
H^{-1}(x_{i_{n+1}}^{\top}\wedge\cdots\wedge x_{i_{d}}^{\top})=\left(
-1\right)  ^{\rho}\cdot x_{i_{1}}\wedge\cdots\wedge x_{i_{n}}\text{.}%
\]
\end{exercise}

Henceforth, unless stated otherwise, we suppose that we are employing the $H$
based on a fixed underlying basis $\mathcal{B}=\left\{  x_{1},\ldots
,x_{d}\right\}  $ with a fixed assignment of subscript labels to basis vectors.

\begin{exercise}
Use $H$ and $H^{-1}$ to show that any linear combination of a set of $\left(
d-1\right)  $-blades equals a $\left(  d-1\right)  $-blade or $0$.
\end{exercise}

\begin{exercise}
The hyperplane blade $x_{1}\wedge\cdots\wedge x_{i-1}\wedge x_{i+1}%
\wedge\cdots\wedge x_{d}$ represents the subspace $\left\{  x_{i}^{\top
}\right\}  ^{0}$ and the corresponding subflat, called a
\emph{\textbf{coordinate hyperplane}}, for which the coordinate value
corresponding to the basis vector $x_{i}$ is always $0$.
\end{exercise}

\subsection{A Dual Exterior Product}

A dual to our original ``wedge'' product $\wedge$ on $\bigwedge\Phi_{+}$ is
provided by the ``vee'' product $\vee$ that we define by
\[
\chi\vee\psi=H^{-1}\left(  H\left(  \chi\right)  \wedge H\left(  \psi\right)
\right)  .
\]
We may describe $\vee$ as being obtained by ``pulling back'' $\wedge$ through
$H$.

\begin{exercise}
Verify that $\bigwedge\Phi_{+}$ with the $\vee$ product is a legitimate vector
algebra by showing that it is just an alias of $\bigwedge\Phi_{+}^{\top}$ with
the $\wedge$ product, obtained by relabeling $\zeta$ and $\eta$ of
$\bigwedge\Phi_{+}^{\top}$ as $H^{-1}\left(  \zeta\right)  $ and
$H^{-1}\left(  \eta\right)  $, and writing $\vee$ instead of $\wedge$. What is
the unit element of this new algebra on $\bigwedge\Phi_{+}$?
\end{exercise}

The vector space $\bigwedge\Phi_{+}$ with the $\vee$ product is thus another
exterior algebra on the same space. In its dual relationship with the original
exterior algebra, hyperplane blades play the r\^ole that point blades (nonzero
vectors) play in the original. That is, the r\^oles of $\bigwedge^{1}\Phi_{+}$
and $\bigwedge^{d-1}\Phi_{+}$ become interchanged, and in fact the r\^oles of
$\bigwedge^{n}\Phi_{+}$ and $\bigwedge^{d-n}\Phi_{+}$ all become interchanged.
The $\vee$ product of independent hyperplane blades represents the
intersection of the hyperplanes that those blades represent, just as the
$\wedge$ product of independent vectors represents the subspace sum of the
subspaces that those vectors represent. The $\vee$ product of dependent
hyperplane blades is $0$. The $\vee$ product of two blades is $0$ if and only
if those blades represent subspaces that have a subspace sum smaller than the
whole space $\Phi_{+}$, just as the $\wedge$ product of two blades is $0$ if
and only if those blades represent subspaces that have an intersection larger
than the $0$ subspace. Following custom (and Grassmann's lead), we will call
the original $\wedge$ product \textbf{progressive} and the new $\vee$ product
\textbf{regressive}, terms that reflect the respective relationships of the
products to subspace sum and intersection. In expressions, the progressive
product will be given precedence so that by $u\wedge v\vee w\wedge x$ will be
meant $(u\wedge v)\vee(w\wedge x)$.

\begin{exercise}
Let $\bar x_{i}$ denote a hyperplane blade that represents $\left\{
x_{i}^{\top}\right\}  ^{0}$. Then if $\left\{  i_{k+1},\ldots,i_{d}\right\}
=\left\{  1,\ldots,d\right\}  \smallsetminus\left\{  i_{1},\ldots
,i_{k}\right\}  $, $\left\{  x_{i_{1}}^{\top}\right\}  ^{0}\cap\cdots
\cap\left\{  x_{i_{k}}^{\top}\right\}  ^{0}$ is represented by $\bar x_{i_{1}%
}\vee\cdots\vee\bar x_{i_{k}}=\pm x_{i_{k+1}}\wedge\cdots\wedge x_{i_{d}}$.
\end{exercise}

\subsection{Results Concerning Zero Coordinates}

The basic significance of a zero coordinate is given in the following result.

\begin{proposition}
In the blade $\beta$, the Pl\"ucker coordinate corresponding to the basis
element $\xi=x_{i_{1}}\wedge\cdots\wedge x_{i_{n}}$ is zero if and only if the
flat represented by $\beta$ intersects the flat represented by the blade
$\bar\xi=x_{i_{n+1}}\wedge\cdots\wedge x_{i_{d}}$, where $\overset{\_}%
{I}=\left\{  i_{n+1},\ldots,i_{d}\right\}  $ is the complement of $I=\left\{
i_{1},\ldots,i_{n}\right\}  $ in $D=\left\{  1,..,d\right\}  $.
\end{proposition}

Proof: $\psi\wedge\bar{\xi}=0$ for each basis $n$-blade $\psi=x_{k_{1}}%
\wedge\cdots\wedge x_{k_{n}}$ with the exception of $\xi$. (Each $n$-element
subset of $D$ except $I$ meets $\overset{\_}{I}$ since the only $n$-element
subset of $I$ is $I$ itself). Therefore $\beta\wedge\bar{\xi}$ is zero or not
according as $\beta$'s Pl\"{u}cker coordinate corresponding to the basis
element $\xi$ is zero or not. $\blacksquare$

\begin{exercise}
The proposition above implies, as expected, that the point represented by the
basis vector $x_{1}$ lies in the coordinate hyperplanes represented by
$\left\{  x_{i}^{\top}\right\}  ^{0},$ $i=2,\ldots,d$.
\end{exercise}

The previous result serves as a prelude to the following generalization that
treats the case where there is a particular family of zero coordinates.

\begin{proposition}
Let $1\leqslant k\leqslant n\leqslant d$. The flat represented by the subspace
$\mathcal{X}=\left\{  x_{i_{1}}^{\top}\right\}  ^{0}\cap\cdots\cap\left\{
x_{i_{k}}^{\top}\right\}  ^{0}$ intersects the flat represented by the
$n$-dimensional subspace $\mathcal{Y}$ in a flat of dimension at least $n-k$
if and only if the flat represented by $\mathcal{Y}$ has a zero Pl\"{u}cker
coordinate corresponding to each of the basis elements $x_{j_{1}}\wedge
\cdots\wedge x_{j_{n}}$ such that $\left\{  i_{1},\ldots,i_{k}\right\}
\subset\left\{  j_{1},\ldots,j_{n}\right\}  $.
\end{proposition}

Proof: Suppose first that the flat represented by $\mathcal{Y}$ has a zero
Pl\"{u}cker coordinate corresponding to each of the basis elements $x_{j_{1}%
}\wedge\cdots\wedge x_{j_{n}}$ such that $\left\{  i_{1},\ldots,i_{k}\right\}
\subset\left\{  j_{1},\ldots,j_{n}\right\}  $. The flats represented by
$\mathcal{X}$ and $\mathcal{Y}$ do intersect since by the previous proposition
the flat represented by $\mathcal{Y}$ intersects (several, perhaps, but one is
enough) a flat represented by a subspace $\left\{  x_{j_{1}}^{\top}\right\}
^{0}\cap\cdots\cap\left\{  x_{j_{n}}^{\top}\right\}  ^{0}$ such that $\left\{
i_{1},\ldots,i_{k}\right\}  \subset\left\{  j_{1},\ldots,j_{n}\right\}  $ and
must therefore intersect the flat represented by its factor $\mathcal{X}$.
Denote by $\bar{x}_{m}$ the basis element that represents the hyperplane
$\left\{  x_{m}^{\top}\right\}  ^{0}$. Then $\mathcal{X}$ is represented by
the regressive product $\xi=\bar{x}_{i_{1}}\vee\cdots\vee\bar{x}_{i_{k}}$, and
each basis $n$-blade $x_{j_{1}}\wedge\cdots\wedge x_{j_{n}}$ is a regressive
product of the $d-n$ hyperplane blades $\bar{x}_{m_{n+1}},\ldots,\bar
{x}_{m_{d}}$ such that $m_{n+1},\ldots,m_{d}\in\left\{  1,\ldots,d\right\}
\smallsetminus\left\{  j_{1},\ldots,j_{n}\right\}  $. Hence the regressive
product $\eta\vee\xi$, where $\eta$ is a blade that represents the flat that
$\mathcal{Y}$ represents, consists of a sum of terms of the form
\[
p_{j_{1},\ldots,j_{n}}\cdot\bar{x}_{m_{n+1}}\vee\cdots\vee\bar{x}_{m_{d}}%
\vee\bar{x}_{i_{1}}\vee\cdots\vee\bar{x}_{i_{k}}%
\]
and upon careful scrutiny of these terms we find that each is zero, but only
because the flat represented by $\mathcal{Y}$ has a zero Pl\"{u}cker
coordinate corresponding to each of the basis elements $x_{j_{1}}\wedge
\cdots\wedge x_{j_{n}}$ such that $\left\{  i_{1},\ldots,i_{k}\right\}
\subset\left\{  j_{1},\ldots,j_{n}\right\}  $. The terms such that $\left\{
j_{1},\ldots,j_{n}\right\}  $ excludes at least one element of $\left\{
i_{1},\ldots,i_{k}\right\}  $ have each excluded $i$ appearing as an $m$,
while the terms such that $\left\{  i_{1},\ldots,i_{k}\right\}  \subset
\left\{  j_{1},\ldots,j_{n}\right\}  $ have no $i$ appearing as an $m$. Hence
$\eta\vee\xi=0$ if and only if the flat represented by $\mathcal{Y}$ has a
zero Pl\"{u}cker coordinate corresponding to each of the basis elements
$x_{j_{1}}\wedge\cdots\wedge x_{j_{n}}$ such that $\left\{  i_{1},\ldots
,i_{k}\right\}  \subset\left\{  j_{1},\ldots,j_{n}\right\}  $. But $\eta
\vee\xi=0$ if and only if the dimension of $\mathcal{X}+\mathcal{Y}$ is
strictly less than $d.$ The dimension of $\mathcal{X}$ is $d-k$ and the
dimension of $\mathcal{Y}$ is $n$ so that by Grassmann's Relation (Corollary
\ref{GrassmannRel})
\begin{equation}
\left(  d-k\right)  +n-\dim\left(  \mathcal{X}\cap\mathcal{Y}\right)
=\dim\left(  \mathcal{X}+\mathcal{Y}\right)  <d\tag{$\ast$}%
\end{equation}
and hence
\[
n-k<\dim\left(  \mathcal{X}\cap\mathcal{Y}\right)  \text{.}%
\]
The flat represented by $\mathcal{X}\cap\mathcal{Y}$ thus has dimension at
least $n-k$.

Suppose on the other hand that $\mathcal{X}$ and $\mathcal{Y}$ intersect in a
flat of dimension at least $n-k$. Then reversing the steps above, we recover
Equation $\left(  \ast\right)  $ so that $\eta\vee\xi=0$ which can hold only
if the flat represented by $\mathcal{Y}$ has a zero Pl\"{u}cker coordinate
corresponding to each of the basis elements $x_{j_{1}}\wedge\cdots\wedge
x_{j_{n}}$ such that $\left\{  i_{1},\ldots,i_{k}\right\}  \subset\left\{
j_{1},\ldots,j_{n}\right\}  $. $\blacksquare$

\begin{exercise}
From Grassmann's Relation alone, show that the subspaces $\mathcal{X}$ and
$\mathcal{Y}$ of the proposition above always intersect in a subspace of
dimension at least $n-k$. Hence, without requiring any zero Pl\"ucker
coordinates, we get the result that the flats represented by $\mathcal{X}$ and
$\mathcal{Y}$ must intersect in a flat of dimension at least $n-k-1$.
\end{exercise}

\subsection{Some Examples in Space\label{Examples}}

Let our flat of interest be $\Phi=\mathcal{F}^{3}$ for some field
$\mathcal{F}$ so that $\Phi_{+}$ is a vector space of dimension $d=4$. The
basis $\mathcal{B}$ will be $\left\{  x_{1},x_{2},x_{3},x_{4}\right\}  $ where
$\left\{  x_{2},x_{3},x_{4}\right\}  $ is a basis for $\Phi_{0}$ and $\Phi
_{1}=x_{1}+\Phi_{0}$ is $\Phi$ in its inflated context. In the following table
we show how the elements of $\mathcal{B}^{\wedge}$ and $\mathcal{B}%
^{\top\wedge}$ are connected by the isomorphism $H\mathrm{,}$ where we have
only shown the subscript sequences (with $\varnothing$ indicating the empty
subscript sequence corresponding to the scalar $1$).
\[%
\begin{tabular}
[c]{lllllllllll}%
$\varnothing$ & $\rightarrow$ & $+1234$ & \ \ \ \ \ \ \ \  & $12$ &
$\rightarrow$ & $+34$ & \ \ \ \ \ \ \ \  & $123$ & $\rightarrow$ & $+4$\\
&  &  &  & $13$ & $\rightarrow$ & $-24$ &  & $124$ & $\rightarrow$ & $-3$\\
$1$ & $\rightarrow$ & $+234$ &  & $14$ & $\rightarrow$ & $+23$ &  & $134$ &
$\rightarrow$ & $+2$\\
$2$ & $\rightarrow$ & $-134$ &  & $23$ & $\rightarrow$ & $+14$ &  & $234$ &
$\rightarrow$ & $-1$\\
$3$ & $\rightarrow$ & $+124$ &  & $24$ & $\rightarrow$ & $-13$ &  &  &  & \\
$4$ & $\rightarrow$ & $-123$ &  & $34$ & $\rightarrow$ & $+12$ &  & $1234$ &
$\rightarrow$ & $+\varnothing$%
\end{tabular}
\]
We will generally use subscript sequences to describe the blades made up of
basis vectors. Thus $3\cdot143$ will indicate either $3\cdot x_{1}\wedge
x_{4}\wedge x_{3}$ or $3\cdot x_{1}^{\top}\wedge x_{4}^{\top}\wedge
x_{3}^{\top}$ as dictated by context. In progressive products we will
generally omit the $\wedge$ sign, and in regressive products we will generally
replace the $\vee$ sign with a ``$.$'', so that $(1+2)(1+3)=13+23$ and
$234.123=H^{-1}\left(  -14\right)  =-23$. We indicate the plane blade obtained
by omitting any individual factor from $1234$ by using the subscript of the
omitted factor with a bar over it as in $234=\overline{1}$, $134=\overline{2}%
$, etc. The elements of $\mathcal{B}^{\wedge}$ may be expressed as regressive
products of such plane blades as follows, where the empty regressive product
of such plane blades is indicated by $\overline{\varnothing}$.
\[%
\begin{tabular}
[c]{lllllllllll}%
$1234$ & $=$ & $\overline{\varnothing}$ & \ \ \ \ \ \ \ \  & $34$ & $=$ &
$-\overline{1}.\overline{2}$ & \ \ \ \ \ \ \ \  & $4$ & $=$ & $-\overline
{1}.\overline{2}.\overline{3}$\\
&  &  &  & $24$ & $=$ & $-\overline{1}.\overline{3}$ &  & $3$ & $=$ &
$-\overline{1}.\overline{2}.\overline{4}$\\
$234$ & $=$ & $\overline{1}$ &  & $23$ & $=$ & $-\overline{1}.\overline{4}$ &
& $2$ & $=$ & $-\overline{1}.\overline{3}.\overline{4}$\\
$134$ & $=$ & $\overline{2}$ &  & $14$ & $=$ & $-\overline{2}.\overline{3}$ &
& $1$ & $=$ & $-\overline{2}.\overline{3}.\overline{4}$\\
$124$ & $=$ & $\overline{3}$ &  & $13$ & $=$ & $-\overline{2}.\overline{4}$ &
&  &  & \\
$123$ & $=$ & $\overline{4}$ &  & $12$ & $=$ & $-\overline{3}.\overline{4}$ &
& $\varnothing$ & $=$ & $\overline{1}.\overline{2}.\overline{3}.\overline{4}$%
\end{tabular}
\]
We shall often find it convenient to engage in a harmless abuse of language by
referring to blades as the flats they represent.

We now present our first example. Consider the plane
\[
\left(  1+2\right)  \left(  1+3\right)  \left(  1+4\right)
=134+214+231+234=\overline{2}-\overline{3}+\overline{4}+\overline{1}.
\]
We wish to compute its intersection with the plane $\overline{2}$, which is
then the line
\[
\left(  \overline{2}-\overline{3}+\overline{4}+\overline{1}\right)
.\overline{2}=-\overline{3}.\overline{2}+\overline{4}.\overline{2}%
+\overline{1}.\overline{2}=-14+13-34.
\]
We now have an answer in the form of a line blade, but we may wish also to
know at least two points that determine the line. A point on a line can be
obtained by taking its regressive product with any coordinate plane that does
not contain it. The result of taking the regressive product of our line with
each coordinate plane is as follows.
\[%
\begin{tabular}
[c]{l}%
$\left(  -\overline{3}+\overline{4}+\overline{1}\right)  .\overline
{2}.\overline{1}=-\overline{3}.\overline{2}.\overline{1}+\overline
{4}.\overline{2}.\overline{1}=4-3$\\
$\left(  -\overline{3}+\overline{4}+\overline{1}\right)  .\overline
{2}.\overline{2}=0$\\
$\left(  -\overline{3}+\overline{4}+\overline{1}\right)  .\overline
{2}.\overline{3}=\overline{4}.\overline{2}.\overline{3}+\overline{1}%
.\overline{2}.\overline{3}=-1-4$\\
$\left(  -\overline{3}+\overline{4}+\overline{1}\right)  .\overline
{2}.\overline{4}=-\overline{3}.\overline{2}.\overline{4}+\overline
{1}.\overline{2}.\overline{4}=-1-3$.
\end{tabular}
\]
Thus the line pierces $\overline{1}$ (the plane at $\infty$) at $4-3$, it
pierces $\overline{3}$ at $-1-4$, and it pierces $\overline{4}$ at $-1-3$, but
it is wholly contained in $\overline{2}$.

Let us now consider two lines in space. Let one be $\lambda=\left(  1\right)
(1+2+3+4)$ and the other be $\mu=\left(  1+2\right)  \left(  1+3+4\right)  .$
They intersect because their progressive product is
\[
\left(  12+13+14\right)  \left(  13+14+21+23+24\right)  =1324+1423=0\text{.}%
\]
We can compute their affine sum by considering progressive products of one
with the individual progressive point factors of the other. Multiplying $1$
and $\mu$ gives $123+124=\overline{4}+\overline{3}$, the affine sum we seek.
Multiplying $1+2+3+4$ and $\mu$ also gives $123+124$. The regressive product
of the two lines is $0$, confirming that their affine sum fails to be all of
space. We can compute their intersection point by considering regressive
products of one with the individual regressive plane factors of the other. We
factor $\lambda$ into plane blades:
\[
\lambda=\left(  12+13+14\right)  =-\left(  \overline{3}.\overline{4}%
+\overline{2}.\overline{4}+\overline{2}.\overline{3}\right)  =-\left(
\overline{3}+\overline{2}\right)  .\left(  \overline{3}+\overline{4}\right)
\text{.}%
\]
$\mu$ may be written as $\overline{2}.\overline{4}+\overline{2}.\overline
{3}-\overline{3}.\overline{4}+\overline{1}.\overline{4}+\overline{1}%
.\overline{3}$ and taking its regressive product with $\overline{3}%
+\overline{2}$ gives
\[
\overline{2}.\overline{4}.\overline{3}+\overline{1}.\overline{4}.\overline
{3}-\overline{3}.\overline{4}.\overline{2}+\overline{1}.\overline{4}%
.\overline{2}+\overline{1}.\overline{3}.\overline{2}=2\cdot1+2+3+4\text{,}%
\]
whereas taking its regressive product with $\overline{3}+\overline{4}$ gives
\[
\overline{2}.\overline{4}.\overline{3}+\overline{1}.\overline{4}.\overline
{3}+\overline{2}.\overline{3}.\overline{4}+\overline{1}.\overline{3}%
.\overline{4}=0.
\]
$\mu$ is therefore contained in the plane $\overline{3}+\overline{4}$ (as we
already know since we found above that this same plane is the affine sum of
$\lambda$ and $\mu$), but pierces the plane $\overline{3}+\overline{2}$ at the
point $2\cdot1+2+3+4$ which must also be the point of intersection of the two lines.

\begin{exercise}
Find the intersection and the affine sum of the two lines $\left(
\bar1\right)  .(\bar1+\bar2+\bar3+\bar4)$ and $\left(  \bar1+\bar2\right)
.\left(  \bar1+\bar3+\bar4\right)  $.
\end{exercise}

\begin{exercise}
Suppose that the regressive product of two blades equals a nonzero scalar.
What is the meaning of this both for the corresponding subspaces of $\Phi_{+}$
and for the corresponding flats of $\Phi$? Give an example where the blades
both are lines in space.
\end{exercise}

\subsection{Factoring a Blade}

As illustrated by the examples of the previous section, there are instances
where we have a blade expressed in expanded form and we instead want to
express it, up to a nonzero scalar factor, as the exterior product of vectors.
That is, we seek an independent set of vectors that span the subspace
represented by the blade. It is always possible to find such a factorization.
One method for doing so will now be presented. We start by introducing the
\textbf{extended coordinate array} of a blade based on a given set of
Pl\"ucker coordinates for it. This is a full skew-symmetric array of values
$P_{j_{1},\ldots,j_{n}}$ which contains the given set of Pl\"ucker coordinates
along with either $\pm$ duplicates of them or zeroes. This $d\times
\cdots\times d$ ($n$ factors of $d$) array is the generalization of a
skew-symmetric $d\times d$ matrix. Given the Pl\"ucker coordinate
$p_{i_{1},\ldots,i_{n}}$ related to the basis $n$-blade $x_{i_{1}}\wedge
\cdots\wedge x_{i_{n}}$, the values $P_{j_{1},\ldots,j_{n}}$ at the array
positions that correspond to this coordinate (the positions with subscripts
that are the permutations of the subscripts of the given coordinate) are given
by
\[
P_{\sigma\left(  i_{1}\right)  ,\ldots,\sigma\left(  i_{n}\right)  }=\left(
-1\right)  ^{\sigma}p_{i_{1},\ldots,i_{n}}%
\]
for each permutation $\sigma$ of $i_{1},\ldots,i_{n}$. Entries in positions
that have any equal subscripts are set to $0$. The $d$-tuple of array values
obtained by varying a given subscript in order from $1$ to $d$, leaving all
other subscripts fixed at chosen constants, is called a \textbf{file},
generalizing the concept of a row or column of a matrix. Any entry in the
array is contained in exactly $n$ files. Here now is the factorization result.

\begin{proposition}
[Blade Factorization]Given a set of Pl\"ucker coordinates for an $n$-blade,
the $n$ files that contain a given nonzero entry of the corresponding extended
coordinate array constitute a factorization when we regard these files as
vector coordinate $d $-tuples.
\end{proposition}

Proof: We may assume that the given Pl\"{u}cker coordinates were obtained as
minor determinants in the manner described in Section \ref{coords} above.
Using the notation of that section, each extended coordinate array value may
be written in determinant form as
\[
P_{j_{1},\ldots,j_{n}}=\left|
\begin{array}
[c]{ccc}%
a_{j_{1},1} & \cdots & a_{j_{1},n}\\
\vdots & \cdots & \vdots\\
a_{j_{n},1} & \cdots & a_{j_{n},n}%
\end{array}
\right|  \text{ .}%
\]
We assume without loss of generality that the given nonzero entry is
$P_{1,\ldots,n}$. We will first show that each of the $n$ files containing
$P_{1,\ldots,n}$ is in the span of the columns $A_{j}$ of the matrix
$A=[a_{i,j}]$, so that each is a coordinate $d$-tuple of a vector in the
subspace represented by the blade with the given Pl\"{u}cker coordinates. Let
the files containing $P_{1,\ldots,n}$ form the columns $B_{j}$ of the $d\times
n$ matrix $B=[b_{i,j}]$ so that the file forming the column $B_{j}$ has the
elements
\[
b_{i,j}=P_{1,\ldots,j-1,i,j+1,\ldots,n}=\left|
\begin{array}
[c]{lllllll}%
a_{1,1} &  &  & \cdots &  &  & a_{1,n}\\
\vdots &  &  &  &  &  & \vdots\\
a_{j-1,1} &  &  & \cdots &  &  & a_{j-1,n}\\
a_{i,1} &  &  & \cdots &  &  & a_{i,n}\\
a_{j+1,1} &  &  & \cdots &  &  & a_{j+1,n}\\
\vdots &  &  &  &  &  & \vdots\\
a_{n,1} &  &  & \cdots &  &  & a_{n,n}%
\end{array}
\right|  =\sum_{k=1}^{n}c_{j,k}\,a_{i,k}\text{,}%
\]
where the last expression follows upon expanding the determinant about row $j$
containing the $a_{i,k}$. Note that the coefficients $c_{j,k}$ are independent
of $i$, and hence
\[
B_{j}=\sum_{k=1}^{n}c_{j,k}\,A_{k}%
\]
showing that the file forming the column $B_{j}$ is in the span of the columns
of $A$.

To verify the independence of the vectors for which the $B_{j}$ are coordinate
$d$-tuples, let us examine the matrix $B$. Observe that the top $n$ rows of
$B$ are $P_{1,\ldots,n}$ times the $n\times n$ identity matrix. Hence $B$
contains an $n\times n$ minor with a nonzero determinant, and the vectors
$w_{j}$ for which the $B_{j}$ are coordinate $d$-tuples must therefore be
independent, since it follows from what we found in Section \ref{coords} that
$w_{1}\wedge\cdots\wedge w_{n}$ is not zero. $\blacksquare$

\begin{exercise}
For each nonzero extended coordinate, apply the method of the proposition
above to the Pl\"ucker coordinates of the line obtained in the first example
of the previous section. Characterize the sets of nonzero extended coordinates
that are sufficient to use to yield all essentially different obtainable
factorizations. Compare your results to the points that in the example were
obtained by intersecting with coordinate planes.
\end{exercise}

\begin{exercise}
Use the method of the proposition above to factor the line $\lambda$ of the
second example of the previous section into the regressive product of plane
blades. Do this first by factoring $H\left(  \lambda\right)  $ in
$\bigwedge\Phi_{+}^{\top}$ and then carrying the result back to $\bigwedge
\Phi_{+}$. Then also do it by factoring $\lambda$ directly in the exterior
algebra on $\bigwedge\Phi_{+}$ that the $\vee$ product induces \emph{(}where
the ``vectors'' are the plane blades\emph{)}.
\end{exercise}

\subsection{Algorithm for Sum or Intersection}

Given the blades $\alpha$ and $\beta$, $\alpha\wedge\beta$ only represents
their affine sum when the flats represented have an empty intersection.
Otherwise the non-blade $0$ is the result. Dually, $\alpha\vee\beta$ only
represents their intersection when the flats represented have a full affine
sum (all of $\Phi,$ that is). Otherwise the non-blade $0$ is again the result.
Blades that yield a $0$ result when multiplied contain redundancy and pose for
us the extra challenge of coping with that redundancy. As we have just seen in
the previous section, a factorization for a blade is always readily
obtainable, and we can exploit this to construct a resultant blade that omits
redundant factors. We give an algorithm for doing this in the case where the
affine sum is sought, and the same algorithm can be used on $H\left(
\alpha\right)  $ and $H\left(  \beta\right)  $ (or the similar $\vee$-based
algorithm on $\alpha$ and $\beta$ expressed in terms of basic hyperplane
blades) when the intersection is sought. Suppose that $\beta$ has been
factored as
\[
\beta=v_{1}\wedge\cdots\wedge v_{n}%
\]
and start off with $\gamma_{0}=\alpha$. Then calculate successively
$\gamma_{1},\cdots,\gamma_{n}$ according to
\[
\gamma_{i}=\left\{
\begin{tabular}
[c]{l}%
$\gamma_{i-1}\wedge v_{i}$ if this is $\neq0$,\\
$\gamma_{i-1}$ otherwise.
\end{tabular}
\ \right.
\]
Then $\gamma_{n}$ is the desired blade that represents the affine sum. This
algorithm is the obvious one based on omitting vectors that can already be
expressed in terms of the previously included ones that couldn't be so
expressed. Which vectors are omitted, but not their number, depends on their
labeling, of course. There are similarities here to the problems of finding
the least common multiple and greatest common divisor of positive whole
numbers, including that an algorithm seems to be required.

\subsection{Adapted Coordinates and $\vee$ Products}

We continue to employ the fixed underlying basis $\mathcal{B}=\left\{
x_{1},\ldots,x_{d}\right\}  $ for the purposes of defining the isomorphism
$H:\bigwedge\Phi_{+}\rightarrow\bigwedge\Phi_{+}^{\top}$ and the resulting
regressive product that it produces. However, simultaneously using other bases
can sometimes prove useful, particularly in the case where we are dealing with
a subspace of $\Phi_{+}$ and we use coordinates adapted to a basis for that
subspace. We will find that doing this will lead to a useful new view of the
regressive product.

We start by considering finding coordinates for an $n$-blade $\gamma$ that
represents a subspace $\mathcal{W}$ that is itself a subspace of the subspace
$\mathcal{U}$ of $\Phi_{+}$. We want these coordinates adapted to the basis
$\left\{  u_{1},\ldots,u_{l}\right\}  $ for the $l$-dimensional subspace
$\mathcal{U}$. The way we will describe these coordinates involves some new notation.

Let $I=\left\{  i_{1},\ldots,i_{n}\right\}  \subset\left\{  1,\ldots
,l\right\}  $ and let $\overline{I}=\left\{  i_{n+1},\ldots,i_{l}\right\}  $
denote its complement in $\left\{  1,\ldots,l\right\}  $. Note that this
implies that $i_{j}=i_{k}$ only if $j=k$. Note also that the elements of both
$I$ and $\overline{I}$ are shown as having subscript labels that imply choices
have been made of an order for each in which their elements may be presented
by following subscript numerical order. Thus we will treat $I$ and
$\overline{I} $ as labeled sets that split $\left\{  1,\ldots,l\right\}  $
into two pieces (one of which could be empty, however). The order choices
$i_{1}<\cdots<i_{n}$ and $i_{n+1}<\cdots<i_{l}$ could be specified if one
wishes. However, it will only be important that the relevant $n$-element
subsets and their complements be given labelings that remain fixed throughout
a particular discussion. The particular labelings will be unimportant. Based
on this foundation, we now describe our new notation. We will denote the two
permutations
\[
\left(
\begin{array}
[c]{ccc}%
1 & \cdots & l\\
i_{1} & \cdots & i_{l}%
\end{array}
\right)  \hspace{0.25in}\text{and}\hspace{0.25in}\left(
\begin{tabular}
[c]{llllll}%
$1$ & $\cdots$ & $l-n$ & $l-n+1$ & $\cdots$ & $l$\\
$i_{n+1}$ & $\cdots$ & $i_{l}$ & $i_{1}$ & $\cdots$ & $i_{n}$%
\end{tabular}
\right)
\]
by $I\overline{I}$ and $\overline{I}I$ respectively. Also we will denote
$u_{i_{1}}\wedge\cdots\wedge u_{i_{n}}$ by $u_{I}$, and $u_{i_{n+1}}%
\wedge\cdots\wedge u_{i_{l}}$ by $u_{\overline{I}}$.

Since $\gamma$ is an $n$-blade that represents a subspace of the vector space
$\mathcal{U}$ that has the basis $\left\{  u_{1},\ldots,u_{l}\right\}  $, we
can expand $\gamma$ in terms of the coordinates $\left(  -1\right)
^{\overline{I}I}a_{I}$ adapted to $\left\{  u_{1},\ldots,u_{l}\right\}  $ as
\[
\gamma=\sum_{I}\left(  -1\right)  ^{\overline{I}I}a_{I}\cdot u_{I}%
\]
where the sum is over all subsets $I\subset\left\{  1,\ldots,l\right\}  $ such
that $\left|  I\right|  =n$, where $\left|  I\right|  $ denotes the number of
elements of the set $I.$ Each $a_{I}$ is determined by
\[
a_{I}\cdot u_{1}\wedge\cdots\wedge u_{l}=u_{\overline{I}}\wedge\gamma
\]
since if $J$ is a particular one of the $I$s in the sum above for $\gamma$,
\[
u_{\overline{J}}\wedge\gamma=\sum_{I}\left(  -1\right)  ^{\overline{I}I}%
a_{I}\cdot u_{\overline{J}}\wedge u_{I}=\left(  -1\right)  ^{\overline{J}%
J}a_{J}\cdot u_{\overline{J}}\wedge u_{J}=a_{J}\cdot u_{1}\wedge\cdots\wedge
u_{l}\text{.}%
\]
A similar result is obtained by using $\left(  -1\right)  ^{I\overline{I}}$
instead of $\left(  -1\right)  ^{\overline{I}I}$ in the expansion for $\gamma$.

We are now ready to begin to apply this to the intersection of a pair of
subspaces in the case where that intersection corresponds to a nonzero
regressive product. Thus let $\mathcal{U}=\left\langle \left\{  u_{1}%
,\ldots,u_{l}\right\}  \right\rangle $ and $\mathcal{V}=\left\langle \left\{
v_{1},\ldots,v_{m}\right\}  \right\rangle $ be subspaces of $\Phi_{+}$ of
dimension $l$ and $m$ respectively. Let $\mathcal{W}=\mathcal{U}%
\cap\mathcal{V}$ and let $\mathcal{U}+\mathcal{V}=\Phi_{+}$. Then by
Grassmann's Relation (Corollary \ref{GrassmannRel}), $\mathcal{W}$ has
dimension $n=l+m-d$ (so that, of course, $l+m\geqslant d$ since $n\geqslant0
$). Let $\gamma$ be a blade that represents $\mathcal{W}$. Then the $n$-blade
$\gamma$ has the expansions
\[
\gamma=\sum_{I\subset\left\{  1,\ldots,l\right\}  }\left(  -1\right)
^{\overline{I}I}a_{I}\cdot u_{I}=\sum_{I\subset\left\{  1,\ldots,m\right\}
}\left(  -1\right)  ^{I\overline{I}}b_{I}\cdot v_{I}%
\]
where the one adapted to the $v$s is intentionally set up using $\left(
-1\right)  ^{I\overline{I}}$ rather than $\left(  -1\right)  ^{\overline{I}I}%
$, and the sums are over only those subsets $I$ such that $\left|  I\right|
=n$. The $a_{I}$ and $b_{I}$ are determined by
\[
\left.  a_{I}\cdot u_{1}\wedge\cdots\wedge u_{l}=u_{\overline{I}}\wedge
\gamma\right.  \hspace{0.25in}\text{and}\hspace{0.25in}\left.  b_{I}\cdot
v_{1}\wedge\cdots\wedge v_{m}=\gamma\wedge v_{\overline{I}}\right.  \text{.}%
\]
We can get some useful expressions for the $a_{I}$ and $b_{I}$ by first
creating some blades that represent $\mathcal{U}$ and $\mathcal{V}$ and
contain $\gamma$ as a factor. To simplify our writing, denote $u_{1}%
\wedge\cdots\wedge u_{l}$ by $\alpha$ and $v_{1}\wedge\cdots\wedge v_{m}$ by
$\beta$. Supposing that $\gamma=w_{1}\wedge\cdots\wedge w_{n}$, the
Replacement Theorem (Theorem \ref{ReplacementTheorem}) guarantees us blades
$\widetilde{\alpha}=u_{j_{1}}\wedge\cdots\wedge u_{j_{l-n}}$ and
$\widetilde{\beta}=v_{k_{1}}\wedge\cdots\wedge v_{k_{m-n}}$ such that for some
nonzero scalars $a$ and $b$ we have $a\cdot\alpha=\widetilde{\alpha}%
\wedge\gamma$ and $b\cdot\beta=\gamma\wedge\widetilde{\beta}$. Note that
$\widetilde{\alpha}\wedge\gamma\wedge\widetilde{\beta}$ is a $d$-blade because
it is the product of vectors that span $\mathcal{U}+\mathcal{V}=\Phi_{+}$. We
find then that
\[
\left.  a_{I}\cdot\alpha\wedge\widetilde{\beta}=b\cdot u_{\overline{I}}%
\wedge\beta\right.  \hspace{0.25in}\text{and}\hspace{0.25in}\left.  b_{I}%
\cdot\widetilde{\alpha}\wedge\beta=a\cdot\alpha\wedge v_{\overline{I}%
}\right.
\]
and since
\[
a\cdot\alpha\wedge\widetilde{\beta}=\widetilde{\alpha}\wedge\gamma
\wedge\widetilde{\beta}=b\cdot\widetilde{\alpha}\wedge\beta
\]
we then have
\[
\left.  \frac{a_{I}}{ab}\cdot\widetilde{\alpha}\wedge\gamma\wedge
\widetilde{\beta}=u_{\overline{I}}\wedge\beta\right.  \hspace{0.25in}%
\text{and}\hspace{0.25in}\left.  \frac{b_{I}}{ab}\cdot\widetilde{\alpha}%
\wedge\gamma\wedge\widetilde{\beta}=\alpha\wedge v_{\overline{I}}\right.
\text{.}%
\]
$\widetilde{\alpha}\wedge\gamma\wedge\widetilde{\beta}$, $u_{\overline{I}%
}\wedge\beta$, and $\alpha\wedge v_{\overline{I}}$ are each the progressive
product of $d$ vectors and therefore are scalar multiples of the progressive
product of all the vectors of any basis for $\Phi_{+}$ such as the fixed
underlying basis $\mathcal{B}=\left\{  x_{1},\ldots,x_{d}\right\}  $ that we
use for the purposes of defining the isomorphism $H:\bigwedge\Phi
_{+}\rightarrow\bigwedge\Phi_{+}^{\top}$ and the resulting regressive product
that it produces. We might as well also then use $\mathcal{B}$ for extracting
scalar multipliers in the above equations. We therefore define \textbf{[}%
$\left.  \xi\right.  $\textbf{]}, the \textbf{bracket} of the progressive
product $\xi$ of $d$ vectors, via the equation $\xi=$\textbf{[}$\left.
\xi\right.  $\textbf{]}$\cdot x_{1}\wedge\cdots\wedge x_{d}$. \textbf{[}%
$\left.  \xi\right.  $\textbf{]} is zero when $\xi=0$ and is nonzero otherwise
(when $\xi$ is a $d$-blade). Thus we have
\[
\left.  a_{I}=\frac{ab}{\text{\textbf{[}}\left.  \widetilde{\alpha}%
\wedge\gamma\wedge\widetilde{\beta}\right.  \text{\textbf{]}}}\cdot
\text{\textbf{[}}\left.  u_{\overline{I}}\wedge\beta\right.  \text{\textbf{]}%
}\right.  \hspace{0.25in}\text{and}\hspace{0.25in}\left.  b_{I}=\frac{ab}%
{\text{\textbf{[}}\left.  \widetilde{\alpha}\wedge\gamma\wedge\widetilde
{\beta}\right.  \text{\textbf{]}}}\cdot\text{\textbf{[}}\left.  \alpha\wedge
v_{\overline{I}}\right.  \text{\textbf{]}}\right.  \text{.}%
\]
Dropping the factor that is independent of $I$, we find that we have the two
equal expressions%
\[
\sum_{\substack{I\subset\left\{  1,\ldots,l\right\}  \\\left|  I\right|  =n
}}\left(  -1\right)  ^{\overline{I}I}\mathbf{[}\left.  u_{\overline{I}}%
\wedge\beta\right.  \mathbf{]}\cdot u_{I}=\sum_{\substack{I\subset\left\{
1,\ldots,m\right\}  \\\left|  I\right|  =n}}\left(  -1\right)  ^{I\overline
{I}}\mathbf{[}\left.  \alpha\wedge v_{\overline{I}}\right.  \mathbf{]}\cdot
v_{I}\;\;\left(  \ast\ast\right)
\]
for an $n$-blade that represents $\mathcal{U}\cap\mathcal{V}$ and which
therefore must be proportional to $\alpha\vee\beta$. The above development of
($\ast\ast$) follows Garry Helzer's online class notes for \emph{Math 431:
Geometry for Computer Graphics} at the University of Maryland (see
\underline{www.math.umd.edu/\symbol{126}gah/}).

We note that, in agreement with what $\alpha\vee\beta$ gives, when the
subspaces represented by $\alpha$ and $\beta$ do not sum to all of $\Phi_{+}$,
the two expressions of ($**$) both give $0$. This is because the totality of
the vectors that make up $\alpha$ and $\beta$ then do not span $\Phi_{+}$ and
therefore neither do any $d$ of them. Let us denote these equal expressions by
$\mu(\alpha,\beta)$ in any case where the sum of the degrees of the blades
$\alpha$ and $\beta$ is at least $d$.

We know that $\mu(\alpha,\beta)=c\cdot\alpha\vee\beta$ but it remains to be
seen how $c$ depends on $\alpha$ and $\beta$. Let us see what happens in a
simple example where $\alpha$ and $\beta$ are made up of $x$s drawn from the
underlying basis $\mathcal{B}$. Let
\[
\alpha=x_{j_{1}}\wedge\cdots\wedge x_{j_{l}}\hspace{0.25in}\text{and}%
\hspace{0.25in}\beta=x_{j_{1}}\wedge\cdots\wedge x_{j_{n}}\wedge x_{j_{l+1}%
}\wedge\cdots\wedge x_{j_{d}}%
\]
where $\left\{  x_{j_{1}},\ldots,x_{j_{d}}\right\}  =\mathcal{B}$. Using the
right-hand expression of ($\ast\ast$), we find that it gives the single
nonzero term $\left(  -1\right)  ^{\sigma}\cdot x_{j_{1}}\wedge\cdots\wedge
x_{j_{n}}$ for $\mu(\alpha,\beta)$, where
\[
\sigma=\left(
\begin{tabular}
[c]{lll}%
$1$ & $\cdots$ & $d$\\
$j_{1}$ & $\cdots$ & $j_{d}$%
\end{tabular}
\ \right)  \text{.}%
\]
We also find that
\[
H\left(  \alpha\right)  =\left(  -1\right)  ^{\sigma}\cdot x_{j_{l+1}}^{\top
}\wedge\cdots\wedge x_{j_{d}}^{\top}\hspace{0.25in}\text{and}\hspace
{0.25in}H\left(  \beta\right)  =\left(  -1\right)  ^{\rho}\cdot x_{j_{n+1}%
}^{\top}\wedge\cdots\wedge x_{j_{l}}^{\top}%
\]
where
\[
\rho=\left(
\begin{tabular}
[c]{lllllllll}%
$1$ & $\cdots$ & $n$ & $n+1$ & $\cdots$ & $n+d-l$ & $n+d-l+1$ & $\cdots$ &
$d$\\
$j_{1}$ & $\cdots$ & $j_{n}$ & $j_{l+1}$ & $\cdots$ & $j_{d}$ & $j_{n+1}$ &
$\cdots$ & $j_{l}$%
\end{tabular}
\ \right)
\]
so that
\begin{align*}
H\left(  \alpha\right)  \wedge H\left(  \beta\right)   &  =\left(  -1\right)
^{\sigma}\left(  -1\right)  ^{\rho}\cdot x_{j_{l+1}}^{\top}\wedge\cdots\wedge
x_{j_{d}}^{\top}\wedge x_{j_{n+1}}^{\top}\wedge\cdots\wedge x_{j_{l}}^{\top}\\
&  =\left(  -1\right)  ^{\sigma}\cdot H\left(  x_{j_{1}}\wedge\cdots\wedge
x_{j_{n}}\right)
\end{align*}
and therefore $\alpha\vee\beta=\left(  -1\right)  ^{\sigma}\cdot x_{j_{1}%
}\wedge\cdots\wedge x_{j_{n}}=\mu(\alpha,\beta)$. This example leads us to
conjecture that it is always the case that $\mu(\alpha,\beta)=\alpha\vee\beta
$, a result that we will soon prove.

\begin{exercise}
Verify the results of the example above in the case where $d=4$, $u_{1}%
=x_{3},u_{2}=x_{2},u_{3}=x_{1}$, and $v_{1}=x_{3},v_{2}=x_{2},v_{3}=x_{4}$.
\emph{(}In evaluating the second expression of \emph{(}$**$\emph{)}, be
careful to base $I$ and $\overline{I}$ on the subscripts of the $v$\emph{s},
not on the subscripts of the $x$\emph{s}, and notice that $\sigma$ comes from
the bracket.\emph{)}
\end{exercise}

While $\vee$ is a bilinear function defined on all of $\bigwedge\Phi_{+}%
\times\bigwedge\Phi_{+}$, $\mu$ is only defined for pairs of blades that have
degrees that sum to at least $d$. Extending $\mu$ as a bilinear function
defined on all of $\bigwedge\Phi_{+}\times\bigwedge\Phi_{+}$ will allow us
prove equality with $\vee$ by considering only what happens to the elements of
$\mathcal{B}^{\wedge}\times\mathcal{B}^{\wedge}$. We thus define the bilinear
function $\widetilde{\mu}:\bigwedge\Phi_{+}\times\bigwedge\Phi_{+}%
\rightarrow\bigwedge\Phi_{+}$ by defining it on $\mathcal{B}^{\wedge}%
\times\mathcal{B}^{\wedge}$ as
\[
\widetilde{\mu}\left(  x_{J},x_{K}\right)  =\left\{
\begin{tabular}
[c]{l}%
$\mu\left(  x_{J},x_{K}\right)  $ if $\left|  J\right|  +\left|  K\right|
\geqslant d$,\\
$0$ otherwise.
\end{tabular}
\ \right.
\]
where $J,K\subset\left\{  1,\ldots,d\right\}  $. We now show that when
$\alpha$ and $\beta$ are blades that have degrees that sum to at least $d$,
$\widetilde{\mu}\left(  \alpha,\beta\right)  =\mu\left(  \alpha,\beta\right)
$, so that $\widetilde{\mu}$ extends $\mu$. Suppose that
\[
\alpha=\sum_{J}a_{J}\cdot x_{J}\hspace{0.25in}\text{and}\hspace{0.25in}%
\beta=\sum_{K}b_{K}\cdot x_{K}\text{.}%
\]
Then, since the bracket is clearly linear, using in turn each of the
expressions of ($\ast\ast$) which define $\mu$, we find that\smallskip\
\begin{align*}
\mu(\alpha,\beta)  &  =\sum_{I}\left(  -1\right)  ^{\overline{I}%
I}\text{\textbf{[}}\left.  u_{\overline{I}}\wedge\sum_{K}b_{K}\cdot
x_{K}\right.  \text{\textbf{]}}\cdot u_{I}=\sum_{K}b_{K}\cdot\mu(\alpha
,x_{K})=\\
\  &  =\sum_{K}b_{K}\cdot\sum_{J}a_{J}\cdot\mu(x_{J},x_{K})=\widetilde{\mu
}\left(  \alpha,\beta\right)
\end{align*}
and $\widetilde{\mu}$ indeed extends $\mu$. We now are ready to formally state
and prove the conjectured result.

\begin{theorem}
For any $\eta,\zeta\in\bigwedge\Phi_{+}$, $\widetilde{\mu}(\eta,\zeta
)=\eta\vee\zeta$. Hence for blades $\alpha=u_{1}\wedge\cdots\wedge u_{l}$ and
$\beta=v_{1}\wedge\cdots\wedge v_{m}$, of respective degrees $l $ and $m$ such
that $l+m\geqslant d$, we have%
\[
\alpha\vee\beta=\sum_{\substack{I\subset\left\{  1,\ldots,l\right\}  \\\left|
I\right|  =n}}\left(  -1\right)  ^{\overline{I}I}\emph{[}\left.
u_{\overline{I}}\wedge\beta\right.  \emph{]}\cdot u_{I}=\sum
_{\substack{I\subset\left\{  1,\ldots,m\right\}  \\\left|  I\right|
=n}}\left(  -1\right)  ^{I\overline{I}}\emph{[}\left.  \alpha\wedge
v_{\overline{I}}\right.  \emph{]}\cdot v_{I}\hspace{0.25in}\left(  \ast
\vee\ast\right)
\]
where $n=l+m-d$.
\end{theorem}

Proof: We need only verify that $\widetilde{\mu}(\eta,\zeta)=\eta\vee\zeta$
for $\left(  \eta,\zeta\right)  \in\mathcal{B}^{\wedge}\times\mathcal{B}%
^{\wedge}$. Thus we assume that $\left(  \eta,\zeta\right)  \in\mathcal{B}%
^{\wedge}\times\mathcal{B}^{\wedge}$. Since $\widetilde{\mu}(\eta,\zeta
)=\eta\vee\zeta=0$ when the subspaces represented by $\eta$ and $\zeta$ do not
sum to $\Phi_{+}$ (including the case when the degrees of $\eta$ and $\zeta$
sum to less than $d$), we assume that the subspaces represented by $\eta$ and
$\zeta$ do sum to $\Phi_{+}$ (and therefore that the degrees of $\eta$ and
$\zeta$ sum to at least $d$). The cases where $\eta=\alpha$ and $\zeta=\beta$
where $\alpha$ and $\beta$ are given as
\[
\alpha=x_{j_{1}}\wedge\cdots\wedge x_{j_{l}}\hspace{0.25in}\text{and}%
\hspace{0.25in}\beta=x_{j_{1}}\wedge\cdots\wedge x_{j_{n}}\wedge x_{j_{l+1}%
}\wedge\cdots\wedge x_{j_{d}}%
\]
have been verified in the example above. The forms of $\alpha$ and $\beta$ are
sufficiently general that every case we need to verify is obtainable by
separately permuting factors inside each of $\alpha$ and $\beta$. This merely
results in $\eta=\left(  -1\right)  ^{\tau}\cdot\alpha$ and $\zeta=\left(
-1\right)  ^{\upsilon}\cdot\beta$ so that $\eta\vee\zeta=\left(  -1\right)
^{\tau}\left(  -1\right)  ^{\upsilon}\cdot\alpha\vee\beta=\widetilde{\mu}%
(\eta,\zeta)$ by bilinearity. $\blacksquare$

\smallskip\ 

The formulas $\left(  \ast\vee\ast\right)  $ of the theorem above express
$\alpha\vee\beta$ in all the nontrivial cases, i. e., when the respective
degrees $l$ and $m$ of the blades $\alpha$ and $\beta$ sum to at least $d,$
and they do it in a \emph{coordinate-free} manner. This is not to say that the
idea of coordinates (or bases) is missing from these formulas or their
derivation. But these formulas do not involve \emph{a priori} coordinates for
$\alpha$ and $\beta$, or their factors, in terms of some overall basis.
Coordinates appear in these formulas as outputs, not as inputs. We have
defined the bracket in terms of the overall basis $\mathcal{B}$ (the same one
used to define $H$\textbf{),} and the bracket is a coordinate function. The
derivation of these formulas proceeded by adapting coordinates to the blade
factors, and then bringing in the bracket to express the adapted coordinates
as functions of the blade factors. The result is a formulation of the
regressive product where only the blade factors appear explicitly. These
formulas are known to be important in symbolic algebraic computations, as
opposed to the numeric computations that we have previously carried out using
\emph{a priori} coordinates based on $\mathcal{B}$. They also now render the
following result readily apparent.

\begin{corollary}
Any two regressive products are proportional.
\end{corollary}

Proof: The bracket depends on the basis $\mathcal{B}$ only up to a nonzero
factor. $\blacksquare$

\newpage\ 

\subsection{Problems}

\ 

1. Theorem \ref{GenEnvDesargues} and its proof require that $P$ be distinct
from the other points. Broaden the theorem by giving a similar proof for the
case where $P$ coincides with $A$.\medskip\ \ 

2. What is the effect on $H$ of a change in the assignment of the subscript
labels $1,\ldots,d$ to the vectors of $\mathcal{B}$?\medskip\ \ 

3. Determine the proportionality factor between two regressive products as the
determinant of a vector space map.\bigskip

4. Notating as in Section \ref{Examples}, $12+34$ is not a blade in the
exterior algebra of a 4-dimensional vector space.\bigskip

5. Referring to Theorem \ref{Hodge}, for arbitrary vectors $v_{1},\ldots
,v_{n}$ of $\Phi_{+}$, $n\leqslant d$,%
\[
H\left(  v_{1}\wedge\cdots\wedge v_{n}\right)  =\left|
\begin{array}
[c]{cccccc}%
x_{1}^{\top}\left(  v_{1}\right)  & \cdots & x_{1}^{\top}\left(  v_{n}\right)
& x_{1}^{\top}\left(  \cdot_{1}\right)  & \cdots & x_{1}^{\top}\left(
\cdot_{d-n}\right) \\
\vdots &  & \vdots & \vdots &  & \vdots\\
x_{n}^{\top}\left(  v_{1}\right)  & \cdots & x_{n}^{\top}\left(  v_{n}\right)
& x_{n}^{\top}\left(  \cdot_{1}\right)  & \cdots & x_{n}^{\top}\left(
\cdot_{d-n}\right) \\
x_{n+1}^{\top}\left(  v_{1}\right)  & \cdots & x_{n+1}^{\top}\left(
v_{n}\right)  & x_{n+1}^{\top}\left(  \cdot_{1}\right)  & \cdots &
x_{n+1}^{\top}\left(  \cdot_{d-n}\right) \\
\vdots &  & \vdots & \vdots &  & \vdots\\
x_{d}^{\top}\left(  v_{1}\right)  & \cdots & x_{d}^{\top}\left(  v_{n}\right)
& x_{d}^{\top}\left(  \cdot_{1}\right)  & \cdots & x_{d}^{\top}\left(
\cdot_{d-n}\right)
\end{array}
\right|
\]
where $x_{i}^{\top}\left(  \cdot_{k}\right)  $ and $x_{j}^{\top}\left(
\cdot_{k}\right)  $ denote the unevaluated functionals $x_{i}^{\top}$ and
$x_{j}^{\top}$, each of which must be evaluated at the same $k$-th vector of
some sequence of vectors. The resulting determinant is then a function of
$x_{1}^{\top},\cdots,x_{d}^{\top}$ each potentially evaluated once at each of
d-n arguments, which we identify with a blade in $\bigwedge\Phi_{+}^{\top}$
via the isomorphism of Theorem \ref{DualExtIso}.

\newpage

\section{Vector Projective Geometry}

\subsection{The Projective Structure on $\mathcal{V}$}

Recall that in Chapter \ref{AffineGeo} we introduced the set $\mathsf{V}%
\left(  \mathcal{V}\right)  $ of all subspaces of the vector space
$\mathcal{V} $, defined the affine structure on $\mathcal{V}$ as
$\mathsf{A}\left(  \mathcal{V}\right)  =\mathcal{V}+\mathsf{V}\left(
\mathcal{V}\right)  $, and interpreted the result as an affine geometry based
on $\mathcal{V}$. We now take just $\mathsf{V}\left(  \mathcal{V}\right)  $ by
itself, call it the \textbf{projective structure} on $\mathcal{V}$, and
interpret it as a projective geometry based on $\mathcal{V}$. The
\textbf{projective flats} are thus the subspaces of $\mathcal{V}$. The
\textbf{points} are the one-dimensional subspaces of $\mathcal{V}$, the
\textbf{lines} are the two-dimensional subspaces of $\mathcal{V}$, the
\textbf{planes} are the three-dimensional subspaces of $\mathcal{V}$, the
$(n-1)$-\textbf{planes} are the $n$-dimensional subspaces of $\mathcal{V}$,
and the \textbf{hyperplanes} are the subspaces of $\mathcal{V}$ of codimension
one. When $\dim\mathcal{V}=2$, the projective geometry based on $\mathcal{V}$
is called a \textbf{projective line}, when $\dim\mathcal{V}=3$, it is called a
\textbf{projective plane}, and in the general finite-dimensional case when
$\dim\mathcal{V}=d$, it is called a \textbf{projective }$(d-1)$-\textbf{space}%
. The $n$-dimensional subspaces of $\mathcal{V}$ are said to have
\textbf{projective dimension} $n-1$. Thus we introduce the function
$\operatorname*{pdim}$ defined by
\[
\operatorname*{pdim}\left(  \mathcal{X}\right)  =\dim\left(  \mathcal{X}%
\right)  -1
\]
to give the projective dimension of the flat or subspace $\mathcal{X}$. The
trivial subspace $0=\left\{  0\right\}  $ is called the \textbf{null flat} and
is characterized by $\operatorname*{pdim}\left(  0\right)  =-1$. Grassmann's
Relation (Corollary \ref{GrassmannRel}) may be reread as%
\[
\operatorname*{pdim}\mathcal{X}+\operatorname*{pdim}\mathcal{Y}%
=\operatorname*{pdim}(\mathcal{X}+\mathcal{Y})+\operatorname*{pdim}%
(\mathcal{X}\cap\mathcal{Y})\text{.}%
\]
We also refer to a subspace sum as the \textbf{join} of the corresponding
projective flats$.$ We denote by $\mathsf{P}\left(  \mathcal{V}\right)  $ the
set of all the points of the geometry and we put $\mathsf{P}^{+}\left(
\mathcal{V}\right)  =\mathsf{P}\left(  \mathcal{V}\right)  \cup0$. All flats
are joins of points, with the null flat being the empty join. If $\mathcal{X}
$ is a subspace of $\mathcal{V}$, then $\mathsf{V}\left(  \mathcal{X}\right)
\subset\mathsf{V}\left(  \mathcal{V}\right)  $ gives us a projective geometry
in its own right, a \textbf{subgeometry} of $\mathsf{V}\left(  \mathcal{V}%
\right)  $, with $\mathsf{P}\left(  \mathcal{X}\right)  \subset\mathsf{P}%
\left(  \mathcal{V}\right)  $.

$\mathcal{V}$ is conceptually the same as our $\Phi_{+}$ of the previous
chapter, except that it is stripped of all references to the affine flat
$\Phi$. No points of $\mathsf{P}\left(  \mathcal{V}\right)  $ are special. The
directions that were viewed as special ``points at infinity'' in the viewpoint
of the previous chapter become just ordinary points without any designated
$\Phi_{0}$ for them to inhabit. We will at times find it useful to designate
some hyperplane through $0$ in $\mathcal{V}$ to serve as a $\Phi_{0} $ and
thereby allow $\mathcal{V}$ to be interpreted in a generalized affine manner.
However, independent of any generalized affine interpretation, we have the
useful concept of homogeneous representation of points by vectors which we
continue to exploit.

\subsection{Projective Frames}

\emph{Notice: }Unless otherwise stated, we assume for the rest of this chapter
that $\mathcal{V}$ is a nontrivial vector space of finite dimension $d$ over
the field $\mathcal{F}$.

\smallskip\ 

The points $X_{0},\ldots,X_{n}$ of $\mathsf{P}\left(  \mathcal{V}\right)  $
are said to be in \textbf{general position} if they are distinct and have
vector representatives that are as independent as possible, i. e., either all
of them are independent, or, if $n\geqslant d$, every $d$ of them are
independent. The case $n=d$ is of particular interest: any $d+1$ points in
general position are said to form a \textbf{projective frame} for
$\mathsf{V}\left(  \mathcal{V}\right)  $.

\begin{exercise}
Any $3$ distinct points of a projective line form a projective frame for it.
\end{exercise}

To say that the $d+1$ distinct points $X_{0},\ldots,X_{d}$ are a projective
frame for $\mathsf{V}\left(  \mathcal{V}\right)  $ is the same as to say that
there are respective representing vectors $x_{0},\ldots,x_{d}$ such that
$\left\{  x_{1},\ldots,x_{d}\right\}  $ is an independent set (and therefore a
basis for $\mathcal{V}$) and $x_{0}=a_{1}\cdot x_{1}+\cdots+a_{d}\cdot x_{d}$,
where all of the scalar coefficients $a_{i}$ are nonzero. Given that $\left\{
x_{1},\ldots,x_{d}\right\}  $ is a basis, then with the same nonzero $a_{i}$,
$\left\{  a_{1}\cdot x_{1},\ldots,a_{d}\cdot x_{d}\right\}  $ is also a basis.
Hence the following result.

\begin{proposition}
Points $X_{0},\ldots,X_{d}$ of a projective frame for $\mathsf{V}\left(
\mathcal{V}\right)  $ can be represented respectively by vectors $x_{0}%
,\ldots,x_{d}$ such that $\left\{  x_{1},\ldots,x_{d}\right\}  $ is a basis
for $\mathcal{V}$ and $x_{0}=x_{1}+\cdots+x_{d}$. $\blacksquare$
\end{proposition}

A representation of the points $X_{0},\ldots,X_{d}$ of a projective frame by
respective vectors $x_{0},\ldots,x_{d}$ such that $x_{0}=x_{1}+\cdots+x_{d}$
is said to be \textbf{standardized}. In such a standardized representation,
$X_{0}$ is called the \textbf{unit point}, and $X_{1},\ldots,X_{d}$ are known
by various terms such as \textbf{base points}, \textbf{vertices of the simplex
of reference}, or \textbf{fundamental points of the coordinate system}. In
terms of coordinate $d$-tuples, any $d+1$ points $X_{0},\ldots,X_{d}$ in
general position can thus be assigned the respective coordinate $d$-tuples
$\varepsilon_{0}=\left(  1,1,\ldots,1\right)  ,$ $\varepsilon_{1}=\left(
1,0,\ldots,0\right)  ,$ $\varepsilon_{2}=\left(  0,1,0,\ldots,0\right)
,\ldots,\varepsilon_{d}=\left(  0,\ldots,0,1\right)  $. However, up to the
unavoidable (and ignorable) nonzero overall scalar multiplier, there is only
one basis for $\mathcal{V}$ that permits this assignment.

\begin{proposition}
The basis $\left\{  x_{1},\ldots,x_{d}\right\}  $ of the previous proposition
is unique up to a nonzero overall scalar multiplier.
\end{proposition}

Proof: Suppose that $X_{1},\ldots,X_{d}$ are also represented respectively by
the vectors $y_{1},\ldots,y_{d}$ which together automatically form a basis for
$\mathcal{V}$ due to the general position hypothesis, and $X_{0}$ is
represented by $y_{0}=y_{1}+\cdots+y_{d}$. Then $y_{i}=a_{i}\cdot x_{i}$ for
each $i$ with $a_{i}$ a nonzero scalar. Hence
\[
y_{0}=y_{1}+\cdots+y_{d}=a_{1}\cdot x_{1}+\cdots+a_{d}\cdot x_{d}=a_{0}%
\cdot\left(  x_{1}+\cdots+x_{d}\right)
\]
and since $\left\{  x_{1},\ldots,x_{d}\right\}  $ is a basis for $\mathcal{V}$
we must have $a_{i}=a_{0}$ for all $i>0$. Hence $y_{i}=a_{0}\cdot x_{i}$ for
all $i>0$ which is what was to be shown. $\blacksquare$

\smallskip\ 

Designating a standardized representation of the points of a projective frame
thus uniquely specifies a system of homogeneous coordinates on $\mathsf{P}%
\left(  \mathcal{V}\right)  $. Given that homogeneous coordinate $d$-tuples
for the unit point and the base points are $\varepsilon_{0}$ and
$\varepsilon_{1},\ldots,\varepsilon_{d}$, coordinate $d$-tuples of all the
points of $\mathsf{P}\left(  \mathcal{V}\right)  $ are completely determined
in the homogeneous sense (i. e., up to a nonzero scalar multiplier).

\subsection{Pappus Revisited}

Employing a projective frame can ease the proof of certain theorems of
projective geometry. The Theorem of Pappus is one such. We previously stated
and proved one of its affine versions as Theorem \ref{Pappus}. That proof
employed rather elaborate barycentric calculations. The projective version
that we now state and prove has essentially the same statement in words as
Theorem \ref{Pappus} but now refers to points and lines in a projective plane
over a suitable field.

\begin{theorem}
[Theorem of Pappus]Let $l,l\,^{\prime}$ be distinct coplanar lines with
distinct points $A,B,C\subset l\smallsetminus l\,^{\prime}$ and distinct
points $A\,^{\prime},B\,^{\prime},C\,^{\prime}\subset l\,^{\prime
}\smallsetminus l$. If $BC\,^{\prime}$ meets $B\,^{\prime}C$ in $A\,^{\prime
\prime}$, $AC\,^{\prime}$ meets $A\,^{\prime}C$ in $B\,^{\prime\prime}$, and
$AB\,^{\prime}$ meets $A\,^{\prime}B$ in $C\,^{\prime\prime}$, then
$A\,^{\prime\prime},B\,^{\prime\prime},C\,^{\prime\prime}$ are distinct
collinear points.
\end{theorem}%

\begin{center}
\includegraphics[
trim=0.000000in 0.000000in -0.047249in -0.027293in,
height=2.9992in,
width=3.4999in
]%
{FINALPAP.eps}%
\\
Illustrating the Theorem of Pappus
\end{center}

Proof: We choose the points $B^{\prime},A,B,C^{\prime}$ as our projective
frame. These points are in general position since any three of them always
contain one that is not allowed to be on the line through the other two. We
designate a standardized representation by choosing vectors $x_{1},x_{2}%
,x_{3}$ to respectively represent $A,B,C^{\prime}$ and $x_{0}=x_{1}%
+x_{2}+x_{3}$ to represent the unit point $B^{\prime}$.\medskip\ \ 

$AB=\left\langle \left\{  x_{1},x_{2}\right\}  \right\rangle $ and $C$ on it
can be represented as $x_{1}+c\cdot x_{2}$, or what is the same, as the
coordinate tuple $\left(  1,c,0\right)  $.\medskip\ 

$B^{\prime}C^{\prime}=\left\langle \left\{  x_{1}+x_{2}+x_{3},x_{3}\right\}
\right\rangle $ and $A^{\prime}$ on it can be represented by $\left(
1,1,a\right)  $.\medskip\ \ 

$a\neq0$ because letting $a=0$ in the representation $\left(  1,1,a\right)  $
for $A^{\prime}$ would put $A^{\prime}$ on $AC.$ This is so because
$AC=\left\langle \left\{  x_{1},x_{1}+c\cdot x_{2}\right\}  \right\rangle $
contains $c\cdot\left(  x_{1}+x_{2}\right)  =\left(  c-1\right)  \cdot
x_{1}+\left(  x_{1}+c\cdot x_{2}\right)  $ and $c\neq0$ since $C\neq A$.
Similarly, $c\neq1$ because letting $c=1$ in the representation $\left(
1,c,0\right)  $ of $C$ would put $C$ on $A^{\prime}B^{\prime}$. This is so
because $A^{\prime}B^{\prime}=\left\langle \left\{  x_{1}+x_{2}+a\cdot
x_{3},x_{1}+x_{2}+x_{3}\right\}  \right\rangle $ which contains $\left(
1-a\right)  \cdot\left(  x_{1}+x_{2}\right)  =\left(  x_{1}+x_{2}+a\cdot
x_{3}\right)  -a\cdot\left(  x_{1}+x_{2}+x_{3}\right)  $ and $a\neq1$ since
$A^{\prime}\neq B^{\prime}$.\medskip

$BC^{\prime}=\left\langle \left\{  x_{2},x_{3}\right\}  \right\rangle $,
$B^{\prime}C=\left\langle \left\{  x_{1}+x_{2}+x_{3},x_{1}+c\cdot
x_{2}\right\}  \right\rangle $ and $\left(  1-c\right)  \cdot x_{2}%
+x_{3}=\left(  x_{1}+x_{2}+x_{3}\right)  -\left(  x_{1}+c\cdot x_{2}\right)  $
lies in both so that $A^{\prime\prime}$ is represented by $\left(
0,1-c,1\right)  $.\medskip\ \ 

$C^{\prime}A=\left\langle \left\{  x_{1},x_{3}\right\}  \right\rangle $,
$A^{\prime}C=\left\langle \left\{  x_{1}+x_{2}+a\cdot x_{3},x_{1}+c\cdot
x_{2}\right\}  \right\rangle $ and $\left(  1-c\right)  \cdot x_{1}-ca\cdot
x_{3}=x_{1}+c\cdot x_{2}-c\cdot\left(  x_{1}+x_{2}+a\cdot x_{3}\right)  $ lies
in both so that $B^{\prime\prime}$ is represented by $\left(
1-c,0,-ca\right)  $.\medskip\ \ 

$AB^{\prime}=\left\langle \left\{  x_{1},x_{1}+x_{2}+x_{3}\right\}
\right\rangle $, $A^{\prime}B=\left\langle \left\{  x_{1}+x_{2}+a\cdot
x_{3},x_{2}\right\}  \right\rangle $ and $x_{1}+a\cdot x_{2}+a\cdot
x_{3}=\left(  1-a\right)  \cdot x_{1}+a\cdot\left(  x_{1}+x_{2}+x_{3}\right)
=\left(  x_{1}+x_{2}+a\cdot x_{3}\right)  +\left(  a-1\right)  \cdot x_{2}$
lies in both so that $C^{\prime\prime}$ is represented by $\left(
1,a,a\right)  $.\medskip\ \ 

The dependency relation
\[
a\cdot\left(  0,1-c,1\right)  +\left(  1-c,0,-ca\right)  +\left(  c-1\right)
\cdot\left(  1,a,a\right)  =0
\]
shows that $A^{\prime\prime},B^{\prime\prime},C^{\prime\prime}$ are collinear.
$A^{\prime\prime},B^{\prime\prime},C^{\prime\prime}$ are distinct since
$a\neq0$ and $c\neq1$ (shown earlier in this proof) make it impossible for any
one of $\left(  0,1-c,1\right)  $, $\left(  1-c,0,-ca\right)  $, $\left(
1,a,a\right)  $ to be a multiple of any other. $\blacksquare$

What we have just done in proving the Theorem of Pappus becomes even more
transparent when viewed in the light of a generalized affine interpretation.
Using the same basis as for our chosen standardized representation, and
letting $x_{1}=1$ be our $\Phi_{1}$ just as we did in the previous chapter, we
get the generalized affine interpretation of a rather simple affine theorem.
The flat $\Phi$ that we view as having been inflated to $\Phi_{1}$ has point
$A$ at its origin, and the points $B$ and $C^{\prime}$ have been ``sent to
infinity'' in two different directions. As a result, $A^{\prime\prime}$ is at
infinity, and the conclusion of the theorem is that $B^{\prime\prime}%
C^{\prime\prime}$ is parallel to $B^{\prime}C$, as depicted in the figure
below. The proof above is just a proof of this affine theorem using the
inflated flat method. We prove exactly the same thing when we give an ordinary
proof of the affine theorem, as we now do. In coordinate form, the finite
points are then $A=\left(  0,0\right)  $, $C=\left(  c,0\right)  $,
$B^{\prime}=\left(  1,1\right)  $, $A^{\prime}=\left(  1,a\right)  $,
$C^{\prime\prime}=\left(  a,a\right)  $, $B^{\prime\prime}=\left(  0,ac\left(
c-1\right)  ^{-1}\right)  $. Thus
\[
C^{\prime\prime}-B^{\prime\prime}=\left(  a,a\left(  1-c\right)  ^{-1}\right)
=a\left(  1-c\right)  ^{-1}\cdot\left(  1-c,1\right)  =a\left(  1-c\right)
^{-1}\cdot\left(  B^{\prime}-C\right)
\]
so that $B^{\prime\prime}C^{\prime\prime}$ is parallel to $B^{\prime}C$ as was
to be shown.%

\begin{center}
\includegraphics[
trim=0.000000in 0.097803in -0.005807in 0.004180in,
height=4.1796in,
width=4.8395in
]%
{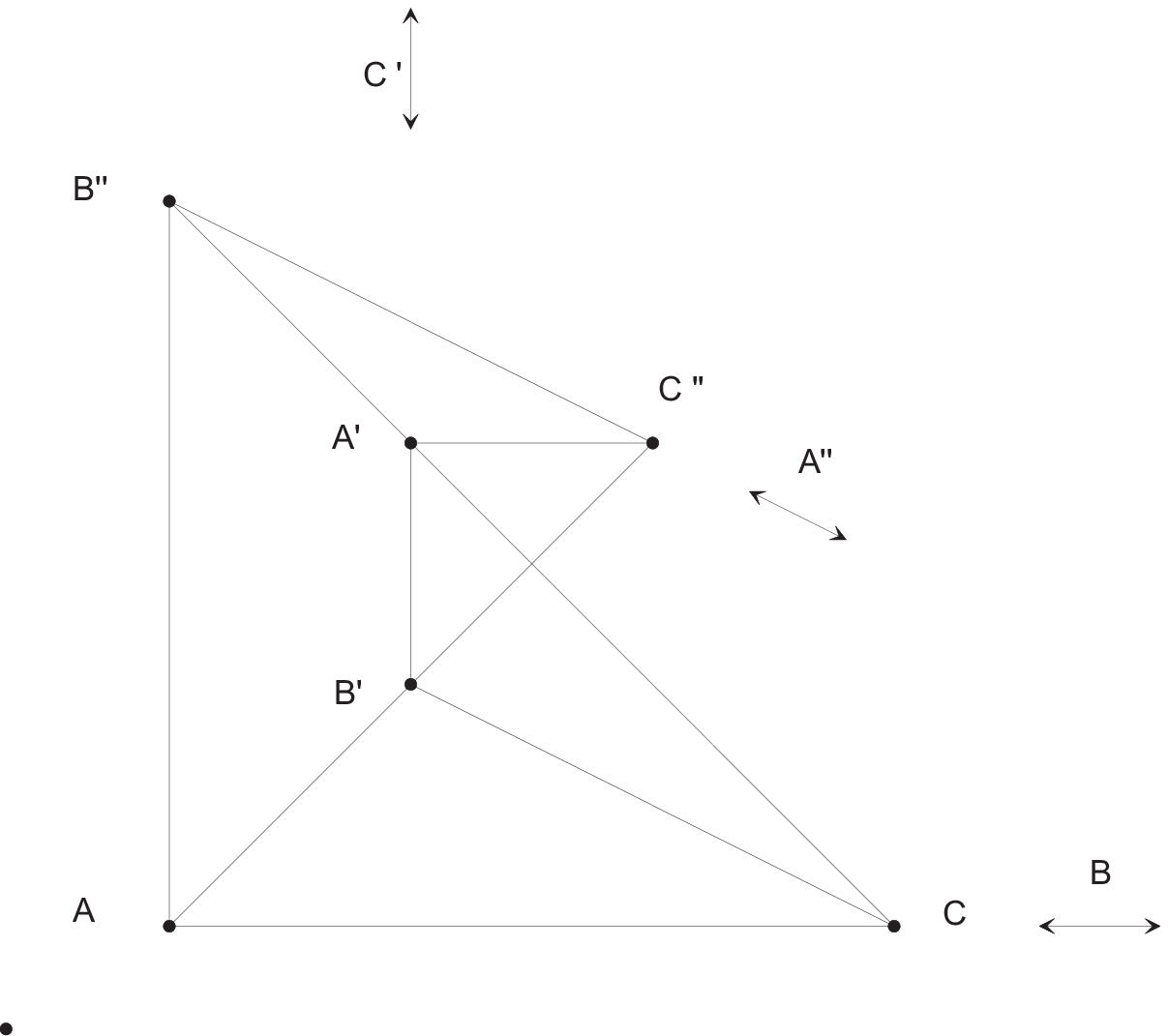}%
\\
Generalized Affine Interpretation
\end{center}

\subsection{Projective Transformations}

Besides bases and frames, certain functions play an important r\^ole in
projective geometry. Let $\mathcal{W}$ be an alias of $\mathcal{V}$. Any
vector space isomorphism $f:\mathcal{V}\rightarrow\mathcal{W}$ induces a
bijection from $\mathsf{P}\left(  \mathcal{V}\right)  $ to $\mathsf{P}\left(
\mathcal{W}\right)  $. By a \textbf{projective transformation} we will mean
any such bijection from $\mathsf{P}\left(  \mathcal{V}\right)  $ to
$\mathsf{P}\left(  \mathcal{W}\right)  $ induced by a vector space isomorphism
from $\mathcal{V}$ to $\mathcal{W}$. The terms \textbf{projectivity} and
\textbf{homography }are also commonly used. For any nonzero scalar $a$, $f$
and $a\cdot f$ clearly induce the same projective transformation. On the other
hand, if $g:\mathcal{V}\rightarrow\mathcal{W}$ is another vector space
isomorphism that induces the same projective transformation as $f,$ then
$g=a\cdot f$ for some nonzero scalar $a$, as we now show. Let $\left\{
x_{1},\ldots,x_{d}\right\}  $ be a basis for $\mathcal{V}$. For each $i$ we
must have $g\left(  x_{i}\right)  =a_{i}\cdot f\left(  x_{i}\right)  $ for
some nonzero scalar $a_{i}$. Also $g\left(  x_{1}+\cdots+x_{d}\right)  =a\cdot
f\left(  x_{1}+\cdots+x_{d}\right)  $ for some nonzero scalar $a$. But then
\[
a_{1}\cdot f\left(  x_{1}\right)  +\cdots+a_{d}\cdot f\left(  x_{d}\right)
=a\cdot f\left(  x_{1}\right)  +\cdots+a\cdot f\left(  x_{d}\right)
\]
and since $\left\{  f\left(  x_{1}\right)  ,\ldots,f\left(  x_{d}\right)
\right\}  $ is an independent set, $a_{i}=a$ for all $i$. Hence $g\left(
x_{i}\right)  =a\cdot f\left(  x_{i}\right)  $ for all $i$ and therefore
$g=a\cdot f$. Thus, similar to the classes of proportional vectors as the
representative classes for points, we have the classes of proportional vector
space isomorphisms as the representative classes for projective
transformations. A vector space isomorphism homogeneously represents a
projective transformation in the same fashion that a vector represents a
point. For the record, we now formally state this result as the following theorem.

\begin{theorem}
\label{ProjTransf}Two isomorphisms between finite-dimensional vector spaces
induce the same projective transformation if and only if these isomorphisms
are proportional. $\blacksquare$
\end{theorem}

It is clear that projective transformations send projective frames to
projective frames. Given an arbitrary projective frame for $\mathsf{V}\left(
\mathcal{V}\right)  $ and another for $\mathsf{V}\left(  \mathcal{W}\right)
$, there is an obvious projective transformation that sends the one to the
other. This projective transformation, in fact, is uniquely determined.

\begin{theorem}
Let $X_{0},\ldots,X_{d}$ and $Y_{0},\ldots,Y_{d}$ be projective frames for
$\mathsf{V}\left(  \mathcal{V}\right)  $ and $\mathsf{V}\left(  \mathcal{W}%
\right)  $, respectively. Then there is a unique projective transformation
from $\mathsf{P}\left(  \mathcal{V}\right)  $ to $\mathsf{P}\left(
\mathcal{W}\right)  $ which for each $i$ sends $X_{i}$ to $Y_{i}$.
\end{theorem}

Proof: Let $x_{1},\ldots,x_{d}$ be representative vectors corresponding to
$X_{1},\ldots,X_{d}$ and such that $x_{0}=x_{1}+\cdots+x_{d}$ represents
$X_{0}$. Similarly, let $y_{1},\ldots,y_{d}$ be representative vectors
corresponding to $Y_{1},\ldots,Y_{d}$ and such that $y_{0}=y_{1}+\cdots+y_{d}$
represents $Y_{0}$. Then the vector space map $f:\mathcal{V}\rightarrow
\mathcal{W}$ that sends $x_{i}$ to $y_{i}$ for $i>0$ induces a projective
transformation that for each $i$ sends $X_{i}$ to $Y_{i}$. Suppose that the
vector space map $g:\mathcal{V}\rightarrow\mathcal{W}$ also induces a
projective transformation that for each $i$ sends $X_{i}$ to $Y_{i}$. Then
there are nonzero scalars $a_{0},\ldots,a_{d}$ such that, for each $i$,
$g\left(  x_{i}\right)  =a_{i}\cdot y_{i}$ so that then
\[
g\left(  x_{0}\right)  =g\left(  x_{1}+\cdots+x_{d}\right)  =a_{1}\cdot
y_{1}+\cdots+a_{d}\cdot y_{d}=a_{0}\cdot\left(  y_{1}+\cdots+y_{d}\right)
\]
and since $\left\{  y_{1},\ldots,y_{d}\right\}  $ is an independent set, it
must be that all the $a_{i}$ are equal to $a_{0}$. Hence $g=a_{0}\cdot f$, so
that $f$ and $g$ induce exactly the same projective transformation.
$\blacksquare$

\subsection{Projective Maps}

Because $\mathsf{P}\left(  \mathcal{W}\right)  $ does not contain the null
flat, only a one-to-one vector space map from $\mathcal{V}$ to some other
vector space $\mathcal{W}$ will induce a function from $\mathsf{P}\left(
\mathcal{V}\right)  $ to $\mathsf{P}\left(  \mathcal{W}\right)  $. However,
any vector space map from $\mathcal{V}$ to $\mathcal{W}$ does induce a
function from $\mathsf{P}^{+}\left(  \mathcal{V}\right)  $ to $\mathsf{P}%
^{+}\left(  \mathcal{W}\right)  $. By a \textbf{projective map} we will mean
any function from $\mathsf{P}^{+}\left(  \mathcal{V}\right)  $ to
$\mathsf{P}^{+}\left(  \mathcal{W}\right)  $ induced by a vector space map
from $\mathcal{V}$ to $\mathcal{W}$. Technically, a projective transformation
is not a projective map because it is from $\mathsf{P}\left(  \mathcal{V}%
\right)  $ to $\mathsf{P}\left(  \mathcal{W}\right)  $, not from
$\mathsf{P}^{+}\left(  \mathcal{V}\right)  $ to $\mathsf{P}^{+}\left(
\mathcal{W}\right)  $. However, each projective transformation clearly extends
to a unique bijective projective map from $\mathsf{P}^{+}\left(
\mathcal{V}\right)  $ to $\mathsf{P}^{+}\left(  \mathcal{V}\right)  $.

\begin{exercise}
The composite of vector space maps induces the composite of the separately
induced projective maps, and the identity induces the identity.
\end{exercise}

A vector space map from $\mathcal{V}$ to $\mathcal{W}$ also induces a function
from $\mathsf{V}\left(  \mathcal{V}\right)  $ to $\mathsf{V}\left(
\mathcal{W}\right)  $ and we could just as well have used this as our
projective map. A function is also induced from $\mathsf{P}\left(
\mathcal{V}\right)  \smallsetminus\mathsf{P}\left(  \mathcal{K}\right)  $ to
$\mathsf{P}\left(  \mathcal{W}\right)  $, where $\mathcal{K}$ is the kernel of
the vector space map, and some authors call this a projective map. These
various concepts of projective map all amount to essentially the same thing
because each has the same origin in terms of classes of vector space maps. For
general vector space maps, just as for the isomorphisms, these classes are the
proportional classes as the following exercise records.

\begin{exercise}
Let $f,g:\mathcal{V}\rightarrow\mathcal{W}$ be vector space maps from the
finite-dimensional vector space $\mathcal{V}$ and suppose that $f$ and $g$
induce the same projective map. Show that $f$ and $g$ are proportional by
considering how they act on the basis vectors $x_{1},\ldots,x_{k}$ of a
complementary subspace of their common kernel, and on $x_{1}+\cdots+x_{k}$.
\end{exercise}

\subsection{Central Projection}

In $\mathsf{V}\left(  \mathcal{V}\right)  $, fix a hyperplane $t$ (the
\textbf{target}), and fix a point $C$ (the \textbf{center}) such that $C\notin
t$. Given the point $P\in\mathsf{P}\left(  \mathcal{V}\right)  \smallsetminus
C$, we can form the line $CP$ through $C$ and $P$. $CP$ will intersect $t$ in
a unique point $Q.$ The assignment of $Q=CP\cap t$ to $P\in\mathsf{P}\left(
\mathcal{V}\right)  \smallsetminus C$ is called \textbf{central projection
}from the point $C$ onto the hyperplane $t$. By also assigning the null flat
of $\mathsf{V}\left(  t\right)  $ both to $C$ and to the null flat of
$\mathsf{V}\left(  \mathcal{V}\right)  $, we get a function from
$\mathsf{P}^{+}\left(  \mathcal{V}\right)  $ to $\mathsf{P}^{+}\left(
t\right)  $ which turns out to be a projective map. Also, by restricting the
domain to the points of any hyperplane not containing $C,$ and restricting the
codomain to the points of$\ t$, we will see that a projective transformation
between hyperplanes results, one that derives from the vector space self-map
of projection onto $t$ along $C$ (Section \ref{Projections}).%

\begin{center}
\includegraphics[
trim=0.000000in 0.099461in 0.000000in 0.178615in,
height=3.915in,
width=5.3549in
]%
{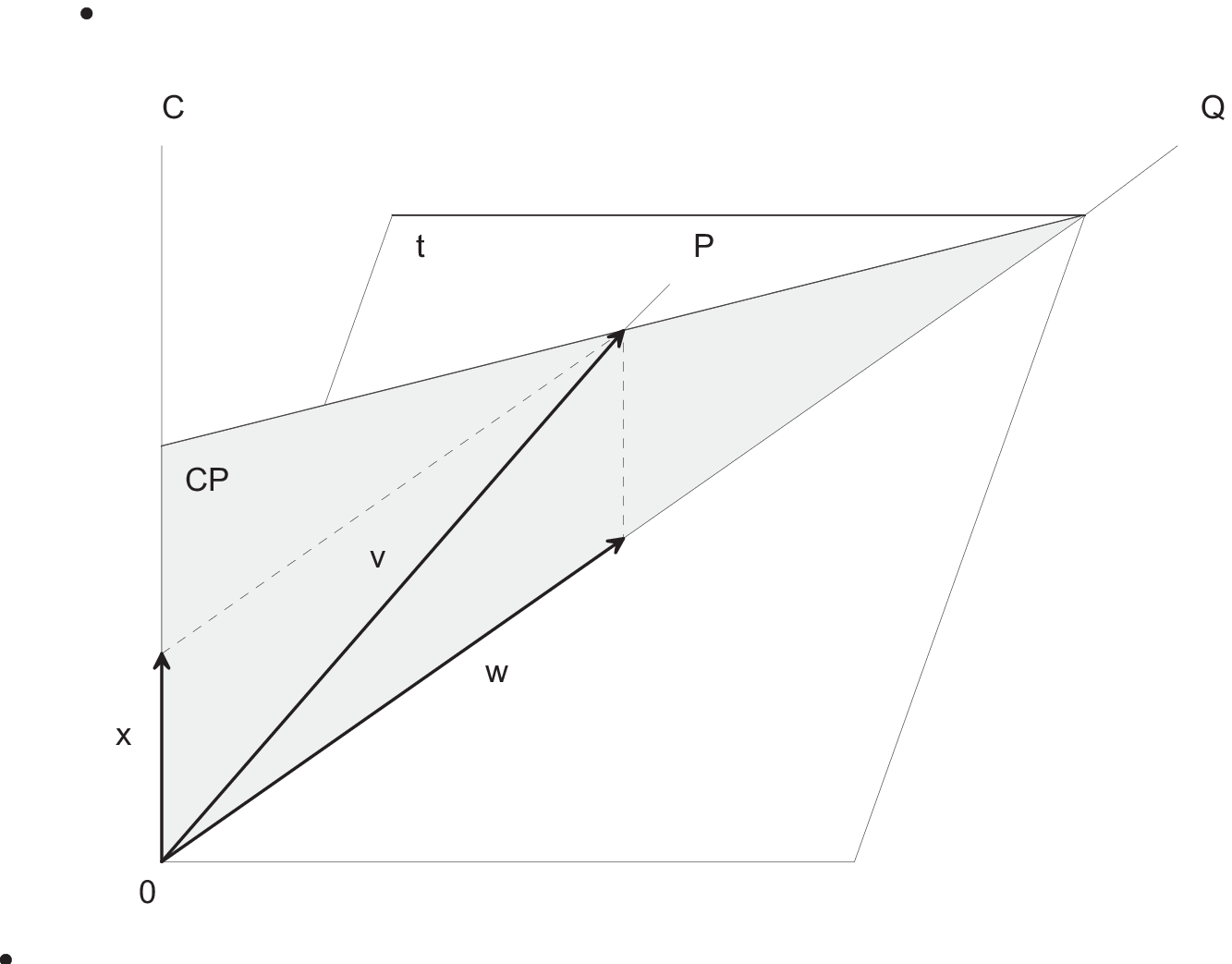}%
\\
Projecting P from C onto t
\end{center}

The figure above, not the usual schematic but rather a full depiction of the
actual subspaces in $\mathcal{V}$, illustrates central projection in the
simple case where $\mathcal{V}$ is of dimension $3$ and the target $t$ is a
projective line (a $2$-dimensional subspace of $\mathcal{V}$). Let $P$ be
represented by the vector $v$. Then we can (uniquely) express $v$ as $v=x+w$
where $x\in C$ and $w\in t$ since $\mathcal{V}=C\oplus t$. But $w=v-x\in CP$,
so $w\in CP\cap t$ and therefore $w$ represents the point $Q$ that we seek.
Thus a vector representing $Q$ may always be obtained by applying the
projection onto $t$ along $C$ to any vector that represents $P$, thereby
making the assignment of $Q$ to $P$ a projective map. $C$ is the kernel of the
projection onto $t$ along $C$, so the points of any projective line $s$ that
does not contain $C$ will be placed in one-to-one correspondence with the
points of $t$. Thus the restriction of the projection onto $t$ along $C$ to a
projective line $s$ that does not contain $C$ is a vector space isomorphism
from $s$ onto $t$, and the induced map from the points of $s$ onto the points
of $t$ is a projective transformation, an example of a \textbf{perspectivity
}between hyperplanes.

Although we have only illustrated here the case where $\mathcal{V}$ has
dimension $3$, the same considerations clearly apply in any $\mathcal{V}$ of
finite dimension $d$. Also, as detailed in the following exercise, we may
readily generalize to where center and target are any complementary subspaces
in $\mathcal{V}$ and obtain a concept of central projection from a center of
any positive dimension $n$ to a target of complementary dimension $d-n$.

\begin{exercise}
Let $\mathcal{C}$ be a subspace of dimension $n>0$ in the vector space
$\mathcal{V}$ of dimension $d$, and let $\mathcal{T}$ be a complement of
$\mathcal{C}$. Let $P$ be a $1$-dimensional subspace of $\mathcal{V}$ such
that $P\cap\mathcal{C}=0$. Then the join $\mathcal{C}P$ of $\mathcal{C}$ and
$P$ intersects $\mathcal{T}$ in a $1$-dimensional subspace $Q.$ If $v$ is a
nonzero vector in $P,$ then the projection of $v$ onto $\mathcal{T}$ along
$\mathcal{C}$ is a nonzero vector in $Q$.
\end{exercise}

\subsection{Problems}

\ 

1. Let $f$ and $g$ be isomorphisms between possibly infinite-dimensional
vector spaces $\mathcal{V}$ and $\mathcal{W}$, and let $f$ and $g$ be such
that for all vectors $v$%
\[
g\left(  v\right)  =a\left(  v\right)  \cdot f\left(  v\right)  \text{ for a
nonzero scalar }a\left(  v\right)  \text{.}%
\]
Then for all vectors $u$ and $v$, $a\left(  u\right)  =a\left(  v\right)  $.
Thus Theorem \ref{ProjTransf} holds in the infinite-dimensional case as well.

(Consider separately the cases where $u$ and $v$ represent the same point in
$\mathsf{P}\left(  \mathcal{V}\right)  $ and when they do not. Write $g\left(
u\right)  $ two different ways in the former case and write $g\left(
u+v\right)  $ in two different ways in the latter case.)

\medskip\ 

2. Suppose that $\mathcal{V}$ is over the finite field of $q$ elements. Then
if $d\geqslant q$, $d+1$ distinct points of $\mathsf{P}\left(  \mathcal{V}%
\right)  $, \emph{but no more}, can be in general position. But if $2\leqslant
d<q$, then at least $d+2$ points of $\mathsf{P}\left(  \mathcal{V}\right)  $
can be in general position.

\newpage

\section{Scalar Product Spaces}

\subsection{Pairings and Their Included Maps}

Let $\mathcal{V}$ and $\mathcal{W}$ be vector spaces over the same field
$\mathcal{F}$. We will refer to a bilinear functional $g:\mathcal{V}%
\times\mathcal{W}\rightarrow\mathcal{F}$ as a \textbf{pairing} of
$\mathcal{V}$ and $\mathcal{W}$. Consider the pairing $g:\mathcal{V}%
\times\mathcal{W}\rightarrow\mathcal{F}$. With the second argument of $g$ held
fixed at $w$, allowing the first argument to vary produces the values of a
linear functional $g_{\_w}\in\mathcal{V}^{\top}$. Similarly, the values of a
linear functional $g_{v\_}\in$ $\mathcal{W}^{\top}$ are produced when the
first argument is held fixed at $v$ and the second argument is allowed to
vary. Now, if we take $g_{v\_}$ and let $v$ vary, we get a map $g_{1}%
:\mathcal{V}\rightarrow\mathcal{W}^{\top}$, and if we take $g_{\_w}$ and let
$w$ vary, we get a map $g_{2}:\mathcal{W}\rightarrow\mathcal{V}^{\top}$. Thus
$g$ includes within itself the two maps $g_{1}$ and $g_{2}$ which we will
refer to as the \textbf{included maps }belonging to $g$. The various $g$s are
related by
\[
g_{v\_}(w)=(g_{1}(v))(w)=g(v,w)=(g_{2}(w))(v)=g_{\_w}(v).
\]
For any map $f:\mathcal{V}\rightarrow\mathcal{W}^{\top}$, $(f(v))(w)$ is the
value of a bilinear functional of $v$ and $w$, and hence any such $f$ is the
$g_{1}$ included map of some pairing of $\mathcal{V}$ and $\mathcal{W}$. By
the same token, any map from $\mathcal{W}$ to $\mathcal{V}^{\top}$ can be a
$g_{2}$. However, each included map clearly determines the other, so only one
of them may be specified arbitrarily.

\subsection{Nondegeneracy}

A pairing $g\ $is called \textbf{nondegenerate} if for each nonzero $v$, there
is some $w$ for which $g(v,w)\neq0$, and for each nonzero $w$, there is some
$v$ for which $g(v,w)\neq0$. A nondegenerate pairing $g:\mathcal{V}%
\times\mathcal{W}\rightarrow\mathcal{F}$ is sometimes referred to as a
\textbf{perfect} (or \textbf{duality})\textbf{\ pairing}, and is said to put
$\mathcal{V}$ \textbf{in duality} with $\mathcal{W}$.

\begin{exercise}
\label{Evaluation}For any vector space $\mathcal{V}$ over $\mathcal{F}$, the
natural \textbf{evaluation pairing} $e:\mathcal{V}^{\top}\times\mathcal{V}%
\rightarrow\mathcal{F}$, defined by $e(f,v)$ $=f(v)$ for $f\in\mathcal{V}%
^{\top}$ and $v\in\mathcal{V}$, puts $\mathcal{V}^{\top}$ in duality with
$\mathcal{V}$. \emph{(}Each $v\neq0$ is part of some basis and therefore has a
coordinate function $v^{\top}$.\emph{)}
\end{exercise}

The nice thing about nondegeneracy is that it makes the included maps
one-to-one, so they map different vectors to different functionals. For, if
the distinct vectors $s$ and $t$ of $\mathcal{V}$ both map to the same
functional (so that $g_{t\_}=g_{s\_}$), then $g_{u\_}$, where $u=t-s\neq0$, is
the zero functional on $\mathcal{W}$. Hence for the nonzero vector $u$,
$g_{u\_}(w)=g(u,w)=0$ for all $w\in\mathcal{W}$, and $g$ therefore fails to be nondegenerate.

\begin{exercise}
If the included maps are both one-to-one, the pairing is nondegenerate.
\emph{(}The one-to-one included maps send only the zero vector to the zero functional.\emph{)}
\end{exercise}

If $g$ is nondegenerate and one of the spaces is finite-dimensional, then so
is the other. ($\mathcal{W}$ is isomorphic to a subspace of $\mathcal{V}%
^{\top}$, so if $\mathcal{V}$ is finite-dimensional, so is $\mathcal{W}$.)
Supposing then that a nondegenerate $g$ is a pairing of finite-dimensional
spaces, we have $\dim\mathcal{W}\leqslant\dim\mathcal{V}^{\top}=\dim
\mathcal{V}$ and also $\dim\mathcal{V}\leqslant\dim\mathcal{W}^{\top}%
=\dim\mathcal{W}$. Hence the spaces all have the same finite dimension and the
included maps are therefore isomorphisms, since we know (Theorem \ref{Injmap})
that a one-to-one map into an alias must also be onto.

The included maps belonging to a nondegenerate pairing $g$ of a
finite-dimensional $\mathcal{V}$ with itself (which puts $\mathcal{V}$ in
duality with itself) are isomorphisms between $\mathcal{V}$ and $\mathcal{V}%
^{\top}$, each of which allows any functional in $\mathcal{V}^{\top}$ to be
represented by a vector of $\mathcal{V}$ in a basis-independent manner. Any
$\varphi\in\mathcal{V}^{\top}$ is represented by $v_{\varphi}=g_{1}%
^{-1}(\varphi)$, and also by $w_{\varphi}=g_{2}^{-1}(\varphi)$, so that for
any $w\in\mathcal{V}$%
\[
\varphi(w)=(g_{1}(g_{1}^{-1}(\varphi))(w)=g(g_{1}^{-1}(\varphi
),w)=g(v_{\varphi},w)
\]
and for any $v\in\mathcal{V}$%
\[
g(v,w_{\varphi})=g(v,g_{2}^{-1}(\varphi))=(g_{2}(g_{2}^{-1}(\varphi
))(v)=\varphi(v).
\]

\subsection{Orthogonality, Scalar Products}

Designating an appropriate pairing $g$ of $\mathcal{V}$ with itself will give
us the means to define \emph{orthogonality} on a vector space by saying that
the vectors $u$ and $v$ are \textbf{orthogonal }(written $u\bot v$) if
$g(u,v)=0$. Requiring nondegeneracy will insure that no nonzero vector is
orthogonal to all nonzero vectors. However, we also want to rule out the
undesirable possibility that $g(u,v)=0$ while $g(v,u)\neq0$. Thus we will also
insist that $g$ be \textbf{reflexive}: $g(u,v)=0$ if and only if $g(v,u)=0$.
$g$ will be reflexive if it is \emph{symmetric} ($g(u,v)=g(v,u)$ for all
$u,v$), or if it is \emph{alternating} ($g(v,v)=0$ for all $v$, which implies
$g(u,v)=-g(v,u)$ for all $u,v$). Here we will be confining our attention to
the case where $g$ is symmetric and nondegenerate. By designating a specific
symmetric nondegenerate self-pairing as the \textbf{scalar product }on a
vector space, the vector space with its designated scalar product becomes a
\textbf{scalar product space}. The scalar product $g:\mathcal{V}%
\times\mathcal{V}\rightarrow\mathcal{F}$ then provides the orthogonality
relation $\bot$ on $\mathcal{V}$.

\begin{exercise}
\label{Nondegenerate}A symmetric pairing $g$ of a vector space with itself is
nondegenerate if and only if $g(v,w)=0$ for all $w$ implies $v=0$.
\end{exercise}

Notice that the two included maps of a scalar product are identical, and
abusing notation somewhat, we will use the same letter $g$ to refer both to
the included map and to the scalar product itself. $g:\mathcal{V}%
\times\mathcal{V}\rightarrow\mathcal{F}$ and $g:\mathcal{V}\rightarrow
\mathcal{V}^{\top}$ will not be readily confused. When $\mathcal{V}$ is
finite-dimensional, the scalar product's included map $g$ is an isomorphism
that allows each $\varphi\in\mathcal{V}^{\top}$ to be represented by the
vector $g^{-1}(\varphi)$ via
\[
\varphi(w)=g(g^{-1}(\varphi),w),
\]
so that the evaluation of the functional $\varphi$ at any vector $w$ can be
replaced by the evaluation of the scalar product at $(g^{-1}(\varphi),w)$. In
slightly different notation, with $e$ being the evaluation pairing introduced
in Exercise \ref{Evaluation} above, the formula reads%
\[
e(\varphi,w)=g(g^{-1}(\varphi),w),
\]
which shows how the two seemingly different pairings, $e$ and $g$, are
essentially the same when $\mathcal{V}$ is a finite-dimensional scalar product
space. The natural nondegenerate pairing $e$ between $\mathcal{V}^{\top}$ and
$\mathcal{V}$ thus takes on various interpretations on $\mathcal{V}$ depending
on the choice of scalar product on the finite-dimensional space. For instance,
given a particular vector $w$, the functionals $\varphi$ such that
$\varphi(w)=0$ determine by $v=g^{-1}(\varphi)$ the vectors $v$ for which
$v\bot w$. Another example is the interpretation of the dual of a basis for
$\mathcal{V}$ as another basis for $\mathcal{V}$.

\subsection{Reciprocal Basis}

\emph{Notice:} We assume for the remainder of the chapter that we are treating
a scalar product space $\mathcal{V}$ of finite dimension $d$ over the field
$\mathcal{F}$, and that $\mathcal{V}$ has the scalar product $g$.\medskip

Let $\mathcal{B}=\{x_{1},\ldots,x_{d}\}$ be a basis for $\mathcal{V}$ and
$\mathcal{B}^{\top}=\{x_{1}^{\top},\ldots,x_{d}^{\top}\}$ its dual on
$\mathcal{V}^{\top}$. For each $i$, let $x_{i}^{\bot}=g^{-1}(x_{i}^{\top})$.
Then $\mathcal{B}^{\bot}=\{x_{1}^{\bot},\ldots,x_{d}^{\bot}\}$ is another
basis for $\mathcal{V}$, known as the \textbf{reciprocal} of $\mathcal{B}$.
Acting through the scalar product, the elements of $\mathcal{B}^{\bot}$ behave
the same as the coordinate functions from which they derive via $g^{-1}$ and
therefore satisfy the \textbf{biorthogonality} conditions%
\[
g(x_{i}^{\bot},x_{j})=x_{i}^{\top}(x_{j})=\delta_{i,j}=\left\{
\begin{array}
[c]{c}%
1,\text{ if }i=j\text{,}\\
0,\text{ if }i\neq j\text{,}%
\end{array}
\right.
\]
for all $i$ and $j$. On the other hand, if some basis $\{y_{1},\ldots,y_{d}\}$
satisfies $g(y_{i},x_{j})=(g(y_{i}))(x_{j})=\delta_{i,j}$ for all $i$ and $j$,
then $g(y_{i})$ must by definition be the coordinate function $x_{i}^{\top}$,
which then makes $y_{i}=g^{-1}(x_{i}^{\top})=x_{i}^{\bot}$. We therefore
conclude that biorthogonality characterizes the reciprocal. Hence, if we
replace $\mathcal{B}$ with $\mathcal{B}^{\bot}$, the original biorthogonality
formula displayed above tells us that $\left(  \mathcal{B}^{\bot}\right)
^{\bot}=\mathcal{B}$.

For each element $x_{j}$ of the above basis $\mathcal{B}$, $g(x_{j}%
)\in\mathcal{V}^{\top}$ has a (unique) dual basis representation of the form%
\[
g(x_{j})=\underset{i=1}{\overset{d}{\sum}\ }g_{i,j}\cdot x_{i}^{\top}\text{.}%
\]
Applying this representation to $x_{k}$ yields $g_{k,j}=g(x_{j},x_{k})$, so we
have the explicit formula%
\[
g(x_{j})=\underset{i=1}{\overset{d}{\sum}\ }g(x_{i},x_{j})\cdot x_{i}^{\top
}\text{.}%
\]
Applying $g^{-1}$ to this result, we get the following formula that expresses
each vector of $\mathcal{B}$ in terms of the vectors of $\mathcal{B}^{\bot}$
as%
\[
x_{j}=\overset{d}{\underset{i=1}{\sum}\ }g(x_{i},x_{j})\cdot x_{i}^{\bot
}\text{.}%
\]
Each element of $\mathcal{B}^{\bot}$ thus is given in terms of the elements of
$\mathcal{B}$ by%
\[
x_{j}^{\bot}=\overset{d}{\underset{i=1}{\sum}\ }g^{i,j}\cdot x_{i}\text{,}%
\]
where $g^{i,j}$ is the element in the $i$th row and $j$th column of the
\emph{inverse} of%
\[
\left[  g(x_{i},x_{j})\right]  =\left[  g_{i,j}\right]  =\left[
\begin{array}
[c]{ccc}%
g(x_{1},x_{1}) & \cdots & g(x_{1},x_{d})\\
\vdots & \cdots & \vdots\\
g(x_{d},x_{1}) & \cdots & g(x_{d},x_{d})
\end{array}
\right]  \text{.}%
\]
This we may verify by noting that since the $g^{i,j}$ satisfy%
\[
\overset{d}{\underset{k=1}{\sum}\ }g_{i,k}g^{k,j}=\delta_{i}^{j}=\left\{
\begin{array}
[c]{l}%
1\text{ if }i=j\text{,}\\
0\text{ otherwise,}%
\end{array}
\right.
\]
we have%
\[
\overset{d}{\underset{i=1}{\sum}}g^{i,j}\cdot x_{i}=\overset{d}{\underset
{i=1}{\sum}}g^{i,j}\underset{k=1}{\sum\ }g_{k,i}\cdot x_{k}^{\bot}=\overset
{d}{\underset{k=1}{\sum}}\overset{d}{\underset{i=1}{\sum}\ }g_{k,i}%
g^{i,j}\cdot x_{k}^{\bot}=\overset{d}{\underset{k=1}{\sum}}\delta_{k}^{j}\cdot
x_{k}^{\bot}=x_{j}^{\bot}\text{.}%
\]

\begin{exercise}
$x_{1}\wedge\cdots\wedge x_{d}=G\cdot x_{1}^{\bot}\wedge\cdots\wedge
x_{d}^{\bot}$ where $G=\det[g(x_{i},x_{j})]$.
\end{exercise}

\subsection{Related Scalar Product for Dual Space}

Given the scalar product $g$ on $\mathcal{V}$, we define the \textbf{related}
scalar product $\widetilde{g}$ on $\mathcal{V}^{\top}$. This we do in a
natural way by mapping the elements of $\mathcal{V}^{\top}$ back to
$\mathcal{V}$ via $g^{-1}$, and then declaring the scalar product of the
mapped elements to be the scalar product of the elements themselves. Thus for
$\varphi,\psi\in\mathcal{V}^{\top}$, we define $\widetilde{g}:\mathcal{V}%
^{\top}\times\mathcal{V}^{\top}\rightarrow\mathcal{F}$ by $\widetilde
{g}(\varphi,\psi)=g(g^{-1}(\varphi),g^{-1}(\psi))$. It is easy to see that
this $\widetilde{g}$ is indeed a scalar product. Note that $g(v,w)=\widetilde
{g}(g(v),g(w))$, so the included map $g:\mathcal{V}\rightarrow\mathcal{V}%
^{\top}$ is more than just an isomorphism $-$ it also preserves scalar
product. We call such a map an \textbf{isometry}. It makes$\mathcal{\ V}$ and
$\mathcal{V}^{\top}$ into \textbf{isometric} scalar product spaces which are
then scalar product space \textbf{aliases}, and we have justification to view
them as two versions of the same scalar product space by identifying each
$v\in\mathcal{V}$ with $g(v)\in\mathcal{V}^{\top}$.

For the elements of the dual basis $\mathcal{B}^{\top}=\{x_{1}^{\top}%
,\ldots,x_{d}^{\top}\}$ of the basis $\mathcal{B}=\{x_{1},\ldots,x_{d}\}$ for
$\mathcal{V}$, we then have
\[
\widetilde{g}(x_{1}^{\top},x_{j}^{\top})=g\left(  g^{-1}(x_{i}^{\top}%
),g^{-1}(x_{j}^{\top})\right)  =g\left(  x_{i}^{\bot},x_{j}^{\bot}\right)
\text{,}%
\]
and using the result we found in the previous section, the new included map
$\widetilde{g}:\mathcal{V}^{\top}\rightarrow\mathcal{V}$ has the (unique)
basis representation%
\[
\widetilde{g}(x_{j}^{\top})=\sum_{i=1}^{d}\ g\left(  x_{i}^{\bot},x_{j}^{\bot
}\right)  \cdot x_{i}\text{.}%
\]
Swapping $\mathcal{B}$ and $\mathcal{B}^{\bot}$, a formula derived in the
previous section becomes the equally valid formula%
\[
x_{j}^{\bot}=\overset{d}{\underset{i=1}{\sum}\ }g(x_{i}^{\bot},x_{j}^{\bot
})\cdot x_{i}%
\]
which has the same right hand side as the preceding basis representation
formula, so we see that $\widetilde{g}(x_{j}^{\top})=x_{j}^{\bot}$. The
included map $\widetilde{g}$ must therefore be $g^{-1}$. Applying
$g=\widetilde{g}^{-1}$ to both sides of the basis representation formula for
$\widetilde{g}(x_{j}^{\top})$, we get%
\[
x_{j}^{\top}=\sum_{i=1}^{d}\ g\left(  x_{i}^{\bot},x_{j}^{\bot}\right)  \cdot
g(x_{i})\text{,}%
\]
or%
\[
x_{j}^{\top}=\sum_{i=1}^{d}\ g\left(  x_{i}^{\bot},x_{j}^{\bot}\right)
\cdot\left(  x_{i}^{\top}\right)  ^{\bot}\text{.}%
\]
As we already know from the previous section, the inverse relationship is%
\[
\left(  x_{i}^{\top}\right)  ^{\bot}=\widetilde{g}^{-1}(x_{j})=g(x_{j}%
)=\sum_{i=1}^{d}\ g\left(  x_{i},x_{j}\right)  \cdot x_{i}^{\top}\text{.}%
\]

\begin{exercise}
\label{grecip}$\left[  g\left(  x_{i}^{\bot},x_{j}^{\bot}\right)  \right]
=\left[  g^{i,j}\right]  $, i. e., $\left[  g\left(  x_{i}^{\bot},x_{j}^{\bot
}\right)  \right]  =\left[  g\left(  x_{i},x_{j}\right)  \right]
^{-1}=\left[  g_{i,j}\right]  ^{-1}$.
\end{exercise}

\subsection{Some Notational Conventions}

Notational conventions abound in practice, and we note some of them now. It is
not unreasonable to view $\widetilde{g}$ as extending $g$ to $\mathcal{V}%
^{\top}$, and it would not be unreasonable to just write $g$ for both
(remembering that the included map of the $g$ on $\mathcal{V}^{\top}$ is the
inverse of the included map of the original $g$ on $\mathcal{V}$). Unless
explicitly stated otherwise, we will assume henceforth that the related scalar
product, no tilde required, is being used on $\mathcal{V}^{\top}$. In the
literature, either scalar product applied to a pair of vector or functional
arguments is often seen written $(v,w)$ or $(\varphi,\psi)$, entirely omitting
$g$ at the front, but we will avoid this practice here.

Relative to a fixed basis $\mathcal{B}=\{x_{1},\ldots,x_{d}\}$ for
$\mathcal{V}$ and its dual for $\mathcal{V}^{\top}$, the matrix $\left[
g(x_{i},x_{j})\right]  $ of the included map of the original $g$ may be
abbreviated as $\left[  g_{ij}\right]  $ (without any comma between the
separate subscripts $i$ and $j$), and the matrix $\left[  g(x_{i}%
,x_{j})\right]  ^{-1}$ of the related scalar product on $\mathcal{V}^{\top}$
may similarly be abbreviated as $\left[  g^{ij}\right]  $. The typical element
$g_{ij}$ of the matrix $\left[  g_{ij}\right]  $ is often used to designate
the $g$ on $\mathcal{V}$, and similarly $g^{ij}$ to designate the related $g$
on $\mathcal{V}^{\top}$. By the same token, $v^{i}$ designates the vector
$v=v^{1}\cdot x_{1}+\cdots+v^{d}\cdot x_{d}$ and $\varphi_{i}$ designates the
functional $\varphi=\varphi_{1}\cdot x_{1}^{\top}+\cdots+\varphi_{d}\cdot
x_{d}^{\top}$. The detailed evaluation of $g(v,w)$ is usually then written as
$g_{ij}v^{i}w^{j}$ which, since $i$ appears twice (only), as does $j$, is
taken to mean that summation over both $i$ and $j$ has tacitly been performed
over their known ranges of $1$ to $d$. This convention, called the
\emph{summation convention} (attributed to Albert Einstein), can save a lot of
ink by omitting many $\sum$ signs. Similarly, the detailed evaluation of the
scalar product $g(\varphi,\psi)=\widetilde{g}(\varphi,\psi)$ of the two
functionals $\varphi$ and $\psi$ appears as $g^{ij}\varphi_{i}\psi_{j}$. Other
examples are the evaluation of $\varphi=g(v)$ as $v_{i}=g_{ij}v^{j}$ (the
vector $v$ is converted to a functional $\varphi$ usually designated $v_{i}$
rather than $\varphi_{i}$) and the reverse $v=g^{-1}(\varphi)$ as
$v^{i}=g^{ij}v_{j}$. Note that $v_{i}w^{i}=g_{ij}v^{j}w^{i}=g(v)(w)=g(v,w)$
and $g^{ij}g_{jk}v^{k}=\delta_{k}^{i}v^{k}=v^{i}$.

This component-style notation, including the summation convention, is usually
referred to as \emph{tensor notation}. It has been a very popular notation in
the past, and is still widely used today in those applications where it has
gained a traditional foothold. Some concepts are easier to express using it,
but it can also bury concepts in forests of subscripts and superscripts. The
above samples of tensor notation used in scalar product calculations were
included to show how scalar product concepts are often handled notationally in
the literature, and not because we are going to adopt this notation here.

\subsection{Standard Scalar Products}

We will now define \emph{standard scalar products} on the finite-dimensional
vector space $\mathcal{V}$ over the field $\mathcal{F}$. As we know, we may
determine a bilinear function by giving its value on each pair $(x_{i},x_{j})
$ of basis vectors of whatever basis $\{x_{1},\ldots,x_{d}\}$ we might choose.
The bilinear functional $g:\mathcal{V}\times\mathcal{V}\rightarrow\mathcal{F}$
will be called a \textbf{standard} \textbf{scalar} \textbf{product} if for the
chosen basis $\{x_{1},\ldots,x_{d}\}$ we have $g(x_{i},x_{j})=0$ whenever
$i\neq j$ and each $g(x_{i},x_{i})=\eta_{i}$ is either $+1$ or $-1$. A
standard scalar product will be called \textbf{definite} if the $\eta_{i}$ are
all equal, and any other standard scalar product will be called
\textbf{indefinite}. The standard scalar product for which $\eta_{i}=+1$ for
all $i$ will be referred to as the \textbf{positive }standard scalar product
or as the \textbf{usual }standard scalar product. With the \emph{standard
basis} $\{(1,0,\ldots,0),\ldots,(0,\ldots,0,1)\}$ as the chosen basis, the
usual standard scalar product on the real $d$-dimensional space $\mathbb{R}%
^{d}$ is the well-known \emph{Euclidean inner product}, or \emph{dot product}.
A standard scalar product with chosen basis $\{x_{1},\ldots,x_{d}\}$ always
makes $x_{i}\bot x_{j}$ whenever $i\neq j$, so the chosen basis is always an
\textbf{orthogonal basis} under a standard scalar product, and because the
scalar product of a chosen basis vector with itself is $\pm1$, the chosen
basis for a standard scalar product is more specifically an
\textbf{orthonormal basis}.

\begin{exercise}
A standard scalar product is indeed a scalar product. \emph{(}For\emph{\ }the
nondegeneracy, Exercise \ref{Nondegenerate} above may be applied using each
vector from the chosen basis as a $w$.\emph{)}
\end{exercise}

\begin{exercise}
On real $2$-dimensional space $\mathbb{R}^{2}$, there is a standard scalar
product that makes some nonzero vector orthogonal to itself. However,
$\mathbb{R}^{d}$ with the Euclidean inner product as its scalar product has no
nonzero vector that is orthogonal to itself.
\end{exercise}

With a standard scalar product, the vectors of the reciprocal of the chosen
basis can be found in terms of those of the chosen basis without computation.

\begin{exercise}
The included map $g$ for a standard scalar product satisfies $g(x_{i}%
)=\eta_{i}\cdot x_{i}^{\top}$ for each basis vector $x_{i}$ of the chosen
basis, and also then $x_{i}^{\bot}=g^{-1}(x_{i}^{\top})=\eta_{i}\cdot x_{i}$
for each $i$.
\end{exercise}

A scalar product that makes a particular chosen basis orthonormal can always
be imposed on a vector space, of course. However, it is not necessarily true
that an orthonormal basis exists for each possible scalar product that can be
specified on a vector space. As long as $1+1\neq0$ in $\mathcal{F}$, an
orthogonal basis can always be found, but due to a lack of square roots in
$\mathcal{F}$, it may be impossible to find any orthogonal basis that can be
normalized. However, over the real numbers, normalization is always possible,
and in fact, over the reals, the number of $g(x_{i},x_{i})$ that equal $-1$ is
always the same for every orthonormal basis of a given scalar product space
due to the well-known \emph{Sylvester's Law of Inertia} (named for James
Joseph Sylvester who published a proof in 1852).

\begin{exercise}
Let $\mathcal{F}=\{0,1\}$ and let $\mathcal{V}=\mathcal{F}^{2}$. Let $x_{1}$
and $x_{2}$ be the standard basis vectors and let $[g(x_{i},x_{j}]=\left[
\begin{array}
[c]{rr}%
0 & 1\\
1 & 0
\end{array}
\right]  .$ Then for this $g$, $\mathcal{V}$ has no orthogonal basis.
\end{exercise}

\subsection{Orthogonal Subspaces and Hodge Star}

Corresponding to each subspace $\mathcal{X}\vartriangleleft\mathcal{V}$ is its
\textbf{orthogonal subspace} $\mathcal{X}^{\bot}$ consisting of all the
vectors from $\mathcal{V}$ that are orthogonal to every vector of
$\mathcal{X}$. A vector $v$ is thus in $\mathcal{X}^{\bot}$ if and only if
$g(v,x)=g(v)(x)=0$ for all $x\in\mathcal{X}$ or hence if and only if $g(v)$ is
in the \emph{annihilator} of $\mathcal{X}$ (Section \ref{Annihilator}). Thus
$\mathcal{X}^{\bot}=g^{-1}(\mathcal{X}^{0})$, where $\mathcal{X}%
^{0}\vartriangleleft\mathcal{V}^{\top}$ is the annihilator of $\mathcal{X}%
\vartriangleleft\mathcal{V}$, and $\mathcal{X}^{\bot}$ is therefore just
$\mathcal{X}^{0}$ being interpreted in $\mathcal{V}$ by using the scalar
product's ability to represent linear functionals as vectors.

Suppose now that we fix a basis $\mathcal{B}_{H}=\{x_{1},\ldots,x_{d}\}$ and
use it, as we did in Theorem \ref{Hodge}, in defining the map $H:\bigwedge
\mathcal{V}\rightarrow\bigwedge\mathcal{V}^{\top}$ that maps each blade to a
corresponding annihilator blade. Then replacing each $x_{i}^{\top}$ in the
corresponding annihilator blade in $\bigwedge\mathcal{V}^{\top}$ by
$x_{i}^{\bot}$ changes it into a corresponding \textbf{orthogonal blade} in
$\bigwedge\mathcal{V}$. The process of replacing each $x_{i}^{\top}$ by
$x_{i}^{\bot}$ may be viewed as extending $g^{-1}$ to $\bigwedge
\mathcal{V}^{\top}$and the resulting isomorphism, (which, incidentally, is
easily seen to be independent of the basis choice) is denoted $\bigwedge
g^{-1}:\bigwedge\mathcal{V}^{\top}\rightarrow\bigwedge\mathcal{V}$. The
composite $\bigwedge g^{-1}\circ H$ will be denoted by $\ast:\bigwedge
\mathcal{V}\rightarrow\bigwedge\mathcal{V}$ and will be called the
\textbf{Hodge star}, for Scottish geometer William Vallance Douglas Hodge
(1903 -- 1975), although it is essentially the same as an operator that
Grassmann used. Thus the ``annihilator blade map'' $H$ is reinterpreted as the
``orthogonal blade map'' $\ast$ for blades and the subspaces they represent.
Applying $\ast$ to exterior products of elements of $\mathcal{B}_{H}$ gives%
\[
\ast\left(  x_{i_{1}}\wedge\cdots\wedge x_{i_{n}}\right)  =\left(  -1\right)
^{\rho}\cdot x_{i_{n+1}}^{\bot}\wedge\cdots\wedge x_{i_{d}}^{\bot}\text{,}%
\]
where $\rho$ is the permutation $i_{1},\ldots,i_{d}$ of $1,\ldots,d$ and
$\left(  -1\right)  ^{\rho}=+1$ or $-1$ according as $\rho$ is even or odd.
Employing the usual standard scalar product with $\mathcal{B}_{H}$ as its
chosen basis, we then have, for example, $\ast(x_{1}\wedge\cdots\wedge
x_{n})=x_{n+1}\wedge\cdots\wedge x_{d}$.

\begin{exercise}
In $\mathbb{R}^{2}$, with $\mathcal{B}_{H}$ the standard basis $x_{1}=(1,0)$,
$x_{2}=(0,1)$, compute $\ast((a,b))$ and check the orthogonality for each of
these scalar products $g$ with matrix $\left[
\begin{array}
[c]{ll}%
g(x_{1},x_{1}) & g(x_{1},x_{2})\\
g(x_{2},x_{1}) & g(x_{2},x_{2})
\end{array}
\right]  $\emph{=}%
\[
a)\ \left[
\begin{array}
[c]{ll}%
1 & \ \ 0\\
0 & \ \ 1
\end{array}
\right]  ,\ b)\ \left[
\begin{array}
[c]{rr}%
1 & 0\\
0 & -1
\end{array}
\right]  ,\ c)\ \left[
\begin{array}
[c]{ll}%
0 & \ \ 1\\
1 & \ \ 0
\end{array}
\right]  ,\ d)\ \left[
\begin{array}
[c]{ll}%
\frac{5}{3} & \ \ \frac{4}{3}\\
\frac{4}{3} & \ \ \frac{5}{3}%
\end{array}
\right]  \text{.}%
\]
\end{exercise}

\begin{exercise}
For $\mathcal{X}\vartriangleleft\mathcal{V}$, $\dim\mathcal{X}^{\bot}%
=\dim\mathcal{V}-\dim\mathcal{X}$.
\end{exercise}

\begin{exercise}
For $\mathcal{X}\vartriangleleft\mathcal{V}$, $\mathcal{V}=\mathcal{X}%
\oplus\mathcal{X}^{\bot}$ if and only if no nonzero vector of $\mathcal{X}$ is
orthogonal to itself.
\end{exercise}

\begin{exercise}
For $\mathcal{X}\vartriangleleft\mathcal{V}$, $\mathcal{X}^{\bot\bot}=\left(
\mathcal{X}^{\bot}\right)  ^{\bot}=\mathcal{X}$ since $\mathcal{X}%
\vartriangleleft\mathcal{X}^{\bot\bot}$ and $\dim\mathcal{X}^{\bot\bot}%
=\dim\mathcal{X}$.
\end{exercise}

The Hodge star was defined above by the simple formula $\bigwedge g^{-1}\circ
H$. This formula could additionally be scaled, as some authors do, in order to
meet some particular normalization criterion, such as making a blade and its
star in some sense represent the same geometric content. For example,
$\bigwedge g^{-1}\circ H$ is sometimes scaled up by the factor $\sqrt{\left|
\det\left[  g(x_{i},x_{j})\right]  \right|  }$ when $\mathcal{F}$ is the field
of real numbers. However, for the time being at least, $\mathcal{F}$ will be
kept general, and the unscaled simple formula will continue to be our definition.

\subsection{Scalar Product on Exterior Powers\label{ExtScalarProd}}

We will now show that the scalar product $g$ on $\mathcal{V}$ may be extended
to each $\bigwedge^{p}\mathcal{V}$ as the bilinear function $g$ that on
$p$-blades has the value%
\[
g(v_{1}\wedge\cdots\wedge v_{p},w_{1}\wedge\cdots\wedge w_{p})=\det\left[
g(v_{i},w_{j})\right]  \text{.}%
\]
Thus we need to show that there is a unique bilinear function $g:\bigwedge
^{p}\mathcal{V}\times\bigwedge^{p}\mathcal{V}\rightarrow\mathcal{F}$ that
satisfies the above specification, and that this $g$ is symmetric and
nondegenerate. (For $p=0$, by convention we put $g(1,1)=1$ so that
$g(a\cdot1,b\cdot1)=ab$. The cases $p=0$ or $1$ need no further attention, so
we assume $p>1$ in the following treatment.)

For fixed $(v_{1},\ldots,v_{p})$, $\det\left[  g(v_{i},w_{j})\right]  $ is an
alternating $p$-linear function of $(w_{1},\ldots,w_{p})$ and by Theorem
\ref{UnivExt} there is a unique $g_{v_{1},\ldots,v_{p}}\in\left(
\bigwedge^{p}\mathcal{V}\right)  ^{\top}$ such that $g_{v_{1},\ldots,v_{p}%
}(w_{1}\wedge\cdots\wedge w_{p})=\det\left[  g(v_{i},w_{j})\right]  $. Now the
assignment of $g_{v_{1},\ldots,v_{p}}$ to $(v_{1},\ldots,v_{p})$ is an
alternating $p$-linear function and hence there is a unique map $g_{\wedge
}:\bigwedge^{p}\mathcal{V}\rightarrow\left(  \bigwedge^{p}\mathcal{V}\right)
^{\top}$such that $g_{\wedge}(v_{1}\wedge\cdots\wedge v_{p})=g_{v_{1}%
,\ldots,v_{p}}$ and therefore such that $g_{\wedge}(v_{1}\wedge\cdots\wedge
v_{p})(w_{1}\wedge\cdots\wedge w_{p})=\det\left[  g(v_{i},w_{j})\right]  $.
Evidently this $g_{\wedge}$ is an included map of the uniquely determined
bilinear functional $g:\bigwedge^{p}\mathcal{V}\times\bigwedge^{p}%
\mathcal{V}\rightarrow\mathcal{F}$ given by $g(s,t)=g_{\wedge}(s)(t)$ for any
elements $s,t$ of $\bigwedge^{p}\mathcal{V} $.

The symmetry of the $g$ we have just defined on all of $\bigwedge
^{p}\mathcal{V}$ follows readily from the symmetry of $g$ on $\mathcal{V}$,
and the equality of the determinant of a matrix with that of its transpose.

We will show nondegeneracy by showing that the included map is one-to-one. Let
$\mathcal{V}$ have the basis $\mathcal{B}=\{x_{1},\ldots,x_{d}\}$ from which
we make a basis $\mathcal{B}^{\wedge^{p}}$ for $\bigwedge^{p}\mathcal{V}$ in
the usual manner by taking exterior products of the elements of each subset of
$p$ elements of $\mathcal{B}$. Elements of $\mathcal{B}^{\wedge^{p}}$ will
denoted $x_{I}=x_{i_{1}}\wedge\cdots\wedge x_{i_{p}}$, $x_{J}=x_{j_{1}}%
\wedge\cdots\wedge x_{j_{p}}$, etc. using capital-letter subscripts that
denote the multi-index subscript sets $I=\{i_{1},\ldots,i_{p}\}$ and
$J=\{j_{1},\ldots,j_{p}\}$ of the $p$ elements of $\mathcal{B} $ that are
involved. Indexing now with these multi-indices, our familiar explicit dual
basis representation formula for the effect of the included map (now also
called $g)$ on basis vectors from $\mathcal{B}^{\wedge^{p}}$ appears as%
\[
g(x_{J})=\underset{\left|  I\right|  =p}{\overset{}{\sum}\ }g(x_{I}%
,x_{J})\cdot(x_{I})^{\top}\text{,}%
\]
with the scalar product on basis vector pairs being%
\[
g(x_{I},x_{J})=g(x_{i_{1}}\wedge\cdots\wedge x_{i_{p}},x_{j_{1}}\wedge
\cdots\wedge x_{j_{p}})=\det\left[  g(x_{i_{k}},x_{j_{l}})\right]
\]
of course.

The concept of the $p$th exterior power of a map, introduced by Exercise
\ref{MapExtPower}, will now find employment. Consider the $p$th exterior power
$\bigwedge^{p}g:\bigwedge^{p}\mathcal{V}\rightarrow\bigwedge^{p}%
\mathcal{V}^{\top}$ of the included map $g$ of the original scalar product on
$\mathcal{V}$. We have%
\[%
{\textstyle\bigwedge\nolimits^{p}}
g(x_{J})=g(x_{j_{1}})\wedge\cdots\wedge g(x_{j_{p}})=\left(  \underset
{i_{1}=1}{\overset{d}{\sum}\ }g(x_{i_{1}},x_{j_{1}})\cdot x_{i_{1}}^{\top
}\right)  \wedge\cdots\wedge\left(  \underset{i_{p}=1}{\overset{d}{\sum}%
\ }g(x_{i_{p}},x_{j_{p}})\cdot x_{i_{p}}^{\top}\right)
\]
from which it follows that%
\[%
{\textstyle\bigwedge\nolimits^{p}}
g(x_{J})=\sum_{\left|  I\right|  =p}\det\left[  g(x_{i_{k}},x_{j_{l}})\right]
\cdot x_{I}^{\top}\text{.}%
\]
So we see that the new included map has a basis representation
\[
g(x_{J})=\sum_{\left|  I\right|  =p}g_{I,J}\cdot(x_{I})^{\top}\text{,}%
\]
and the $p$th exterior power of the original included map has a basis
representation
\[%
{\textstyle\bigwedge\nolimits^{p}}
g(x_{J})=\sum_{\left|  I\right|  =p}g_{I,J}\cdot x_{I}^{\top}\text{,}%
\]
both with the same coefficients $g_{I,J}$, namely%
\[
g_{I,J}=\det\left[  g(x_{i_{k}},x_{j_{l}})\right]  \text{.}%
\]
Therefore, applying Proposition \ref{Matrixdet}, if one of these maps is
invertible, then so is the other one. Now, the invertibility of the original
included map $g:\mathcal{V}\rightarrow\mathcal{V}^{\top}$ makes its $p$th
exterior power invertible as well, because the $g(x_{i})$ being a basis for
$\mathcal{V}^{\top}$make the $\bigwedge^{p}g(x_{I})$ a basis for
$\bigwedge^{p}\mathcal{V}^{\top}$. Therefore the new included map is also
invertible and our new symmetric bilinear $g$ is indeed nondegenerate as we
wished to show. What is going on here is rather transparent: the $p$th
exterior power of the included map $g:\mathcal{V}\rightarrow\mathcal{V}^{\top
}$ followed by the isomorphism$:\widehat{\Phi}:\bigwedge^{p}\mathcal{V}^{\top
}\rightarrow\left(  \bigwedge^{p}\mathcal{V}\right)  ^{\top}$ of Theorem
\ref{DualExtIso}, evidently is the included map of the extension. That is, we
have the formula $g(s,t)=((\widehat{\Phi}\circ\bigwedge^{p}g)(s))(t)$, valid
for all $(s,t)\in\bigwedge^{p}\mathcal{V}\times\bigwedge^{p}\mathcal{V}$. If
we agree to identify $\left(  \bigwedge^{p}\mathcal{V}\right)  ^{\top}$ with
$\bigwedge^{p}\mathcal{V}^{\top}$ through the isomorphism $\widehat{\Phi}$,
then $\bigwedge^{p}g$ ``is'' the included map of our scalar product extension
to $\bigwedge^{p}\mathcal{V}$.

The Hodge star provides an isomorphic mapping of $\bigwedge^{p}\mathcal{V}$
with $\bigwedge^{d-p}\mathcal{V}$ that is also an isometry (except for scale),
as the following proposition records.

\begin{proposition}
\label{StarIso}Let $s,r\in\bigwedge^{p}\mathcal{V}$, let $\left\{
x_{1},\ldots,x_{d}\right\}  $ be the basis for $\mathcal{V}$ used in defining
$\ast$, and let $G=\det\left[  g(x_{i},x_{j})\right]  $. Then%
\[
g(r,s)=Gg(\ast r,\ast s)\text{.}%
\]
\end{proposition}

Proof: It suffices to prove the result for $r=x_{I}$ and $s=x_{J}$. Now%
\[
g(x_{I},x_{J})=\det\left[
\begin{array}
[c]{ccc}%
g_{i_{1},j_{1}} & \cdots & g_{i_{1},j_{p}}\\
\vdots & \cdots & \vdots\\
g_{i_{p},j_{1}} & \cdots & g_{i_{p},j_{p}}%
\end{array}
\right]
\]
and in light of Exercise \ref{grecip}%
\[
g(\ast x_{I},\ast x_{J})=\left(  -1\right)  ^{\rho}\left(  -1\right)
^{\sigma}\det\left[
\begin{array}
[c]{ccc}%
g^{i_{p+1},j_{p+1}} & \cdots & g^{i_{p+1},j_{d}}\\
\vdots & \cdots & \vdots\\
g^{i_{d},j_{p+1}} & \cdots & g^{i_{d},j_{d}}%
\end{array}
\right]
\]
where $\rho$ and $\sigma$ are the respective permutations $i_{1},\ldots,i_{d}$
and $j_{1},\ldots,j_{d}$ of $\{1,\ldots,d\}$. By Jacobi's Determinant Identity
(Lemma \ref{Jacobi})%
\[
\det\left[
\begin{array}
[c]{ccc}%
g_{i_{1},j_{1}} & \cdots & g_{i_{1},j_{p}}\\
\vdots & \cdots & \vdots\\
g_{i_{p},j_{1}} & \cdots & g_{i_{p},j_{p}}%
\end{array}
\right]  =\det\left[
\begin{array}
[c]{ccc}%
g_{i_{1},j_{1}} & \cdots & g_{i_{1},j_{d}}\\
\vdots & \cdots & \vdots\\
g_{i_{d},j_{1}} & \cdots & g_{i_{d},j_{d}}%
\end{array}
\right]  \det\left[
\begin{array}
[c]{ccc}%
g^{i_{p+1},j_{p+1}} & \cdots & g^{i_{p+1},j_{d}}\\
\vdots & \cdots & \vdots\\
g^{i_{d},j_{p+1}} & \cdots & g^{i_{d},j_{d}}%
\end{array}
\right]
\]
from which it follows that $g(x_{I},x_{J})=Gg(\ast x_{I},\ast x_{J})$ once we
have attributed the $\left(  -1\right)  ^{\rho}\left(  -1\right)  ^{\sigma}$
to the rearranged rows and columns of the big determinant. $\blacksquare$

\subsection{Another Way to Define $\ast\label{other*}$}

A version of the Hodge star may be defined using the scalar products defined
on the $\bigwedge^{p}\mathcal{V}$ in the previous section. Let $\mathcal{B}%
=\{x_{1},\ldots,x_{d}\}$ be a designated basis for $\mathcal{V}$ and consider
the exterior product $s\wedge t$ where $s$ is a fixed element of
$\bigwedge^{p}\mathcal{V}$ and $t$ is allowed to vary over $\bigwedge
^{d-p}\mathcal{V}$. Each such $s\wedge t$ is of the form $c_{s}\left(
t\right)  \cdot x_{1}\wedge\cdots\wedge x_{d}$ where $c_{s}\in\left(
\bigwedge^{d-p}\mathcal{V}\right)  ^{\top}$. Hence there is a unique element
$\ast s\in\bigwedge^{d-p}\mathcal{V}$ such that $g\left(  \ast s,t\right)
=c_{s}\left(  t\right)  =\left(  x_{1}\wedge\cdots\wedge x_{d}\right)  ^{\top
}\left(  s\wedge t\right)  $ and this unique element $\ast s$ for which%
\[
s\wedge t=g\left(  \ast s,t\right)  \cdot x_{1}\wedge\cdots\wedge x_{d}%
\]
will be the new version of the Hodge star on $\bigwedge^{p}\mathcal{V}$ that
we will now scrutinize and compare with the original version. We will focus on
elements $s$ of degree $p$ (meaning those in $\bigwedge^{p}\mathcal{V)}$, but
we do intend this $\ast$ to be applicable to arbitrary elements of
$\bigwedge\mathcal{V}$ by applying it separately to the components of each
degree and summing the results.

If we have an element $r\in\bigwedge^{d-p}\mathcal{V}$ that we allege is our
new $\ast s$ for some $s\in^{{}}\bigwedge^{p}\mathcal{V}$, we can verify that
by showing that $g\left(  r,t\right)  $ is the coefficient of $x_{1}%
\wedge\cdots\wedge x_{d}$ in $s\wedge t$ for all $t$ in a basis for
$\bigwedge^{d-p}\mathcal{V}$. If we show this for such an $r$ corresponding to
each $s$ that is in a basis for $\bigwedge^{p}\mathcal{V}$, then we will have
completely determined the new $\ast$ on $\bigwedge^{p}\mathcal{V}$.

Let us see what happens when we choose the designated basis $\mathcal{B}$ to
be the basis $\mathcal{B}_{H}$ used in defining the annihilator blade map $H$,
with exactly the same assignment of subscript labels to the basis vectors
$x_{i}$. With $s$ and $t$ equal to the respective basis monomials
$x_{I}=x_{i_{1}}\wedge\cdots\wedge x_{i_{p}}$ and $x_{J}=x_{j_{p+1}}%
\wedge\cdots\wedge x_{j_{d}}$, we have%
\[
s\wedge t=x_{I}\wedge x_{J}=\varepsilon_{I,J}\cdot x_{1}\wedge\cdots\wedge
x_{d}%
\]
where $\varepsilon_{I,J}=0$ if $\left\{  j_{p+1},\ldots,j_{d}\right\}  $ is
not complementary to $\left\{  i_{1},\ldots,i_{p}\right\}  $ in $\{1,\ldots
,d\}$, and otherwise $\varepsilon_{I,J}=\left(  -1\right)  ^{\sigma}$ with
$\sigma$ the permutation $i_{1},\ldots,i_{p},j_{p+1},\ldots,j_{d}$ of
$1,\ldots,d$. The original Hodge star gives%
\[
\ast x_{I}=\ast\left(  x_{i_{1}}\wedge\cdots\wedge x_{i_{p}}\right)  =\left(
-1\right)  ^{\rho}\cdot x_{i_{p+1}}^{\bot}\wedge\cdots\wedge x_{i_{d}}^{\bot}%
\]
where $\rho$ is the permutation $i_{1},\ldots,i_{d}$ of $1,\ldots,d$, so for
this original $\ast x_{I}$,%
\begin{align*}
g\left(  \ast x_{I},x_{J}\right)   &  =\left(  -1\right)  ^{\rho}\det\left[
\begin{array}
[c]{ccc}%
g(x_{i_{p+1}}^{\bot},x_{j_{p+1}}) & \cdots & g(x_{i_{p+1}}^{\bot},x_{j_{d}})\\
\vdots & \cdots & \vdots\\
g(x_{i_{d}}^{\bot},x_{j_{p+1}}) & \cdots & g(x_{i_{d}}^{\bot},x_{j_{d}})
\end{array}
\right] \\
&  =\left(  -1\right)  ^{\rho}\det\left[
\begin{array}
[c]{ccc}%
x_{i_{p+1}}^{\top}(x_{j_{p+1}}) & \cdots & x_{i_{p+1}}^{\top}(x_{j_{d}})\\
\vdots & \cdots & \vdots\\
x_{i_{d}}^{\top}(x_{j_{p+1}}) & \cdots & x_{i_{d}}^{\top}(x_{j_{d}})
\end{array}
\right]  \text{.}%
\end{align*}
If $J$ is not complementary to $I$ in $\{1,\ldots,d\}$ then $j_{p+1}%
,\ldots,j_{d}$ is not a permutation of $i_{p+1},\ldots,i_{d}$ and some
$j_{p+k}$ equals none of $i_{p+1},\ldots,i_{d}$ so that the $k$th column is
all zeroes and the latter determinant vanishes. On the other hand, if $J$ is
complementary to $I$ in $\{1,\ldots,d\}$ then $J=\left\{  i_{p+1},\ldots
,i_{d}\right\}  $ and there is no reason why we cannot assume that
$j_{p+1}=i_{p+1},\ldots,j_{d}=i_{d}$, which then makes the determinant equal
to $1$ and makes $\sigma$ equal to $\rho$. Hence using the same basis ordered
the same way, our new version of the Hodge star is exactly the same as the
original version.

\begin{exercise}
For a given scalar product space, using $\mathcal{B}^{\prime}=\{x_{1}^{\prime
},\ldots,x_{d}^{\prime}\}$ instead of $\mathcal{B}=$ $\{x_{1},\ldots,x_{d}\}$
to define the Hodge star gives $h\cdot\ast$ instead of $\ast$, where $h$ is
the nonzero factor, independent of $p$, such that $h\cdot x_{1}^{\prime}%
\wedge\cdots\wedge x_{d}^{\prime}=x_{1}\wedge\cdots\wedge x_{d}$. Thus, no
matter which definition is used for either, or what basis is used in defining
either, for any two of our Hodge stars, over all of $\bigwedge\mathcal{V}$ the
values of one are the same constant scalar multiple of the values of the
other. That is, ignoring scale, for a given scalar product space $\mathcal{V}$
all of our Hodge stars on $\bigwedge\mathcal{V}$ are identical, and the scale
depends only on the basis choice \emph{(}including labeling\emph{) }for
$\mathcal{V}$ used in each definition.
\end{exercise}

\begin{exercise}
\emph{(}Continuation\emph{)} Suppose that for the scalar product $g$, the
bases $\mathcal{B}$ and $\mathcal{B}^{\prime}$ have the same \textbf{Gram
determinant}, i.e., $\det\left[  g(x_{i},x_{j})\right]  =\det\left[
g(x_{i}^{\prime},x_{j}^{\prime})\right]  $, or equivalently, $g(x_{1}%
\wedge\cdots\wedge x_{d},x_{1}\wedge\cdots\wedge x_{d})=g(x_{1}^{\prime}%
\wedge\cdots\wedge x_{d}^{\prime},x_{1}^{\prime}\wedge\cdots\wedge
x_{d}^{\prime})$. Taking it as a given that in a field, $1$ has only itself
and $-1$ as square roots, $\mathcal{B}$ and $\mathcal{B}^{\prime}$ then
produce the same Hodge star \emph{up to sign}.
\end{exercise}

The following result is now apparent.

\begin{proposition}
Two bases of a given scalar product space yield the same Hodge star up to sign
if and only if they have the same Gram determinant. $\blacksquare$
\end{proposition}

\medskip

Finally, we obtain some more formulas of interest, and some conclusions based
on them. Putting $t=\ast r$ in the defining formula above gives%
\[
s\wedge\ast r=g\left(  \ast s,\ast r\right)  \cdot x_{1}\wedge\cdots\wedge
x_{d}\text{.}%
\]
Similarly, we find%
\[
r\wedge\ast s=g\left(  \ast r,\ast s\right)  \cdot x_{1}\wedge\cdots\wedge
x_{d}=g\left(  \ast s,\ast r\right)  \cdot x_{1}\wedge\cdots\wedge
x_{d}\text{,}%
\]
and therefore for all $r,s\in\bigwedge^{p}\mathcal{V}$ we have%
\[
r\wedge\ast s=s\wedge\ast r\text{.}%
\]
Applying Proposition \ref{StarIso}, we also get%
\[
s\wedge\ast r=G^{-1}g\left(  r,s\right)  \cdot x_{1}\wedge\cdots\wedge
x_{d}\text{,}%
\]
where $G=\det\left[  g(x_{i},x_{j})\right]  $. The results of the next two
exercises then follow readily.

\begin{exercise}
For all $r,s\in\bigwedge^{p}\mathcal{V}$,
\[
g\left(  r,s\right)  =G\cdot\ast\left(  r\wedge\ast s\right)  =G\cdot
\ast\left(  s\wedge\ast r\right)  \text{.}%
\]
\end{exercise}

\begin{exercise}
\label{*expansion}For any $r\in\bigwedge^{p}\mathcal{V}$,%
\[
\ast r=G^{-1}\sum_{|I|=p}\left(  -1\right)  ^{\rho}g(r,x_{I})\cdot
x_{\,\overline{I\,}}\text{,}%
\]
where $I=\{i_{1},\ldots,i_{p}\}$, $\overline{I\,}=\{i_{p+1},\ldots,i_{d}\}$,
$I\cup\overline{I\,}=\left\{  1,\ldots,d\right\}  $, and%
\[
\rho=\left(
\begin{array}
[c]{ccc}%
1 & \cdots & d\\
i_{1} & \cdots & i_{d}%
\end{array}
\right)  \text{.}%
\]
Thus, not only do we have a nice expansion formula for $\ast r$, but we may
also conclude that $\ast r$ has yet another definition as the unique $t\in$
$\bigwedge^{d-p}\mathcal{V}$ for which $s\wedge t=G^{-1}g\left(  r,s\right)
\cdot x_{1}\wedge\cdots\wedge x_{d}$ for all $s\in\bigwedge^{p}\mathcal{V} $.
\end{exercise}

\subsection{Hodge Stars for the Dual Space}

We will now define a related Hodge star $\widetilde{\ast}$ on $\bigwedge
\mathcal{V}^{\top}$ in a way similar to how we defined $\widetilde{g}$ on
$\mathcal{V}^{\top}$. We will use $\bigwedge g^{-1}$ to map a blade from
$\bigwedge\mathcal{V}^{\top}$ to $\bigwedge\mathcal{V}$, use $\ast$ to get an
orthogonal blade for the result and map that back to $\bigwedge\mathcal{V}%
^{\top}$ using the inverse of $\bigwedge g^{-1}$. So we will define
$\widetilde{\ast}:\bigwedge\mathcal{V}^{\top}\rightarrow\bigwedge
\mathcal{V}^{\top}$ as $\left(  \bigwedge g^{-1}\right)  ^{-1}\circ\ast
\circ\bigwedge g^{-1}$, which works out to be $H\circ\bigwedge g^{-1}$ when
the unscaled $\bigwedge g^{-1}\circ H$ is substituted for $\ast$. Thus,
applying $\widetilde{\ast}$ to an exterior product of vectors from the dual
basis $\mathcal{B}_{H}^{\top}$ gives
\[
\widetilde{\ast}\left(  x_{i_{1}}^{\top}\wedge\cdots\wedge x_{i_{n}}^{\top
}\right)  =H\left(  x_{i_{1}}^{\bot}\wedge\cdots\wedge x_{i_{n}}^{\bot
}\right)  \text{.}%
\]
Employing the usual standard scalar product with $\mathcal{B}_{H}$ as its
chosen basis, we then have, for example, $\widetilde{\ast}(x_{1}^{\top}%
\wedge\cdots\wedge x_{n}^{\top})=x_{n+1}^{\top}\wedge\cdots\wedge x_{d}^{\top
}$.

\begin{exercise}
Using $\widetilde{g}$ to define orthogonality for $\mathcal{V}^{\top}$,
$\widetilde{\ast}\beta$ is an orthogonal blade of the blade $\beta$ in
$\bigwedge\mathcal{V}^{\top}$.
\end{exercise}

As an alternative to the related Hodge star $\widetilde{\ast}$ that we have
just defined, we could define a Hodge star on $\bigwedge\mathcal{V}^{\top}$ in
the manner of the previous section, using the related scalar product
$\widetilde{g}$ on $\mathcal{V}^{\top}$ extended to each $\bigwedge
^{p}\mathcal{V}^{\top}$.

\begin{exercise}
Compare $\widetilde{\ast}$ with the Hodge star defined on $\bigwedge
\mathcal{V}^{\top}$ in the manner of the previous section, using the related
scalar product $\widetilde{g}$ on $\mathcal{V}^{\top}$.
\end{exercise}

\newpage

\subsection{Problems}

\begin{enumerate}
\item For finite-dimensional spaces identified with their double duals, the
included maps belonging to a pairing are the duals of each other.

\item A pairing of a vector space with itself is reflexive if and only if it
is either symmetric or alternating.

\item Give an example of vectors $u,v,w$ such that $u\bot v$ and $v\bot w$ but
it is not the case that $u\bot w$.

\item The pairing $g$ of two $d$-dimensional vector spaces with respective
bases $\{x_{1},\ldots,x_{d}\}$ and $\{y_{1},\ldots,y_{d}\}$ is nondegenerate
if and only if%
\[
\det\left[
\begin{array}
[c]{ccc}%
g(x_{1},y_{1}) & \cdots & g(x_{1},y_{d})\\
\vdots & \cdots & \vdots\\
g(x_{d},y_{1}) & \cdots & g(x_{d},y_{d})
\end{array}
\right]  \neq0.
\]

\item A pairing $g$ of two finite-dimensional vector spaces $\mathcal{V}$ and
$\mathcal{W}$ with respective bases $\{x_{1},\ldots,x_{d}\}$ and
$\{y_{1},\ldots,y_{e}\}$ is an element of the functional product space
$\mathcal{V}^{\top}\mathcal{W}^{\top}$ and may be expressed as $g=\sum
_{i,j}g(x_{i},y_{j})\,\,x_{i}^{\top}\,y_{j}^{\top}$.

\item When $\mathcal{F}$ is the field with just two elements, so that $1+1=0$,
any scalar product on $\mathcal{F}^{2}$ will make some nonzero vector
orthogonal to itself.

\item For a scalar product space over a field where $1+1\neq0$, the scalar
product $g$ satisfies $g(w-v,w-v)=g(v,v)+g(w,w)$ if and only if $v\bot w$.

\item What, if anything, is the Hodge star on $\bigwedge\mathcal{V}$ when
$\dim\mathcal{V}=1$?

\item With the same setup as Exercise \ref{*expansion},%
\[
r=\sum_{|I|=p}(-1)^{\rho}g(\ast r,x_{\overline{I}})\cdot x_{I}=\sum
_{|J|=d-p}(-1)^{\overline{J}J}g(\ast r,x_{J})\cdot x_{\overline{J}},
\]
while%
\[
\ast(\ast r)=G^{-1}\sum_{|J|=d-p}(-1)^{J\overline{J}}g(\ast r,x_{J})\cdot
x_{\overline{J}}.
\]
where $J=\{j_{1},\ldots,j_{d-p}\}$, $\overline{J\,}=\{j_{d-p+1},\ldots
,j_{d}\}$, and $J\cup\overline{J\,}=\{1,\ldots,d\}$. Hence%
\[
\ast(\ast r)=G^{-1}(-1)^{p(d-1)}\cdot r=G^{-1}(-1)^{p(d-p)}\cdot r.
\]

\item $\ast^{-1}\left(  \ast s\wedge\ast t\right)  =H^{-1}\left(  H\left(
s\right)  \wedge H\left(  t\right)  \right)  =s\vee t$ for $s,t\in
\bigwedge\mathcal{V}$. On the other hand, $\widetilde{\ast}\left(
\widetilde{\ast}^{-1}\sigma\wedge\widetilde{\ast}^{-1}\tau\right)  =H\left(
H^{-1}\left(  \sigma\right)  \wedge H^{-1}\left(  \tau\right)  \right)
=\sigma\widetilde{\vee}\tau$ for $\sigma,\tau\in\bigwedge\mathcal{V}^{\top}$.
How does $\widetilde{\vee}$ compare to the regressive product defined on
$\bigwedge\mathcal{V}^{\top}$ using the dual of the basis used for
$\bigwedge\mathcal{V}$? How do $\ast\left(  \ast s\wedge\ast t\right)  $ and
$\widetilde{\ast}\left(  \widetilde{\ast}\sigma\wedge\widetilde{\ast}%
\tau\right)  $ relate to the same respective regressive products?

\item In $\mathbb{R}^{3}$ with the standard basis used for defining $H$ and as
the chosen basis for the usual standard scalar product, the familiar cross
product $u\times v$ of two vectors $u$ and $v$ is the same as $\ast(u\wedge
v)$.

\item Using multi-index notation as in Section \ref{ExtScalarProd} above,%
\begin{align*}
x_{I}\wedge\ast x_{I}  &  =g\left(  x_{I},x_{I}\right)  \cdot x_{1}^{\bot
}\wedge\cdots\wedge x_{d}^{\bot}\\
&  =g\left(  x_{I},x_{I}\right)  \cdot\left(  G^{-1}\cdot x_{1}\wedge
\cdots\wedge x_{d}\right)
\end{align*}
where $\{x_{1},\ldots,x_{d}\}$ is the basis used in defining $\ast$, and
$G=\det\left[  g(x_{i},x_{j})\right]  $.

\item Using $\{x_{1},\ldots,x_{d}\}$ as the basis in defining $\ast$, and
using multi-index notation as in Section \ref{ExtScalarProd} above, then%
\[
\ast x_{J}=G^{-1}\sum_{\left|  I\right|  =\left|  J\right|  }\left(
-1\right)  ^{\rho}g\left(  x_{I},x_{J}\right)  \cdot x_{\,\overline{I\,}%
}\ \text{,}%
\]
but%
\[
\widetilde{\ast\,}x_{J}^{\top}=\sum_{\left|  I\right|  =\left|  J\right|
}\left(  -1\right)  ^{\rho}\widetilde{g}\left(  x_{I},x_{J}\right)  \cdot
x_{\,\overline{I\,}}^{\top}\ \text{,}%
\]
where $I=\{i_{1},\ldots,i_{p}\}$, $J=\{j_{1},\ldots,j_{p}\}$, $\overline
{I\,}=\{i_{p+1},\ldots,i_{d}\}$ is such that $I\cup\overline{I\,}=\left\{
1,\ldots,d\right\}  $, $\rho=\left(
\begin{array}
[c]{ccc}%
1 & \cdots & d\\
i_{1} & \cdots & i_{d}%
\end{array}
\right)  $ and $G=\det\left[  g(x_{i},x_{j})\right]  $.

\item Let $\{x_{1},\ldots,x_{d}\}$ be the basis used in defining $\ast$, let
$\mathcal{F}$ be the real numbers, and let $d=\dim\mathcal{V}>1$. The
symmetric bilinear functional $g$ on $\mathcal{V}$ such that
\[
g(x_{i},x_{j})=\left\{
\begin{array}
[c]{r}%
1\text{, when }i\neq j\\
0\text{, otherwise }%
\end{array}
\right.
\]
is a scalar product that makes $x_{j}\wedge\ast x_{j}=0$ for each $j$. (The
nondegeneracy of $g$ follows readily upon observing that multiplying the
matrix $\left[  g(x_{i},x_{j})\right]  $ on the right by the matrix $\left[
h_{i,j}\right]  $, where
\[
h_{i,j}=\left\{
\begin{array}
[c]{r}%
-1\text{, when }i\neq j\\
d-2\text{, otherwise }%
\end{array}
\right.  \text{,}%
\]
gives $-(d-1)$ times the identity matrix.)

\item If two bases have the same Gram determinant with respect to a given
scalar product, then they have the same Gram determinant with respect to any
scalar product.

\item All chosen bases that produce the same positive standard scalar product
also produce the same Hodge star up to sign.
\end{enumerate}
\end{document}